\documentclass{article}
\usepackage{amsmath}
\usepackage{amsfonts}
\usepackage{amssymb}
\newcommand{\bm}[1]{\mbox{\boldmath{$#1$}}}
\newtheorem{theorem}{Theorem}[section]
\newtheorem{lemma}[theorem]{Lemma}
\newtheorem{proposition}[theorem]{Proposition}
\newtheorem{corollary}[theorem]{Corollary}
\newtheorem{problem}[theorem]{Problem}
\newtheorem{remark}[theorem]{Remark}
\newtheorem{definition}[theorem]{Definition}

\begin{document}

\title{Existence and Uniqueness Results for Double-Free-Boundary Bernoulli Problems in Fluid Dynamics\footnote{ \vspace{.1cm}\newline Some of these results (especially Thm. \ref{thm1.2}) were presented in the author's talk entitled Operator methods for double-free-boundary fluid problems in annular domains. At the 7 th. AIMS International Conference on Dynamical Systems, Differential Equations and Applications, University of Texas at Arlington, 17-21 May, 2008, Special Session on Partial Differential Equations in Fluid Mechanics and Mathematical Physics, organized by Ning Ju, Xiaming Wang, and Jiahang Wu.}}
\author{Andrew F Acker III\\Department of Mathematics and Statistics\\Wichita State University\\Wichita, KS 67208-0033\\USA\\e-mail: acker@math.wichita.edu}
\date{April 18, 2016}

\maketitle

\begin{abstract}
\noindent
We prove existence theorems for the double-free-boundary Bernoulli problem in two space dimensions, in which, given the strictly-positive, smooth "flow-speed" functions $a_i(p):\Re^2\rightarrow\Re_+$, $i=1,2$, one seeks an ideal fluid flow in an annular flow-domain $\Omega$ (also to be determined), whose boundary components $\Gamma_i$, $i=1,2$, are such that $|{\nabla}U(p)|=a_i(p)$ on $\Gamma_i$, $i=1,2$, where $U(p):\Omega\rightarrow\Re$ denotes the stream-function in $\Omega$. 
The existence result (Thms. \ref{thm1.2} and \ref{thm2.1.0}) states that given a strict inner-solution pair $(\tilde{\Gamma}_1,\tilde{\Gamma}_2)$ (defined such that $(-1)^i(a_i(p)-|{\nabla} \tilde{U}(p)|)<0$ on $\tilde{\Gamma}_i$, $i=1,2$) and a strict outer-solution pair $(\hat{\Gamma}_1,\hat{\Gamma}_2)$ (defined such that $(-1)^i(a_i(p)-|{\nabla} \hat{U}(p)|)>0$ on $\hat{\Gamma}_i$, $i=1,2$), such that the inner-solution pair lies component-wise inside the outer-solution pair, there exists a classical solution pair $(\dot{\Gamma}_1,\dot{\Gamma}_2)$, which lies component-wise outside the inner-solution pair and inside the outer-solution pair.
\vspace{.1in}

\noindent
Consider Prob. \ref{prob 1.1} in the special case where the functions $a_1(p)$, $a_2(p)$ both coincide with a single strictly logarithmically subharmonic flow-speed function $a(p):G\rightarrow\Re_+$ defined in a finite annular domain $G$. We show that there cannot exist more than one classical solution.
\vspace{.1in}

\noindent
Finally, we apply Thm. \ref{thm1.2} to construct continuously-varying and suitably monotone families of solution-pairs for the double-free-boundary Bernoulli problem corresponding to various parameter-pairs $(\lambda_1,\lambda_2)\in\Re^2_+$ and the related function-pairs $(a_1(p),a_2(p))$ defined such that $a_1(p):=\lambda_1\,a(p)$ and $a_2(p):=\lambda_2\,a(p)$.
\end{abstract}

\tableofcontents
\section{Introduction}
\label{section 1}
\subsection{Existence}
\label{subsection 1.1}
Consider an ideal fluid, flowing in a channel in two-dimensional space under the influence of a potential-energy function of the two space coordinates. In the simplest approximation of such a flow, called the "narrow-stream limit", the flow-channel is modeled by a smooth simple closed curve chosen to locally minimize the arc-length integral of the flow-speed, which is given by a positive $C^1$-function $a(p)$ of the space coordinates $p=(x,y)$. The Euler equation for the total-flow-speed minimization problem states that
\begin{equation}
\label{eqn 1.0}
K(p)=\big(a_{\boldsymbol\nu}(p)/a(p)\big)
\end{equation}
along the flow-curve $\Gamma$, where $\boldsymbol{\nu}$ and $K(p)$ denote the left-hand normal to $\Gamma$ at $p\in\Gamma$ and the corresponding counter-clockwise-oriented curvature at $p$. For example, for a flow governed by Bernoulli's law (at constant pressure and fluid density), the speed $a(p)$ of the flow at any point $p\in\Re^2$ would be related to the potential energy at that same point in such a way that the sum of twice the potential energy density with the square of the flow-speed equals a positive constant. Thus the path of the flow through a potential-energy terrain in $\Re^2$ would tend to follow the high-potential-energy ridges (corresponding to valleys of the flow-speed function), while avoiding low-potential-energy regions, where the flow-speed is relatively high.
\vspace{.1in}

\noindent
In the much more detailed 2-dimensional stream-model to be studied in the present paper, we assume a flow-channel (or stream-bed) of finite width, in which there is a flow of ideal fluid. This has been a widely accepted flow model in free-boundary studies of jets, wakes, cavitation, ocean waves, etc. The ideal flow is by definition an incompressible, irrotational, inviscid, steady-state, two-dimensional flow in a (therefore) two-dimensional flow-region (or stream-bed).
We also assume a finite, annular flow-channel. Thus, we have invoked the accepted mathematical idealization of the stream which flows around a finite closed loop back into itself, and thus doesn't begin or end anywhere (or flow to or from infinity). The two boundary components of the annular flow-channel, as well as the (harmonic) stream function for the ideal fluid-flow in the flow-channel, are all determined together as elements of the solution of a double-free-boundary problem, in which the flow-speed along the free boundaries is again determined point-wise by a positive flow-speed function $a(p)$, or by two independent positive flow-speed functions $a_1(p)$, $a_2(p)$, independently governing the flow-speeds along the two boundary components of the stream. In either case, the flow-speed functions are a direct input into the stream-model, which applies to a large class of flow-speed functions. Our flow-model leads to the following flow-problem.

\begin{problem} {(Double-free-boundary Bernoulli problem for annular  domains in  $\Re^2$: the case of two independent flow-speed functions (for the equivalent formulation involving periodic strip-like domains, see Prob. \ref{prob2.1.1}))}
\label{prob 1.1}
In $\Re^2$ (where $p=(x,y)$), let ${\rm X}$ denote the family of all simple closed curves $\Gamma$. Let ${\bf X}$ denote the family of all ordered pairs ${\bf \Gamma}=(\Gamma_1,\Gamma_2)$ of elements of ${\rm X}$ such that $\Gamma_1<\Gamma_2$ in the sense that ${\rm Cl}(D(\Gamma_1))\subset D(\Gamma_2)$, where we use $D(\Gamma)$ (resp. $E(\Gamma))$ to denote the interior (resp. exterior) complement of $\Gamma\in {\rm X}$. Given the strictly-positive $C^{2}$-functions $a_i(p)=a_i(x,y):\Re^2\rightarrow\Re_+$, $i=1,2$, which we call "flow-speed functions", we seek a smooth curve-pair ${\bf\Gamma}=(\Gamma_1,\Gamma_2)\in{\bf X}$ (the flow boundaries) such that 
\begin{equation}
\label{eqn 1.1}
|{\nabla} U(p)|=\,a_i(p)\,\,{\rm on}\,\,\Gamma_i,
\end{equation}
$i=1,2$, where the region $\Omega=\Omega({\bf\Gamma}):=D(\Gamma_2)\cap E(\Gamma_1)$ is the annular (i.e. doubly connected) stream bed, and $U(p):=U({\bf\Gamma};p)$ denotes the stream function (or the capacitary potential) in ${\rm Cl}(\Omega({\bf\Gamma}))$, which solves the Dirichlet boundary value problem:
\begin{equation}
\Delta U=0\,\, {\rm in}\,\, \Omega:=\Omega({\bf\Gamma}), U(\Gamma_1)=0, U(\Gamma_2)=1.
\label{eqn 1.2}
\end{equation}
(Given any vector ${\bm\lambda}=(\lambda_1,\lambda_2)\in\Re_+^2 $, we call ${\bf\Gamma}_{\boldsymbol\lambda}=(\Gamma_{{\boldsymbol\lambda},1},\Gamma_{{\boldsymbol\lambda},2})\in{\bf X}$ a solution of Prob. \ref{prob 1.1} at ${\bm\lambda}$ if the corresponding capacitary potential $U_{\boldsymbol\lambda}(p):=U({\bm\Gamma}_{\boldsymbol\lambda};p)$ satisfies (\ref{eqn 1.1}) with the new flow-speed-functions $a_{{\boldsymbol\lambda},i}(p)=\lambda_i\,a_i(p):\Re^2\rightarrow\Re_+$, $i=1,2$.)
\end{problem}

\noindent
In studying Prob. \ref{prob 1.1}, it is very helpful to  keep in mind the following pair of closely related single-free-boundary Bernoulli problems, which have received much more attention in the literature:
\begin{problem} {(Exterior (resp. interior) Bernoulli free-boundary Problem for annular domains in $\Re^2$)}
\label{prob 1.2}
Given a curve $\Gamma^*\in {\rm X}$ and a strictly-positive continuous function $a(p):\Re^2\rightarrow\Re_+$, we seek a smooth curve $\Gamma\in{\rm X}$ such that $\Gamma>(<)\,\Gamma^*$ and such that $|{\nabla}U(p)|=a(p)$ on $\Gamma$, where $U(p)=U(\Gamma;p)$ denotes the stream function in the annular (i.e. doubly-connected) domain $\Omega=\Omega(\Gamma)$ between $\Gamma^*$ and $\Gamma$. 
\end{problem}

\begin{remark} 
\label{rem 0.0}
(Generality of the stream model) Probs. \ref{prob 1.1} and \ref{prob 1.2} are based on very general flow-models, in the sense that the flow-speed functions $a_1(p),a_2(p):\Re^2\rightarrow\Re_+$, which indirectly represent the terrain, are completely unrestricted apart from their positivity and smoothness. Indeed, if one were to remove the stated requirement that the solutions be annular, then there would be many possible types of solutions, in fact any sufficiently-smooth, simply or multiply-connected, bounded open set $\Omega$ would be a solution of Prob. \ref{prob 1.1} for numerous pairs of strictly-positive and arbitrarily smooth flow-speed functions $a_1(p),a_2(p):\Re^2\rightarrow\Re_+$. To show this, just choose any such region $\Omega$ and partition $\partial\Omega$ into two parts, each of which is a union of smooth simple closed boundary curves. Then choose $a_i(p):=|\nabla U(p)|$ on $\Gamma_i$, $i=1,2,$ where $U$ denotes the harmonic function in $\Omega$ satisfying the boundary conditions: $U=i-1$ on $\Gamma_i$, $i=1,2$. Finally, continue the functions $a_i(p)$, $i=1,2,$ into the rest of $\Re^2$ in any convenient smooth way. Similarly, in the context of Prob. \ref{prob 1.2}, any given bounded, open, multiply-connected, smooth set $\Omega$ having $\Gamma^*$ as a boundary component, would be a solution for a suitable function $a(p):\Re^2\rightarrow\Re_+$ chosen such that $a(p)=|\nabla U(p)|$ on $\Gamma:=(\partial\Omega)\setminus \Gamma^*$, where $U$ is the harmonic function such that $U(\Gamma)=0$ and $U(\Gamma^*)=1$. 
\end{remark} 

\noindent
For the case of prob. \ref{prob 1.2}, we point out the following classical existence theorem due to Arne Beurling:
\begin{theorem} {(Beurling's Theorem \cite{AB1}(1957), \cite{AB2}(1989))} 
\label{thm1.1}
Let $\boldsymbol{\Gamma}^-\in{\rm X}$ and $\Gamma^+\in{\rm X}$ (with $\Gamma^-<\Gamma^+<\Gamma^*$) denote respective super and sub solutions of the interior Bernoulli Problem (Prob. \ref{prob 1.2}), in the sense that $|{\nabla} U^+(p)|>a(p)$ on $\Gamma^+$ and $|{\nabla} U^-(p)|<a(p)$ on $\Gamma^-$, where $U^\pm(p)=U(\Gamma^\pm,\Gamma^*;p)$ in the closure of $\Omega^\pm:=\Omega(\Gamma^\pm,\Gamma^*)$. Then there exists a classical solution $\Gamma\in{\rm X}$ such that $\Gamma^-<\Gamma<\Gamma^+<\Gamma^*$. 
\vspace{.1in}

\noindent
Similarly, if $\Gamma^-\in {\rm X}$ and $\Gamma^+\in{\rm X}$ (with $\Gamma^*<\Gamma^-<\Gamma^+)$ are inner and outer solutions of the exterior problem, in the sense that $|\nabla U^-(p)|>a(p)$ (resp. $|\nabla U^+(p)|<a(p)$ on $\Gamma^-$ (resp. $\Gamma^+$), where $U^\pm(p)=U(\Gamma^*,\Gamma^\pm;\,p)$ in the closure of $\Omega^\pm:=\Omega(\Gamma^*;\Gamma^\pm)$. Then there exists a classical solution $\Gamma\in {\rm X}$ such that $\Gamma^*<\Gamma^-<\Gamma<\Gamma^+$.
\end{theorem}

\noindent
The main objective of Chapters \ref{section 2} and \ref{section 3} of the present paper is to generalize Beurling's Theorem to Prob. \ref{prob 1.1} by proving the following result: 

\begin{theorem} 
{(Existence of a classical solution-pair between inner and outer solution-pairs (see Thms. \ref{thm2.1.0}, \ref{thm 2.6.2}, and \ref{thm 2.7.1}))}
\label{thm1.2}
Let be given the ordered curve-pairs ${\bf\tilde{\Gamma}}=(\tilde{\Gamma}_1,\tilde{\Gamma}_2)\in{\bf X}$ and ${\bf\hat{\Gamma}}=(\hat{\Gamma}_1,\hat{\Gamma}_2)\in{\bf X}$, each composed of $C^2$ curves. Assume that ${\bf\tilde{\Gamma}}$ and ${\bf\hat{\Gamma}}$ are respective strict inner and strict outer solutions (i.e. strict sub- and super-solutions) of Prob. \ref{prob 1.1}, in the sense that
\begin{equation}
\label{eqn 1.3}
|{\nabla} \tilde{U}(p)|<a_1(p)\,\, {\rm on}\,\, \tilde{\Gamma}_1\,;\,\, 
|{\nabla} \tilde{U}(p)|>a_2(p)\,\, {\rm on}\,\, \tilde{\Gamma}_2,
\end{equation}
\begin{equation}
\label{eqn 1.4}
|{\nabla} \hat{U}(p)|>a_1(p)\,\, {\rm on}\,\, \hat{\Gamma}_1\,;\,\, 
|{\nabla} \hat{U}(p)|<a_2(p)\,\, {\rm on}\,\, \hat{\Gamma}_2,
\end{equation}
where we set $\tilde{U}(p)=U({\bf\tilde{\Gamma}};p)$ and $\hat{U}(p)=U({\bf\hat{\Gamma}};p)$. If ${\bf\tilde{\Gamma}}<{\bf\hat{\Gamma}}$ (in the sense that $\tilde{\Gamma}_i<\hat{\Gamma}_i$, $i=1,2$), then there exists a curve-pair ${\bf\dot{\Gamma}}=(\dot{\Gamma}_1,\dot{\Gamma}_2)\in{\bf X}$, such that ${\bf\tilde{\Gamma}}<{\bf\dot{\Gamma}}<{\bf\hat{\Gamma}}$ and ${\bf\dot{\Gamma}}$ is a classical solution of Prob. \ref{prob 1.1}, in the sense that   
\begin{equation}
\label{eqn 1.5}
|{\nabla} \dot{U}(p)|=a_1(p)\,\, {\rm on}\,\, \dot{\Gamma}_1\,;\,\, 
|{\nabla} \dot{U}(p)|=a_2(p)\,\, {\rm on}\,\, \dot{\Gamma}_2,
\end{equation}
where we set $\dot{U}(p)=U({\bf\dot{\Gamma}};p)$. In other words, if a smooth strict inner solution-pair ${\bf\tilde{\Gamma}}$ (characterized by (\ref{eqn 1.3})) lies component-wise inside a smooth strict outer solution-pair ${\bf\hat{\Gamma}}$ (characterized by (\ref{eqn 1.4})), then there exists a classical solution-pair ${\bf\dot{\Gamma}}$ (satisfying (\ref{eqn 1.5}) classically) such that the latter lies component-wise between the inner and outer solution-pairs.
\end{theorem}
\begin{remark}
\label{rem 0.1}

(a) The author first studied the double-free-boundary Bernoulli free-boundary problem (Prob. \ref{prob 1.1}) in [A3](1978). By comparison to the present existence and uniqueness results for Prob. \ref{prob 1.1}, the corresponding results in [A3] are less general, but also easier to prove, both due to additional assumptions concerning the given flow-speed functions $a_1(p), a_2(p):\Re^2\rightarrow\Re_+$.
\vspace{.1in}

\noindent
(b) The double-free-boundary soft-terrain problem (Prob. \ref{prob 1.1}) is related to the 
double-free-boundary hard-barrier problem, which was the main topic studied by Arne Beurling in his Mittag-Leffler lecture series in 1977-78, later reconstructed from notes and published in Beurling's collected works [AB2](1989). The present author, who first learned of these lectures in 1989, presented an independent, although less general study of single and double free-boundary hard-barrier problems at Oberwohlfach in 1977 (see [A3]). Similar problems were also studied by Tepper (1984)
\vspace{.1in}

\noindent
(c) Although Thm. \ref{thm1.2} seems to be the most natural generalization of Beurling's Theorem to Prob. \ref{prob 1.1}, the author's proof of Thm. \ref{thm1.2} (in Chapters \ref{section 2} and \ref{section 3}) is not based on Beurling's methods in his proof of Beurling's theorem. Instead, we define a one-parameter family of operators ${\bf T}_\varepsilon:{\bf X}\rightarrow{\bf X}$, $\varepsilon\in(0,\varepsilon_0)$, such that (${\bf i}$) the operator fixed point problem has a solution between any suitably ordered inner and outer solutions, and (${\bf ii}$) a convergent sequence of operator fixed-points $\big({\bf\Gamma}_n\big)_{n=1}^\infty$ corresponding to a null-sequence $\big(\varepsilon_n\big)_{n=1}^\infty$ defines a weak solution of Prob. \ref{prob 1.1}. This is a geometrically very general version of the Operator Method (see Acker [A5], [A7]), which was was first developed by the present author in in 1978, and whose previous applications have been largely restricted to the convex and starlike cases, although far greater geometric generality is possible, as we will see. 
\vspace{.1in}

\noindent
(d) We point out the author's counterexample in [A3](1980), which seems to show that without assumptions not present in the stream-model under discussion, the method of variational inequalities (see Alt and Caffarelli [AC](1981)) is not ideally suited to proving the existence of annular solutions of Prob. \ref{prob 1.2}. Namely, in [A3], it is shown in the context of the interior Bernoulli problem that if one chooses the simple closed curve $\Gamma^*$ and the constant $c>0$ suitably (where we set $a(p)=c$ everywhere), then either the minimizing configuration is not annular, or else it does not satisfy the Bernoulli condition on the free boundary. On the other hand, if a simple closed curve $\Gamma$ solving Prob. \ref{prob 1.2} (in either case) is a member of a continuously-varying, elliptically-ordered family of simple closed curves $\Gamma_\lambda$, $\lambda\in I$, such that $\Gamma_\lambda$ solves Prob. \ref{prob 1.2} (in the same case) at $\lambda$ for each $\lambda\in I$, then $\Gamma$ is the (at least local) minimizer of the variational-inequalities functional (see [A1], [A2](1978)).
\vspace{.1in}
\end{remark}

\subsection{Global uniqueness}
\label{subsection 1.2}
Returning again to the "narrow-stream-limit" model in the presence of only one flow-speed function $a(p)$ (introduced at the beginning of Section \ref{subsection 1.1}), one can easily visualize that if the terrain of the flow-speed function $a(p)$ consists of a suitable system of valleys separated by ridges, then the path-integral of the function $a(p)$, defined on any path $\Gamma\in{\rm X}$, can have many local relative minimizers, each satisfying the Euler equation (\ref{eqn 1.0}). To eliminate these obvious causes of non-uniqueness of the solution of (\ref{eqn 1.0}), it is natural to restrict our attention to the case of logarithmically-subharmonic flow-speed functions $a(p)$, since these functions don't have any high-altitude ridges to separate valleys. In connection with this restriction, it is also necessary to restrict the size of the domain of the function $a(p)$, since logarithmic subharmonicity is an increasingly restrictive property for uniformly bounded functions in increasingly large domains. 
\vspace{.1in}

\noindent
Assume that the positive function $a(p)$ is strictly logarithmically subharmonic relative to an annular domain $G\subset\Re^2$. Then the solution $\Gamma\in{\rm X}$ of (\ref{eqn 1.0}) is unique relative to $G$ in the following sense: Let $\Gamma_1, \Gamma_2\in{\rm X}$ denote two solutions of (\ref{eqn 1.0}), which intersect at two successive points $p_1, p_2$, so that the arc-segment $\gamma_1$ (resp. $\gamma_2$) of $\Gamma_1$ (resp. $\Gamma_2$) initiating at $p_1$ and terminating at $p_2$ forms the left (right) boundary of a bounded, simply-connected region $\omega\subset G$. Then for continuous argument functions $\theta_i(p):\Gamma_i\rightarrow\Re$, $i=1,2$, we have $\theta_2(p_2)-\theta_2(p_1)\geq \theta_1(p_2)-\theta_1(p_1)$, from which it follows by (\ref{eqn 1.0}), the divergence theorem, and the strict subharmonicity of the function ${\rm ln}\big(a(p)\big)$ in $G$ that  
\begin{equation}
\label{eqn 2.1.1b}
0\leq\int_{\gamma_2} K(p)\,ds-\int_{\gamma_1}K(p)\,ds=\int_{\gamma_2}\big({\rm ln}(a(p))\big)_{\boldsymbol{\nu}}\,ds-\int_{\gamma_1}\big({\rm ln}(a(p))\big)_{\boldsymbol{\nu}}\,ds
\end{equation}
$$=\int\int_{\omega}\Delta\,{\rm ln}\big(a(p)\big)\,dxdy\,<0.$$
For the case of the thin-stream limit, the asserted uniqueness follows from this contradiction, as well as a similar contradiction which occurs in the case where the solutions $\Gamma_1, \Gamma_2$ are disjoint.
\vspace{.1in}

\noindent
In Chapter \ref{section 4}, we explore the uniqueness of solutions of the finite-width stream problem (Prob. \ref{prob 1.1}) in the case of a single strictly logarithmically subharmonic flow-speed function $a(p):G\rightarrow\Re_+$ in an annular domain $G$ such that $a_1(p)=a(p)=a_2(p)$ in $G$. We study the uniqueness question from two perspectives. First (in Thms. \ref{thm 4I.bib} and \ref{thm 4B.bjb}), we show that if ${\bf\Gamma}_1=(\Gamma_{1,1}, \Gamma_{1,2})\in{\bf X}$ and ${\bf\Gamma}_2=(\Gamma_{2,1},\Gamma_{2,2})\in{\bf X}$ denote any two distinct classical solutions with corresponding stream beds $\Omega_1:=\Omega({\bf\Gamma}_1)$ and $\Omega_2:=\Omega({\bf\Gamma}_2)$, then the stream beds are similar enough so that there exists a simple closed curve which follows both channels while encircling the interior complement of $G$. Using this, we show that one of the streams follows the channel of the other. Finally, we show that this leads to a contradiction. Alternatively, under about the same assumptions, we apply Thm. \ref{thm1.2} to construct continuously and monotonically-varying parametrized solution families for Prob. \ref{prob 1.1} (see Def. \ref{def 4.1.1} and Thm. \ref{thm 4.1.2}). As we will show, these monotone solution families can be made the basis of a uniqueness proof (see Thm. \ref{thm 4.1.1} and Cor. \ref{cor 4.1.2}).
\vspace{.1in}

\noindent
These results are related to the author's uniqueness studies in [A6],[A7](1989) for (mostly) the interior Bernoulli problem (see Thms. 2.1 and 3.1 in [A7]). Other uniqueness studies for the interior Bernoulli Problem are due to Liu [YL](1995) and to Cardaliaguet and Tahraoui, [CT](2002). The primary uniqueness result for the exterior Bernoulli Problem in the starlike case is known as the Lavrentiev Principle [LS](1967), pp. 413-434.

\section{Operator fixed-points and weak solutions}
\label{section 2} 
\subsection{Existence results for fixed-point problems}
\label{subsection 2.1} 
\begin{definition}
\label{def 2.1.1} 
{(Preliminary definitions and notation)} The results of this paper are in $\Re^2$. We denote the typical point by $p=(x,y)$. We use $B_r(p)$ and $B(r;p)$ to denote the open ball of radius $r>0$ and center-point $p\in\Re^2$, while $\overline{B}_r(p)=\overline{B}(p;r)$ denotes their closures. Given a value $P>0$, we call a set $M$ (or a function $U(x,y):M\rightarrow\Re$) $P$-periodic in $x$ (or just $P$-periodic) if $(x,y)\in M\iff(x+P,y)\in M$ (resp. $U(x,y)=U(x+P,y)$ in $M$). Let ${\rm X}$ denote the family of all infinite, $P$-periodic (in $x$) directed arcs $\Gamma$, with positive direction from $x=-\infty$ to $x=\infty$, such that the winding number of $\Gamma$ about any point $p\notin\Gamma$ is $\pm 1/2$. For any arc $\Gamma\in{\rm X}$, we use $D_1(\Gamma)$ (resp. $D_2(\Gamma)$) to denote the lower (resp. upper) complement of $\Gamma$ in $\Re^2$, defined to be the set of points $p\in\Re^2\setminus\Gamma$ such that $W(\Gamma;p)=-1/2$ (resp. $W(\Gamma;p)=1/2$), where $W(\Gamma;p)$ denotes the winding number of $\Gamma$ about $p$. Then we write $\Gamma_1\leq\Gamma_2$ in ${\rm X}$ if $D_1(\Gamma_1)\subset D_1(\Gamma_2)$ or if $D_2(\Gamma_2)\subset D_2(\Gamma_1)$, and we write $\Gamma_1<\Gamma_2$ in ${\rm X}$ if ${\rm Cl}\big(D_1(\Gamma_1))\subset D_1(\Gamma_2)$ or ${\rm Cl}\big(D_2(\Gamma_2))\subset D_2(\Gamma_1)$. We use ${\bf X}$ to denote the family of all ordered pairs ${\bf \Gamma}=(\Gamma_1,\Gamma_2)\in{\rm X}\times{\rm X}$ such that $\Gamma_1<\Gamma_2$. For ${\bf \Gamma}_1,{\bf\Gamma}_2\in {\bf X}$, the inequalities ${\bf\Gamma}_1\leq{\bf\Gamma}_2$ and ${\bf\Gamma}_1<{\bf\Gamma}_2$ are defined component-wise. For ${\bf \Gamma}=(\Gamma_1,\Gamma_2)$ in ${\bf X}$, we define the $P$-periodic (in $x$) strip-like region $\Omega({\bf \Gamma}):=D_1(\Gamma_2)\cap D_2(\Gamma_1)$, and we define the capacitary potential $U({\bf \Gamma};p):{\rm Cl}(\Omega({\bf \Gamma}))\rightarrow\Re$ to be the continuous $P$-periodic (in $x$) function in ${\rm Cl}(\Omega({\bf\Gamma}))$ which is harmonic in $\Omega({\bf \Gamma})$ and satisfies the boundary conditions: $U({\bf \Gamma};p)=0$ on $\Gamma_1$ and $U({\bf\Gamma};p)=1$ for $\Gamma_2$. (For convenience, we also define $U_i({\bm\Gamma};p)$, $i=1,2$, such that $U_1({\bm\Gamma};p)=U({\bm\Gamma};p)$ and $U_2({\bm\Gamma};p):=1-U({\bm\Gamma};p)$.) For any finite arc $\gamma$ or region $\omega$, we use $||\gamma||$ and $K(\gamma)$ to denote, respectively, the Euclidean arc-length and total curvature of $\gamma$ or, while $||\omega||$ denotes the  Euclidean area of $\omega$. For $P$-periodic (in $x$) arcs and regions, the same notation refers to the length, total curvature, or area of the restriction of $\gamma$ or $\omega$ to one $P$-period. 
\vspace{.1in}

\noindent
Given the constants $0<\underline{A}\leq \overline{A}$ and $A_1, A_2>0$, we use $\boldsymbol{{\cal A}}$ to denote the class of all strictly-positive, twice-continuously-differentiable functions $a(p):\Re^2\rightarrow\Re_+$ such that $\underline{A}\leq a(p)\leq \overline{A}$ and the absolute values of all first (resp. second) order partial derivatives of $a(p)$ are uniformly bounded above by $A_1$ (resp. $A_2$). Given $\varrho\in(0,1]$, let $\boldsymbol{\cal{A}}_{\varrho}$ denote the family of all functions in $\boldsymbol{\cal{A}}$ which are also in $C^{\,3,\varrho}(\Re^2)$. (We remark that $A_1\leq \sqrt{2\delta A_2}$, where $\delta:=\overline{A}-\underline{A}$.) We make frequent use of the  constant $\varepsilon_0:=\min\{1/2, \underline{A}^2/2A_1\}$.
\end{definition}
\begin{problem}
\label{prob2.1.1}
{(Double free-boundary Bernoulli problem with two independent flow-speed functions; the strip-like, $P$-periodic case)} 
Given the strictly-positive, $P$-periodic (in $x$), $C^2$-functions $a_1(p),$ $a_2(p):\Re^2\rightarrow\Re_+$, we seek a pair ${\bf \Gamma}=(\Gamma_1,\Gamma_2)\in{\bf X}$ such that $|{\nabla} U({\bf\Gamma};p)|=a_1(p)$ on $\Gamma_1$ and $|{\nabla} U({\bf\Gamma};p)|=a_2(p)$ on $\Gamma_2$.
\end{problem}
\begin{definition}
\label{def 2.1.2} 
{(Lower and upper solutions of the double free-boundary Ber-noulli problem)} 
A pair of directed $C^{2}$-arcs ${\bf \Gamma}=(\Gamma_1,\Gamma_2)\in{\bf X}$ is called a lower (resp. upper) solution of Prob. \ref{prob2.1.1} if $|{\nabla} U({\bf\Gamma};p)|\leq(\geq)\,a_1(p)$ on $\Gamma_1$ and $|{\nabla} U({\bf\Gamma};p)|\geq$ $(\leq)\,a_2(p)$ on $\Gamma_2$. It is a strictly lower (resp. strictly upper) solution of Prob. \ref{prob2.1.1} if $|{\nabla} U({\bf\Gamma};p)|<(>)\,a_1(p)$ on $\Gamma_1$ and $|{\nabla} U({\bf\Gamma};p)|>(<)\,a_2(p)$ on $\Gamma_2$. 
\end{definition}
\begin{theorem} 
{(Existence of solutions of Prob. \ref{prob2.1.1})}
\label{thm2.1.0} In the context of Prob. \ref{prob2.1.1}, given the strict upper and lower solutions, ${\bf\Gamma}^+$ and ${\bf\Gamma}^-$, both in ${\bf X}\cap C^{2}$, and such that ${\bf\Gamma}^-<{\bf\Gamma}^+$, there exists a solution ${\bf\Gamma}\in{\bf X}\cap C^{1,1}$ such that ${\bf\Gamma}^-<{\bf\Gamma}<{\bf\Gamma}^+$.
\end{theorem}
\begin{definition} 
{(Solution at {\bm\lambda}, upper and lower solutions at {\bm\lambda})}
\label{def 2.1.3}
For any pair ${\bm \lambda}=(\lambda_1,\lambda_2)\in\Re_+^2$, a curve-pair ${\bf\Gamma}=(\Gamma_1,\Gamma_2)\in{\bf X}\cap C^{2}$ is called a solution (resp. lower solution, upper solution) of Prob. \ref{prob2.1.1} at ${\bm\lambda}$ if ${\bf\Gamma}$ is a solution (lower solution, upper solution) of Prob. \ref{prob2.1.1} with the positive $C^{2}$-functions $a_1(p),a_2(p):\Re^2\rightarrow\Re_+$ replaced by the functions $\lambda_1a_1(p),\lambda_2a_2(p):\Re^2\rightarrow\Re_+$. 
\end{definition}
\begin{remark} 
{(Equivalence of the annular case to the periodic case)}
\label{rem 2.1.1}
Prob. \ref{prob2.1.1} is equivalent to Prob. \ref{prob 1.1} in the following sense: Let the simple closed curves $\hat{\Gamma}_1$ and $\hat{\Gamma}_2$ (expressed in polar coordinates $r,\theta$ relative to an origin encircled by $\hat{\Gamma}_1$) be equivalent to the $2\pi$-periodic (in $x$) arcs $\Gamma_1$ and $\Gamma_2$ (expressed in Cartesian coordinates $(x,y)$ under the mapping $r={\rm exp}(-y)$, $\theta=x$. Similarly, let the polar-coordinate functions $\hat{a}_1(r,\theta),\hat{a}_2(r,\theta)$ correspond to the $2\pi$-periodic Cartesian-coordinate functions $a_1(x,y),a_2(x,y)$ such that $r\,\hat{a}_1(r,\theta)=a_2(x,y)$ and $r\,\hat{a}_2(r,\theta)=a_1(x,y)$, where $r={\rm exp}(-y)$ and $\theta=x$. Then the curve-pair $(\hat{\Gamma}_1,\hat{\Gamma}_2)$ solves Prob. \ref{prob 1.1} relative to the function-pair $(\hat{a}_1,\hat{a}_2)$ if and only if the curve-pair $(\Gamma_1,\Gamma_2)$ solves Prob. \ref{prob2.1.1} relative to the function-pair $(a_1,a_2)$. The proof follows essentially from the fact that if $\hat{U}(r,\theta)$ and $U(x,y)$ are corresponding capacitary potentials in the domains $\hat{\Omega}$ and $\Omega$ such that $\partial\hat{\Omega}=\hat{\Gamma}_1\cup\hat{\Gamma}_2$ and $\partial\Omega=\Gamma_1\cup\Gamma_2$, then $\hat{U}({\rm exp}(jz))=U(z)$ in $\Omega$ for $j=\sqrt{-1}$, whence $r|{\nabla} U(r,\theta)|=|{\nabla} U(x,y)|$ by differentiation.      
\end{remark}
\begin{remark} (a) {(The one-dimensional model of Prob. \ref{prob2.1.1})}
\label{rem 2.1.2}
Some insight into the general existence and uniqueness of solutions of Prob. \ref{prob2.1.1} can be gained by studying the simpler one-dimensional case, corresponding to the cross-section of a straight 2-dimensional flow. Here, given positive, smooth real-valued functions $a_1(x),a_2(x):\Re\rightarrow\Re_+$ (with reciprocals $b_i(x)=(1/a_i(x))$, $i=1,2$), one seeks a pair $(x_1,x_2)$ with $x_1<x_2$ such that $a_1(x_1)=(1/(x_2-x_1))=a_2(x_2)$, or, more conveniently, such that $({\bf i})$: $\,b_1(x_1)=(x_2-x_1)=b_2(x_2)$. The pair $(x_1,x_2)$ (with $x_1<x_2$) is a strict lower (resp. strict upper) solution of Prob. $({\bf i})$ if and only if $({\bf ii})$: $b_1(x_1)<(>)(x_2-x_1)<(>)b_2(x_2)$. Clearly, there is no solution of $({\bf i})$ in the case where the functions $a_1(x)$ and $a_2(x)$ are distinct constants (for which there are also no weak upper and lower solutions). Every pair $(x_1,x_1+(1/C))$, $x_1\in\Re$, is a solution in the case where $a_1(x)=C=a_2(x)$ for a constant $C$. If we assume that $a_1(x),a_2(x)\geq\underline{A}>0$ and $|a_1'(x)|,|a_2'(x)|<A_1$ for all $x\in\Re$, then $(x_1,x_2)$ cannot be a solution if $|a_2(x)-a_1(x)|<2(A_1/\underline{A})$ for any point $x\in[x_1,x_2]$. The natural place for a solution is near a point where $(a_2(x)-a_1(x))$ changes sign. 
\vspace{.1in}

\noindent
(b) {\it (One-dimensional model in the case where $a_1(x)=a(x)=a_2(x)$)}
In the interesting special case where $a_1$ and $a_2$ reduce to the same function $a(x)$, there are no solutions in any interval in which $a(x)$ is strictly monotone. The natural location to seek solutions is near local extrema of the function $a(x)$. Consider the case where the function $a(x)$ is strictly decreasing (resp. increasing) to the left (right) of a single minimum at $x_0\in\Re$ and $a(x)\rightarrow+\infty$ as $x\rightarrow\pm\infty$. For any value $\alpha>a(x_0)$, we have $a(x_1)=a(x_2)=\alpha$ for unique values $x_1=x_1(\alpha)<x_0$ and $x_2=x_2(\alpha)>x_0$. Also, the function $\phi(\alpha):=
|x_2(\alpha)-x_1(\alpha)|$ increases continuously from $0$ to $+\infty$ as $\alpha$ increases from $a(x_0)$ to $\infty$. Therefore there exists a unique root $\alpha=\alpha_0$ of the equation $\phi(\alpha)=(1/\alpha)$. The pair $(x_1(\alpha_0), x_2(\alpha_0))$ uniquely solves the problem, as is clear from the construction. In the alternate case where the positive function $a(x)$ is strictly increasing (decreasing) to the left (right) of a single maximum at a point $x_0$ and $a(x)\rightarrow 0$ as $x\rightarrow\pm\infty$, there are again unique points $x_1(\alpha)<x_0<x_2(\alpha)$ such that $a(x_1(\alpha))=\alpha=a(x_2(\alpha))$ for any $\alpha\in(0,a(x_0))$, but the function $\phi(\alpha):=|x_2(\alpha)-x_1(\alpha)|$ decreases continuously from $+\infty$ to $0$ as $\alpha$ increases from $0$ to $a(x_0)$. The solutions are again all of the form $(x_1(\alpha_0),x_2(\alpha_0))$ corresponding to roots of the equation $\phi(\alpha)=(1/\alpha)$, but there may be no root or many roots, even infinitely many, depending on the exact details of the given function $a(x)$. Whatever solution-pairs  exist can all be strictly ordered by inclusion of the corresponding intervals (which all contain $x_0$), and therefore cannot be ordered by components.
\vspace{.1in}

\noindent
(c) {\it (Proof of Thm. \ref{thm2.1.0} in the one-dimensional model)}
Let ${\it X}$ denote the set of all pairs of real numbers $(x_1,x_2)$ such that $x_1<x_2$. In ${\it X}$, we define inequality notation (for example $(x_1,x_2)<(\leq)(y_1,y_2)$) componentwise. We choose $\eta_0=\min\{1/2,1/2B_1\}$, where $b_i(p)\geq\underline{B}>0$ and $|b_i'(x)|\leq B_1$ for $i=1,2$ and all $x\in\Re$. For any $\varepsilon\in(0,\eta_0)$, we define the continuous operator $T_\varepsilon: {\it X}\rightarrow {\it X}$ such that for any pair $(x_1,x_2)\in {\it X}$, we have
$T_{\varepsilon}(x_1,x_2)=(x_{\varepsilon,1}^*,x_{\varepsilon,2}^*)\in {\it X}$, where the values $x_{\varepsilon,1}^*\in(-\infty,
(1-\varepsilon)x_1+\varepsilon\,x_2)$ and $x_{\varepsilon,2}^*\in(\varepsilon\,x_1+(1-\varepsilon)\,x_2,\infty)$ are uniquely determined by the requirements that $({\bf iii})$: $x_{\varepsilon,1}^*+\varepsilon\,b_1(x_{\varepsilon,1}^*)=(1-\varepsilon)x_1+\varepsilon x_2$ and $x_{\varepsilon,2}^*-\varepsilon\,b_2(x_{\varepsilon,2}^*)=\varepsilon x_1+(1-\varepsilon)x_2$ (from which it follows that $x_{\varepsilon,2}^*-x_{\varepsilon,1}^*=(1-2\varepsilon)(x_2-x_1)+\varepsilon\,(b_2(x_{\varepsilon,2}^*)+b_1(x_{\varepsilon,1}^*))\geq 2\underline{B}\varepsilon$). A pair $(x_1,x_2)\in {\it X}$ is called an operator "fixed point" at $\varepsilon$ if $({\bf iv})$: $T_\varepsilon(x_1,x_2)=(x_1,x_2)$, or, equivalently: $b_1(x_{\varepsilon,1})=(x_2-x_1)=b_2(x_{\varepsilon,2})$, and it is called 
a weak lower (resp. upper) solution of the fixed-point problem $({\bf iv})$ at $\varepsilon$ if $({\bf v})$: $T_\varepsilon(x_1,x_2)\geq (\leq)(x_1,x_2)$, or, equivalently, $b_1(x_{\varepsilon,1}^*)\leq (\geq)\,(x_2-x_1)\leq(\geq)\,b_2(x_{\varepsilon,2}^*)$. 
For $(x_1,x_2), (y_1,y_2)\in {\it X}$, and $\varepsilon\in(0,\eta_0)$, it follows from the definitions $({\bf iii})$ that $({\bf vi})$: if $(x_1,x_2)\leq(y_1,y_2)$, then $T_\varepsilon(x_1,x_2)\leq T_\varepsilon(y_1,y_2)$.
Also, it follows from $({\bf i})$, $({\bf iv})$, and $({\bf v})$ that $({\bf vi})$: if, given any positive null-sequence $\big(\varepsilon_k\big)$, a corresponding sequence of fixed points (of $T_\varepsilon$ at $\varepsilon=\varepsilon_k$) has a convergent subsequence, then the limit-pair solves $({\bf i})$, and, finally $({\bf vii})$: any strict lower (resp. upper) solution of $({\bf i})$ is also a strict lower (upper) solution of $({\bf iv})$ at $\varepsilon$ if $\varepsilon\in(0,\eta_0)$ is sufficiently small. 
\vspace{.1in}

\noindent 
In the context of Prob. $({\bf i})$, let be given a strict lower solution $(x_1^-,x_2^-)$ and a strict upper solution $(x_1^+,x_2^+)$, both in ${\it X}$, such that $(x_1^-,x_2^-)<(x_1^+,x_2^+)$. We then slightly increase (resp. decrease) the components of $(x_1^-,x_2^-)$ (resp. $(x_1^+,x_2^+)$) while preserving these properties. It follows that $({\bf viii})$:$(x_1^-,x_2^-)<T_\varepsilon(x_1^-,x_2^-)\leq\cdots\leq T_\varepsilon^n(x_1^-,x_2^-)\leq T_\varepsilon^n(x_1^+,x_2^+)\leq\cdots\leq T_\varepsilon(x_1^+,x_2^+)<(x_1^+,x_2^+)$ for any $n\in N$ and any sufficiently small value $\varepsilon\in(0,\eta_0)$, where the first and last inequalities in $({\bf viii})$ both follow from $({\bf vii})$, and then all the remaining inequalities follow by multiple applications of $({\bf vi})$. By continuity, the pairs $(\tilde{x}_{\varepsilon,1}^\pm,\tilde{x}_{\varepsilon,2}^\pm):=\lim_{n\rightarrow\infty}\,T_\varepsilon^n(x_1^\pm,x_2^\pm)$ (which may be identical) are both "fixed points" of the operator $T_\varepsilon$ such that $({\bf ix})$: $(x_1^-,x_2^-)<(\tilde{x}_{\varepsilon,1}^-,\tilde{x}_{\varepsilon,2}^-)\leq(\tilde{x}_{\varepsilon,1}^+,\tilde{x}_{\varepsilon,2}^+)<(x_1^+,x_2^+)$. It follows by $({\bf vi})$ and a compactness argument that there exists a solution $(\tilde{x}_1,\tilde{x}_2)$ of Prob. $({\bf i})$ such that $(x_1^-,x_2^-)\leq(\tilde{x}_{1},\tilde{x}_{2})\leq(x_1^+,x_2^+)$. At this point, we complete the proof by returning to the original components of the pairs $(x_1^\pm,x_2^\pm)$.
\end{remark}
\begin{definition}
{($\boldsymbol{{\cal A}}$-operator definitions)} (See Lemmas \ref{lem 2.1.4} and \ref{lem 2.1.4*}.)
\label{def 2.1.4}
(a) Given the functions $a_1(p), a_2(p)\in\boldsymbol{{\cal A}}$, we define their reciprocal functions $b_1(p),$ $b_2(p):\Re^2\rightarrow\Re_+$, such that $b_i(p):=\big(1\big/a_i(p)\big)$, and the related functions $\phi_i(p):{\rm Cl}\big(D_i(\Gamma)\big)\rightarrow[0,\infty)$, $i=1,2$, such that $\phi_i(p)=a_i(p)\,{\rm dist}\big(p,\Gamma\big)$.
\vspace{.1in}

\noindent
(b) Let $\varepsilon_0:=\min\{1/2, \underline{A}^2/2A_1\}$.
For any $\varepsilon\in(0,\varepsilon_0)$, $\Gamma\in{\rm X}$, and $i=1,2,$ we define $\dot{G}_{\varepsilon,i}(\Gamma)$ to be the (open) set of all points $p\in D_i(\Gamma)$ such that $\phi_i(p)>\varepsilon$, and we define $G_{\varepsilon,i}(\Gamma)$ to be the (open) set of all points $p\in\dot{G}_{\varepsilon,i}(\Gamma)$ such that $p$ is joined to $\{y=(-1)^i\,\infty\}$ by an arc $\gamma\subset\dot{G}_{\varepsilon,i}(\Gamma)$. Finally, we define the open set $\hat{G}_{\varepsilon,i}(\Gamma)$ to be the union of the collection of all the open balls $B\big(p\,;\,\varepsilon\,b_i(p)\big)$ whose center-points $p$ are elements of the set $G_{\varepsilon,i}(\Gamma)$.
\vspace{.1in}

\noindent
(c) Equivalently, we define $\hat{G}_{\varepsilon,i}(\Gamma)$ to be the union of the set of all closed balls $\overline{B}\big(p\,;\varepsilon\,b_i(p)\big)\subset D_i(\Gamma)$, whose center-points $p$ are joined to $\{y=(-1)^i\,\infty\}$ by directed arcs $\gamma\subset D_i(\Gamma)$ such that $\overline{B}\big(q\,;\varepsilon\,b_i(q)\big)\subset D_i(\Gamma)$ for all points $q\in\gamma$ (including the initial point $p\in\gamma$), and we then define $G_{\varepsilon,i}(\Gamma)$ to be the set of the center-points of all the closed balls $\overline{B}\big(p\,;\varepsilon\,b_i(p)\big)\subset \hat{G}_{\varepsilon,i}(\Gamma)$. Alternatively, one can define $G_{\varepsilon,i}(\Gamma):=\big\{p\in\hat{G}_{\varepsilon,i}(\Gamma):\phi_i(p)>\varepsilon\big\}.$ We remark that $G_{\varepsilon,i}(\Gamma)$ can be characterized as the set of all points $p\in G_{\varepsilon,i}(\Gamma)$ which are joined to $\{y=(-1)^i\,\infty\}$ by arcs $\gamma\subset G_{\varepsilon,i}(\Gamma)$. 
\vspace{.1in}

\noindent
(d) Similarly, we define the closed set $\dot{H}_{\varepsilon,i}(\Gamma):=\big\{p\in {\rm Cl}\big(D_i(\Gamma)\big):\,\phi_i(p)\geq\varepsilon\big\}$, and we use ${H}_{\varepsilon,i}(\Gamma)$ to denote the closed set of all points in the set $\dot{H}_{\varepsilon,i}(\Gamma)$ which are joined to $\{y=(-1)^i\,\infty\}$ by directed arcs $\gamma$ lying entirely within the same set $\dot{H}_{\varepsilon,i}(\Gamma)$.
The set $H_{\varepsilon,i}(\Gamma)$ can be characterized as the set of all points $p\in H_{\varepsilon,i}(\Gamma)$ such that $p$ is joined to $\{y=(-1)^i\,\infty\}$ by an arc $\gamma\subset H_{\varepsilon,i}(\Gamma)$. We also use $\hat{H}_{\varepsilon,i}(\Gamma)$ to denote the closure of the union of all the balls $B\big(p\,;\varepsilon\,b_i(p)\big)$ whose center-points $p$ are elements of the set ${H}_{\varepsilon,i}(\Gamma)$. Equivalently, we define $\hat{H}_{\varepsilon,i}(\Gamma)$ to be the union of the set of all closed balls $\overline{B}\big(p\,;\varepsilon\,b_i(p)\big)\subset D_i(\Gamma)$, whose center-points $p$ can be joined to $\{y=(-1)^i\,\infty\}$ by directed arcs $\gamma\subset D_i(\Gamma)$ such that $\overline{B}\big(q\,;\varepsilon\,b_i(q)\big)\subset D_i(\Gamma)$ for all points $q\in\gamma$ (including the initial point $p\in\gamma$), and we then define $G_{\varepsilon,i}(\Gamma)$ to be the set of the center-points of all the closed balls $\overline{B}\big(p\,;\varepsilon\,b_i(p)\big)\subset \hat{G}_{\varepsilon,i}(\Gamma)$. 
\end{definition}
\begin{definition}
\label{def 2.1.5}
{(Operator definitions for Prob. \ref{prob2.1.1})} 
We define the capacitary-potential operators 
\begin{equation}
\label{eqn 2.1.1}
{\bm\Phi}_\varepsilon({\bm\Gamma})=\big(\Phi_{\varepsilon,1}({\bm\Gamma}),\Phi_{\varepsilon,2}({\bm\Gamma})\big):{\bf X}\rightarrow{\tilde{\bf X}}_\varepsilon,\,\,\,\varepsilon\in(0,1/2),
\end{equation}
(where $\tilde{\bf X}_\varepsilon:={\bf \Phi}_\varepsilon({\bf X})=$ the range of ${\bf \Phi}_\varepsilon$) component-wise such that 
\begin{equation}
\label{eqn 2.1.2}
\Phi_{\varepsilon,i}({\bf \Gamma})=\big\{U_i({\bf \Gamma};p)=\varepsilon\big\},\,\,\, i=1,2.
\end{equation}
\noindent
In other words, ${\bf \Phi}_\varepsilon({\bf\Gamma})$ denotes the level curve at altitude $\varepsilon$ of one of the capacit-ary-potential functions $U_i({\bf\Gamma};p)$ associated with the strip-like domain $\Omega({\bf\Gamma})$ between the components of ${\bf\Gamma}$. 
\vspace{.1in}

\noindent
In the context of Def. \ref{def 2.1.4}, we define the parametrized families of $\boldsymbol{{\cal A}}$-operators:

\begin{equation}
\label{eqn 2.1.3}
{\bf\Psi}_{\varepsilon}^\pm({\bm\Gamma})=(\Psi_{\varepsilon,1}^\pm({\bm\Gamma}),\Psi_{\varepsilon,2}^\pm({\bm\Gamma                                                                                                                                                                              })):{\tilde{\bf X}}_\varepsilon\rightarrow{\bf X}, \,\,\,\varepsilon\in(0,\varepsilon_0),
\end{equation}
such that $\Psi_{\varepsilon,i}^\pm({\bf\Gamma}):=\Psi_{\varepsilon,i}^\pm({\Gamma}_i)$ for $i=1,2$, where 
\begin{equation}
\label{eqn 2.1.4}
\Psi_{\varepsilon,1}^+(\Gamma_1):=\partial {H}_{\varepsilon,1}(\Gamma_1);\,\,\Psi_{\varepsilon,2}^+(\Gamma_2):=\partial {G}_{\varepsilon,2}(\Gamma_2),
\end{equation}
\begin{equation}
\label{eqn 2.1.5}
\Psi_{\varepsilon,1}^-(\Gamma_1):=\partial {G}_{\varepsilon,1}(\Gamma_1);\,\,\Psi_{\varepsilon,2}^-(\Gamma_2):=\partial {H}_{\varepsilon,2}(\Gamma_2),
\end{equation}
both for all ${\bf\Gamma}=(\Gamma_1, \Gamma_2)\in \tilde{\bf X}_{\varepsilon}$ (or both for all $\Gamma_1, \Gamma_2\in\tilde{\rm X}_\varepsilon$, where $\tilde{\rm X}_\varepsilon$ denotes the set of all first components and all second components of pairs in $\tilde{\bf X}_\varepsilon$). Here, for each $\varepsilon\in(0,\varepsilon_0)$, ${\bf\Gamma}=(\Gamma_1,\Gamma_2)\in\tilde{{\bf X}}_\varepsilon$, and $i=1,2$, the above sets $\partial G_{\varepsilon,i}(\Gamma_i)$ and $\partial H_{\varepsilon,i}(\Gamma_i)$ are interpreted (by Lems. \ref{lem 2.1.4} and \ref{lem 2.1.4*}) to be $P$-periodic (in $x$) directed arcs in ${\rm X}$ such that $\partial G_{\varepsilon,i}(\Gamma_i)$ are double-point free, while the arcs $\partial H_{\varepsilon,i}(\Gamma_i)$ do not cross themselves. 
\vspace{.1in}

\noindent 
Finally, in terms of Eqs. (\ref{eqn 2.1.1})-(\ref{eqn 2.1.5}), we define the family of composite operators
\begin{equation}
\label{eqn 2.1.6}
{\bf T}_\varepsilon^{\pm}({\bf \Gamma})=\big(T_{\varepsilon,1}^{\pm}({\bf \Gamma}),T_{\varepsilon,2}^{\pm}({\bf \Gamma})\big):{\bf X}\rightarrow {\bf X},\,\, 0<\varepsilon<\varepsilon_0,
\end{equation}
such that 
\begin{equation}
\label{eqn 2.1.7}
{\bf T}_\varepsilon^\pm({\bm\Gamma})={\bm \Psi}_{\varepsilon}^\pm\big({\bm\Phi}_\varepsilon({\bm\Gamma})\big)\,\,{\rm for}\,\,{\rm all}\,\,{\bm\Gamma}\in{\bf X}.
\end{equation}

\noindent
We remark that for any $\varepsilon\in(0,\varepsilon_0)$ and ${\bm\Gamma}\in{\bf X}$, we have ${\bf T}_{\varepsilon}({\bm \Gamma})\in{\bf X}$ by Lem. \ref{lem 2.1.4*}(c).
\end{definition}
\begin{problem}
\label{prob 2.1.2} 
{(The operator "fixed point" problem)}
Given $\varepsilon\in(0,\varepsilon_0)$, the positive $C^2$-functions $a_1(p),a_2(p):\Re^2\rightarrow\Re_+$, and the operators ${\bf T}_\varepsilon^\pm:{\bf X}\rightarrow{\bf X}$ defined in Def. \ref{def 2.1.5}, we seek pairs ${\bf \Gamma}_{\varepsilon}^\pm=(\Gamma_{\varepsilon,1}^{\pm},\Gamma_{\varepsilon,2}^{\pm})\in{\bf X}$ (called "fixed points" at $\varepsilon$) such that
\begin{equation}
\label{eqn 2.1.8}
{\bf T}_\varepsilon^\pm({\bf \Gamma}_\varepsilon^\pm)={\bf \Gamma}_\varepsilon^\pm.
\end{equation}
\end{problem}
\begin{definition}
\label{def 2.1.6}
{(Family of operator "fixed points")} Given the positive $C^2$-functions $a_1(p),a_2(p):\Re\rightarrow\Re_+$, for any $\varepsilon\in(0,\varepsilon_0)$, we use ${\bm{\cal F}}_{\varepsilon}$ to denote the family of all "fixed points" ${\bf \Gamma}_\varepsilon^\pm\in{\bf X}$ of the operators ${\bf T}_\varepsilon^\pm$ (i.e. solutions of Prob. \ref{prob 2.1.2}). 
\end{definition}
\begin{theorem} {(Characterization of operator fixed points)}
\label{thm 2.1.1}
\noindent
In view of the definitions of the operators ${\bf T}_\varepsilon^\pm:{\bf X}\rightarrow{\bf X}$, $\varepsilon\in(0,\varepsilon_0)$ (in terms of Def. \ref{def 2.1.4} and Def. \ref{def 2.1.5}, Eqs. (\ref{eqn 2.1.2})-(\ref{eqn 2.1.5})), any pair of "fixed points" (solution pairs) ${\bm\Gamma}_\varepsilon^\pm=(\Gamma_{\varepsilon,1}^\pm,\Gamma_{\varepsilon,2}^\pm)\in{\bm{\cal F}}_\varepsilon$ (which exist by Thms. \ref{thm 2.1.2} and \ref{thm 2.1.3}), must satisfy the equations: 
\begin{equation}
\label{eqn 2.1.10}
\Gamma_{\varepsilon,1}^+=\partial{H}_{\varepsilon,1}(\tilde{\Gamma}_{\varepsilon,1}^+)\subset\big\{p\in D_1\big(\tilde{\Gamma}_{\varepsilon,1}^+\big):a_1(p)\,{\rm dist}\big(p,\tilde{\Gamma}_{\varepsilon,1}^+\big)=\varepsilon\big\},
\end{equation}
\begin{equation}
\label{eqn 2.1.10a}
\Gamma_{\varepsilon,2}^+=\partial{G}_{\varepsilon,2}(\tilde{\Gamma}_{\varepsilon,2}^+)\subset\big\{p\in D_2\big(\tilde{\Gamma}_{\varepsilon,2}^+\big):a_2(p)\,{\rm dist}\,\big(p,\tilde{\Gamma}_{\varepsilon,2}^+\big)=\varepsilon\big\},
\end{equation}
\begin{equation}
\label{eqn 2.1.10b}
\Gamma_{\varepsilon,1}^-=\partial{G}_{\varepsilon,1}(\tilde{\Gamma}_{\varepsilon,1}^-)\subset\big\{p\in D_1\big(\tilde{\Gamma}_{\varepsilon,1}^-\big):a_1(p)\,{\rm dist}\big(p,\tilde{\Gamma}_{\varepsilon,1}^-\big)=\varepsilon\big\},
\end{equation}
\begin{equation}
\label{eqn 2.1.10c}
\Gamma_{\varepsilon,2}^-=\partial{H}_{\varepsilon,2}(\tilde{\Gamma}_{\varepsilon,2}^-)\subset\big\{p\in D_2\big(\tilde{\Gamma}_{\varepsilon,2}^-\big):a_2(p)\,{\rm dist}\,\big(p,\tilde{\Gamma}_{\varepsilon,2}^-\big)=\varepsilon\big\},
\end{equation} 
where $\tilde{\Gamma}_{\varepsilon,i}^\pm:=\Phi_{\varepsilon,i}^\pm({\bm\Gamma}_{\varepsilon,i}^\pm)=\{U_{\varepsilon,i}^\pm(p)=\varepsilon\}$, and  where we define $U_{\varepsilon,i}^\pm(p):=U_i({\bm \Gamma}_{\varepsilon}^\pm; p)$ in the closure of the domain $\Omega_{\varepsilon}^\pm:=\Omega({\bm\Gamma}_\varepsilon^\pm)$. 
It follows from (\ref{eqn 2.1.10})-(\ref{eqn 2.1.10c}) that for any $i=1,2$ and $\varepsilon\in(0,\varepsilon_0)$, we have 
\begin{equation}
\label{eqn 2.1.13}
a_i(p)\,{\rm dist}\big(p, \tilde{\Gamma}_{\varepsilon,i}^\pm\big)=\varepsilon
\,\,{\rm for}\,\,{\rm all}\,\,{\rm points}\,\,p\in\Gamma_{\varepsilon,i}^\pm.
\end{equation}
Therefore, any solution ${\bf\Gamma}_\varepsilon=(\Gamma_{\varepsilon,1},\Gamma_{\varepsilon,2})\in{\bm{\cal F}}_\varepsilon$ of Prob. \ref{prob 2.1.2} (representing either ${\bm\Gamma}_\varepsilon^+$ or ${\bm\Gamma}_\varepsilon^-$) can be characterized as a pair ${\bm\Gamma}_\varepsilon\in{\bf X}$ such that $a_i(p)\,{\rm dist}(p,\tilde{\Gamma}_{\varepsilon,i})=\varepsilon$ for all points  $p\in\Gamma_{\varepsilon,i}$, $i=1,2$, where $\tilde{\Gamma}_{\varepsilon,i}:=\Phi_{\varepsilon,i}({\bm\Gamma}_\varepsilon)$. 
\end{theorem}

\begin{definition}
\label{def 2.1.7} 
{(Lower and upper solutions of operator fixed-point problems)}
For $\varepsilon\in(0,\varepsilon_0)$, let ${\bf T}_\varepsilon$ denote one of the operators ${\bf T}_\varepsilon^\pm$. Then a curve-pair ${\bf\Gamma}_\varepsilon\in{\bf X}$ is a lower (resp.  strict lower, upper, strict upper) solution of the operator "fixed point" problem (Prob. \ref{prob 2.1.2}) at $\varepsilon$ if ${\bf T}_\varepsilon({\bf\Gamma}_\varepsilon)\geq$ (resp. $>$, $\leq$, or $<$) ${\bf\Gamma}_\varepsilon$
\end{definition}
\begin{lemma}
\label{lem 2.1.1} 
{(Inequalities for operators (See Rem. \ref{rem 2.1.2}))}
In the context of Defs. \ref{def 2.1.4} and \ref{def 2.1.5}:
\vspace{.1in}

\noindent
(a) We have $\Gamma_1\leq\Phi_{\varepsilon,1}({\bf \Gamma})<\Phi_{\varepsilon,2}({\bf\Gamma})\leq\Gamma_2$
for every ${\bf \Gamma}=(\Gamma_1,\Gamma_2)\in{\bf X}$ and $\varepsilon\in(0,1/2)$.
\vspace{.1in}

\noindent
(b) For any $\varepsilon\in(0,1/2)$ and ${\bf \Gamma}_1,{\bf \Gamma}_2\in{\bf X}$ such that ${\bf \Gamma}_1\leq{\bf \Gamma}_2$, we have ${\bf \Phi}_\varepsilon({\bf \Gamma}_1)\leq{\bf \Phi}_\varepsilon({\bf \Gamma}_2)$, where ${\bf \Phi}_\varepsilon=(\Phi_{\varepsilon,1},\Phi_{\varepsilon,2})$.
\vspace{.1in}

\noindent
(c) We have $\Psi_{\varepsilon,1}^{-}(\Gamma)\leq\Psi_{\varepsilon,1}^{+}(\Gamma)\leq\Gamma\leq\Psi_{\varepsilon,2}^{-}(\Gamma)\leq\Psi_{\varepsilon,2}^{+}(\Gamma)$ for every $\Gamma\in \tilde{\rm X}_\varepsilon$ and $\varepsilon\in(0,\varepsilon_0)$. 
\vspace{.1in}

\noindent
(d) We have ${\rm dist}(\Psi_{1,\varepsilon}^\pm(\Gamma_1),\Psi_{2,\varepsilon}^\pm(\Gamma_2))\geq (2\varepsilon/\,\overline{A})$ for any ${\bf\Gamma}=(\Gamma_1,\Gamma_2)\in{\bf X}$ and any $\varepsilon\in(0,\varepsilon_0)$.
\vspace{.1in}

\noindent
(e) If $\Gamma_1\leq\Gamma_2$ in ${\rm X}$, then
$\Psi_{\varepsilon,1}^{\pm}(\Gamma_1)\leq\Psi_{\varepsilon,1}^{\pm}(\Gamma_2)$
and $\Psi_{\varepsilon,2}^{\pm}(\Gamma_1)\leq\Psi_{\varepsilon,2}^{\pm}(\Gamma_2)$,
both for all $\varepsilon\in(0,\varepsilon_0)$.
\vspace{.1in}

\noindent
(f) If $\varepsilon\in(0,\varepsilon_0)$ and ${\bf \Gamma}_1\leq{\bf \Gamma}_2$ in ${\bf X}$, then ${\bf T}_\varepsilon^{\pm}({\bf \Gamma}_1)\leq{\bf T}_\varepsilon^{\pm}({\bf \Gamma}_2)$. Therefore ${\bf T}_\varepsilon^-({\bf \Gamma}_1)\leq{\bf T}_\varepsilon^+({\bf \Gamma}_2)$.
\end{lemma}

\noindent
{\bf Proofs.}
Part (a) follows from the comparison principle for capacitary potentials, and Part (b) follows component-wise from the same.
\vspace{.1in}

\noindent
Turning to Part (c), for fixed $\varepsilon\in(0,\varepsilon_0)$ and  $\Gamma\in\tilde{\rm X}_\varepsilon$, we have $\dot{G}_{\varepsilon,i}(\Gamma):=\big\{p\in D_i(\Gamma):a_i(p)\,{\rm dist}\big(p,\Gamma\big)>\varepsilon\big\}\subset\dot{H}_{\varepsilon,i}(\Gamma):=\big\{p\in D_i(\Gamma):
a_i(p)\,{\rm dist}\big(p,\Gamma)\geq\varepsilon\}\subset D_i(\Gamma)$, for $i=1,2$, from which it follows that $G_{\varepsilon,i}(\Gamma)\subset H_{\varepsilon,i}(\Gamma)\subset D_i(\Gamma)$ for $i=1,2$, and therefore that $\partial{G}_{\varepsilon,1}(\Gamma)\leq\partial{H}_{\varepsilon,1}(\Gamma)\leq\Gamma\leq\partial{H}_{\varepsilon,2}(\Gamma)\leq\partial{G}_{\varepsilon,2}(\Gamma)$. In view of this, the claim (c) follows from Eqs. (\ref{eqn 2.1.4}) and (\ref{eqn 2.1.5}).
\vspace{.1in}

\noindent
Concerning Part (d), we have ${\rm dist}\,(\Gamma_1,\Psi_{\varepsilon,1}^\pm(\Gamma_1))\geq(\varepsilon/\overline{A})$ and ${\rm dist}\,(\Gamma_2,\Psi_{\varepsilon,2}^\pm(\Gamma_2))$ $\geq(\varepsilon\big/\,\overline{A})$, where $\Psi_{\varepsilon,1}^\pm(\Gamma_1)\leq\Gamma_1\leq\Gamma_2\leq\Psi_{\varepsilon,2}^\pm(\Gamma_2)$.
\vspace{.1in}

\noindent
Turning to Part (e), we assume that $\Gamma_1\leq\Gamma_2$ in $\tilde{\rm X}_\varepsilon$. Then for all $p\in \partial G_{\varepsilon,i}(\Gamma_i)$, we have that $\varepsilon=a_i(p)\,{\rm dist}\big(p,\Gamma_i\big)\leq a_i(p)\,{\rm dist}\,\big(p,\Gamma_{3-i}\big)$, from which it follows that $\partial G_{\varepsilon,i}(\Gamma_i)\subset {\rm Cl}\big(G_{\varepsilon,i}(\Gamma_{3-i})\big)={\rm Cl}\big(D_i(\partial G_{\varepsilon,i}(\Gamma_{3-i}))\big)$, $i=1,2$. It follows from this by (\ref{eqn 2.1.4}) and (\ref{eqn 2.1.5}) that $\Psi_{\varepsilon}^+(\Gamma_2)=\partial G_{\varepsilon,2}(\Gamma_2)\geq\partial G_{\varepsilon,2}(\Gamma_1)=\Psi_{\varepsilon}^+(\Gamma_1)$ and $\Psi_\varepsilon^-(\Gamma_1)=\partial G_{\varepsilon,1}(\Gamma_1)\leq\partial G_{\varepsilon,1}(\Gamma_2)=\Psi_\varepsilon^-(\Gamma_2)$. The remaining two inequalities follow from an analogous argument based on the sets $H_{\varepsilon,i}(\Gamma_i).$
\vspace{.1in}

\noindent
Finally, concerning Part (f), the first assertion follows from Parts (b) and (e), and the second assertion follows from Parts (b) and (c).

\begin{lemma}
\label{lem 2.1.2} Given $\varepsilon\in(0,1)$, let ${\bf\Gamma}_\varepsilon\in{\bf X}$ denote a lower (resp. upper) solution of Prob. \ref{prob 2.1.2} at  $\varepsilon$. Then ${\bf \Gamma}_{\varepsilon,n}:={\bf T}_\varepsilon^n({\bf\Gamma}_\varepsilon)\in{\bf X}$ is a lower (resp. upper) solution of Prob. \ref{prob 2.1.2} at $\varepsilon$ such that ${\bf\Gamma}_{\varepsilon,1}\geq\,(\leq)\,{\bf\Gamma}_\varepsilon$.
\end{lemma}
 
\begin{theorem} 
{(Existence of operator "fixed points" between upper and lower solutions of Prob. \ref{prob 2.1.2} (See Rem. \ref{rem 2.1.2}))}
\label{thm 2.1.2}
Given $\varepsilon\in(0,\varepsilon_0)$, let ${\bf\Gamma}_{*,\varepsilon}^\pm\in {\bf X}$ denote two curve-pairs such that
\begin{equation}
\label{eqn 2.1.13*}
{\bf\Gamma}_{*,\varepsilon}^-\leq{\bf\Gamma}_{*,\varepsilon}^+,
\end{equation}
\begin{equation}
\label{eqn 2.1.14}
{\bf\Gamma}_{*,\varepsilon}^-\leq {\bf T}_\varepsilon^-({\bf\Gamma}_{*,\varepsilon}^-);\,\,  {\bf T}_\varepsilon^+({\bf\Gamma}_{*,\varepsilon}^+)\leq{\bf\Gamma}_{*,\varepsilon}^+.
\end{equation}
Then: (a) There exist (not necessarily distinct) "fixed points" ${\bf \Gamma}_\varepsilon^\pm\in {\bf X}$ of the double-free-boundary operators ${\bf T}_\varepsilon^\pm$, respectively, such that
\begin{equation}
\label{eqn 2.1.15}
{\bf\Gamma}_{*,\varepsilon}^-\leq {\bf\Gamma}_\varepsilon^-\leq{\bf\Gamma}_\varepsilon^+\leq {\bf\Gamma}_{*,\varepsilon}^+,
\end{equation}
and such that ${\bf\Gamma}_\varepsilon^-\leq{\bf\Gamma}_\varepsilon\leq{\bf\Gamma}_\varepsilon^+$ for any other "fixed point" $\bf{\Gamma}_\varepsilon$ of either operator such that ${\bf\Gamma}_{*,\varepsilon}^-\leq {\bf\Gamma}_\varepsilon\leq{\bf\Gamma}_{*,\varepsilon}^+$. 
\vspace{.1in}

\noindent
(b) In fact the maximal (resp. minimal) fixed point ${\bf\Gamma}_\varepsilon^+$ (resp. ${\bf\Gamma}_\varepsilon^-$) of the operator ${\bf T}_\varepsilon^+$ (resp. ${\bf T}_\varepsilon^-$) is the limit of the weakly monotone decreasing (resp. weakly monotone increasing) sequence of upper solutions ${\bf\Gamma}_{\varepsilon,n}^+:=({\bf T}_\varepsilon^+)^n({\bf\Gamma}_{*,\varepsilon}^+)\in {\bf X}$ (resp. lower solutions  ${\bf\Gamma}_{\varepsilon,n}^-:=({\bf T}_\varepsilon^-)^n({\bf\Gamma}_{*,\varepsilon}^-)\in {\bf X}$).
\end{theorem}

\begin{definition}
\label{def 2.1.8}
{(An invariant set ${\bf Y})$}
Let $\tilde{\bf \Gamma}=(\tilde{\Gamma}_1,\tilde{\Gamma}_2)\in{\bf X}\cap C^2$ and $\hat{{\bf \Gamma}}=(\hat{\Gamma}_1,\hat{\Gamma}_2)\in{\bf X}\cap C^2$ be respective lower and upper solutions of Prob. \ref{prob2.1.1} (see Def. \ref{def 2.1.2}) such that $\tilde{{\bf\Gamma}}<\hat{{\bf\Gamma}}$. We define ${\bf Y}:=\{{\bf \Gamma}\in{\bf X}: \tilde{{\bf\Gamma}}\leq{\bf\Gamma}\leq\hat{{\bf\Gamma}}\}.$ By the definitions of strict lower and upper solutions, there exists a value $\varepsilon_1=\varepsilon_1(\tilde{{\bm\Gamma}},\hat{{\bm\Gamma}})\in(0,\varepsilon_0)$ so small that for any value $\varepsilon\in(0,\varepsilon_1]$, we have 
$$a_1(p)\,{\rm dist}(p,\{\tilde{U}_1=\varepsilon\})\geq\varepsilon\,\, {\rm on}\,\, \tilde{\Gamma}_1;\,\,a_2(p)\,{\rm dist}(p,\{\tilde{U}_2=\varepsilon\})\leq\varepsilon\,\,{\rm on}\,\,\tilde{\Gamma}_2,$$
$$a_1(p)\,{\rm dist}(p,\{\hat{U}_1=\varepsilon\})\leq\varepsilon\,\,{\rm on}\,\,\hat{\Gamma}_1;\,\, a_2(p)\,{\rm dist}(p,\{\hat{U}_2=\varepsilon\})\geq\varepsilon\,\,{\rm on}\,\,\hat{\Gamma}_2,$$
\noindent
where $\tilde{U}_i(p):=U_i(\tilde{{\bf\Gamma}};p)$ and $\hat{U}_i(p):=U_i(\hat{{\bf\Gamma}};p)$, both for $i=1,2$. It easily follows from this that ${\bf T}_\varepsilon^\pm(\tilde{{\bf \Gamma}})\geq\tilde{{\bf\Gamma}}$
and ${\bf T}_\varepsilon^\pm(\hat{{\bf \Gamma}})\leq\hat{{\bf\Gamma}}$ whenever $\varepsilon\in(0,\varepsilon_1]$. In view of the monotonicity properties of the operators (Lem. \ref{lem 2.1.1}(f)), we have
$${\bf T}_\varepsilon^\pm({\bf\Gamma}):{\bf Y}\rightarrow{\bf Y},\,\, \varepsilon\in (0,\varepsilon_1).$$
\end{definition}

\begin{theorem} 
{(Existence of "Fixed points" in an invariant set)}
\label{thm 2.1.3}
We let ${\bf\tilde{\Gamma}}\in{\bf X}\cap C^2$ and ${\bf\hat{\Gamma}}\in{\bf X}\cap C^2$) denote respective strict lower and upper solutions of Prob. \ref{prob2.1.1} (see Def. \ref{def 2.1.2}) such that ${\bf\tilde{\bf\Gamma}}<{\bf\hat{\bf\Gamma}}$. Then:
\vspace{.1in}

\noindent
(a) For any $\varepsilon\in(0,\varepsilon_1]$ (with $\varepsilon_1:=\varepsilon_1(\tilde{\Gamma},\hat{\Gamma})$ as in Def. \ref{def 2.1.8}), all the assertions of Thm. \ref{thm 2.1.2} hold, where one defines ${\bf \Gamma}_{*,\varepsilon}^-:={\bf\tilde{\Gamma}}$ and ${\bf \Gamma}_{*,\varepsilon}^+:={\bf\hat{\Gamma}}$. In particular, there exists at least one fixed point ${\bm\Gamma}_\varepsilon^\pm\in{\bf Y}$ of either operator ${\bf T}_\varepsilon^\pm$. 
\vspace{.1in}

\noindent
(b) There exists a constant $M$ such that
\begin{equation}
\label{eqn 2.1.16}
K(\Gamma_{\varepsilon,1}), K(\Gamma_{\varepsilon,2}), ||\Gamma_{\varepsilon,1}||, ||\Gamma_{\varepsilon,2}||\leq M,
\end{equation}
uniformly for all sufficiently small $\varepsilon\in(0,\varepsilon_1]$, and all fixed points ${\bf\Gamma}_\varepsilon\in{\bf Y}$ of either of the operators ${\bf T}_\varepsilon^\pm$ at $\varepsilon$, such that $K(\Gamma_{\varepsilon,1}),K(\Gamma_{\varepsilon,2}), ||\Gamma_{\varepsilon,1}||, ||\Gamma_{\varepsilon,2}||<\infty$. 
\end{theorem}
\begin{corollary}
\label{cor 2.1.1}
For sufficiently small $\varepsilon\in(0,\varepsilon_1]$, the estimate (\ref{eqn 2.1.16}) in Thm. \ref{thm 2.1.3}(b) applies to the "fixed points" ${\bf\Gamma}_\varepsilon^\pm\in{\bf Y}$ of the operators ${\bf T}_\varepsilon^\pm$, since they are both limits of monotone sequences of operator iterates (see the proofs of Thms. \ref{thm 2.1.1} and \ref{thm 3.1.2}(b)).
\end{corollary}

\noindent
{\bf  Proof of Thm. \ref{thm 2.1.2}.}  (See Rem. \ref{rem 2.1.2}(b).) For any fixed value $\varepsilon\in(0,\varepsilon_0)$, it follows from (\ref{eqn 2.1.13}) and (\ref{eqn 2.1.14}), by multiple application of Lems. \ref{lem 2.1.1}(f) and \ref{lem 2.1.2}, that
\begin{equation}
\label{eqn 2.1.17}
{\bf\Gamma}_{*,\varepsilon}^-\leq{\bf \Gamma}_{\varepsilon,1}^-\leq{\bf\Gamma}_{\varepsilon,2}^-\leq \cdots\leq{\bf\Gamma}_{\varepsilon,n}^-\leq
{\bf\Gamma}_{\varepsilon,n}^+\leq\cdots\leq{\bf\Gamma}_{\varepsilon,2}^+\leq{\bf\Gamma}_{\varepsilon,1}^+\leq{\bf\Gamma}_{*,\varepsilon}^+
\end{equation}
\noindent
in ${\bf X}$ for all $n\in N$, where we define ${\bf\Gamma}_{\varepsilon,k}^\pm:=({\bf T}_\varepsilon^\pm\big)^k({\bf\Gamma}_{*,\varepsilon}^\pm)\in {\bf X}$ for all $k\in N$, from which it follows that 
\begin{equation}
\label{eqn 2.1.17a}
{\bf\Gamma}_{\varepsilon,n+1}^\pm={\bf T}_{\varepsilon}^\pm\big({\bf\Gamma}_{\varepsilon,n}^\pm\big),
\end{equation} 
for all $n\in N$. We also have that 
\begin{equation}
\label{eqn 2.1.18}
{\rm dist}(\Gamma_{\varepsilon,k,1}^\pm,\Gamma_{\varepsilon,k,2}^\pm)\geq (\varepsilon/\,\overline{A}\,)
\end{equation}
for $k\in N$, due to Lem. \ref{lem 2.1.1}(d). For $\varepsilon\in(0,\varepsilon_0)$, we let $\Gamma_{\varepsilon,i}^{+}$, $i=1,2$, denote the boundary of the union $D_{\varepsilon,i}^{+}$ of the weakly increasing sequence (under set inclusion) of the upper complements $D_{\varepsilon,n,i}^{+}$ of the curves $\Gamma_{\varepsilon,n,i}^{+}$, $n\in N$. Similarly, we define $\Gamma_{\varepsilon,i}^{-}$, $i=1,2$, to be the boundary of the the union $D_{\varepsilon,i}^{-}$ of the weakly increasing sequence (under set inclusion) of the lower complements $D_{\varepsilon,n,i}^{-}$ of the curves $\Gamma_{\varepsilon,n,i}^{-}$, $n\in N$. Then ${\bf\Gamma}_\varepsilon^\pm:=(\Gamma_{\varepsilon,1}^\pm,\Gamma_{\varepsilon,2}^\pm)\in{\bf X}$, since ${\rm dist}(\Gamma_{\varepsilon,1}^\pm,\Gamma_{\varepsilon,2}^\pm)\geq(\varepsilon/\,\overline{A}\,)$ by (\ref{eqn 2.1.18}). Since, for the above definitions of ${\bf\Gamma}_\varepsilon^\pm$, we have ${\bf\Gamma}_{\varepsilon,n}^+\downarrow{\bf\Gamma}_\varepsilon^+$ and ${\bf \Gamma}_{\varepsilon,n}^-\uparrow{\bf\Gamma}_\varepsilon^-$, both as $n\rightarrow\infty$, it follows by (\ref{eqn 2.1.17a}) and continuity-properties of the operators (see Lems. \ref{lem 2.1.5} and \ref{lem 2.1.6}) that 
\begin{equation}
\label{eqn 2.1.19}
{\bf\Gamma}_\varepsilon^\pm=\lim_{n\rightarrow\infty}{\bf\Gamma}_{\varepsilon,n+1}^\pm\lim_{n\rightarrow\infty}{\bf T}_{\varepsilon}^\pm\big({\bf\Gamma}_{\varepsilon,n}^\pm\big)=
{\bf T}_\varepsilon^\pm(\lim_{n\rightarrow\infty}{\bf\Gamma}_{\varepsilon,n}^\pm\big)=
{\bf T}_\varepsilon^\pm\big({\bf\Gamma}_\varepsilon^\pm\big).
\end{equation}
\noindent
Therefore, the pairs of directed arcs ${\bf\Gamma}_\varepsilon^+\in{\bf X}$ (resp. ${\bf\Gamma}_\varepsilon^-\in{\bf X}$), which are the respective "fixed points" of the operators ${\bf T}_\varepsilon^+$ (resp. ${\bf T}_\varepsilon^-$), can be obtained as limits of weakly decreasing (increasing) sequences of upper (lower) solutions, as asserted. The remaining assertions in Part (a) easily follow from (\ref{eqn 2.1.17}).
\vspace{.1in}

\noindent
{\bf Proof of Thm. \ref{thm 2.1.3}.} Part (a) follows from Thm. \ref{thm 2.1.2}, while Part (b) follows from Thm. \ref{thm 3.1.2}(c) and Lem. \ref{lem 2.2.5}.
\vspace{.1in}

\noindent
{\bf Summary of the proof in Chapters \ref{section 2} and \ref{section 3} of Thms. \ref{thm1.2} and \ref{thm2.1.0}.} By Rem. \ref{rem 2.1.1}, it suffices to prove Thm. \ref{thm2.1.0} in the context of Prob. \ref{prob2.1.1}. Toward the proof of Thm. \ref{thm2.1.0} following the plan developed in Rem. \ref{rem 2.1.2}(c), we have thus far defined two (closely related) one-parameter families of monotone operators ${\bf T}_\varepsilon^\pm$, $\varepsilon\in(0,\varepsilon_0)$ (see Defs. \ref{def 2.1.4} and \ref{def 2.1.5}), whose properties, stated in Lem. \ref{lem 2.1.1} and in Section \ref{subsection 2.2}, are the basis for the proof in this section of Thm. \ref{thm 2.1.2}, and therefore of Thm. \ref{thm 2.1.3}(a). At this point, we invoke the estimates in Section \ref{subsection 5.1} to assert the existence of uniform bounds (independent of sufficiently small $\varepsilon>0$) for the arc-length and total curvature of the fixed points (see Thm. \ref{thm 2.1.3}(b)). It remains to prove, based on Thm. \ref{thm 2.1.3} and Section \ref{subsection 2.3} (properties of operator fixed points), that there exist suitable weak solutions Prob. \ref{prob2.1.1} satisfying the same uniform upper bound on their total curvature (see Thm. \ref{thm 2.3.1}). Following this, Chapter \ref{section 3} will be devoted to proof that $|\nabla U_i(p)|=a_i(p)$ on $\Gamma_i$, in a strong sense for $i=1,2$, and that $\Gamma_i$ has further regularity properties (see Thms. \ref{thm 2.6.2} and \ref{thm 2.7.1}).

\subsection{General qualitative properties of the $\boldsymbol{{\cal A}}$-operators}
\label{subsection 2.2}

\begin{lemma} 
\label{lem 2.1.4}
{(Qualitative geometry of the sets $G_{\varepsilon,i}(\Gamma)$ and $H_{\varepsilon,i}(\Gamma)$)} In the context of Def. \ref{def 2.1.4}, given $\varepsilon\in(0,\varepsilon_0)$, an arc $\Gamma\in \tilde{\rm X}_\varepsilon$, and functions $a_i\in\boldsymbol{{\cal A}}$, $i=1,2,$ we set $\phi_i(p):=r(p)\,a_i(p)$ in ${\rm Cl}\big(D_i(\Gamma)\big)$ for $i=1,2$, where $r(p):={\rm dist}\big(p,\Gamma\big)$. Then:
\vspace{.1in}

\noindent
(a) For any given points $p,q$ such that $q\in\Gamma$ and $|p-q|=r(p)>0$, we have $\phi_i(p-\delta\nu)<\phi_i(p)$ for $0<\delta<r(p)$, where $\nu$ is the unit vector pointing from $q$ to $p$.
\vspace{.1in}

\noindent
(b) The continuous function $\phi_i(p):{\rm Cl}\big(D_i(\Gamma)\big)\rightarrow[0,\infty)$ cannot have a weak local minimum at any point $p_0\in D_i(\Gamma)$ at which $\phi_i(p_0)<\varepsilon_0$.
\vspace{.1in}

\noindent
(c) Assume for a given non-empty, bounded open set $\omega_i\subset D_i(\Gamma)$, that $\phi_i(p)\geq\varepsilon$ on $\partial\omega_i$ for some constant $\varepsilon\in(0,\varepsilon_0)$. Then $\phi_i(p)>{\varepsilon}$ throughout $\omega_i$.
\vspace{.1in}

\noindent
(d) For any fixed value $\varepsilon\in(0,\varepsilon_0)$, and any fixed arc $\Gamma\in\tilde{\rm X}_\varepsilon$, the $P$-periodic (in $x$) sets $G_i:=G_{\varepsilon,i}(\Gamma)$, $\overline{G}_i:={\rm Cl}\big(G_{\varepsilon,i}(\Gamma)\big)$, and $H_i:=H_{\varepsilon,i}(\Gamma)\,$, $i=1,2$ (such that $G_i\subset\overline{G}_i\subset H_i$, $i=1,2$), are all uniquely determined by Def. \ref{def 2.1.4}. 
Moreover, for $\eta>0$ sufficiently large, each one of these sets contains the strip $R_i(\eta):=\{|y|>\eta\}\cap D_i(\Gamma)$, and in fact it follows from Def. \ref{def 2.1.4} that every point of the set $G_i$ is connected to $R_i(\eta)$ by a closed arc $\gamma\subset D_i(\Gamma)$ such that $\overline{B}\big(q\,;\,\varepsilon\,b_i(q)\big)\subset D_i(\Gamma)$ (i.e. $\phi_i(q)>\varepsilon$) for all $q\in\gamma$, whereas every point of $\,\overline{G}_i$ and $H_i$ is joined to $R_i(\eta)$ by a closed arc $\gamma\subset D_i(\Gamma)$ such that ${B}\big(q\,;\,\varepsilon\,b_i(q)\big)\subset D_i(\Gamma)$ (i.e. $\phi_i(q)\geq {\varepsilon}$) for all $q\in\gamma$.  
\vspace{.1in}

\noindent
(e) In the context of Part (d), the connected sets $G_i$, $H_i$, and $\overline{G}_i:={\rm Cl}(G_i)$, $i=1,2$, are all simply-connected 
in the sense of having no "holes" (as discussed in the proof). Therefore, the arcs $\partial G_i$, $i=1,2$, must be double-point free. Also, in the notation of Def. \ref{def 2.1.1}, we have $G_i=D_i(\partial G_i)$, $\overline{G}_i={\rm Cl}\big(D_i(\partial G_i)\big)$, and $H_i=\partial H_i\cup D_i(\partial H_i)$, all for $i=1,2$. 
\vspace{.1in}

\noindent
(f) Each of the sets in Part (d) satisfies the condition: ${\rm dist}\big(p,\Gamma\big)>\varepsilon\,b_i(p)$ (i.e. $\phi_i(p)>\varepsilon$) at all its interior points, while also satisfying the condition: ${\rm dist}\big(p,\Gamma\big)=\varepsilon\,b_i(p)$ (i.e. $\phi_i(p)=\varepsilon$) at all its boundary points. 
\vspace{.1in}

\noindent
(g) For fixed $\varepsilon\in(0,\varepsilon_0)$, $\Gamma\in\tilde{\rm X}_\varepsilon$, $i=1,2$, and for any sufficiently large value $\eta>0$, the $P$-periodic (in $x$) sets $\hat{G}_i=\hat{G}_{\varepsilon,i}(\Gamma)$ and $\hat{H}_i=\hat{H}_{\varepsilon,i}(\Gamma)$ (see Def. \ref{def 2.1.4}) are simply connected sets containing $R_i(\eta)$. Therefore $\hat{G}_i=D_i(\partial\hat{G}_i)$ and ${\rm Int}\big(\hat{H}_i\big)=D_i(\partial\hat{H}_i)$. 
\vspace{.1in}
\end{lemma}
\noindent
{\bf Proof of Part (a).} 
We observe that $|r(p)-r(q)|\leq|p-q|$ for any $p,q\in\Re^2$, so that $\phi_i(p)$ is Lipschitz continuous. Also, it easily follows from the definition of the distance function that $r(p-\delta\boldsymbol{\nu})\leq r(p)-\delta$ for any pair of points $(p,q)\in \Re^2\times\Gamma$, any unit vector $\boldsymbol{\nu}$ such that $|p-q|=r(p)$ and $p-q=|p-q|\boldsymbol{\nu}$, and any value $0\leq\delta<r(p)$. Therefore, we have 
\begin{equation}
\label{eqn 2.1.31}
\phi_i(p-\delta\boldsymbol{\nu})=r(p-\delta\boldsymbol{\nu})\,a_i(p-\delta\boldsymbol{\nu})\leq\big(r(p)-\delta\big)\,\big(a_i(p)+|\nabla a_i(p_\delta^*)|\,\delta\big)
\end{equation}
$$\leq\phi_i(p)+r(p)\,|\nabla a_i(p_\delta^*)|\delta-a_i(p)\,\delta$$
$$\leq\phi_i(p)-\big(a^2_i(p)-\phi_i(p)\,|\nabla a_i(p_\delta^*)|\big)\big(\delta\big/a_i(p)\big),$$
where $p_\delta^*$ denotes a point on the line-segment joining $p$ to $p-\delta\nu$. It follows that if $\phi_i(p)<\varepsilon_0$, so that $$2\,\phi_i(p)\,|\nabla a_i(p_\delta^*)|<\,2\,\varepsilon_0\,|\nabla a_i(p_\delta^*)|\leq 2\varepsilon_0\,A_1\leq\underline{A}^2\leq a^2_i(p),$$ 
\noindent
then
\begin{equation}
\label{eqn 2.a.33}
\phi_i(p-\delta\boldsymbol{\nu})<\phi_i(p)-(a_i(p)/2)\,\delta
\end{equation}
for $0\leq\delta<r(p)$, which is impossible if $p$ is a local minimum (at positive altitude below $\varepsilon_0$) of the function $\phi_i$, and if $\delta>0$ is sufficiently small. 
\vspace{.1in}

\noindent
{\bf Proof of Parts (b) and (c).} Part (b) follows directly from Part (a). Concerning Part (c), we have $\phi_i(p)\geq\varepsilon$ in ${\rm Cl}(\omega_i)$, since otherwise the continuous function $\phi_i(p):{\rm Cl}(\omega_i)\rightarrow\Re$ has a local interior minimum contradicting Part (b). With this established, we can repeat the argument for Part (b) to obtain a contradiction under the assumption that $\phi_i(p)=\varepsilon$ for any point $p\in\omega_i$.
\vspace{.1in}

\noindent
{\bf Proof of Part (d).} We first observe that
\begin{equation}
\label{eqn 2.1.21**}
\inf\big\{\phi_i(x,y):x\in\Re\big\}\rightarrow\,+\infty\,\,\,{\rm as}\,\,\,y\rightarrow(-1)^i\,\infty.
\end{equation}
We apply (\ref{eqn 2.1.21**}) to choose a value $\eta>0$ so large that $\inf\{\phi_i(p):p\in R_i(\eta)\}>2\,\varepsilon$, and therefore that $R_i(\eta)\subset G_i\subset\overline{G}_i\subset H_i\subset\hat{G}_i\subset\hat{H}_i$, all for $i=1,2$, 
where we set $R_i(\eta):=\{|y|>\eta\}\cap D_i(\Gamma)$. The remaining assertions in Part (d) follow from this.
\vspace{.1in}

\noindent
{\bf Proof of Part (e).} 
For any fixed $\Gamma\in\tilde{{\rm X}}_{\varepsilon}$ and $i=1,2$, let $\dot{G}_i\,({\rm resp.}\,\dot{H}_i$) denote the open (closed) subset of $D_i(\Gamma)$ in which $\phi_i(p):=\phi_i(\Gamma,p)>\varepsilon$ (resp. $\phi_i(p)\geq\varepsilon$), and let $G_i$ (resp. $H_i$) denote the open (closed) subset of $\dot{G}_i$ (resp. $\dot{H}_i$) consisting of those points which are joined to $\{y=(-1)^i\,\infty\}$ by arcs $\gamma\subset\dot{G}_i$ (resp. $\gamma\subset\dot{H}_i$), and let $O_i\big({\rm Cl}(\dot{G}_i)\big)$ (resp. $O_i\big({\rm Cl}(G_i)\big)$, $O_i(\dot{H}_i)$, $O_i(H_i)$) denote the family of all bounded, arc-wise connected, open sets $\omega_i\subset\Re^2$ such that $\partial\omega_i\subset{\rm Cl}(\dot{G}_i)$ (resp. $\partial\omega_i\subset{\rm Cl}(G_i)\big)$, $\partial\omega_i\subset \dot{H}_i$, $\partial\omega_i\subset H_i$). It follows from Part (c) that (${\bf i}$): $\omega_i\subset\dot{G}_i$ (resp. $\omega_i\subset\dot{H}_i$) for any open set $\omega_i\in O_i\big({\rm Cl}(\dot{G}_i)\big)$ (resp. $\omega_i\in O_i(\dot{H}_i)$), and it follows from (${\bf i}$) that (${\bf ii}$): $\omega_i\subset {G}_i$ (resp. $\omega_i\subset{H}_i$) for any open set $\omega_i\in O_i\big({\rm Cl}({G}_i)\big)$ (resp. $\omega_i\in O_i({H}_i)$), due to the fact that any point $p\in\omega_i$ is joined by the shortest possible straight line-segment $\gamma_1\subset{\rm Cl}(\omega_i)$ to a point $q\in\partial\omega_i$, which is in turn joined by an arc $\gamma_2\subset{\rm Cl}\big(\dot{G}_i\big)$ (resp. $\gamma_2\subset\dot{H}_i$) to $\{y=(-1)^i\,\infty\}$, thus joining $p$ to $\{y=(-1)^i\,\infty\}$ by a composite arc $\gamma=\gamma_1+\gamma_2$ through ${\rm Cl}\big(\dot{G}_i\big)$ (resp. $\dot{H}_i$). It follows that (${\bf iii}$): the boundary $\partial G_i$ has no double points, since, given the connectedness of $G_i$, the boundary $\partial G_i$ can have a double-point only if there exists at least one bounded, non-empty, open set $\omega_i\in O_i\big({\rm Cl}(G_i)\big)$ such that $\omega_i\cap G_i=\emptyset$. But, in view of (${\bf ii}$), this leads to the contradiction that $\omega_i\subset G_i$ and $\omega_i\cap G_i=\emptyset$, proving (${\bf iii}$). Finally, it also follows from (${\bf ii}$) that (${\bf iv}$): $H_i$ has no "holes" $\omega_i$.  
\vspace{.1in}

\noindent
Toward an alternate perspective on Part (e), we remark that (${\bf v}$): for any fixed open set $\omega_i\in O_i(\dot{H}_i)$, any point $p\in \omega_i\cap{H}_i$, and any value $\delta=\delta(p)>0$ which is small enough so that $B_{\delta}(p)\subset\omega_i$, we have that $B_{\delta}(p)\subset {H}_i$ (due to (${\bf i}$) and a variant of the above composite-arc argument). It follows from (${\bf v}$) that (${\bf vi}$): we have $\omega_i\subset H_i$ for $i\in\{1,2\}$ and any open set $\omega_i\in O_i\big(\dot{H}_i\big)$ such that $\omega_i\cap{H}_i\not=\emptyset$. Therefore, we have that (${\bf vii}$): ${\rm Int}\big(\dot{H}_i\big)\subset H_i\subset \dot{H}_i$, from which it follows that where . proving the simple-connectedness assertion for the sets $H_i$, $i=1,2$. 
\vspace{.1in}

\noindent
Also, for $i=1,2$, let $S(\dot{G}_i)$ (resp. $S\big({\rm Cl}(\dot{G}_i)\big)$) denote the family of all bounded, arc-wise connected, open sets $\omega_i\subset\Re^2$ such that $\partial\omega_i\subset \dot{G}_i$ (resp. $\partial\omega_i\subset {\rm Cl}\big(\dot{G}_i\big)$).
Analogous reasoning based on Parts (b) and (c) shows that we have $\omega_i\subset G_i$ for any $\omega_i\in S(G_i)$ (resp. $\omega_i\in S\big(\,\overline{G}_i\big)$) such that $\omega_i\cap G_i\not=\emptyset$.
\vspace{.1in}

\noindent
{\bf Proof of Parts (f) and (g).} We have that $\phi_i(p)\geq\varepsilon$ in the sets ${\rm Cl}\big(G_i\big)$ and $H_i$, by definition. Therefore, we must have that $\phi_i(p)>\varepsilon$ throughout the interiors of the smallest simply-connected sets containing $G_i$ and $H_i$, which we denote by $G_i^*$ and $H_i^*$, respectively. In fact for any bounded open set $\omega_i\in O_i\big({\rm Cl}(G_i)\big)$ or $\omega_i\in O_i\big(H_i\big)$ (see Part (e)), we have $\phi_i(p)>\varepsilon$ throughout $\omega_i$ by Parts (b) and (c). Finally, in view of the continuity of the functions 
$\phi_i(p):{\rm Cl}\big(D_i(\Gamma)\big)\rightarrow[0,\infty)$, the proof that $\phi_i\big(\partial G_i\big)=\phi_i\big(\partial H_i\big)=\varepsilon$ follows from the maximality of the sets $G_i$, $H_i$ under set inclusion subject to the conditions that $\phi_i(p)\geq\varepsilon$ in the sets ${\rm Cl}\big(G_i\big)$ and $H_i$. For example, if $\phi_i(p)>\varepsilon$ at some point $p\in\partial G_i$ (resp. $p\in\partial H_i$), then the domain $G_i$ (resp. $H_i$) can be enlarged to contain a closed ball $\overline{B}(p\,;r)$ centered at $p$ and having a sufficiently small radius $r>0$ such that $\phi_i(p)>0$ throughout $\overline{B}(p\,;r)$. Also, every point $q\in \overline{B}(p\,;r)$ can be joined to $p$ by a radial line-segment $\gamma$ of length $|p-q|\leq r$. Thus, any point $q\in \overline{B}(p\,;r)$ such that $q\notin G_i$ (resp. $q\notin H_i$) can be joined to any point $q'\in B(p\,;r)\cap G_i$ (resp. $q'\in B(p\,;r)\cap H_i$) by two radial arcs in sequence, the first joining $p$ to $q$ and the second joining $p$ to $q'$. Since also $q'$ can be joined to $\{y=(-1)^i\,\infty\}$ by an arc $\gamma$ through $G_i$, we conclude that $q$ can be joined to $\{y=(-1)^i\,\infty\}$ by an arc $\gamma$ passing through $G_i\cup B(p\,;r)$ (resp. $H_i\cup B(p\,;r)$) on which $\phi_i(p)\geq\varepsilon$, completing the proof of Part (f). Finally, concerning Part (g), one shows that the sets $\hat{G}_i$ and $\hat{H}_i$, whose points are obviously all connected to $R_i(\eta)$, are in fact simply-connected, because, by an argument given in the proof of Part (e), any point $p$ in a bounded, connected, open subset $\omega_i$ of $\hat{G}_i$ (resp.  $\hat{H}_i$) can be joined by an arc through $\hat{G}_i$ (resp. $\hat{H}_i$) to $R_i(\eta)$ provided that all boundary points $q\in\partial \omega_i$ have the same property.  
\vspace{.1in}

\noindent
\begin{lemma}
\label{lem 2.1.4*}
{(Alternative approach to geometry of $\boldsymbol{{\cal A}}$-operators)} In the context of Def. \ref{def 2.1.4} and Lem. \ref{lem 2.1.4}:
\vspace{.1in}

\noindent
(a) The boundaries $\partial\hat{G}_i$ and $\partial\hat{H}_i$ of the $P$-periodic (in $x$) simply-connected sets $\hat{G}_i=\hat{G}_{\varepsilon,i}(\Gamma)$, $\hat{H}_i=\hat{H}_{\varepsilon,i}(\Gamma)$, $i=1,2,$ are double-point-free directed $C^{1,1}$-arcs such that $\hat{G}_i(\Gamma)=D_i(\partial \hat{G}_i)$ and $\hat{H}_i(\Gamma)={\rm Cl}\big(D_i(\partial \hat{H}_i)\big)$.
\vspace{.1in}

\noindent
(b) For any fixed $\varepsilon\in(0,\varepsilon_0)$, $\Gamma\in \tilde{\rm X}_\varepsilon$, and $i=1,2$, we call the points $p\in\partial\hat{H}_i:=\partial\hat{H}_{\varepsilon,i}(\Gamma)$ and $q\in\partial H_i:=\partial H_{\varepsilon,i}(\Gamma)$ "related by (\ref{eqn 2.1.21*})" if and only if 
\begin{equation}
\label{eqn 2.1.21*}
q=p+\varepsilon\,b_i(q)\,\hat{N}_i(p),
\end{equation}  
where we set $b_i(p):=\big(1\big/a_i(p)\big)$, and use $\hat{N}_i(p)$ to denote the (Lipschitz continuously varying) left-hand unit normal to the arc $\partial\hat{H}_{i}$ at any point $p\in\partial\hat{H}_i$. Then:
\vspace{.1in}

\noindent
(b1) For every point $p\in\partial\hat{H}_i$, there exists at least one point $q\in\partial H_i$ such that the points $p$ and $q$ are related by (\ref{eqn 2.1.21*}). 
\vspace{.1in}

\noindent
(b2) For every $q\in\partial H_i$, there exists a point $p\in\partial\hat{H}_i$ such that $p$ and $q$ are related by (\ref{eqn 2.1.21*}).
\vspace{.1in}

\noindent
(b3) The point $q$ is uniquely and Lipschitz-continuously determined by the point $p$ relative to the set of all ordered pairs $(p,q)\in\partial\hat{H}_i\times\partial H_i$ solving (\ref{eqn 2.1.21*}).
\vspace{.1in}

\noindent
(b4) In view of (b1), (b2), and (b3), the set of all ordered pairs of points $(p,q)\in\partial\hat{H}_i\times\partial H_i$ satisfying the relation (\ref{eqn 2.1.21*}) actually constitutes the graph of a Lipschitz-continuous function 
\begin{equation}
\label{eqn 2.1.22}
q=f_i(p):\partial\hat{H}_i\rightarrow\partial H_i,
\end{equation}
mapping the arc $\partial\hat{H}_i$ onto the arc $\partial H_i$. (This part also holds with $H_i$ and $\hat{H}_i$ replaced by $G_i$ and $\hat{G}_i$.)
\vspace{.1in}

\noindent
(c) For any fixed $\varepsilon\in(0,\varepsilon_0)$ and $\Gamma\in \tilde{\rm X}_\varepsilon$, the $P$-periodic (in $x$) set boundaries $\partial G_i:=\partial G_{\varepsilon,i}(\Gamma)$, $i=1,2,$ coincide with the positively-oriented, double-point-free directed arcs $\Gamma_i\in{\rm X}$, defined such that $G_1$ (resp. $G_2$) is the lower (upper) complement of $\Gamma_1$ (resp. $\Gamma_2$). Also, the $P$-periodic (in $x$) set-boundaries $\partial\,H_i=\partial\,H_{\varepsilon,i}(\Gamma)$, $i=1,2,$ coincide with $P$-periodic (in $x$) directed arcs $\Gamma_i^*\in{\rm X}$ which are positively-oriented by the requirement that the complement of $H_1$ (resp. $H_2$) be the upper (lower) complement of $\Gamma_1^*$ (resp. $\Gamma_2^*$), which always lies locally to the left (right) of $\Gamma_1^*$ (resp. $\Gamma_2^*$). 
\end{lemma}
\noindent
{\bf Proof of Part (a).} 
Concerning the smoothness of $\partial\hat{G}$, $\partial\hat{H}$, for any $\varepsilon\in(0,1/2)$, there exists a value $\delta=\delta(\varepsilon)>0$ such that any arc $\Gamma\in\tilde{\rm X}_\varepsilon$ has tangent balls of radius $\delta$ at all arc-points and on both sides of the arc. Also, for any arc $\Gamma\in\tilde{\rm X}_\varepsilon$, the corresponding arc $\partial\hat{G}_{\varepsilon,i}(\Gamma)$ has an interior tangent ball of the form $B\big(p\,;\varepsilon\,b_i(p)\big)$, $p\in\Re^2$, located in the interior of the region $\hat{G}_{\varepsilon,i}(\Gamma)$). Using these properties, one can show that for each arc $\Gamma\in\tilde{\rm X}_\varepsilon$, the corresponding arc $\partial\hat{G}_{\varepsilon,i}(\Gamma)$ has both an exterior tangent ball $B(p\,;\delta)$ and an interior tangent ball $B\big(p\,;\varepsilon\,b_i(p)\big)$ at every arc point. Therefore, the arc $\partial\hat{G}_{\varepsilon,i}(\Gamma)$ has uniformly bounded curvature at all points, and is therefore the smooth, double-point-free image of a $C^{1,1}$ arc-length parametrization $\hat{p}_{\varepsilon,i}(t):\Re\rightarrow\partial\hat{G}_{\varepsilon,i}(\Gamma)$. Also, the mapping $\hat{N}_{\varepsilon,i}(p):\partial\hat{G}_{\varepsilon,i}(\Gamma)\rightarrow\partial B(0,1)$ is Lipschitz continuous. (Obviously, the uniform, absolute curvature-bound increases as $\varepsilon$ decreases for $\varepsilon\in(0,\varepsilon_0)$.) The remaining assertions follows in the first case from the fact that $\varepsilon=a_i(p)\,{\rm dist}\big(p,\Gamma\big)=a_i(p)|p-q|$ for any points $p\in\partial G_{\varepsilon,i}(\Gamma)$ and $q\in\partial\hat{G}_{\varepsilon,i}(\Gamma)$.  
\vspace{.1in}

\noindent
{\bf Proof of Part (b1).} In fact the unknown point $q\in\partial H_i$ can be expressed in the general form: $q=q(\alpha):=p+\hat{N}_ip)\,\alpha$, in which $\alpha\in\Re$. By substitution, we see that the point $q(\alpha)$ solves (\ref{eqn 2.1.21*}) for given $p\in\partial\hat{H}_i$ if and only if the value $\alpha\in\Re\,$ solves the related equation $g_i(p, \alpha)=0$, where we define the continuous function $g_i(p, \alpha):\Re\rightarrow\Re$  such that
\begin{equation}
\label{eqn 2.1.23}
g_i(p,\alpha):=b_i\big(p+\hat{N}_i(p)\,\alpha\big)-\alpha.
\end{equation}
By definition, we have that $g_i(p, 0)=b_i(p)>0$, and that $g_i(p, \alpha)\rightarrow-\infty$ as $\alpha\rightarrow\infty$. Therefore, by the intermediate value theorem, there exists at least one value $\alpha\in(0,\infty)$ such that $g_i(p,\alpha)=0$, and therefore such that the ordered pair $(p,q(\alpha))$ satisfies (\ref{eqn 2.1.23}).
\vspace{.1in}

\noindent
{\bf Proof of Part (b2).} By Lem. \ref{lem 2.1.4}(d), we have ${\rm dist}\big(q,                                                           \Gamma\big)=\phi_i(q)\,b_i(q)>\varepsilon\,b_i(q)$ for all $q\in{\rm Int}\big(H_{\varepsilon,i}(\Gamma)\big)$, from which it follows that $\overline{B}\big(q\,;\varepsilon\/ b_i(q)\,\,\big)$ $\subset D_i(\Gamma)$ for all $q\in{\rm Int}\big(H_{\varepsilon,i}(\Gamma)\big)$. On the other hand, we have ${\rm dist}\big(q,\Gamma\big)=\phi_i(q)\,b_i(q)=\varepsilon\,b_i(q)$ for all $q\in\partial H_{\varepsilon,i}(\Gamma)$ (see Lem. \ref{lem 2.1.4}(d)). Therefore, among all points $q\in H_{\varepsilon,i}(\Gamma)$, we have that
\begin{equation}
{\rm dist}\big(q,\partial\hat{H}_{\varepsilon,i}(\Gamma)\big)={\rm dist}\big(q,\Gamma\big)=\varepsilon\,b_i(q)\iff q\in\partial H_{\varepsilon,i}(\Gamma).
\end{equation}
It follows that for any point $q\in\partial H_i$, there exists a corresponding point $p\in\partial\hat{H}_i$ such that $|q-p|={\rm dist}\big(q,\partial\hat{H}_i\big)=\varepsilon\,b_i(q)$. Clearly, $p\in\partial\hat{H}_i$ is such that $(q-p)$ is parallel to $\boldsymbol{\nu}(p)$ and $(p,q)$ satisfies the requirement of Part (b2).  
\vspace{.1in}

\noindent
{\bf Proof of Part (b3).} 
For fixed $\varepsilon\in(0,\varepsilon_0)$ and  $i\in\{1,2\}$, and for any two solution-pairs $(p,q)$, $(p',q')$ of (\ref{eqn 2.1.21*}), we have the equations:
$q=p+\varepsilon\,b(q)\,\hat{N}(p)$ and $q'=p'+\varepsilon\,b(q')\,\hat{N}(p')$, from which it follows by applying mean-value estimates to the various differences that 
\begin{equation}
\label{eqn 2.1.24}
|q-q'|\leq|p-p'|+\varepsilon\,|{\nabla} b(q^*)|\,|q-q'|+\varepsilon\,b(q)\,|\hat{K}(p^*)|\,|p-p'|,
\end{equation}
where $p^*$ lies between $p$ and $p'$. Therefore, we have
\begin{equation}
\label{eqn 2.1.25}
\big|q-q'\big|\leq\Big(\frac{1+\varepsilon\,\sup\{b(p)|\hat{K}(p)|\}}{1-\varepsilon\sup\{|{\nabla
}\,b(p)|\}}\Big)\big|p-p'\big|
\end{equation}
(where $\hat{K}(p)$ denotes the signed curvature of the arc $\hat{\Gamma}_{\varepsilon,i}$ at the point $p\in\hat{\Gamma}_{\varepsilon,i}$), provided only that $\varepsilon\sup\{|{\nabla} b_i(p)|:p\in\Re^2\}<1$ for $i=1,2$, as always occurs if $\varepsilon\in(0,\varepsilon_0)$. 
\vspace{.1in}

\noindent
{\bf Proof of Part (c).} The first claim follows from the simple-connectedness of the sets $G_i$ and ${\rm Cl}\big(G_i\big)$, as treated in Lem. \ref{lem 2.1.4}(e). Turning to the second claim, it follows from Part (b4), Eq. (\ref{eqn 2.1.22}), that the uniquely-determined boundary sets $\partial H_i$, $i=1,2$, of Lem. \ref{lem 2.1.4}(d) must coincide with the images $\Gamma_i^*$, $i=1,2$, of the positively-oriented, Lipschitz-continuous, composite arc-parametrizations of the form $q=q_i(t)=f_i(p_i(t)):\Re\rightarrow\partial H_i$, where $p=p_i(t):\Re\rightarrow\partial\hat{H}_i$ denotes a double-point-free arc-length parametrization of the smooth arc $\partial \hat{H}_i$. Clearly, the directed arcs $\Gamma_i^*$, $i=1,2$, retrace the set boundaries $\partial H_{\varepsilon,i}(\Gamma)$, $i=1,2$, as was asserted. We remark that the arcs $\Gamma_i^*$ may have double-points, but it is impossible for either one to "cross over itself", due to Lem. \ref{lem 2.1.4}(e),(f).
\begin{remark} {(Comments on Lem. \ref{lem 2.1.4})}
\label{rem 2.1.5}
(a) Given $\varepsilon\in(0,\varepsilon_0)$ and $\Gamma\in\tilde{{\rm X}}_\varepsilon$, let $\Gamma_{\varepsilon,i}:=\big\{p\in D_i(\Gamma):\phi_i(\Gamma;\,p)=\varepsilon\big\}$, and let $R_{\varepsilon,i}$ denote the set generated by all the straight line-segments $L_{\varepsilon,i}(p)$ joining points $p\in\Gamma_{\varepsilon,i}$ to corresponding points $q\in\Gamma$ such that $|p-q|={\rm dist}\big(p,\Gamma\big)$. Then $\phi_i(p)<\varepsilon$ throughout $R_{\varepsilon,i}$. Moreover, if $\varepsilon$ is less than the minimum radius of curvature of $\Gamma$, then the lines $L_{\varepsilon,i}(p)$ depend continuously on $p\in\Gamma_{\varepsilon,i}$, and therefore $R_{\varepsilon,i}$ is the annular domain such that $\partial R_{\varepsilon,i}=\Gamma\cup\Gamma_{\varepsilon,i}$. 
\vspace{.1in}

\noindent
(b) If we assume each point $q\in\partial H_i$ is related to one and only one point $p\in\Gamma$ under Eq. \ref{eqn 2.1.21*}, then $\Gamma=\partial\hat{H}_i$ and the mapping preserves the natural order of points $q\in\partial H_i$ and $p\in \Gamma=\partial\hat{H}_i$. However, if one point $q\in\partial H_i$ is related to each of two distinct points $p_1,p_2\in\Gamma$ (where we assume $p_1$ (resp. $p_2$) is minimal (maximal) relative to the natural ordering in $\Gamma$), then of course $p_1, p_2\in\Gamma$ both satisfy Eq.  (\ref{eqn 2.1.21*}), as before, but the remaining points $p\in\Gamma$ such that $p_1<p<p_2$ (in terms of the natural ordering in $\Gamma$) are not likely to be related to any point $q\in\partial H_i$, since it is possible to have $|p-q|>>\varepsilon\,b_i(q)$ for some points $p\in\Gamma$ such that $p_1<p<p_2$.
\end{remark}
\begin{lemma} {(One-sided continuity of the sets $G_{\varepsilon,i}(\Gamma)$ and $H_{\varepsilon,i}(\Gamma)$ as functions of the arcs $\Gamma\in\tilde{\rm X}_\varepsilon$)}
\label{lem 2.1.5}
(a) Let be given a weakly monotone sequence of arcs $\big(\Gamma_n\big)_{n=1}^\infty$, all in $\tilde{\rm X}_\varepsilon$ for some fixed value $\varepsilon\in(0,\varepsilon_0)$. Choose $i\in\{1,2\}$ such that the open sets $D_i(\Gamma_n)$ form a weakly increasing sequence (as ordered by set inclusion). Thus either $i=1$ and the arc-sequence $\big(\Gamma_n\big)_{n=1}^\infty$ is weakly increasing or else $i=2$ and the arc-sequence $\big(\Gamma_n\big)_{n=1}^\infty$ is weakly decreasing. Then we have $D^*:=\bigcup_{n=1}^\infty D_i(\Gamma_n)=D_i(\Gamma)$, where we define $\Gamma:=\partial{D}^*\in\tilde{\rm X}_\varepsilon$. Also, for the same values $\varepsilon\in(0,\varepsilon_0)$ and $i\in\{1,2\}$, we have (in terms of Def. \ref{def 2.1.4} and Lem. \ref{lem 2.1.4}) that $G_{\varepsilon,i}(\Gamma)=\bigcup_{n=1}^\infty G_{\varepsilon,i}(\Gamma_n)$, where $G_{\varepsilon,i}(\Gamma_n)=D_i(\Gamma_{\varepsilon,n,i})$
and $G_{\varepsilon,i}(\Gamma)=D_i(\Gamma_{\varepsilon,i})$, in which $\Gamma_{\varepsilon,n,i}:=\partial G_{\varepsilon,i}(\Gamma_n)\in{\rm X}$ and $\Gamma_{\varepsilon,i}:=\partial G_{\varepsilon,i}(\Gamma)\in{\rm X}$.
\vspace{.1in}

\noindent
(b) Again let be given a weakly monotone sequence of arcs $\big(\Gamma_n\big)_{n=1}^\infty$, all in $\tilde{\rm X}_\varepsilon$ for some fixed value $\varepsilon\in(0,\varepsilon_0)$. Choose $i\in\{1,2\}$ such that the closed sets $E_i(\Gamma_n):={\rm Cl}\big(D_i(\Gamma_n)\big)$ form a weakly-decreasing sequence (ordered by set inclusion). (Thus $i=1$ (resp. $i=2$) if the arc-sequence  $\big(\Gamma_n\big)_{n=1}^\infty$ is weakly decreasing (increasing.) 
Then $E^*:=\bigcap_{n=1}^\infty E_i(\Gamma_n)=E_i(\Gamma)$ for some arc $\Gamma:=\partial E^*\in\tilde{\rm X}_\varepsilon$. Also assume the arc-length per $P$-period (in $x$) of the arcs $\Gamma_n$ is uniformly bounded above over all $n\in N$. Then (in terms of Def. \ref{def 2.1.4}) we have $\bigcap_{n=1}^\infty H_{\varepsilon,i}(\Gamma_n)={H}_{\varepsilon,i}(\Gamma)$, where $H_{\varepsilon,i}(\Gamma_n)={\rm Cl}\big(D_i(\Gamma_{\varepsilon,n,i})\big)$ and $H_{\varepsilon,i}(\Gamma)={\rm Cl}\big(D_i(\Gamma_{\varepsilon,i})\big)$, in which $\Gamma_{\varepsilon,n,i}:=\partial H_{\varepsilon,i}(\Gamma_n)\in{\rm X}$ and $\Gamma_{\varepsilon,i}:=\partial H_{\varepsilon,i}(\Gamma)\in {\rm X}$.
\vspace{.1in}

\noindent
(c) For any fixed value $\varepsilon\in(0,\varepsilon_0)$, it follows directly from Parts (a) and (b) that $\partial G_{\varepsilon,1}(\Gamma_n)\uparrow \partial G_{\varepsilon,1}(\Gamma)$ and $\partial H_{\varepsilon,2}(\Gamma_n)\uparrow \partial H_{\varepsilon,2}(\Gamma)$, both as  $\Gamma_n\uparrow\Gamma$ in $\tilde{\rm X}_\varepsilon$, and that $\partial G_{\varepsilon,2}(\Gamma_n)\downarrow \partial G_{\varepsilon,2}(\Gamma)$ and $\partial H_{\varepsilon,1}(\Gamma_n)\downarrow \partial H_{\varepsilon,1}(\Gamma)$, both as  $\Gamma_n\downarrow\Gamma$ in $\tilde{\rm X}_\varepsilon$. In view of Def. \ref{def 2.1.5}, Eqs. (\ref{eqn 2.1.4}) and (\ref{eqn 2.1.5}), it follows that $\Psi_{\varepsilon,i}^+(\Gamma_{n,i})\downarrow\Psi_{\varepsilon,i}^+(\Gamma_i)$ for $\Gamma_{n,i}\downarrow\Gamma_i$ in $\tilde{\rm X}_\varepsilon$, and $\Psi_{\varepsilon,i}^-(\Gamma_{n,i})\uparrow\Psi_{\varepsilon,i}^-(\Gamma_i)$, for $\Gamma_{n,i}\uparrow\Gamma_i$ in $\tilde{\rm X}_\varepsilon$, both as $n\rightarrow\infty$ for either $i=1$ or $i=2$, where $\uparrow$ and $\downarrow$ refer, respectively, to weakly monotone increasing or weakly monotone decreasing uniform convergence. 
\end{lemma}

\noindent
{\bf Proof of Part (a).}
For $i=1,2$, we have $D_i(\Gamma_n)\subset D_i(\Gamma_{n+1})\subset D_i(\Gamma)$ for all $n\in N$. Therefore, $\dot{G}_{\varepsilon,i}(\Gamma_n)\subset\dot{G}_{\varepsilon,i}(\Gamma_{n+1})\subset\dot{G}_{\varepsilon,i}(\Gamma)$, and we have $\gamma\subset\dot{G}_{\varepsilon,i}(\Gamma_{n+1})$ (resp. $\gamma\subset\dot{G}_{\varepsilon,i}(\Gamma)$) for any directed arc $\gamma\subset\dot{G}_{\varepsilon,i}(\Gamma_n)$ (resp. $\gamma\subset\dot{G}_{\varepsilon,i}(\Gamma_{n+1})$) joining any point $p\in\dot{G}_{\varepsilon,i}(\Gamma_n)$ (resp. $p\in\dot{G}_{\varepsilon,i}(\Gamma_{n+1})$) to $\{y=(-1)^i\,\infty\}$. Therefore, we have ${G}_{\varepsilon,i}(\Gamma_n)\subset{G}_{\varepsilon,i}(\Gamma_{n+1})\subset{G}_{\varepsilon,i}(\Gamma)$, for $i=1,2$, and all $n\in N$, from which it follows that (${\bf i}$):  $\bigcup_{n=1}^\infty {G}_{\varepsilon,i}(\Gamma_n)\subset {G}_{\varepsilon,i}(\Gamma)$. 
\vspace{.1in}

\noindent
In view of (${\bf i}$), it suffices to show that (${\bf ii}$): ${G}_{\varepsilon,i}(\Gamma)\subset\bigcup_{n=1}^\infty{G}_{\varepsilon,i}(\Gamma_n)$, for $i=1,2$ and all $n\in N.$ To this end, we let $p_{\varepsilon,i}$ denote any particular point in the open set ${G}_{\varepsilon,i}(\Gamma)$. Then $p_{\varepsilon,i}\in\dot{G}_{\varepsilon,i}(\Gamma)$ and $p_{\varepsilon,i}$ can be joined to $\{y=(-1)^i\,\infty\}$ by a suitable directed arc $\gamma_{\varepsilon,i}$ contained in the open set $\dot{G}_{\varepsilon,i}(\Gamma)$. Given a suitable arc $\gamma_{\varepsilon,i}$ and a value $\kappa>0$ so large that ${\rm dist}\big(p,\Gamma_n)\big)\geq \varepsilon\,b_i(p)$ for all $n\in N$ and all points $p\notin R_\kappa:=\Re\times[-\kappa,\kappa]$, we can replace the directed arc $\gamma_{\varepsilon,i}$ by a new directed arc $\gamma_{\varepsilon,i}^*:=\gamma_{\varepsilon,i}^{(1)}+\gamma_{\varepsilon,i}^{(2)}$, consisting of an initial arc-segment $\gamma_{\varepsilon,i}^{(1)}$, which coincides with the arc $\gamma_{\varepsilon,i}$ up to the first point, called $q_{\varepsilon,i}$, at which $\gamma_{\varepsilon,i}$ touches the line $\{y=(-1)^i\,\kappa\}$, followed by a second arc-segment $\gamma_{\varepsilon,i}^{(2)}$, which follows a straight vertical path from the point $q_{\varepsilon,i}$ to $\{y=(-1)^i\,\infty\}$. Since $\dot{G}_{\varepsilon,i}(\Gamma)=\bigcup_{n=1}^\infty \dot{G}_{\varepsilon,i}(\Gamma_n)$, it follows that the set-family  $\big(\dot{G}_{\varepsilon,i}(\Gamma_n)\big)_{n=1}^\infty$ is an open covering of the compact set $\gamma_{\varepsilon,i}^{(1)}$, for which there exists a finite sub-covering. Since the sets of the open covering are weakly monotone increasing under set inclusion, it follows that $\gamma_{\varepsilon,i}^{(1)}\subset\dot{G}_{\varepsilon,i}(\Gamma_n)$, and therefore that $\gamma_{\varepsilon,i}^*\subset \dot{G}_{\varepsilon,i}(\Gamma_n)$, and therefore again that $p_{\varepsilon,i}\in G_{\varepsilon,i}(\Gamma)$ for all $n\in N$ which are sufficiently large depending in each case on the point $p_{\varepsilon,i}\in{G}_{\varepsilon,i}(\Gamma)$. Therefore $p_{\varepsilon,i}\in\bigcup_{n=1}^\infty G_{\varepsilon,i}(\Gamma_n)$ for any point $p_{\varepsilon,i}\in{G}_{\varepsilon,i}(\Gamma)$, from which (${\bf ii}$) follows.  
\vspace{.1in}

\noindent
{\bf Proof of Part (b).} Since $E_i(\Gamma)\subset E_i(\Gamma_{n+1})\subset E_i(\Gamma_n)$ for all $n\in N$, we have $\dot{H}_{\varepsilon,i}(\Gamma)\subset\dot{H}_{\varepsilon,i}(\Gamma_{n+1})\subset \dot{H}_{\varepsilon,i}(\Gamma_n)$, from which it follows that $\gamma\subset\dot{H}_{\varepsilon,i}(\Gamma_{n+1})$ (resp. $\gamma\subset\dot{H}_{\varepsilon,i}(\Gamma_n)$) for any directed arc $\gamma\subset\dot{H}_{\varepsilon,i}(\Gamma)$ (resp. $\gamma\subset\dot{H}_{\varepsilon,i}(\Gamma_{n+1})$) joining any point $p\in\dot{H}_{\varepsilon,i}(\Gamma)$ (resp. $p\in\dot{H}_{\varepsilon,i}(\Gamma_{n+1})$) to $\{y=(-1)^i\,\infty\}$. Therefore, 
we have ${H}_{\varepsilon,i}(\Gamma)\subset{H}_{\varepsilon,i}(\Gamma_{n+1})\subset{H}_{\varepsilon,i}(\Gamma_n)$, both for $i=1,2$, and all $n\in N$, from which it follows that 
(${\bf i}$): ${H}_{\varepsilon,i}(\Gamma)\subset\bigcap_{n=1}^\infty {H}_{\varepsilon,i}(\Gamma_n)$.
\vspace{.1in}

\noindent
In view of (${\bf i}$), it suffices to prove that (${\bf ii}$): $\bigcap_{n=1}^\infty {H}_{\varepsilon,i}(\Gamma_n)\subset{H}_{\varepsilon,i}(\Gamma)$. Toward this end, we observe that for any given point $p_{\varepsilon,i}\in\bigcap_{n=1}^\infty{H}_{\varepsilon,i}(\Gamma_n)$, there exists a sequence $\big(\gamma_{\varepsilon,n,i}\big)_{n=1}^\infty$ of directed arcs such that (${\bf iii}$): for each $n\in N$, the point $p_{\varepsilon,i}$ is joined to $\{y=(-1)^i\,\infty\}$ by the arc $\gamma_{\varepsilon,n,i}$, which is located entirely inside the closed set $\dot{H}_{\varepsilon,i}(\Gamma_n):=\big\{p\in E_i(\Gamma_n):a_i(p)\,{\rm dist}\big(p,\Gamma_n\big)\geq\varepsilon\big\}$. Therefore, we have (${\bf iv}$): $\gamma_{\varepsilon,n,i}\subset \dot{H}_{\varepsilon,i}(\Gamma_m)$ for any $m,n\in N$ such that $m\leq n$. Assuming that (${\bf v}$): $\kappa_{n,1}=\kappa_{n,1}(\varepsilon)$ (resp. $\kappa_{n,2}=\kappa_{n,2}(\varepsilon))$ is maximum (resp. minimum) subject to the requirement that $\partial H_{\varepsilon,1}(\Gamma_n), \partial H_{\varepsilon,2}(\Gamma_n)\subset\Re\times[\kappa_{n,1},\kappa_{n,2}]$), we let $\gamma_{\varepsilon,n,i}^*=\gamma_{\varepsilon,n,i}^{(1)}+\gamma_{\varepsilon,n,i}^{(2)}$ for all $n\in N$, where, for $i=1,2$, $\,\,\gamma_{\varepsilon,n,i}^{(1)}$ denotes the initial segment of the arc $\gamma_{\varepsilon,n,i}$, up to its first point of contact, called $r_{\varepsilon,n,i}$, with the line $\{y=\kappa_{n,i}\}$, and $\gamma_{\varepsilon,n,i}^{(2)}$ denotes the vertical path from $r_{\varepsilon,n,i}$ to $\{y=(-1)^i\,\infty\}$. As an application of (${\bf v}$), we have that $\gamma_{\varepsilon,n,i}^{(2)}\subset \dot{H}_{\varepsilon,n,i}(\Gamma_n)$, and therefore (${\bf vi}$): $\gamma_{\varepsilon,n,i}^*\in\dot{H}_{\varepsilon,n,i}(\Gamma_n)$. We assume that (${\bf vii}$): the corresponding arc-length sequence $\big(||\gamma_{\varepsilon,n,i}^{(1)}||\big)_{n=1}^\infty$ is uniformly bounded, so that the arc-sequence $\big(\gamma_{\varepsilon,n,i}^{(1)}\big)_{n=1}^\infty$ has a convergent subsequence. It follows that (${\bf viii}$): the related arc-sequence $\big(\gamma_{\varepsilon,n,i}^*\big)_{n=1}^\infty$ also has a convergent subsequence, still indexed by $n\in N$, such that the arcs $\gamma_{\varepsilon,n,i}^*$ converge uniformly as $n\rightarrow\infty$ to a limit-arc $\gamma_{\varepsilon,i}$ joining the point $p_{\varepsilon,i}$ to $\{y=(-1)^i\,\infty\}$. Then by (${\bf iv}$) and (${\bf viii}$), we have (${\bf ix}$): $\gamma_{\varepsilon,i}\subset \dot{H}_{\varepsilon,i}(\Gamma_m)$ for any fixed $m\in N$. Since $m\in N$ can be chosen arbitrarily large in (${\bf ix}$), it follows that $\gamma_{\varepsilon,i}\subset\bigcap_{m=1}^\infty \dot{H}_{\varepsilon,i}(\Gamma_m)$, where one easily sees that $\bigcap_{m=1}^\infty\dot{H}_{\varepsilon,i}(\Gamma_m)=\bigcap_{m=1}^\infty\big\{p\in E_i(\Gamma_m):a_i(p)\,{\rm dist}\big(p,\Gamma_m\big)\geq\varepsilon\big\}=\big\{p\in E_i(\Gamma):a_i(p)\,{\rm dist}\big(p,\Gamma\big)\geq\varepsilon\big\}=\dot{H}_{\varepsilon,i}(\Gamma)$, and where $\Gamma\in{\rm X}$ is chosen such that $\bigcap_{n=1}^\infty E_i(\Gamma_n)=E_i(\Gamma)$). It follows that $p_{\varepsilon,i}\in{H}_{\varepsilon,i}(\Gamma)$, completing the proof of (${\bf ii}$) subject to Assumption (${\bf vii}$).
\vspace{.1in}

\noindent
It remains to prove that among the possible arc-sequences $\big(\gamma_{\varepsilon,n,i}\big)_{n=1}^\infty$ satisfying (${\bf iii}$) for each $n\in N$, there exists at least one sequence such that $\big(||\gamma_{\varepsilon,n,i}^{(1)}||\big)_{n=1}^\infty$ is uniformly bounded (Assumption (${\bf vii}$)). Given any point $p_{\varepsilon,n,i}\in\big(\Re\times[\kappa_{n,1},\kappa_{n,2}]\big)\cap{H}_{\varepsilon,i}(\Gamma_n)$, we use $\gamma_{\varepsilon,n,i}^{(1)}$ to denote a directed arc which first follows the straight-line-segment from the initial point $p_{\varepsilon,n,i}$ to a point $q_{\varepsilon,n,i}$ which is closest to $p_{\varepsilon,n,i}$ relative to the set $\partial H_{\varepsilon,i}(\Gamma_n)$, then proceeds from $q_{\varepsilon,n,i}$ along the arc $\partial{H}_{\varepsilon,i}(\Gamma_n)$ up to its terminal point $r_{\varepsilon,n,i}\in(\partial H_{\varepsilon,i}(\Gamma_n)\cap\{y=\kappa_{n,i}\})$, which is the first point of contact of the arc $\partial H_{\varepsilon,i}(\Gamma_n)$ with the line $\{y=\kappa_{n,i}\}$. We also use 
$\gamma_{\varepsilon,n,i}^{(2)}$ to denote the vertical line-segment joining its initial point $r_{\varepsilon,n,i}$ to $\{y=(-1)^i\,\infty\}$, and we use $\gamma_{\varepsilon,n,i}$ to denote the composite directed arc joining $p_{\varepsilon,n,i}$ to $\{y=(-1)^i\,\infty\}$ which first follows the arc $\gamma_{\varepsilon,n,i}^{(1)}$, then the arc $\gamma_{\varepsilon,n,i}^{(2)}$. Clearly for each $n\in N$, the composite directed arc $\gamma_{\varepsilon,n,i}$ satisfies (${\bf iii}$), and the length of the arc-segment $\gamma_{\varepsilon,n,i}^{(1)}$ cannot exceed the length of the arc $\partial H_{\varepsilon,i}(\Gamma_n)$ relative to one $P$-period (in $x$). The latter is uniformly bounded as $n\rightarrow\infty$, as we show in Lem. \ref{lem 3.5.2}, Eq. (\ref{eqn 3.5.4}). In particular, we have
\begin{equation}
\label{eqn 2.1.26a}
||\gamma_{\varepsilon,n,i}^{(1)}||\leq \big({\rm dist}\big(p_{\varepsilon,i},\partial{H}_{\varepsilon,i}(\Gamma_n)\big)+||\partial{H}_{\varepsilon,i}(\Gamma_n)||\big)
\end{equation}
$$\leq 2\,||\partial{H}_{\varepsilon,i}(\Gamma_n)||\leq 2C\,||\Gamma_n||\leq 2CL,$$
uniformly for all $n\in N$, where $L:=\sup_{n\in N}\,\{||\Gamma_n|\}$ and $C=C(\varepsilon)$ denotes a uniform upper bound for the Lipschitz constants in (\ref{eqn 2.1.25}). 
\begin{lemma} 
{(One-sided continuity of the operators ${\bf T}^\pm_{\varepsilon}({\bf\Gamma}):{\bf X}\rightarrow{\bf X}$)}
\label{lem 2.1.6}
Let be given $\varepsilon\in(0,\varepsilon_0)$, a weakly monotone sequence $\big({\bf\Gamma}_{\varepsilon,n}\big)_{n=1}^\infty$ of arc-pairs ${\bf\Gamma}_{\varepsilon,n}=(\Gamma_{\varepsilon,n,1},\Gamma_{\varepsilon,n,2})\in {\bf X}$, 
and their limit ${\bf\Gamma}_{\varepsilon}=(\Gamma_{\varepsilon,1},\Gamma_{\varepsilon,2})\in{\bf X}$.
Assume there exists a uniform upper bound for the sequences $\big(||\Gamma_{\varepsilon,n,i}||\big)_{n=1}^\infty$, $i=1,2,$ of arc-lengths per $P$-period (in $x$). Then we have: 
\begin{equation}
\label{eqn 2.1.26}
{\bf\Phi}_\varepsilon({\bf\Gamma}_{\varepsilon,n})\,\uparrow\,(\downarrow)\, {\bf\Phi}_\varepsilon({\bf\Gamma}_{\varepsilon})\,\,{\rm if}\,\,{\bf \Gamma}_{\varepsilon,n}\,\uparrow\,(\downarrow)\,{\bf\Gamma}_{\varepsilon},
\end{equation}
\begin{equation}
\label{eqn 2.1.27}
\Psi_{\varepsilon,i}^+\Gamma_{n,i})\downarrow\Psi_{\varepsilon,i}^+(\Gamma_i)\,\,{\rm if}\,\, \Gamma_{n,i}\downarrow\Gamma_i,\,\,{\rm for}\,\, i=1,2, 
\end{equation}
\begin{equation}
\label{eqn 2.1.28}
\Psi_{\varepsilon,i}^-(\Gamma_{n,i})\uparrow\Psi_{\varepsilon,i}^-(\Gamma_i)\,\,{\rm if}\,\, 
\Gamma_{n,i}\uparrow\Gamma_i,\,\,{\rm for}\,\,i=1,2, 
\end{equation}
all monotonically as $n\rightarrow\infty$. For the parametrized operator family ${\bf T}_{\varepsilon}^\pm(\Gamma):={\bf\Psi}_{\varepsilon}^\pm\big({\bf\Phi}_{\varepsilon}({\bf\Gamma})\big):{\bf X}\rightarrow{\bf X}$, it follows that
\begin{equation}
\label{eqn 2.1.29}
{\bf T}_{\varepsilon}^+({\bf\Gamma}_{\varepsilon,n})\,\downarrow\, {\bf T}_\varepsilon^+({\bf\Gamma}_{\varepsilon})\,\,{\rm for}\,\,{\bf\Gamma}_{\varepsilon,n}\,\downarrow\,{\bf\Gamma}_\varepsilon, 
\end{equation}
\begin{equation}
\label{eqn 2.1.30}
{\bf T}_{\varepsilon}^-({\bf\Gamma}_{\varepsilon,n})\,\uparrow\, {\bf T}_\varepsilon^-({\bf\Gamma}_{\varepsilon})\,\,{\rm for}\,\,{\bf\Gamma}_{\varepsilon,n}\uparrow{\bf\Gamma}_\varepsilon, 
\end{equation}
both monotonically as $n\rightarrow\infty$.
\end{lemma}

\noindent
{\bf Proof.} We have $$\Omega_{\varepsilon,n}:=\Omega(\Gamma_{\varepsilon,n,1},\Gamma_{\varepsilon,n,2})\uparrow\,(\downarrow)\,\Omega_\varepsilon:=\Omega(\Gamma_{\varepsilon,1},\Gamma_{\varepsilon,2})$$
monotonically as $n\rightarrow\infty$, in the sense that $\Gamma_{\varepsilon,n,1}\uparrow\,(\downarrow)\,\Gamma_{\varepsilon,1}$ and $\Gamma_{\varepsilon,n,2}\,\uparrow\,(\downarrow)\,\Gamma_{\varepsilon,2}$. Therefore $$U(\Gamma_{\varepsilon,n,1},\Gamma_{\varepsilon,n,2};p)\,\downarrow\,(\uparrow)\, U(\Gamma_{\varepsilon,1},\Gamma_{\varepsilon,2};p),$$ monotonically and uniformly in compact subsets of $\Omega_\varepsilon$ as $n\rightarrow\infty$, from which (\ref{eqn 2.1.26}) follows. Also, the assertions  (\ref{eqn 2.1.27}) and (\ref{eqn 2.1.28}) follow from Lem. \ref{lem 2.1.4}(c), and  the assertions (\ref{eqn 2.1.29}) and (\ref{eqn 2.1.30}) both follow from (\ref{eqn 2.1.26}), (\ref{eqn 2.1.27}), and (\ref{eqn 2.1.28}).

\subsection{Properties of operator fixed points}
\label{subsection 2.3}
\begin{lemma}
\label{lem 2.2.2}
{(Alternative characterization of the operator fixed points)}
For any $\varepsilon\in(0,\varepsilon_0)$, any solution ${\bm \Gamma}_\varepsilon=(\Gamma_{\varepsilon,1},\Gamma_{\varepsilon,2})\in{\bm{\cal F}}_\varepsilon$ of Prob. \ref{prob 2.1.2} at $\varepsilon$, and either $i=1$ or $i=2$, let $\tilde{\omega}_{\varepsilon,i}:=\{0<U_{\varepsilon,i}(p)<\varepsilon\}$ and $\tilde{\Gamma}_{\varepsilon,i}:=\{U_{\varepsilon,i}(p)=\varepsilon\}$, where $U_{\varepsilon,i}(p):=U_{i}({\bm\Gamma}_{\varepsilon};p)$ in the closure of $\Omega_{\varepsilon}:=\Omega({\bm\Gamma}_\varepsilon)$. Also, let $\hat{\Gamma}_{\varepsilon,i}:=(\partial \hat{\omega}_{\varepsilon,i})\cap\Omega_\varepsilon$, where $\hat{\omega}_{\varepsilon,i}$ denotes the intersection with ${\rm Cl}(\Omega_\varepsilon)$ of the union of the collection of all balls $G_{\varepsilon,i}(p):=\{q\in\Re^2:a_i(p)|q-p|<\varepsilon\}=B_{(\varepsilon/a_i(p))}(p)$ having center-points $p\in\Gamma_{\varepsilon,i}$. Here, we have that $\hat{\omega}_{\varepsilon,i}\subset\tilde{\omega}_{\varepsilon,i}$, due to the fact that
$$a_i(p)\,|p-q|<\varepsilon=a_i(p)\,{\rm dist}(p,\tilde{\Gamma}_{\varepsilon,i})$$ for all $p\in\Gamma_{\varepsilon,i}$ and $q\in G_{\varepsilon,i}(p)$, by Thm. \ref{thm 2.1.1}, Eq. (\ref{eqn 2.1.13}). 
Then for any point $p_{\varepsilon,i}\in\Gamma_{\varepsilon,i}$ and any corresponding point $\hat{q}_{\varepsilon,i}\in\hat{\Gamma}_{\varepsilon,i}\cap \partial G_{\varepsilon,i}(p_{\varepsilon,i})$, 
the generalized ball $H_{\varepsilon,i}(\hat{q}_{\varepsilon,i}):=\big\{p\in\Re^2: a_i(p)|p-\hat{q}_{\varepsilon,i}|<\varepsilon\big\}\subset \Omega_\varepsilon$ "centered" at the point $\hat{q}_{\varepsilon,i}$ is such that $p_{\varepsilon,i}\in\partial H_{\varepsilon,i}(\hat{q}_{\varepsilon,i})\cap\Gamma_{\varepsilon,i}$. 
\end{lemma}
\vspace{.1in}

\noindent
{\bf Proof.} For any particular point $q\in\hat{\Gamma}_{\varepsilon,i}$, we have that $q\notin G_{\varepsilon,i}(p)$ for any point $p\in\Gamma_{\varepsilon,i}$. Therefore, we have that $a_i(p)|q-p|\geq\varepsilon$ for all points $p\in\Gamma_{\varepsilon,i}$ and $q\in\hat{\Gamma}_{\varepsilon,i}$. Therefore, for any point $p\in\Gamma_{\varepsilon,i}$, we have that $p\notin H_{\varepsilon,i}(q)$ 
for any point $q\in\hat{\Gamma}_{\varepsilon,i}$. Since $q$ denotes any point in $\hat{\Gamma}_{\varepsilon,i}$, we conclude that (${\bf i}$):  $p\notin\bigcup_{q\in\hat{\Gamma}_{\varepsilon,i}}H_{\varepsilon,i}(q)$ for all points $p\in\Gamma_{\varepsilon,i}$.
\vspace{.1in}

\noindent
Also, for any point $p\in\Gamma_{\varepsilon,i}$, and any corresponding point $q=q(p)$ such that $q\in\hat{\Gamma}_{\varepsilon,i}\cap\partial G_{\varepsilon,i}(p)$, we have that $a_i(p)|p-q|\leq \varepsilon$ and therefore $p\in\partial H_{\varepsilon,i}(q)$. 
In view of (${\bf i}$), it follows that for any point $p\in\Gamma_{\varepsilon,i}$, we have $p\in \partial H_{\varepsilon,i}(q)$ for any generalized ball $H_{\varepsilon,i}(q)\subset\Omega_\varepsilon$ whose "center-point" $q$ is such that $q\in\partial G_{\varepsilon,i}(p)\cap\hat{\Gamma}_{\varepsilon,i}$. 
\begin{lemma}
\label{lem 2.2.1}
{(Preliminary estimates for "fixed points")} (a) For any $\varepsilon\in(0,\varepsilon_0)$ and any fixed point ${\bf\Gamma}_\varepsilon=(\Gamma_{\varepsilon,1},\Gamma_{\varepsilon,2})\in{\bm{\cal F}}_\varepsilon$, we have $|{\nabla} U_{\varepsilon,i}(p)|\leq 2\overline{A}$ throughout $\Omega_\varepsilon\setminus(\tilde{\omega}_{\varepsilon,1}\cup\tilde{\omega}_{\varepsilon,2})$, and even throughout the larger region $\Omega_\varepsilon\setminus(\dot{\omega}_{\varepsilon,1}\cup\dot{\omega}_{\varepsilon,2})$, $i=1,2$, where $U_{\varepsilon,i}(p):=U_i({\bf\Gamma}_\varepsilon;p)$ in the closure of $\Omega_\varepsilon:=\Omega({\bm\Gamma}_\varepsilon)$ and $\tilde{\omega}_{\varepsilon,i}:=\{p\in\Omega_\varepsilon: U_{\varepsilon,i}(p)<\varepsilon\}\supset\dot{\omega}_{\varepsilon,i}:=\{p\in\Omega_\varepsilon:{\rm dist}(p,\Gamma_{\varepsilon,i})<\big(\varepsilon\big/2\overline{A}\big)\}$.  
\vspace{.1in}

\noindent
(b) It follows from Part (a) that ${\rm dist}(p,\partial\Omega_\varepsilon)\geq((\alpha-\varepsilon)/2\overline{A}))$ for any $\alpha\in(\varepsilon,1/2]$ and any point $p\in\Omega_\varepsilon$ such that $\alpha\leq U_{\varepsilon,i}(p)\leq (1-\alpha)$, $i=1,2$, and therefore that ${\rm dist}(\Gamma_{\varepsilon,1},\Gamma_{\varepsilon,2})\geq((1-2\varepsilon)/\,2\overline{A})$ for $\varepsilon\in(0,\varepsilon_0)$.
\vspace{.1in}

\noindent
(c) It also follows that for any $\varepsilon\in(0,\varepsilon_0)$ and $\alpha\in(0,1)$, we have that $U_{\varepsilon,i}(p)\geq\alpha\varepsilon$ in the neighborhood $N_d(\Omega_\varepsilon\setminus\tilde{\omega}_{\varepsilon,i})$, where $d:=((1-\alpha)\varepsilon/2\overline{A})$. 
\end{lemma}

\noindent
{\bf Proof.} Let $\dot{\Gamma}_{\varepsilon,i}$ (resp. $\tilde{\Gamma}_{\varepsilon,i}$)) denotes the boundary of $\dot{\omega}_{\varepsilon,i}$ (resp. $\tilde{\omega}_{\varepsilon,i}$) relative to $\Omega_\varepsilon$.
We have ${\rm dist}(p,\Omega_\varepsilon\setminus\tilde{\omega}_{\varepsilon,i})=\big(\varepsilon\big/a_i(p)\big)\geq\big(\varepsilon\big/\,\overline{A}\,\big)$ for every point $p\in\Gamma_{\varepsilon,i}$, from which it follows that ${\rm dist}(p,\Gamma_{\varepsilon,i}\cup\tilde{\Gamma}_{\varepsilon,i})\geq(\varepsilon\big/2\overline{A}\,)$ for all $p\in\dot{\Gamma}_{\varepsilon,i}$. Since $0\leq U_{\varepsilon,i}\leq\varepsilon$ in $\tilde{\omega}_{\varepsilon,i}$, $i=1,2$, it follows by a well-known gradient estimate that $|{\nabla} U_{\varepsilon,i}(p)|\leq 2\overline{A}$ on $\dot{\Gamma}_{\varepsilon,i}$, for $i=1,2$, from which it follows by maximum principle for sub-harmonic functions that $|{\nabla} U_{\varepsilon}(p)|\leq 2\overline{A}$ for all points $p$ in $\Omega_\varepsilon\setminus(\dot{\omega}_{\varepsilon,1}\cup\,\dot{\omega}_{\varepsilon,2})$. 
\vspace{.1in}

\noindent
At this point, the estimates in Part (b) follows from Part (a) by integrating $|{\nabla} U_{\varepsilon,i}|$ on the straight line-segment $\gamma_{\varepsilon,i}$ joining any point $p\in\Omega_\varepsilon$ such that $U_{\varepsilon,i}(p)\geq\alpha$ to the closest point $q\in\Gamma_{\varepsilon,i}$. Turning to the estimate (c), it follows from (a) by integrating $|{\nabla} U_{\varepsilon,i}|$ on the straight line-segment $\gamma_{\varepsilon,i}$ joining any point $p\in {\rm Cl}(\tilde{\omega}_{\varepsilon,i})\cap N_d(\Omega_\varepsilon\setminus\tilde{\omega}_{\varepsilon,i})$ to a point $q\in \tilde{\Gamma}_{\varepsilon,i}$ that $U_{\varepsilon,i}(p)\geq U_{\varepsilon,i}(q)-\int_{\gamma_{\varepsilon,i}}|{\nabla} U_{\varepsilon,i}|\,ds\geq \varepsilon-2\overline{A}\,d\geq \varepsilon -(1-\alpha)\,\varepsilon=\alpha\,\varepsilon$, from which the assertion follows.
\begin{lemma}
\label{lem 2.2.2a} {(Convexity of level curves of the distance function)} 
In $\Re^2$, given a bounded convex set $H$, let the function $\phi(p):{\rm Cl}(H)\rightarrow\Re$ be defined such that $\phi(p)={\rm dist}\big(p,\partial H\big)$. Then for any value $0<\alpha\leq \phi_0:=\sup\{\phi(p):p\in H\}$, the subset $H_\alpha:=\{\phi(p)\geq\alpha\}$ of $H$ is convex.
\end{lemma}

\noindent
{\bf Proof.} Given values $\alpha\in(0,\phi_0)$, $\lambda\in(0,1)$, and points $p_0,p_1\in {\rm Cl}(H_\alpha)$, let $p_\lambda:=(1-\lambda)p_0+\lambda p_1\in {\rm Cl}(H)$, and choose the point $q_\lambda\in\partial H$ and the straight line $L_\lambda$ such that $|p_\lambda-q_\lambda|=\phi(p_\lambda)$, $q_\lambda\in L_\lambda$, and $H\cap L_\lambda=\emptyset$. Then $|p_\lambda-q_\lambda|=(1-\lambda)|p_0-q_0|+\lambda|p_1-q_1|$, where $q_0, q_1$ denote respectively the orthogonal projections of $p_0, p_1$ onto the straight line $L_\lambda$. Since $|p_i-q_i|\geq\phi(p_i)$, $i=0,1$, it follows that $\phi(p_\lambda)=|p_\lambda-q_\lambda|\geq(1-\lambda)\phi(p_0)+\lambda\phi(p_1)\geq\big((1-\lambda)\alpha+\lambda\alpha\big)=\alpha$.
Therefore, we have that $p_\lambda\in H_\alpha$ for all $\lambda\in(0,1)$, from which it follows that $H_\alpha$ is convex.
\begin{lemma}
\label{lem 2.2.3}
There exist constants $C_1^*,C_2^*>0$ and $\varepsilon_1\in(0,\varepsilon_0)$ such that for any value $\varepsilon\in(0,\varepsilon_1]$, any $i\in\{1,2\}$, any generalized ball $H_{\varepsilon,i}=H_{\varepsilon,i}(Q_0):=\{p\in\Re^2:a_i(p)|p-Q_0|<\varepsilon\}$ with "center-point" $Q_0\in\Re^2$, and any boundary point $p_{\varepsilon,i}\in\partial H_{\varepsilon,i}$, there exists a Euclidean ball $B\big(\hat{q}_{\varepsilon,i}\,; \hat{R}_{\varepsilon,i}\big)$
with center-point $\hat{q}_{\varepsilon,i}$ and radius $\hat{R}_{\varepsilon,i}:=\hat{R}_{\varepsilon,i}(Q_0):=\big(\varepsilon/a_i(Q_0)\big)(1-C_1^*\varepsilon)$ such that $B\big(\hat{q}_{\varepsilon,i}\,; \hat{R}_{\varepsilon,i}\big)\subset H_{\varepsilon,i}$, $p_{\varepsilon,i}\in\partial B\big(\hat{q}_{\varepsilon,i}\,; \hat{R}_{\varepsilon,i}\big)$, and $|\hat{q}_{\varepsilon,i}-Q_0|\leq C_2^*\varepsilon^2$.
\end{lemma}

\noindent
{\bf Proof.} We suppress the subscript $i$ (and sometimes $\varepsilon$). Let $K_{\varepsilon}(p)$ denotes the counter-clockwise oriented curvature at the point $p\in\partial H_\varepsilon$ of the arc $\partial H_{\varepsilon}$ corresponding to a counter-clockwise-oriented arc-length parametrization. By (\ref{eqn 3.3.12}) and ({\ref{eqn 3.3.20}), we have that $K_{\varepsilon}(p)\geq 0$ for all $p\in \partial H_\varepsilon$, and therefore that (${\bf i}$): the set $H_\varepsilon$ is convex, provided that the value $\varepsilon\in(0,\varepsilon_0)$ satisfies the additional requirement that: $C_5\varepsilon<1$, where the constant $C_5$ is defined in (\ref{eqn 3.3.20a}). Thus, we choose $\varepsilon_1:={\rm\min}\{\varepsilon_0,(1/C_5)\}=\min\{1/2,(\underline{A}^2/2A_1),$ $(1/C_5)\}$. It follows under assumptions (\ref{eqn 3.3.12}) and (\ref{eqn 3.3.20}) that (${\bf ii}$): 
$$\big(\varepsilon\,K_{\varepsilon}(p)\big/\big(a(Q_0)+(A_1/\underline{A})\varepsilon)\big)\leq\big(\varepsilon K_{\varepsilon}(p)\big/a(p)\big)=|p-Q_0|\,K_{\varepsilon}(p)\leq\, 1+C_5\varepsilon$$
for $\varepsilon\in(0,\varepsilon_1]$ and $p\in \partial H_{\varepsilon}$. A calculation based on (${\bf ii}$) shows that (${\bf iii}$): for any $\varepsilon\in(0,\varepsilon_1]$, the (point-wise) radius of counter-clockwise curvature of the (directed) boundary-arc $\partial H_{\varepsilon}$ is uniformly bounded from below by the positive value $\hat{R}_{\varepsilon}:=\big((\varepsilon\big/a(Q_0)\big)(1-C_1^*\varepsilon)>0$, where $C_1^*=\big((A_1/\underline{A}^2)+C_5\big)$). Then, in view of (${\bf i}$) and Lem. \ref{lem 2.2.2a}, the region $H_{\varepsilon}$ coincides with the $\hat{R}_{\varepsilon}$-neighborhood of the compact, convex set $\hat{H}_{\varepsilon}:=\big\{p\in H_{\varepsilon}:{\rm dist}\big(p,\partial
H_{\varepsilon}\big)\geq \hat{R}_{\varepsilon}\big\}$, and in fact (${\bf iv}$): for each point $p\in\partial H_{\varepsilon}$, there exists a point $\hat{q}\in\partial\hat{H}_{\varepsilon}$ such that $B\big(\hat{q}\,; \hat{R}_\varepsilon\big)\subset H_\varepsilon$ and $p\in\partial B\big(\hat{q}\,;\hat{R}_{\varepsilon}\big)$. 
Also, by (\ref{eqn 3.3.13}), we have that (${\bf v}$): $B\big(Q_0; {R}_{\varepsilon}^-\big)\subset H_{\varepsilon}\subset B\big(Q_0; {R}_{\varepsilon}^+\big)$ for any $\varepsilon\in(0,\varepsilon_0)$, where we set $R_{\varepsilon}^\pm:=\big(\varepsilon/a(Q_0)\big)(1\pm C_3^*\varepsilon)$ and $C_3^*:=(A_1/\underline{A}^2)$. By (${\bf iv}$) and (${\bf v}$), we have 
$B\big(\hat{q}\,; \hat{R}_\varepsilon\big)\subset H_\varepsilon\subset B\big(Q_0; R^+_\varepsilon\big)$, from which it follows for any point $\hat{q}\in\partial\hat{H}_\varepsilon$ (corresponding to a point $p\in\partial H_\varepsilon$) that
$\big|\hat{q}-Q_0\big|\leq\big(R_\varepsilon^+-\hat{R}_\varepsilon\big)\leq\big(\varepsilon\big/a(Q_0)\big)\big(C_1^*+C_3^*\big)\varepsilon$.
The assertions follow, where we set $C_2^*:=\big((C_1^*+C_3^*)/\underline{A}\big)$.
\begin{corollary}
\label{cor 2.2.2}
{(Uniform interior tangent balls)} There exist constants $C_i^*$, $i=1,2$ and $\varepsilon_1\in(0,\varepsilon_0)$ such that for any value $\varepsilon\in(0,\varepsilon_1]$, any solution ${\bm \Gamma}_\varepsilon=(\Gamma_{\varepsilon,1},\Gamma_{\varepsilon,2})\in{\bm{\cal F}}_\varepsilon$ of Prob. \ref{prob 2.1.2} at $\varepsilon$, and any point $p_{\varepsilon,i}\in\Gamma_{\varepsilon,i}$, there exist Euclidian balls $B_{\varepsilon,i}^\pm:=B\big(\hat{q}_{\varepsilon,i}^\pm\,; \hat{R}_{\varepsilon,i}^\pm\big)\subset \Omega_\varepsilon$ with radii $\hat{R}_{\varepsilon,i}^\pm:=\big(\varepsilon\big/a_i(Q_{\varepsilon,i}^\pm)\big)\big(1-C_1^*\varepsilon\big)$ and "center-points" $\hat{q}_{\varepsilon,i}^\pm$ such that 
$p_{\varepsilon,i}\in \partial B_{\varepsilon,i}^\pm$ and $|\hat{q}_{\varepsilon,i}^\pm-Q_{\varepsilon,i}^\pm|\leq C_2^*\varepsilon^2$, where $Q_{\varepsilon,i}^\pm\in\tilde{\Gamma}_{\varepsilon,i}$ denote the two extremal elements of the set $\tilde{\Gamma}_{\varepsilon,i}\cap\partial G_{\varepsilon,i}(p_{\varepsilon,i})$, relative to the natural ordering of the elements of the arc $\tilde{\Gamma}_{\varepsilon,i}$.
\end{corollary}
\begin{lemma}
\label{lem 2.2.4} 
{(Uniform lower bound for the gradient of the capacitary potent-ial in the "fixed point" problem)} There exists a constant $C:=\big(\,\,\underline{A}\,\big/\big(2\,{\rm ln}(8\overline{A}\big/\underline{A})\big)$ $>0$ such that $|{\nabla} U_{\varepsilon,i}(p)|$ $\geq C$, $i=1,2$, uniformly for all sufficiently small $\varepsilon\in(0,\varepsilon_0)$, all "fixed points" ${\bf\Gamma}_\varepsilon=(\Gamma_{\varepsilon,1},\Gamma_{\varepsilon,2})\in{\bm{\cal F}}_\varepsilon$, and all points $p\in\Omega_\varepsilon:=\Omega({\bf\Gamma}_\varepsilon)$, where $U_{\varepsilon,i}(p):=U_i({\bf\Gamma}_\varepsilon;p)$.
\end{lemma}

\noindent
{\bf Proof.} In the context of Lemmas \ref{lem 2.2.2} and \ref{lem 2.2.3} and Cor. \ref{cor 2.2.2}, given any points $p_{\varepsilon,i}\in\Gamma_{\varepsilon,i}$, for $\varepsilon\in(0,\varepsilon_1]$ and $i=1,2$, we define the harmonic capacitary potentials 
\begin{equation}
\label{eqn 2.2.X1}
V_{\varepsilon,i}^\pm(q):=\frac{\varepsilon\,\,{\rm ln}\big(\hat{R}_{\varepsilon,i}
^\pm\big/|q-\hat{q}_{\varepsilon,i}^\pm|\big)}{2\,\,{\rm ln}(1/\lambda)} 
\end{equation}
in the annuli $A_{\varepsilon,i}^\pm:=\big\{\lambda\hat{R}_{\varepsilon,i}^\pm<|q-\hat{q}_{\varepsilon,i}^\pm|<\hat{R}_{\varepsilon,i}^\pm\big\}\subset\Omega_\varepsilon$ tangent to $\Gamma_{\varepsilon,i}$ at $p_{\varepsilon,i}$, where we set $\lambda:=\big(\,\underline{A}\,\big/8\overline{A}\,\big)$. Then, for any $\varepsilon\in(0,\varepsilon_1]$ such that $8\overline{A}C_2^*\varepsilon\leq 1$, and for any point $q\in {\rm Cl}\big(B\big(\hat{q}_{\varepsilon,i}\,; \lambda\hat{R}_{\varepsilon,i}^\pm\big)\big)$, we have 
\begin{equation}
\label{eqn 2.2.X2}
{\rm dist}(q,\tilde{\Gamma}_{\varepsilon,i})\leq|q-\hat{q}_{\varepsilon,i}^\pm|+|\hat{q}_{\varepsilon,i}^\pm-Q_{\varepsilon,i}^\pm|\leq C_2^*\varepsilon^2+\lambda\hat{R}_{\varepsilon_i}^\pm
\end{equation}
$$\leq
\big(C_2^*\varepsilon^2+\big(\lambda\varepsilon\big/\underline{A}\,\big)\big)\leq
\big(C_2^*\varepsilon^2+\big(\varepsilon/8\overline{A})\big)\leq\big(\varepsilon\big/4\,\overline{A}\,\big),$$
from which it follows by Lem. \ref{lem 2.2.1}(c) (with $\alpha=1/2$) that $U_{\varepsilon,i}(p)\geq(\varepsilon/2)=V_{\varepsilon,i}^\pm(p)$ in $\partial B(\hat{q}_{\varepsilon,i}^\pm\,; \lambda\hat{R}_{\varepsilon,i}^\pm\big)\big)$. Since also $V_{\varepsilon,i}^\pm(q)=0\leq U_{\varepsilon,i}(q)$ for all $q\in\partial 
B(\hat{q}_{\varepsilon,i}^\pm\,;$ $\hat{R}_{\varepsilon,i}^\pm)$, it follows by the maximum principle for harmonic functions that $U_{\varepsilon,i}(q)\geq V_{\varepsilon,i}^\pm(q)$ throughout the annulus $A_{\varepsilon,i}^\pm$. Therefore, we have 
\begin{equation}
\label{eqn 2.2.X3}
|{\nabla} U_{\varepsilon,i}(p_{\varepsilon,i})|\geq|{\nabla} V_{\varepsilon,i}^\pm(p_{\varepsilon,i})|=\big(\varepsilon\big/2\,{\rm ln}(1/\lambda)\,\hat{R}_{\varepsilon,i}^\pm\big)\geq\big(\,\underline{A}\,\big/2\,{\rm ln}(1/\lambda)\big)
\end{equation}
at any point $p_{\varepsilon,i}\in\Gamma_{\varepsilon,i}$ at which the interior normal derivative $|{\nabla} U_{\varepsilon,i}(p_{\varepsilon,i})|$ exists. Under the assumption that  $\Omega_\varepsilon$ is sufficiently smooth, so that the function $\phi_{\varepsilon,i}(p):={\rm ln}\big(|{\nabla} U_{\varepsilon,i}(p)|\big):{\rm Cl}\big(\Omega_{\varepsilon}\big)\rightarrow\Re$ is continuous, it follows by the maximum principle for harmonic functions that $\phi_{\varepsilon,i}(p)\geq{\rm ln}\big(\,\underline{A}\,\big/2\,{\rm ln}(1/\lambda)\big)$ throughout $\Omega_\varepsilon$, as asserted. 
\vspace{.1in}

\noindent
To prove the above result without the above additional smoothness assumption, we define the periodic arc-pair ${\bm\Gamma}_{\varepsilon,\delta}:=(\Gamma_{\varepsilon,\delta,1},\Gamma_{\varepsilon,\delta,2})\in{\bf X}$ for any values $\varepsilon\in(0, \varepsilon_1]$ and $\delta\in(0,\varepsilon]$, where $\Gamma_{\varepsilon,\delta,1}$ (resp. $\Gamma_{\varepsilon,\delta,2}$) denotes the boundary of the union of all the discs of radius $\delta$ which are contained in the lower (resp. upper) complement $D^-(\Gamma_{\varepsilon,1})$ (resp. $D^+(\Gamma_{\varepsilon,2})$) of $\Gamma_{\varepsilon,1}$ (resp. $\Gamma_{\varepsilon,2}$). We also define the capacitary potentials $U_{\varepsilon,\delta,i}(p):=U_i({\bm\Gamma}_{\varepsilon,\delta};p)$ in the closures of the respective regions $\Omega_{\varepsilon,\delta}:=\Omega({\bm\Gamma}_{\varepsilon,\delta})$ for $\delta\in(0,\varepsilon]$ and $\varepsilon\in(0,\varepsilon_1]$, observing that: (${\bf i}$) if $\Omega_\varepsilon$ is the union of all its interior balls of fixed radius $\alpha>0$, then the same is true of $\Omega_{\varepsilon,\delta}$, (${\bf ii}$) the functions ${\nabla} U_{\varepsilon,\delta}(p):{\rm Cl}(\Omega_{\varepsilon,\delta})\rightarrow\Re^2$ are all continuous, and (${\bf iii}$) 
${\nabla} U_{\varepsilon,\delta}(p)\rightarrow{\nabla} U_{\varepsilon}(p)$ (in any compact subset of $\Omega_\varepsilon$) and $\Gamma_{\varepsilon,\delta}\rightarrow\Gamma_{\varepsilon}$ (in the sense of Hadamard distance between two sets), both as $\delta\rightarrow 0+$.
\vspace{.1in}

\noindent
In view of the properties (${\bf i}$), (${\bf ii}$), and (${\bf iii}$), one can replace ${\bm\Gamma}_\varepsilon$ and $U_\varepsilon$ by ${\bm\Gamma}_{\varepsilon,\delta}$ and $U_{\varepsilon,\delta}$ in the above argument to show first that $$\phi_{\varepsilon,\delta}(p):={\rm ln}\big(|{\nabla} U_{\varepsilon,\delta}(p)|\big)\geq{\rm ln}\big(\,\underline{A}\,\big/2\,{\rm ln}(1/\lambda)\big)-z(\delta),$$ first on $\Gamma_{\varepsilon,\delta,1}\cup\Gamma_{\varepsilon,\delta,2}$, then uniformly throughout $\Omega_{\varepsilon,\delta}$, where $z(\delta):\Re_+\rightarrow\Re_+$ denotes a fixed but arbitrary function such that $z(\delta)\rightarrow 0+$ as $\delta\rightarrow 0+$. The assertion then follows in the limit as $\delta\rightarrow 0+$. 
\begin{corollary}
\label{cor 2.2.1}
In the context of Lemma \ref{lem 2.2.4}, we have 
$$||\gamma_{\varepsilon,i}||\leq(|U_{\varepsilon,i}(q)-U_{\varepsilon,i}(p)|/C)$$
\vspace{.1in}

\noindent
for all sufficiently small $\varepsilon\in(0,\varepsilon_0)$ and for all arcs $\gamma_{\varepsilon,i}$ of steepest ascent of $U_{\varepsilon,i}$ such that $U_{\varepsilon,i}\leq (1/2)$ at both endpoints $p,q\in{\rm Cl}\big(\Omega_\varepsilon\big)$. 
\end{corollary}

\noindent
{\bf Proof.} By Lem. \ref{lem 2.2.4}, we have $C\,||\gamma_{\varepsilon,i}||\leq\int_{\gamma_{\varepsilon,i}}|{\nabla} U_{\varepsilon,i}(p)|ds=|U_{\varepsilon,i}(q)-U_{\varepsilon,i}(p)|$. 
\begin{lemma} 
\label{lem 2.2.5} 
Let ${\bf Y}:=\{{\bf\Gamma}\in{\bf X}:\tilde{{\bf \Gamma}}\leq{\bf\Gamma}\leq\hat{{\bf \Gamma}}\}$, where $\tilde{{\bf \Gamma}},\hat{{\bf \Gamma}}\in{\bf X}\cap C^{\,2}$ denote respective strict lower and upper solutions of Prob. \ref{prob2.1.1}. For any $\varepsilon\in(0,\varepsilon_1]$ (where  $\varepsilon_1=\varepsilon_1\big(\tilde{{\bf\Gamma}},\hat{{\bf\Gamma}})$, see Thm. \ref{thm 2.1.3}), let ${\bf \Gamma}_{\varepsilon}\in{\bf Y}\cap{\bm{\cal F}}_\varepsilon$ denote a "fixed point" of one of the operators ${\bf T}_\varepsilon:={\bf T}_\varepsilon^\pm$. Then, for any fixed $\alpha\in(\varepsilon,1/2]$, the arc-lengths $||\Gamma_{\varepsilon,\alpha,i}||$ and total curvatures $K(\Gamma_{\varepsilon,\alpha,i})$ per $P$-period of the curves $\Gamma_{\varepsilon,\alpha,i}:=\Phi_{\alpha,i}({\bf \Gamma}_\varepsilon)$, $i=1,2$, are uniformly bounded from above over all sufficiently small $\varepsilon\in(0,\varepsilon_1]$. Also, ${\rm capacity}(\Omega_\varepsilon)$, the capacity per $P$-period of the domain $\Omega_\varepsilon:=\Omega({\bf\Gamma}_\varepsilon)$, is uniformly bounded from above, independent of small $\varepsilon\in(0,\varepsilon_1]$.
\end{lemma}

\noindent
{\bf Proof.} For any fixed point ${\bf \Gamma}_\varepsilon\in {\bf Y}$, $\varepsilon\in(0,\varepsilon_1]$, we have $\Omega_\varepsilon:=\Omega({\bf\Gamma}_\varepsilon)\subset\Omega(\tilde{\bf\Gamma})\cup\Omega(\hat{\bf\Gamma})$, whence $||\Omega_\varepsilon||\leq M:=||\Omega(\tilde{\bf\Gamma})\cup\Omega(\hat{\bf\Gamma})||$, where $||\cdot||$ refers to the area of one $P$-period (in $x$) of a $P$-periodic region. On the other hand, by Lem. \ref{lem 2.2.1}(b), we have that $||\Omega_\varepsilon||\geq((\alpha-\varepsilon)/2\overline{A})||\Gamma_{\varepsilon,\alpha,i}||$ for $\alpha\in (\varepsilon,1/2]$, where $||\cdot||$ refers to the length of one $P$-period (in $x$) of a $P$-periodic arc. Therefore 
\begin{equation}
\label{eqn 2.2.X}
||\Gamma_{\varepsilon,\alpha,i}||\leq\big(2\overline{A}\,M\big/(\alpha-\varepsilon)\big),
\end{equation}
$${\rm capacity}(\Omega_\varepsilon)\leq\overline{A}||\Gamma_{\varepsilon,1/2,i}||=\big(4M\overline{A}^2\big/(1-2\varepsilon)\big)$$
for $\alpha\in (\varepsilon,1/2]$ and $i=1,2$. Let $\phi_{\varepsilon}(p):={\rm ln}\big(|{\nabla} U_\varepsilon(p)|\big)$. Then, by Lem. \ref{lem 2.2.1}(b) and Lem. \ref{lem 2.2.4}, there exists a constant $C$, independent of small $\varepsilon>0$, such that $|\phi_\varepsilon(p)|\leq C$ in $\tilde{\omega}_\varepsilon:=\{\varepsilon\leq U_\varepsilon(p)\leq 1-\varepsilon\}$. Moreover, we have ${\rm dist}(p,\partial\tilde{\omega}_\varepsilon)\geq ((\alpha-\varepsilon)/2\overline{A})$ for $p\in\tilde{\omega}_\alpha:=\{\alpha\leq U_\varepsilon(p)\leq 1-\alpha\}$, $\alpha\in(\varepsilon,1/2]$, by Lem. \ref{lem 2.2.1}(b). We conclude by the standard derivative estimate that $$|K_\varepsilon(p)|\leq|{\nabla}\phi_\varepsilon(p)|\leq (2C/{\rm dist}(p,\partial\tilde{\omega}_\varepsilon))\leq (4C\overline{A}/(\alpha-\varepsilon))$$ for $p\in\tilde{\omega}_\alpha$, $\alpha\in(\varepsilon,1/2]$.  
Therefore, $$K(\Gamma_{\varepsilon,\alpha,i})\leq ||\Gamma_{\varepsilon,\alpha,i}||\sup\{K_\varepsilon(p)\}\leq 8CM\overline{A}^2/(\alpha-\varepsilon)^2$$
for $\alpha\in(\varepsilon,1/2]$, where the sup is over all $p\in\tilde{\omega}_\alpha$.
\subsection{Weak solutions from fixed points}
\label{subsection 2.4}
\begin{definition}
\label{def 2.3.1} 
{(Weak solutions of Prob. \ref{prob2.1.1})} 
We use ${\boldsymbol{\cal F}}$ to denote the family of all weak solutions of Prob. \ref{prob2.1.1}. Here, a pair of directed arcs ${\bf\Gamma}=(\Gamma_1,\Gamma_2)\in{\bf X}$ is called a weak solution of Prob. \ref{prob2.1.1} if there exists a positive null-sequence $\big(\varepsilon_n\big)_{n=1}^\infty$ of values in $(0,\varepsilon_0)$ and a corresponding sequence of directed-arc-pairs ${\bm \Gamma}_n\in{\bm{\cal F}}_{\varepsilon_n}$ (see Def. \ref{def 2.1.6}(b)) such that as $n\rightarrow\infty$, we have ${\bf\Gamma}_n\rightarrow{\bf\Gamma}$, component-wise, (and $U_n(p)\rightarrow U(p)$ uniformly relative to any compact subset of $\Omega$, where we define $U_n(p):=U({\bm\Gamma}_n;p)$ (resp. $U(p):=U({\bm \Gamma};p)$) in the closure of the domain $\Omega_n:=\Omega({\bm\Gamma})$ (resp. $\Omega:=\Omega({\bm\Gamma})$)).
\end{definition}
\begin{theorem} {(Existence of weak solutions and their properties)}
\label{thm 2.3.1}
Given the fixed positive $C^2$-functions $a_1(p), a_2(p):\Re^2\rightarrow\Re_+$, let be given a null-sequence $\big(\varepsilon_n\big)_{n=1}^\infty$ of values in $(0,\varepsilon_0)$ and a corresponding sequence $\big({\bf\Gamma}_n\big)_{n=1}^\infty$ of solutions ${\bf\Gamma}_n=(\Gamma_{n,1},\Gamma_{n,2})\in{\bm{\cal F}}_{\varepsilon_n}$ of Prob. \ref{prob 2.1.2} at $\varepsilon_n$ (see Thms. \ref{thm 2.1.2} and \ref{thm 2.1.3}). Then:
\vspace{.1in}

\noindent
(a) There exists a weak solution ${\bf\Gamma}\in{\boldsymbol{\cal F}}$ of Prob. \ref{prob2.1.1}, in the form of a directed-arc-pair ${\bf\Gamma}=(\Gamma_1,\Gamma_2)\in{\bf X}$, parametrized component-wise by the Lipschitz-continuous mappings $p_i(t):\Re\rightarrow\Gamma_i$, $i=1,2$, and a suitable subsequence $\big({\bf\Gamma}_{n(k)}\big)_{k=1}^\infty$ of directed arcs such that ${\bf\Gamma}_{n(k)}\rightarrow{\bf\Gamma}$ as $k\rightarrow\infty$.
\vspace{.1in}

\noindent
(b) We have $D^\pm(\Gamma_{n(k),i})\rightarrow D^\pm(\Gamma_i)$, $i=1,2$, and $\Omega({\bf\Gamma}_{n(k)})\rightarrow \Omega({\bf\Gamma}):=D^+(\Gamma_1)\cap D^-(\Gamma_2)$, all as $k\rightarrow\infty$, where $D^\pm(\Gamma_i):=\big\{W(\Gamma_i;p)=\pm 1/2\big\}$. Thus, $U_i({\bf\Gamma}_{n(k)};p)\rightarrow U_i({\bf\Gamma};p)$ in $\Omega({\bf\Gamma})$, $i=1,2$, as $k\rightarrow\infty$.
\vspace{.1in}

\noindent
(c) The pair ${\bf\Gamma}$ is in ${\bf X}$. In particular, for each $i=1,2$, the arc $\Gamma_i$ is $P$-periodic (in $x$). Also ${\rm dist}(\Gamma_1,\Gamma_2)\geq (1/\,\overline{A})$. Finally, for each $i=1,2$, the complement of $\Gamma_i$ can be partitioned into the sets $D^\pm(\Gamma_i)$.
\vspace{.1in}

\noindent
(d) The component-wise arc-length and total curvature per $P$-period of ${\bf\Gamma}$ are bounded from above. The bound depends only on the constants $\overline{A}, \underline{A}, A_1,A_2.$
\vspace{.1in}

\noindent
(e) There exist positive constants $0<C_0\leq C_1$ such that $C_0\leq |\nabla U({\bf\Gamma};p)|\leq C_1$ uniformly in $\Omega:=\Omega({\bf\Gamma})$. 
\vspace{.1in}

\noindent
(f) For $i=1,2,$ and any $t_0\in\Re$, there exists a value $\delta>0$ such that $p_i(t)\not=p_i(t_0)$ for all $t\in(t_0-\delta,t_0+\delta)$
\end{theorem}
\begin{corollary}
\label{cor 2.3.1}
In the context (and under the assumptions) of Def. \ref{def 2.1.5}, in view of Thms. \ref{thm 2.1.3} and \ref{thm 2.3.1}, there exists at least one weak solution ${\bm\Gamma}\in{\bf Y}$ of Prob. \ref{prob2.1.1} with all the properties listed in Thm. \ref{thm 2.3.1}.
\end{corollary}

\noindent
{\bf Proof of Thm. \ref{thm 2.3.1}.}
Concerning Part (a), for each $n\in N$, we let the $P$-periodic (in $x$) directed arcs $\Gamma_{n,i}$, $i=1,2$, be images of the arc-length parametrizations $p_{n,i}(s):\Re\rightarrow\Gamma_{n,i}$, $i=1,2$. Then $|p_{n,i}(s_1)-p_{n,i}(s_2)|\leq|s_1-s_2|$ for all $s_1, s_2\in\Re$, $i=1,2$, and $n\in N$. Each arc $\Gamma_{n,i}$ is also generated by a mapping $q_{n,i}(t):\Re\rightarrow\Gamma_{n,i}$ such that $q_{n,i}(t):=p_{n,i}(L_{n,i}t)$, where $L_{n,i}$ denotes the arc-length of one $P$-period (in $x$) of $\Gamma_{n,i}$. We have $|q_{n,i}(t_1)-q_{n,i}(t_2)|\leq L_{n,i}|t_1-t_2|$ for $t_1,t_2\in\Re$ and $n\in N$. By Thm. \ref{thm 2.1.3}(b) (which follows from Thm. \ref{thm 3.1.2}(c) and Lem. \ref{lem 2.2.5}, Eq. (\ref{eqn 2.2.X})), the sequence $\big(L_{n,i}\big)_{n=1}^\infty$ is uniformly bounded from above. By the Theorem of Ascoli and Arzela, there exist a suitable subsequence (still indexed by $n$), a pair of values $L_i$, $i=1,2$, and a pair of Lipschitz-continuous functions $q_i(t):\Re\rightarrow \Re^2$, $i=1,2$, such that $L_{n,i}\rightarrow L_i$ and ${\rm sup}\{|q_{n,i}(t)-q_i(t)|:t\in\Re\}\rightarrow 0+$ both as $n\rightarrow\infty$, and, therefore, such that $|q_i(t_1)-q_i(t_2)|\leq L_i|t_1-t_2|$ for all $t_1,t_2\in\Re$. Clearly, the mappings $q_i(t):\Re\rightarrow\Re^2$, $i=1,2$, are Lipschitz-continuous parametrizations of $P$-periodic (in $x$) arcs $\Gamma_i$, $i=1,2$. The same directed arcs are also generated by the parametrizations $p_i(s)=q_i(s/L_i):\Re\rightarrow\Re^2$, $i=1,2$, which are such that $|p_i(s_1)-p_i(s_2)|\leq |s_1-s_2|$ for all $s_1,s_2\in\Re$. 
\vspace{.1in}

\noindent
Part (b) follows from Part (a). The properties of ${\bf\Gamma}$ in Part (c) follow in the limit from essentially the same properties of the fixed-points ${\bf\Gamma}_{n(k)}\in{\bf X}$, in particular, the fact that the arcs  $\Gamma_{n(k),i}$ are all $P$-periodic (in $x$), that ${\rm dist}(\Gamma_{n(k),1},\Gamma_{n(k),2})\geq ((1-2\varepsilon_{n(k)})/\,\overline{A})$ (see Lem. \ref{lem 2.2.1}(b)), and finally that the complement of $\Gamma_{i,k}$ can be partitioned into the sets $D^\pm(\Gamma_{n(k),i}):=\big\{p\notin\Gamma_{n(k),i}: W(\Gamma_{n(k),i};p)=\pm 1/2\big\}$. Therefore ${\bf\Gamma}\in {\bf X}$. Part (d) follows in the limit as $k\rightarrow\infty$ from 
Thm. \ref{thm 2.1.3}(b), and Part (e) follows in the limit from Lem. \ref{lem 2.2.4}. Finally, if the assertion in Part (f) is false for some $i\in\{1,2\}$ and $t_0\in\Re$, then there exists a monotone sequence $(t_m)_{m=1}^\infty$ such that $t_m\rightarrow t_0$ as $m\rightarrow\infty$, and such that $p(t_m)=p(t_0)$ and $p(t_m^*)\not=p(t_0)$ both for all $m$, where $t_m^*$ denotes some value between $t_m$ and $t_{m+1}$. It is easily seen that for each $m$, the total curvature of the arc-segment of $\Gamma_i$ corresponding to the interval between $t_m$ and $t_{m+1}$ must exceed $\pi$. Therefore, the total curvature of $\Gamma_i$ is infinite, contradicting Part (d).
\section{Existence of classical solutions}
\label{section 3}
\subsection{Sharp upper bound for $|\nabla U(p)|$ on $\Gamma$}
\label{subsection 3.1}
\begin{theorem} 
{(Sharp point-wise upper bound for the boundary gradient of the capacitary potential)}
\label{thm 2.4.1}
Let $U(p):=U({\bf\Gamma};p)$ denote the capacitary potential in the closure of $\Omega:=\Omega({\bf\Gamma})$, where ${\bf\Gamma}=(\Gamma_1,\Gamma_2)\in{\boldsymbol{\cal F}}$ denotes a weak solution of Prob. \ref{prob2.1.1} (see Def. \ref{def 2.3.1}) Then
\begin{equation}
\label{eqn 2.4.1}
\limsup_{p\rightarrow p_i}|{\nabla} U(p)|\leq a_i(p_i), 
\end{equation}
for $i=1,2$, where $p\in\Omega$ and $p_i$ denotes any point in $\Gamma_i$. 
In fact let $\hat{\Omega}_i:=\{p\in\Omega_i:U_i(p)<1/2\}$ and let $\hat{a}_i(p)$ denotes the unique continuous function in ${\rm Cl}(\hat{\Omega}_i)$ such that $\hat{a}_i(p)=a_i(p)$ on $\partial\hat{\Omega}_i$ and ${\rm ln}(\hat{a}_i(p))$ is harmonic in $\hat{\Omega}_i$, then there exists a constant $C_i:=\sup\{{\rm ln}(|{\nabla} U_i(p)|/\hat{a}_i(p)):p\in\hat{\Gamma}_i\}$ (where $\hat{\Gamma}_i:=\{U_i(p)=1/2\}$), such that
\begin{equation}
\label{eqn 2.4.2}
|{\nabla} U_i(p)|\leq E_i(p):=\hat{a}_i(p)\,{\rm exp}\big(2C_i\,U_i(p)\big),
\end{equation}
uniformly in $\hat{\Omega}_i$ for $i=1,2$. Observe that $E_i(p)$ is a continuous function in ${\rm Cl}\big(\hat{\Omega}_i\big)$ such that $E_i(p)=a_i(p)$ on $\Gamma_i$.
\end{theorem}
\begin{lemma}
\label{lem 2.4.1}
For any $\varepsilon\in(0,\varepsilon_0)$, let ${\bf\Gamma}_\varepsilon=(\Gamma_{\varepsilon,1},\Gamma_{\varepsilon,2})\in{\bm{\cal F}}_\varepsilon$ denote a fixed point of the operator ${\bf T}_\varepsilon$ (by which we denote either ${\bf T}_\varepsilon^+$ or ${\bf T}_\varepsilon^-$). For $i=1,2$, let $U_{\varepsilon,i}(p):=U_i({\bf\Gamma}_\varepsilon;p)$ denote the corresponding capacitary potential in (the closure of) the domain $\Omega_\varepsilon:=\Omega({\bf\Gamma}_\varepsilon)$. Also let $\hat{\Omega}_{\varepsilon,i}:=\{p\in\Omega_{\varepsilon,i}:U_{\varepsilon,i}(p)<1/2\}$, and let $\hat{a}_{\varepsilon,i}(p):{\rm Cl}(\hat{\Omega}_{\varepsilon,i})\rightarrow\Re$, $i=1,2$, denote the unique continuous function such that $\hat{a}_{\varepsilon,i}(p)=a_i(p)$ on $\partial\hat{\Omega}_{\varepsilon,i}$ and ${\rm ln}(\hat{a}_{\varepsilon,i}(p))$ is harmonic in $\hat{\Omega}_{\varepsilon,i}$. Then there exists a uniform constant $C_{\varepsilon,i}$ such that
\begin{equation}
\label{eqn 2.4.3}
(\varepsilon/d_{\varepsilon,i}(p))\,\leq\, E_{\varepsilon,i}(p):=\hat{a}_{\varepsilon,i}(p)\,{\rm exp}\big(2C_{\varepsilon,i}\,U_{\varepsilon,i}(p)\big)
\end{equation}
in $\hat{\Omega}_{\varepsilon,i}$ for $i=1,2$ and for all $\varepsilon\in(0,\varepsilon_0)$, where for each point $p\in\Omega_\varepsilon$ such that $U_{\varepsilon,i}(p)\leq 1-\varepsilon$, where $d_{\varepsilon,i}(p)$ denotes the distance between $p$ and the level curve of the function $U_{\varepsilon,i}$ at altitude $U_{\varepsilon,i}(p)+\varepsilon$.
\end{lemma}

\noindent
{\bf Proof}}. We choose $i\in\{1,2\}$ and suppress the index $i$, so that $\Gamma_{\varepsilon,i}$, $a_{\varepsilon,i}(p)$, $U_{\varepsilon,i}(p)$, etc. become $\Gamma_\varepsilon$, $\hat{a}_\varepsilon(p)$, $U_\varepsilon(p)$, etc. For $j=\sqrt{-1}$, let $z=x+jy=g_\varepsilon(W)=g_\varepsilon(U+jV)$ denote the analytic, periodic mapping of the infinite strip $R:=(0,1)\times\Re$ onto the $P$-periodic domain $\Omega_\varepsilon$, whose inverse is the analytic function $W=f_\varepsilon(z):=U_\varepsilon(z)+jV_\varepsilon(z)$, where $U_\varepsilon(x,y)$ is the capacitary potential in $\Omega_\varepsilon$ and $V_\varepsilon(x,y)$ is the harmonic conjugate of $U_\varepsilon$. For any value $h\in\Re$, the function ${\rm ln}\left(\varepsilon/(g_\varepsilon(W+\varepsilon +jh)-g_\varepsilon(W)\right)$ is an analytic function of $W$ in strip $\hat{R}_\varepsilon:=\{W\in \Re^2:0<U_\varepsilon<1/2\}$. Therefore, ${\rm ln}\big(\varepsilon\big/(g_\varepsilon(f_\varepsilon(z)+\varepsilon+jh)-z)\big)$ is an analytic function of $z$ in $\hat{\Omega}_\varepsilon:=\{z\in\Omega_\varepsilon:U_\varepsilon(z)<1/2\}$. Thus, $E_{\varepsilon,h}(z):={\rm ln}\big(\varepsilon\big/|g_\varepsilon(f_\varepsilon(z)+\varepsilon+jh)-z|\big)$ is a harmonic function of $z$ in $\hat{\Omega}_\varepsilon$. Thus the function $$E_\varepsilon(z):={\rm sup}_{h\in\Re}E_{\varepsilon,h}(z)={\rm ln}\big(\varepsilon\big/({\rm inf}_{h\in\Re}|(g_\varepsilon(f_\varepsilon(z)+\varepsilon+jh)-z|\big)$$
is sub-harmonic in $\hat{\Omega}_\varepsilon$, where, for any value $\mu\in[0,1-\varepsilon]$, and any point $z\in\Gamma_{\varepsilon,\mu}:=\{z\in\Omega_\varepsilon:U_\varepsilon(z)=\mu\}$, we have
$${\rm inf}_{h\in\Re}\left|(g_\varepsilon(f_\varepsilon(z)+\varepsilon+jh)-z\right|={\rm dist}(z,\Gamma_{\varepsilon,\mu+\varepsilon}).$$
In equivalent notation, the function
$E_\varepsilon(p):={\rm ln}\big(\varepsilon\big/d_\varepsilon(p)\big)$ is sub-harmonic in $p=(x,y)\in\hat{\Omega}_\varepsilon$, where $d_\varepsilon(p):={\rm dist}(p,\Gamma_{\varepsilon,\mu+\varepsilon})$ for all $p\in\Gamma_{\varepsilon,\mu}$. Thus, the function $F_\varepsilon(p):={\rm ln}(\varepsilon/\hat{a}_\varepsilon(p)d_\varepsilon(p))$ is sub-harmonic in $\hat{\Omega}_\varepsilon$, and satisfies $F_\varepsilon(\Gamma_\varepsilon)=0$. For $\varepsilon\in(0,\varepsilon_0]$, it follows by the maximum principle that 
\begin{equation}
\label{eqn 2.4.4}
F_\varepsilon(p)\,\leq \,2C_\varepsilon\,U_\varepsilon(p)
\end{equation}
in $\hat{\Omega}_\varepsilon$, where $C_\varepsilon:={\rm sup}\{F_\varepsilon(p):p\in\hat{\Gamma}_\varepsilon\}$, and where $\hat{\Gamma}_\varepsilon:=\{U_\varepsilon(p)=1/2\}$. The claim follows from (\ref{eqn 2.4.4}) by observing that the constants $C_\varepsilon$ have a uniform upper bound independent of small $\varepsilon\in (0,\varepsilon_0)$.
\vspace{.1in}

\noindent
{\bf Proof of Thm. \ref{thm 2.4.1}.} 
Given the weak solution ${\bf\Gamma}\in{\boldsymbol{\cal F}}$, let ${\bf\Gamma}_\varepsilon=(\Gamma_{\varepsilon,1},\Gamma_{\varepsilon,2})\in{\bm{\cal F}}_{\varepsilon}$, $\varepsilon\in S\subset(0,\varepsilon_0),$ be an approximating sequence of fixed points.
For $i=1,2$, let $U_{\varepsilon,i}(p):=U_i({\bf\Gamma}_\varepsilon;p)$ denote the suitable capacitary potential in (the closure of) $\Omega_\varepsilon:=\Omega({\bf\Gamma}_\varepsilon)$. Then $U_{\varepsilon,i}(p)\rightarrow U_i(p)$ and $\hat{a}_{\varepsilon,i}(p)\rightarrow \hat{a}_i(p)$, $i=1,2$, all as $\varepsilon\rightarrow 0+$ relative to $S$, where $\hat{a}_{\varepsilon,i}(p)$ is defined in Lem. \ref{lem 2.4.1}, and where, in each case, the convergence is uniform relative to any compact subset of $\Omega$. The same applies with $U_{\varepsilon,i}$ or $\hat{a}_{\varepsilon,i}$, $i=1,2$, replaced by any fixed derivative of the same function, as follows by a standard derivative estimate for harmonic functions. Of course it also follows that $C_{\varepsilon,i}\rightarrow C_i:=\sup\{{\rm ln}(|{\nabla} U_i(p)|/\hat{a}_i(p)):p\in\hat{\Gamma}_i\}$ and $E_{\varepsilon,i}(p)\rightarrow E_i(p)$, both as $\varepsilon\rightarrow 0+$ relative to $S$, where $\hat{\Gamma}_i:=\{U_i(p)=1/2\}$, and the convergence is uniform in any compact subset of $\Omega$. For any small $\varepsilon\in S$ and any point $p\in\Omega_\varepsilon$ with $U_{\varepsilon,i}(p)\leq 1-\varepsilon$, let $\gamma_{\varepsilon,i}(p)$ denote the unique arc of steepest ascent of the function $U_{\varepsilon,i}$ which joins $p$ to the level curve of $U_{\varepsilon,i}$ at altitude $U_{\varepsilon,i}(p)+\varepsilon$. By Lem. \ref{lem 2.4.1}, we have 
\begin{equation}
\label{eqn 2.4.5}
||\gamma_{\varepsilon,i}(p)||\geq d_{\varepsilon,i}(p)\geq(\varepsilon/E_{\varepsilon,i}(p)),
\end{equation}
in $\{U_{\varepsilon,i}(p)\leq 1/2\}$. In view of this (and the identity $\int_{\gamma_{\varepsilon,i}(p)}|{\nabla} U_{\varepsilon,i}(p)|ds$ $=\varepsilon$), we have
\begin{equation}
\label{eqn 2.4.6}
\inf \{|{\nabla} U_{\varepsilon,i}(q)|:q\in\gamma_{\varepsilon,i}(p)\}\leq\frac{\varepsilon}{||\gamma_{\varepsilon,i}(p)||}\leq E_{\varepsilon,i}(p)
\end{equation}
for $\varepsilon\in S$, where, by Cor. \ref{cor 2.2.1}, we have $||\gamma_{\varepsilon,i}(p)||\leq(\varepsilon/C_i)$ uniformly for $p\in\hat{\Omega}_\varepsilon$ and all small $\varepsilon\in S$. The functions ${\nabla} U_{\varepsilon,i}(p)$, with $\varepsilon$ sufficiently small in $S$, are equi-continuous relative to any compact subset $\mathfrak{K}$ of $\Omega$. In view of this, it follows from (\ref{eqn 2.4.6}) that
\begin{equation}
\label{eqn 2.4.7}
|{\nabla} U_{\varepsilon,i}(p)|\leq E_{\varepsilon,i}(p)+z_\mathfrak{K}(\varepsilon)
\end{equation}
in any compact subset $\mathfrak{K}$ of $\hat{\Omega}_i:=\{U_i(p)<1/2\}$, where $z_\mathfrak{K}(\varepsilon)\rightarrow 0+$ as $\varepsilon\rightarrow 0+$ relative to $S$. The assertion now follows in the limit as $\varepsilon\rightarrow 0+$.

\begin{corollary}
\label{cor 2.4.1aaa}
In the notation of the proof of Lem. \ref{lem 2.4.1}, 
the function $E_{\varepsilon,0}(z):={\rm ln}\big(\varepsilon\big/\big|g_\varepsilon(f_\varepsilon(z)+\varepsilon)-z\big|\big)$ is harmonic for sufficiently small $\varepsilon\in S$. It follows that the continuous function
\begin{equation}
\label{eqn 2.6.11}
\phi_{\varepsilon,i}(p):={\rm ln}\big(\varepsilon\big/|\Pi_{\varepsilon,i}(p)-p|\big):{\rm Cl}(\dot{\Omega}_{\varepsilon,i})\rightarrow\Re
\end{equation}
is harmonic in the domain $\dot{\Omega}_{\varepsilon,i}:=\{p\in\Omega_\varepsilon:
U_{\varepsilon,i}(p)\leq 1-\varepsilon\}$ for any sufficiently small $\varepsilon\in S$, where for each point $p\in\dot{\Omega}_{\varepsilon,i}$, the notation $\Pi_{\varepsilon,i}(p)$ refers to the point of intersection of the arc of steepest ascent through $p$ of the function $U_{\varepsilon,i}$ with the level curve of the function $U_{\varepsilon,i}$ at the altitude $U_{\varepsilon,i}(p)+\varepsilon$.
\end{corollary}
\subsection{Estimates related to total-curvature bounds}
\label{subsection 3.2}
\begin{lemma}
\label{lem 2.5.1}
(a) Let $\Gamma_\lambda,\lambda\in I:=[a,b]$, denote a weakly monotone increasing, Lipschitz-continuously-varying (in $\lambda$), uniformly-$C^{1,1}$ family of arcs in ${\rm X}$ (Lipschitz-continuously varying with Lipschitz constant $M$ in the sense that $|p_\alpha(s)-p_\beta(s)|\leq M\,|\alpha-\beta|$ for all $s\in\Re$ and $\alpha,\beta\in I$, where $p_\lambda(s):\Re\rightarrow\Gamma_\lambda$ denotes the arc-length parametrization of $\Gamma_\lambda$. Assume that $K_\lambda\leq\overline{K}<\infty$ for all $\lambda\in I$, where $K_\lambda$ denotes the total curvature per $P$-period (in $x$) of the arc $\Gamma_\lambda$. Let $V:\Omega\rightarrow\Re$ denote the function such that $V=\lambda$ on $\Gamma_\lambda$, where $\Omega=\cup_{\lambda\in I}\Gamma_\lambda$. Then:
\begin{equation}
\label{eqn 2.5.1}
\big|||\Gamma_{\alpha}^{\,\prime}||-||\Gamma_{\beta}^{\,\prime}||\big|\leq 2M\overline{K}\,|\alpha-\beta|
\end{equation}
for any $\alpha,\beta\in I$, where where $\Gamma_\alpha^{\,\prime}$, $\Gamma_\beta^{\,\prime}$ denote arc-segments of $\Gamma_\alpha$, $\Gamma_\beta$, resp., spanning not more than one $P$-period in $x$, such that $\Gamma_\alpha^{\,\prime}$ projects onto $\Gamma_\beta^{\,\prime}$ (and visa verse) along the arcs of steepest ascent of the function $V:\Omega\rightarrow\Re$, and where, for non-infinite arcs, we use $||\cdot||$ to denote arc-length.
\vspace{.1in}

\noindent
(b) For any $\varepsilon\in(0,\varepsilon_0)$, and any "fixed point" ${\bf\Gamma}_\varepsilon=(\Gamma_{\varepsilon,1},\Gamma_{\varepsilon,2})\in {\bm {\cal F}}_{\varepsilon}$, let $\Gamma_{\varepsilon,\lambda,i}:=\Phi_{\lambda,i}({\bf\Gamma}_{\varepsilon})=\big\{U_{\varepsilon,i}(p)=\lambda\big\}\in X$ for $\lambda\in[0,1]$ and $i=1,2$, where $U_{\varepsilon,i}(p):=U_i({\bf\Gamma}_{\varepsilon};p)$ in the closure of the domain $\Omega_\varepsilon:=\Omega({\bf\Gamma}_\varepsilon)$. Then it follows directly from Part (a) that
\begin{equation}
\label{eqn 2.5.9}
\big|\,||\Gamma_{\varepsilon,\alpha,i}^{\,\prime}||-||\Gamma^{\,\prime}_{\varepsilon,\beta,i}||\big|\,\leq 2M_\varepsilon\overline{K}\,\big|\beta-\alpha\big|,
\end{equation}
for any $\alpha,\beta\in[0,1]$, where $\Gamma^{\,\prime}_{\varepsilon,\alpha,i}$ and $\Gamma^{\,\prime}_{\varepsilon,\beta,i}$ denote any arc-segments of $\Gamma_{\varepsilon,\alpha,i}$ and $\Gamma_{\varepsilon,\beta,i}$, respectively, which span not more than one $P$-period, and which project onto each other along the arcs of steepest descent of the function $U_\varepsilon(p):=U({\bf\Gamma}_\varepsilon;p)$. Here, $M_\varepsilon$ denotes the Lipschitz constant for the variation of the arc $\Gamma_{\varepsilon,\lambda,i}$ as a function of $\lambda$. For small $\varepsilon>0$, one can choose $M_\varepsilon=(2/C)$, where $C>0$ denotes a positive uniform lower bound for $|{\nabla} U_\varepsilon(p)|$ in $\Omega_\varepsilon$ (see Lem. \ref{lem 2.2.4}).
\vspace{.1in}

\noindent
(c) In the context of Part (b), let $\tilde{\omega}_{\varepsilon,i}:=\big\{U_{\varepsilon,i}(p)<\varepsilon\big\}$ (with boundary components $\tilde{\Gamma}_{\varepsilon,i}:=\big\{U_{\varepsilon,i}(p)=\varepsilon\big\}$ and 
$\Gamma_{\varepsilon,i}:=\big\{U_{\varepsilon,i}(p)=0\big\}$).  Also let $\hat{\omega}_{\varepsilon,i}:=\big\{p\in \tilde{\omega}_{\varepsilon,i}:|p-q|<{\rm dist}(q,\tilde{\Gamma}_{\varepsilon,i})\,\,{\rm for}\,\,{\rm some}\,\,{\rm point}\,\,q\in{\Gamma}_{\varepsilon,i}\big\}$, and let $\hat{\Gamma}_{\varepsilon,i}$ denote the boundary of $\hat{\omega}_{\varepsilon,i}$ relative to $\Omega_\varepsilon$ (i.e. $\hat{\Gamma}_{\varepsilon,i}:=\partial\hat{\omega}_{\varepsilon,i}\setminus{\Gamma}_{\varepsilon,i}$). Finally, let $\tilde{\Gamma}_{\varepsilon,i}^{\,\prime}$ and $\hat{\Gamma}_{\varepsilon,i}^{\,\prime}$ denote any arc-segments of $\tilde{\Gamma}_{\varepsilon,i}$ and $\hat{\Gamma}_{\varepsilon,i}$, respectively, such that both arc-segments have the same initial and terminal points. Then:
\begin{equation}
\label{eqn 2.5.9a}
||\hat{\Gamma}_{\varepsilon,i}^{\,\prime}||\leq ||\tilde{\Gamma}_{\varepsilon,i}^{\,\prime}||\,\, {\rm and}\,\, K(\hat{\Gamma}_{\varepsilon,i}^{\,\prime})\leq K(\tilde{\Gamma}_{\varepsilon,i}^{\,\prime}),
\end{equation} 
both for any $\varepsilon\in(0,\varepsilon_0]$ and $i=1,2$. 
\vspace{.1in}

\noindent
(d) In the context of Parts (b) and (c), there exists a constant $C_0>0$ such that
\begin{equation}
\label{eqn 2.5.3}
0\leq||\tilde{\Gamma}_{\varepsilon,i}^{\,\prime}||-||\hat{\Gamma}_{\varepsilon,i}^{\,\prime}||\leq\, n\,C_0\,K\big(\tilde{\Gamma}_{\varepsilon,i})\,\varepsilon,
\end{equation}
\noindent
where $C_0:=\big((2\overline{A}-C_1)\big/C_1\underline{A}\big)$ and $C_1>0$ denotes the uniform lower bound for $|{\nabla} U_\varepsilon(p)|$ in $\tilde{\omega}_{\varepsilon,i}$, and where $n$ denotes the number of whole $P$-periods of $\tilde{\Gamma}_{\varepsilon,i}$ intersected by $\tilde{\Gamma}_{\varepsilon,i}^{\,\prime}$.
\vspace{.1in}

\noindent
(e) Let $p_{\varepsilon,i}(s)$, $\tilde{p}_{\varepsilon,i}(s)$, $p_{\varepsilon,\alpha,i}(s)$, and $\hat{p}_{\varepsilon,i}(s)$ denote respective positively-oriented arc-length parametrizations of the arcs $\Gamma_{\varepsilon,i}$, $\tilde{\Gamma}_{\varepsilon,i}$, $\Gamma_{\varepsilon,\alpha,i}$, and $\hat{\Gamma}_{\varepsilon,i}\in {\rm X}$ introduced above, all chosen such that $\tilde{p}_{\varepsilon,i}(0)$, $p_{\varepsilon,\alpha,i}(0)$, $\hat{p}_{\varepsilon,i}(0)\in{\gamma}_{\varepsilon,i}(0)$, where ${\gamma}_{\varepsilon,i}(s)$, $s\in\Re$, denotes, among all maximal (under set inclusion) arcs of steepest ascent of the function $U_{\varepsilon,i}(p):=U_i({\bm\Gamma}_{\varepsilon};p)$ originating at the point ${p}_{\varepsilon,i}(s)$, the particular arc whose intersection with the arc 
$\tilde{\Gamma}_{\varepsilon,i}$ occurs first for all $s\in\Re$ or else last for all $s\in\Re$. We also re-parametrize the arcs $\tilde{\Gamma}_{\varepsilon,i}$, $\Gamma_{\varepsilon,\alpha,i}$, and $\hat{\Gamma}_{\varepsilon,i}$ in the same order by the functions $\tilde{q}_{\varepsilon,i}(s):=\tilde{p}_{\varepsilon,i}\big((\tilde{L}_{\varepsilon,i}\big/{L}_{\varepsilon,i})\,s\big)$, $q_{\varepsilon,\alpha,i}(s):=p_{\varepsilon,\alpha,i}\big((L_{\varepsilon,\alpha,i}\big/{L}_{\varepsilon,i})\,s\big)$, 
and $\hat{q}_{\varepsilon,i}(s):=\hat{p}_{\varepsilon,i}\big((\hat{L}_{\varepsilon,i}\big/{L}_{\varepsilon,i})\,s\big)$
, where $L_{\varepsilon,i}:=||\Gamma_{\varepsilon,i}||$, $\tilde{L}_{\varepsilon,i}:=||\tilde{\Gamma}_{\varepsilon,i}||$, $L_{\varepsilon,\alpha,i}:=||\Gamma_{\varepsilon,\alpha,i}||$, 
and $\hat{L}_{\varepsilon,i}:=||\hat{\Gamma}_{\varepsilon,i}||$, (in which $||\Gamma||$ denotes the arc-length of one $P$-period (in $x$) of an arc $\Gamma\in{\rm X}$).
In this context, we define the positively-oriented parametrizations $\tilde{r}_{\varepsilon,i}(s)$ and $r_{\varepsilon,\alpha,i}(s)$ of the arcs $\tilde{\Gamma}_{\varepsilon,i}$ and ${\Gamma}_{\varepsilon,\alpha,i}$, respectively, such that for each $s\in\Re$, we have $\tilde{r}_{\varepsilon,i}(s):=\tilde{q}_{\varepsilon,i}(\tilde{\pi}_{\varepsilon,i}(s))$ and $r_{\varepsilon,\alpha,i}(s):=q_{\varepsilon,\alpha,i}(\pi_{\varepsilon,\alpha,i}(s))$, where the values $\tilde{\pi}_{\varepsilon,i}(s)$ and $\pi_{\varepsilon,\alpha,i}(s)$ are uniquely chosen such that the arc ${\gamma}_{\varepsilon,i}(s)$ through the point ${p}_{\varepsilon,i}(s)$ intersects $\tilde{\Gamma}_{\varepsilon,i}$ and $\Gamma_{\varepsilon,\alpha,i}$ at the points  $\tilde{q}_{\varepsilon,i}(\tilde{\pi}_{\varepsilon,i}(s))$ and $q_{\varepsilon,\alpha,i}(\pi_{\varepsilon,\alpha,i}(s))$, respectively, and we define the parametrization $\hat{r}_{\varepsilon,i}(s)$ of the arc $\hat{\Gamma}_{\varepsilon,i}$ such that $\hat{r}_{\varepsilon,i}(s)):=\hat{q}_{\varepsilon,i}(\hat{\pi}_{\varepsilon,i}(s))$, where for each $s\in\Re$, the value $\hat{\pi}_{\varepsilon,i}(s)$ is chosen to be minimal (resp. maximal) subject to the requirement that $|\hat{q}_{\varepsilon,i}(\hat{\pi}_{\varepsilon,i}(s))-p_{\varepsilon,i}(s)|={\rm dist}(p_{\varepsilon,i}(s),\hat{\Gamma}_{\varepsilon,i})$.
By Parts (a)-(d), we have:
\begin{equation}
\label{eqn 2.5.20}
\big|\tilde{L}_{\varepsilon,i}-{L}_{\varepsilon,i}\big|, \big|{L}_{\varepsilon,\alpha,i}-{L}_{\varepsilon,i}\big|, \big|\hat{L}_{\varepsilon,\alpha,i}-{L}_{\varepsilon,i}\big|\leq M_1\varepsilon,
\end{equation}
\begin{equation}
\label{eqn 2.5.21}
\big|\tilde{q}_{\varepsilon,i}(s)-{p}_{\varepsilon,i}(s)\big|, \big|{q}_{\varepsilon,\alpha,i}(s)-{p}_{\varepsilon,i}(s)\big|, \big|\hat{q}_{\varepsilon,i}(s)-{p}_{\varepsilon,i}(s)\big|\leq M_1\varepsilon,
\end{equation}
\begin{equation}
\label{eqn 2.5.22*}
\big|\tilde{\pi}_{\varepsilon,i}(s)-s\big|, \big|{\pi}_{\varepsilon,\alpha,i}(s)-s\big|, \big|\hat{\pi}_{\varepsilon,i}(s)-s\big|\leq M_1\varepsilon,
\end{equation}
uniformly over all $s\in\Re$, all sufficiently small values $\varepsilon\in S$, and all $\alpha\in (0,\varepsilon]$, where the constant $M_1$ depends only on $\overline{A}$, $\underline{A}$, the upper bound  $\overline{K}$ for total curvature of $\Gamma_{\varepsilon,i}$ (independent of small $\varepsilon>0$), and the positive uniform lower bound $C>0$ (independent of small $\varepsilon>0$) for $|\nabla U_{\varepsilon,i}(p)|$ relative to $\tilde{\omega}_{\varepsilon,i}$ (see Lem. \ref{lem 2.2.4}).\end{lemma}

\noindent
{\bf Proof of Part (a).} For any $\alpha\in I$, we let $p_\alpha(s):\Re\rightarrow\Gamma_\alpha$ denote the $C^{1,1}$-arc-length parametrization of the $P$-periodic arc $\Gamma_\alpha\in{\rm X}$. Then there is a constant $C_0$ such that $|T_\alpha(s_1)-T_\alpha(s_2)|\leq C_0|s_1-s_2|$ and $|N_\alpha(s_1)-N_\alpha(s_2)|\leq C_0|s_1-s_2|$ uniformly for all $s_1,s_2\in\Re$ and $\alpha\in I$, where $T_\alpha(s)$ (resp. $N_\alpha(s)$) denotes the forward tangent vector (left-hand normal vector) to $\Gamma_\alpha$ at the point $p_\alpha(s)\in\Gamma_\alpha$. It follows from the bounded curvature of the arc $\Gamma_\alpha$ that 
\begin{equation}
\label{eqn 2.5.4}
N_\alpha(s)-N_\alpha(s_0)+A_\alpha(s_0,s)\,T_\alpha(s_0)(s-s_0)=E_\alpha(s_0,s)
\end{equation} 
for $s_0<s$ in $\Re$, where $A_\alpha(s_0,s)\in\Re$ and $E_\alpha(s_0,s)\in\Re^2$ are chosen such that $|A_\alpha(s_0,s)||s-s_0|$ does not exceed the magnitude of the net curvature of $\Gamma_\alpha$ between $p_\alpha(s_0)$ and $p_\alpha(s)$ (thus $|A_\alpha(s_0,s)|\leq C_0$) and $|E_\alpha|=O(|s-s_0|^2)$ as $s\rightarrow s_0$. For each $\alpha, \delta, s\in\Re$ such that $\alpha, \alpha+\delta \in I$ we map $p_\alpha(s)\in\Gamma_\alpha$ onto the point 
\begin{equation}
\label{eqn 2.5.5}
q_{\alpha+\delta}(s)=p_\alpha(s)+r_{\alpha,\delta}(s)N_\alpha(s)\in\Gamma_{\alpha+\delta},
\end{equation} 
where $|r_{\alpha,\delta}(s)|$ is as small as possible (by assumption, we have $|r_{\alpha,\delta}(s)|\leq M|\delta|$ for small $\delta\in\Re$). By double application of (\ref{eqn 2.5.5}), followed by (\ref{eqn 2.5.4}) and the fact that $|p_\alpha(s)-p_\alpha(s_0)-T_\alpha(s_0)(s-s_0)|\leq O(|s-s_0|^2)$ as $s\rightarrow s_0$, we conclude that
\begin{equation}
\label{eqn 2.5.6}
q_{\alpha+\delta}(s)-q_{\alpha+\delta}(s_0)=p_\alpha(s)-p_\alpha(s_0)
\end{equation}
$$+r_{\alpha,\delta}(s_0)(N_\alpha(s)-N_\alpha(s_0))+N_\alpha(s)(r_{\alpha,\delta}(s)-r_{\alpha,\delta}(s_0))$$
$$=((1-A\,r_{\alpha,\delta}(s)(s-s_0)+e_1)\,T_{\alpha}(s_0))+N_\alpha(s_0)(r_{\alpha,\delta}(s)-r_{\alpha,\delta}(s_0)+e_2)$$ 
\noindent
where $e_{1,\alpha}, e_{2,\alpha}\in\Re$ are such that $|e_{1,\alpha}|,|e_{2,\alpha}|\leq O(|s-s_0|^2)$ as $s\rightarrow s_0$. Thus
$$|q_{\alpha+\delta}(s)-q_{\alpha+\delta}(s_0)|^2=|1-A_\alpha\,r_{\alpha,\delta}(s_0))(s-s_0)+e_{1,\alpha}|^2+|r_{\alpha,\delta}(s)-r_{\alpha,\delta}(s_0)+e_{2,\alpha}|^2.$$
\noindent
It follows that
\begin{equation}
\label{eqn 2.5.7}
|q_{\alpha+\delta}(s)-q_{\alpha+\delta}(s_0)|\geq|(1-A_\alpha(s_0,s)r_{\alpha,\delta}(s_0))
(s-s_0)+e_1|
\end{equation}
$$\geq |s-s_0|-K_\alpha(s_0,s)M|\delta|-O(|s-s_0|^2),$$
\noindent
where $A_\alpha=A_\alpha(s_0,s)$ and $K_\alpha(s_0,s)$ denotes the total curvature of $\Gamma
_\alpha$ corresponding to the parameter interval $[s_0,s]$. Given $\alpha\in[0,1]$ and an arc $\Gamma_\alpha^{\,\prime}:=\{p_\alpha(s):s\in J\}$, we use (\ref{eqn 2.5.5}) to define the arcs $\overline{\Gamma}_{\alpha+\delta}^{\,\prime}:=\{q_{\alpha+\delta}(s):s\in J\}$ for $(\alpha+\delta)\in I$ and $|\delta|$ sufficiently small. By applying the estimate (\ref{eqn 2.5.7}) to a sufficiently refined partition of $J$, we conclude that there exists a value $\delta_0>0$ such
that $||\overline{\Gamma}_{\alpha+\delta}^{\,\prime}||\geq||\Gamma_\alpha^{\,\prime}||-M\overline{K}|\delta|$ for all $\alpha,(\alpha+\delta)\in I$ such that $|\delta|<\delta_0$. It is easily seen that $||\overline{\Gamma}_{\alpha+\delta}^{\,\prime}||\leq||\Gamma_{\alpha+\delta}^{\,\prime}||+O(\delta^2)$ as $\delta\rightarrow 0$, where $\Gamma_{\alpha+\delta}^{\,\prime}$ consists of the points in $\Gamma_{\alpha+\delta}$ which project onto $\Gamma_\alpha^{\,\prime}$ along the arcs of steepest ascent of $U$. Therefore, there exists a value $\delta_1\in(0,\delta)$ such that
\begin{equation}
\label{eqn 2.5.8}
||\Gamma_{\alpha+\delta}^{\,\prime}||\geq||\Gamma_\alpha^{\,\prime}||-2M\overline{K}|\delta|,
\end{equation}
for $\alpha, (\alpha+\delta)\in I$ such that $|\delta|<\delta_1$. Part (a) 
(Eq. (\ref{eqn 2.5.1})) follows from this. 
\vspace{.1in}

\noindent
{\bf Proof of Part (c).} For any given values $i\in\{1,2\}$ and $\varepsilon\in(0,\varepsilon_0)$,
 and for each point $p\in{\Gamma}_{\varepsilon,i}$, we use $\Pi_{\varepsilon,i}(p)$ to denote the set of all points $q\in\tilde{\Gamma}_{\varepsilon,i}$ such that $|q-p|={\rm dist}(p,\tilde{\Gamma}_{\varepsilon,i})$ (see Def. \ref{def 3.2.3}). Assuming that $p_1<p_2$ relative to the natural (positive) ordering of points in ${\Gamma}_{\varepsilon,i}$, it follows that $q_1\leq q_2$, relative to the natural ordering of points in $\tilde{\Gamma}_{\varepsilon,i}$, where we assume that $q_j\in\Pi_{\varepsilon,i}(p_j)$ for $j=1,2$ (see the proof of Lem. \ref{lem 3.2.1}(b)). Let $C_{\varepsilon,i}$ denote the set of points $p$ in ${\Gamma}_{\varepsilon,i}$ such that $\Pi_{\varepsilon,i}(p)$ is not a single point. For each point $p_0\in C_{\varepsilon,i}$, we use $\hat{\gamma}_{\varepsilon,i}(p_0)$ to denote the maximal circular arc centered at $p_0$ and having radius $\rho_0:={\rm dist}(p_0,\tilde{\Gamma}_{\varepsilon,i})$, such that both endpoints of $\hat{\gamma}_{\varepsilon,i}(p_0)$ belong to $\Pi_{\varepsilon,i}(p_0)$. For any points $q\in\hat{\gamma}_{\varepsilon,i}(p_0)$ and $p_1\in\dot{\Gamma}_{\varepsilon,i}$ such that $p_1\not=p_0$, we have $$|q-r|+|r-p_0|>|q-p_0|=\rho_0=|q_0-r|+|r-p_0|,$$ whence $|q-r|>|q_0-r|$, from which it follows that
$$|q-p_1|=|q-r|+|r-p_1|>|q_0-r|+|r-p_1|\geq |q_0-p_1|\geq {\rm dist}(p_1,\tilde{\Gamma}_{\varepsilon,i}),$$ where $q_0$ is a point in $\tilde{\Gamma}_{\varepsilon,i}$ such that $|q_0-p_0|=\rho_0$ and such that the line-segment joining $p_0$ to $q_0$ intersects the line-segment joining $p_1$ to $q$ at some point $r$. It follows from this that for any given values $i\in\{1,2\}$ and $\varepsilon\in(0,\varepsilon_0]$,  the set $\hat{\Gamma}_{\varepsilon,i}\setminus\tilde{\Gamma}_{\varepsilon,i}$ is the disjoint union of all circular arcs $\hat{\gamma}_{\varepsilon,i}(p_0)$ corresponding to points $p_0\in C_{\varepsilon,i}$. Part (c) follows from this.
\vspace{.1in}

\noindent
{\bf Proof of Part (d).} In the context of Part (c), given values $i\in\{1,2\}$, $\varepsilon\in(0,\varepsilon_0)$, a point $p_0\in C_{\varepsilon,i}$, and the maximal (under set inclusion) circular directed arc-segment  $\hat{\gamma}_{\varepsilon,i}(p_0)$ of $\hat{\Gamma}_{\varepsilon,i}$ with initial (resp. terminal) endpoint $q_{\varepsilon,i,1}(p_0)$ (resp. $q_{\varepsilon,i,2}(p_0)$), we use $\tilde{\gamma}_{\varepsilon,i}(p_0)$ to denote the corresponding directed sub-segment of $\tilde{\Gamma}_{\varepsilon,i}:=\big\{U_{\varepsilon,i}(p)=\varepsilon\big\}$ having the same initial and terminal endpoints. We observe that 
\begin{equation}
\label{eqn 2.5.12}
\theta_{\varepsilon,i}(p_0):=K(\hat{\gamma}_{\varepsilon,i}(p_0))\leq K(\tilde{\gamma}_{\varepsilon,i}(p_0)),
\end{equation}
$$\theta_{\varepsilon,i}(p_0)\,r_i(p_0)\,\varepsilon=||\hat{\gamma}_{\varepsilon,i}(p_0)||\leq||\tilde{\gamma}_{\varepsilon,i}(p_0)||$$
(closely related to Eq. (\ref{eqn 2.5.9a})), where $r_i(p_0):=\big(1\big/a_i(p_0)\big)$. We also have have 
\begin{equation}
\label{eqn 2.5.13}
C_1||\tilde{\gamma}_{\varepsilon,i}(p_0)||\leq\int_{\tilde{\gamma}_{\varepsilon,i}(p_0)}|{\nabla} U_\varepsilon(p)|\,ds=\int_{\hat{\gamma}_{\varepsilon,i}(p_0)}\big(\partial U_\varepsilon(p)\big/\partial\boldsymbol{\nu}\big)\,ds
\end{equation}
$$\leq \int_{\hat{\gamma}_{\varepsilon,i}(p_0)}|{\nabla} U_\varepsilon(p)|\,ds\leq C_2||\hat{\gamma}_{\varepsilon,i}(p_0)||$$
for $i=1,2$, sufficiently small $\varepsilon\in(0,\varepsilon_0)$, and $p_0\in C_{\varepsilon,i}$, where $C_1>0$ denotes a uniform lower bound for $|{\nabla} U_\varepsilon(p)|$ independent of $p\in\Omega_\varepsilon$ and of all sufficiently small $\varepsilon\in(0,\varepsilon_0]$ (which exists by Lem. \ref{lem 2.2.4}), and where $C_2$ denotes a uniform upper bound for $|{\nabla} U_\varepsilon(p)|$ relative to all arcs $\hat{\gamma}_{\varepsilon,i}(p_0)$ corresponding to sufficiently small $\varepsilon\in(0,\varepsilon_0]$ (one can choose $C_2=2\overline{A}$ by Lem. \ref{lem 2.2.1}(a)). It follows from (\ref{eqn 2.5.12}) and (\ref{eqn 2.5.13}) that
\begin{equation}
\label{eqn 2.5.14}
0\leq\big(||\tilde{\gamma}_{\varepsilon,i}(p_0)||-||\hat{\gamma}_{\varepsilon,i}(p_0)||\big)\leq\big((C_2-C_1)\big/C_1\big)||\hat{\gamma}_{\varepsilon,i}(p_0)||
\end{equation}
$$=\big((C_2-C_1\big/C_1\big)\,r_i(p_0)\,\theta_{\varepsilon,i}(p_0))\,\varepsilon\leq\big((C_2-C_1)\big/C_1\underline{A}\big)\,K(\tilde{\gamma}_{\varepsilon,i}(p_0))\,\varepsilon,$$
for $i=1,2$, any sufficiently small $\varepsilon\in(0,\varepsilon_0)$, and any circular arc-segment $\hat{\gamma}_{\varepsilon,i}(p_0)$ of $\hat{\Gamma}_{\varepsilon,i}$ corresponding to a center-point $p_0\in C_{\varepsilon,i}$. In view of (\ref{eqn 2.5.9a}), it follows by summing (\ref{eqn 2.5.14}) over the collection of center-points $p_{\varepsilon,1}, p_{\varepsilon,2}, \cdots, p_{\varepsilon,n(\varepsilon,i)}\in C_{\varepsilon,i}$ corresponding to one $P$-period (in $x$) that
\begin{equation}
\label{eqn 2.5.15}
0\leq||\tilde{\Gamma}_{\varepsilon,i}^{\,\prime}||-||\hat{\Gamma}_{\varepsilon,i}^{\,\prime}||=\sum_{j=1}^{n(\varepsilon,i)}\big(||\tilde{\gamma}_{\varepsilon,i}(p_{\varepsilon,j})||-||\hat{\gamma}_{\varepsilon,i}(p_{\varepsilon,j)}||\big)
\end{equation}
$$\leq\Big([C_2-C_1])\big/C_1\underline{A}\Big)\Big(\sum_{j=1}^{n(\varepsilon,i)}K\big(\tilde{\gamma}_{\varepsilon,i}(p_{\varepsilon,j})\big)\Big)\,\varepsilon\leq\,C_0\,K(\tilde{\Gamma}_{\varepsilon,i})\,\varepsilon,$$
from which the assertion (\ref{eqn 2.5.3}) follows, where $C_0:=\big((2\overline{A}-C_1)\big/C_1\underline{A}\big).$
\vspace{.1in}

\noindent
{\bf Proof of Part (e).} 
The functions $\tilde{\pi}_{\varepsilon,i}(s), \pi_{\varepsilon,\alpha,i}(s):\Re\rightarrow\Re$ are weakly increasing, since arcs of steepest ascent of the function $U_{\varepsilon,i}$ can't cross relative to the domain $\Omega_\varepsilon$. Also, the related mappings $\tilde{\phi}_{\varepsilon,i}(s):=\tilde{\pi}_{\varepsilon,i}(s)-s:\Re\rightarrow\Re$ and $\phi_{\varepsilon,\alpha,i}(s):=\pi_{\varepsilon,\alpha,i}(s)-s:\Re\rightarrow\Re$ are $L_{\varepsilon,i}$-periodic. In this context, the first two estimates of each of (\ref{eqn 2.5.20}), (\ref{eqn 2.5.21}), and (\ref{eqn 2.5.22*}) follow from Part (b) in connection with Cor. \ref{cor 2.2.1}. Also, the function $\hat{\pi}_{\varepsilon,i}(s):\Re\rightarrow\Re$ is weakly increasing, due to Lem. \ref{lem 3.2.1}(b), and the related function $\hat{\phi}_{\varepsilon,i}(s):=\hat{\pi}_{\varepsilon,i}(s)-s:\Re\rightarrow\Re$ is $L_{\varepsilon,i}$-periodic by definition. In view of this, the three final estimates in Eqs. (\ref{eqn 2.5.20}), (\ref{eqn 2.5.21}), and (\ref{eqn 2.5.22*}) all follow from Part (d), Cor. \ref{cor 2.2.1}, and the previously proven portions of the estimates (\ref{eqn 2.5.20}), (\ref{eqn 2.5.21}), and (\ref{eqn 2.5.22*}).
\begin{lemma}
\label{lem 2.5.4} 
Given any $L$-periodic, twice (resp. once)-differentiable function $f(t):\Re\rightarrow\Re^2$ with absolutely integrable second (first) derivative, we have that
\begin{equation}
\label{eqn 2.5.16}
\Big(\int_0^L\big|f'(t)\big|\,dt\Big)^2\leq 8\int_0^L\big|f(t)\big|\,dt\int_0^L\big|f^{\prime\prime}(t)\big|\,dt,
\end{equation}
\begin{equation}
\label{eqn 2.5.17}
\int_0^L\sup_{\lambda\in[-h,h]}\big|f(t+\lambda)-f(t)\big|\,dt\leq\int_0^L\int_{t-h}^{t+h}\big|f'(s)\big|\,ds\,dt\leq 2h\int_0^L\big|f'(t)\big|\,dt,
\end{equation}
where (\ref{eqn 2.5.17}) holds for any $h>0$.
\end{lemma}

\noindent
{\bf Proof sketch.} The first estimate follows from the second-order Taylor remainder formula for $f$, while the second is related to the Theorem of the Mean.
\begin{lemma}
\label{lem 2.5.5}
For any $0<\varepsilon<1$, let the interval $I_\varepsilon$ and the integrable functions $g_\varepsilon(s),h_\varepsilon(s),\theta_\varepsilon(s):I_\varepsilon\rightarrow \Re$ be such that $\boldsymbol{|}h_\varepsilon\boldsymbol{|}\leq A\varepsilon$, $\int_{I_\varepsilon}|g_\varepsilon(s)|ds\leq B\varepsilon^{3/2}$, $\int_{I_\varepsilon} |\theta_\varepsilon(s)|ds\leq C\sqrt{\varepsilon}$, and $|h_\varepsilon|\leq|g_\varepsilon|$ for $|\theta_\varepsilon|\leq\pi/4$. Then $\int_{I_\varepsilon}|h_\varepsilon(s)|ds\leq(B+(4AC/\pi))\varepsilon^{3/2}$ for all $0<\varepsilon<1$, where the constants $A,B,C$ are independent of $0<\varepsilon<1$.
\end{lemma}

\noindent
{\bf  Proof.} Let $J_\varepsilon=\{s\in I_\varepsilon:|\theta_\varepsilon(s)|\geq(\pi/4)\}$. Then $(\pi/4)|J_\varepsilon|\leq\int_{J_\varepsilon}|\theta_\varepsilon(s)| ds\leq C\sqrt{\varepsilon}$, from which it follows that $|J_\varepsilon|\leq(4C\sqrt{\varepsilon}/\pi)$. Therefore $\int_{J_\varepsilon}|h_\varepsilon(s)|ds\leq(4AC\varepsilon^{3/2}/\pi)$. We also have $\int_{I_\varepsilon\setminus J_\varepsilon}|h_\varepsilon(s)|ds\leq\int_{I_\varepsilon\setminus J_\varepsilon}|g_\varepsilon(s)|ds\leq B\varepsilon^{3/2}$. The assertion easily follows.

\begin{proposition}
\label{prop 2.5.1}
Given a countable set $S\subset(0,\varepsilon_0)$ having $0$ as it's sole accumulation point, let ${\bf\Gamma}_\varepsilon=(\Gamma_{\varepsilon,1},\Gamma_{\varepsilon,2})\in{\bm{\cal F}}_\varepsilon$, $\varepsilon\in S$, denote a family of "fixed points" in ${\bf X}$ such that the arc-lengths $L_{\varepsilon,i}$, $\tilde{L}_{\varepsilon,i}$ and total curvatures (all per $P$-period in $x$) of the arcs ${\Gamma}_{\varepsilon,i}$ and $\tilde{\Gamma}_{\varepsilon,i}:=\Phi_{\varepsilon,i}({\bm\Gamma}_\varepsilon)=\big\{U_{\varepsilon,i}(p)=\varepsilon\big\}$, resp., are uniformly bounded from above over all $\varepsilon\in S$ and $i=1,2$. For any $\varepsilon\in S$, let ${p}_{\varepsilon,i}(s):\Re\rightarrow{\Gamma}_{\varepsilon,i}$ denote a positively-oriented arc-length parametrization of the arc ${\Gamma}_{\varepsilon,i}$, $i=1,2$, and let $\tilde{p}_{\varepsilon,i}(s):\Re\rightarrow\tilde{\Gamma}_{\varepsilon,i}$ denote the mapping of each $s\in\Re$ into the always maximal or always minimal (in terms of the natural ordering in $\tilde{\Gamma}_{\varepsilon,i}$) point of intersection $\tilde{p}_{\varepsilon,i}(s)$ of the arc $\tilde{\Gamma}_{\varepsilon,i}$ with any arc $\gamma_{\varepsilon,i}(s)$ of steepest ascent of the function $U_{\varepsilon,i}(p):=U_i({\bm\Gamma}_\varepsilon;p)$ originating at the point $p_{\varepsilon,i}(s)$.  
Then, given a constant $C$, there exists a constant $M$ such that 
\begin{equation}
\label{eqn 2.5.1a}
\int_0^{L_{\varepsilon,i}}f_{\varepsilon,i}(s)\,ds\leq M\varepsilon^{3/2},
\end{equation}
uniformly for $i=1,2$, and for all sufficiently small $\varepsilon\in S$, where we define $L_{\varepsilon,i}:=||\Gamma_{\varepsilon,i}||$ and 
\begin{equation}
\label{eqn 2.5.2a}
0\leq f_{\varepsilon,i}(s):=\big|\tilde{p}_{\varepsilon,i}(s)-{p}_{\varepsilon,i}(s)\big|-{\rm dist}\big({p}_{\varepsilon,i}(s),\tilde{\Gamma}_{\varepsilon,i}\big)\leq C\varepsilon.
\end{equation}
\end{proposition}

\noindent
{\bf Proof of Prop. \ref{prop 2.5.1}.} This proof is in the context and notation of the proof of Lem. \ref{lem 2.5.1} (especially Lem. \ref{lem 2.5.1}(e)). Given ${\bm\Gamma}_{\varepsilon}\in{\bm{\cal F}}_\varepsilon$, we use $\hat{\Gamma}_{\varepsilon,i}\in{\rm X}$ to denote the boundary component relative to $\tilde{\omega}_{\varepsilon,i}\cup\tilde{\Gamma}_{\varepsilon,i}$ of the annular domain
\begin{equation}
\label{eqn 2.5.19}
\hat{\omega}_{\varepsilon,i}:=\{q\in\tilde{\omega}_{\varepsilon,i}:|q-p|<{\rm dist}(p,\tilde{\Gamma}_{\varepsilon,i})\,\,{\rm for\,\, some}\,\,p\in{\Gamma}_{\varepsilon,i}\},
\end{equation}
where $\tilde{\omega}_{\varepsilon,i}$ denotes the strip-like domain bounded by ${\Gamma}_{\varepsilon,i}\cup\tilde{\Gamma}_{\varepsilon,i}$. 
By Lem. \ref{lem 2.5.1}(a)-(e), there exists a constant $M_1$ such that the estimates (\ref{eqn 2.5.20}), (\ref{eqn 2.5.21}), and (\ref{eqn 2.5.22*}) hold uniformly over all $s\in\Re$, $\varepsilon\in S$, and $\alpha\in (0,\varepsilon]$. In terms of further notation in Lem. \ref{lem 2.5.1}(e), we study the functions $f_{\varepsilon,i}(s):\Re\rightarrow\Re$ defined by:
\begin{equation}
\label{eqn 2.5.23}
0\leq f_{\varepsilon,i}(s):=\big|\tilde{r}_{\varepsilon,i}(s)-{p}_{\varepsilon,i}(s)\big|-\big|\hat{r}_{\varepsilon,i}(s)-{p}_{\varepsilon,i}(s)\big|\leq C\varepsilon,
\end{equation}
where the mapping $\tilde{r}_{\varepsilon,i}(s):=q_{\varepsilon,i}(\tilde{\pi}_{\varepsilon,i}(s)):\Re\rightarrow\tilde{\Gamma}_{\varepsilon,i}$ coincides with $\tilde{p}_{\varepsilon,i}(s):\Re\rightarrow\tilde{\Gamma}_{\varepsilon,i}$, and the mapping $\hat{r}_{\varepsilon,i}(s):=\hat{q}_{\varepsilon,i}(\hat{\pi}_{\varepsilon,i}(s)):\Re\rightarrow\hat{\Gamma}_{\varepsilon,i}$ is such that $\big|\hat{r}_{\varepsilon,i}(s)-p_{\varepsilon,i}(s)\big|={\rm dist}(p_{\varepsilon,i}(s),\tilde{\Gamma}_{\varepsilon,i})$ for all $s\in\Re$.
We have
\begin{equation}
\label{eqn 2.5.24}
f_{\varepsilon,i}(s)\leq \tilde{f}_{\varepsilon,i}(s)+\hat{f}_{\varepsilon,i}(s),
\end{equation}
where
\begin{equation}
\label{eqn 2.5.25}
\tilde{f}_{\varepsilon,i}(s):={\rm dist}\big(\tilde{r}_{\varepsilon,i}(s),\hat{\cal L}_{\varepsilon,i}(s)\big)=\big|\tilde{r}_{\varepsilon,i}(s)-\tilde{Q}_{\varepsilon,i}(s)\big|,
\end{equation}
\begin{equation}
\label{eqn 2.5.26}
\hat{f}_{\varepsilon,i}(s):=\big|\hat{r}_{\varepsilon,i}(s))-\tilde{Q}_{\varepsilon,i}(s)\big|,
\end{equation}
in which $\hat{{\cal L}}_{\varepsilon,i}(s)$ denotes the straight line passing through the points ${p}_{\varepsilon,i}(s)$ and $\hat{r}_{\varepsilon,i}(s)$, while $\tilde{Q}_{\varepsilon,i}(s)$ denotes the point closest to the point $\tilde{r}_{\varepsilon,i}(s)$ in the line $\hat{{\cal L}}_{\varepsilon,i}(s)$. By assumption, we have that
\begin{equation}
\label{eqn 2.5.27}
0\leq \tilde{f}_{\varepsilon,i}(s), \hat{f}_{\varepsilon,i}(s)\leq C\varepsilon
\end{equation}
for some constant $C$. Also, due to the perpendicularity of the arcs $\tilde{\Gamma}_{\varepsilon,i}$ and $\hat{{\cal L}}_{\varepsilon,i}(s)$ at their intersection point $\hat{r}_{\varepsilon,i}(s)$, it is easily seen that 
\begin{equation}
\label{eqn 2.5.28}
\hat{f}_{\varepsilon,i}(s)\leq \tilde{f}_{\varepsilon,i}(s)\,\,{\rm whenever}\,\,\theta_{\varepsilon,i}(s)\leq\pi/4,
\end{equation}
where $\theta_{\varepsilon,i}(s)$ denotes the maximum absolute variation in the argument of the forward tangent to the arc-segment of the arc $\tilde{\Gamma}_{\varepsilon,i}$ between the points $\hat{r}_{\varepsilon,i}(s)$ and $\tilde{r}_{\varepsilon,i}(s))$. By (\ref{eqn 2.5.22*}), we have $\big|\hat{\pi}_{\varepsilon,i}(s)-\tilde{\pi}_{\varepsilon,i}(s)\big|\leq M_1\varepsilon$, where of course $M_1\varepsilon\leq\sqrt{\varepsilon}$ if $\varepsilon\in S$ is sufficiently small. In view of this, it follows from Lem. \ref{lem 2.5.4} (Eq. (\ref{eqn 2.5.15})), with $f(t):={\rm arg}\,(\tilde{q}'_{\varepsilon,i}(t))$, that
\begin{equation}
\label{eqn 2.5.29} 
\int_0^{L_{\varepsilon,i}} \theta_{\varepsilon,i}(s)\,ds\,\leq M_2\sqrt{\varepsilon}
\end{equation}
\noindent
for a new constant $M_2$. Finally, in view of Lem. \ref{lem 2.5.5}, it follows from (\ref{eqn 2.5.27}), (\ref{eqn 2.5.28}), and (\ref{eqn 2.5.29}) that if the inequality:
\begin{equation}
\label{eqn 2.5.30}
\int_0^{L_{\varepsilon,i}}\,\tilde{f}_{\varepsilon,i}(s)\,ds\leq M_3\varepsilon^{3/2}
\end{equation}
holds for a fixed constant $M_3$ and all sufficiently small $\varepsilon\in S$, then we also have $\int_0^{L_{\varepsilon,i}}\hat{f}_{\varepsilon,i}(s)\,ds\leq M_4\,\varepsilon^{3/2}$ for some constant $M_4$ and for all sufficiently small $\varepsilon\in S$, from which it follows via (\ref{eqn 2.5.24}) that (\ref{eqn 2.5.1a'}) holds under the same conditions. Toward the proof of (\ref{eqn 2.5.30}), 
we define $\tilde{f}_{\varepsilon,i}(\alpha,s):={\rm dist}\big(r_{\varepsilon,\alpha,i}(s), \hat{{\cal L}}_{\varepsilon,i}(s)\big)$ for all $s\in\Re$, $\alpha\in(0,\varepsilon]$, and sufficiently small $\varepsilon\in S$. Then, in view of the fact that $\tilde{f}_{\varepsilon,i}(0,s)=0$, we have that
\begin{equation}
\label{eqn 2.5.30a}
\big|\tilde{f}_{\varepsilon,i}(s)\big|\leq\int_0^{\varepsilon}\left|\frac{\partial \tilde{f}_{\varepsilon,i}(\alpha,s)}{\partial\alpha}\right|\,d\alpha,
\end{equation}
for any $s\in\Re$, where, for each $\alpha\in(0,\varepsilon]$, we have
\begin{equation}
\label{eqn 2.5.31}
\left|\frac{\partial \tilde{f}_{\varepsilon,i}(\alpha,s)}{\partial\alpha}\right|\leq \frac{|\boldsymbol{\nu}_{\varepsilon,i}\cdot\boldsymbol{\tau}_{\varepsilon,i}|}{|{\nabla} U_{\varepsilon,i}|}\leq\frac{|\boldsymbol{\nu}_{\varepsilon,i}-\boldsymbol{\nu}_{\varepsilon,0,i}|}{|{\nabla} U_{\varepsilon,i}|}=\frac{|\boldsymbol{\nu}^\perp_{\varepsilon,i}-\boldsymbol{\nu}_{\varepsilon,0,i}^\perp|}{|{\nabla} U_{\varepsilon,i}|}.
\end{equation}
Here $|{\nabla} U_{\varepsilon,i}|$ and $\boldsymbol{\nu}_{\varepsilon,i}:=({\nabla} U_{\varepsilon,i}/|{\nabla} U_{\varepsilon,i}|)$ are evaluated at the point $r_{\varepsilon,\alpha,i}(s)$, $\boldsymbol{\tau}_{\varepsilon,i}$ is a unit normal to the line $\hat{{\cal L}}_{\varepsilon,i}(s)$, $\boldsymbol{\nu}_{\varepsilon,0,i}$ denotes any vector in the direction of $\hat{{\cal L}}_{\varepsilon,i}(s)$, and $\boldsymbol{\nu}^\perp_{\varepsilon,i}$ and $\boldsymbol{\nu}_{\varepsilon,0,i}^\perp$ denote clockwise 90-degree rotations of $\boldsymbol{\nu}_{\varepsilon,i}$ and $\boldsymbol{\nu}_{\varepsilon,0,i}$, respectively. We also set $\boldsymbol{\nu}_{\varepsilon,0,i}^{\perp}=\hat{q}'_{\varepsilon,i}(\hat{\pi}_{\varepsilon,i}(s))$ for the same $\varepsilon\in S$ and $s\in\Re$. We have by definition (see Lem. \ref{lem 2.5.1}(e)) that 
$$q_{\varepsilon,\alpha,i}(s):=p_{\varepsilon,\alpha,i}\big(\big(L_{\varepsilon,i}\big/L_{\varepsilon,\alpha,i}\big)\,\,s\big),$$ from which it follows that $$q'_{\varepsilon,\alpha,i}(s)=\big(L_{\varepsilon,i}\big/L_{\varepsilon,\alpha,i}\big)\,p'_{\varepsilon,\alpha,i}\big(\big({L}_{\varepsilon,i}\big/L_{\varepsilon,\alpha,i}\big)\,\,s\big)$$ for any $s\in\Re$, where we have $|p'_{\varepsilon,\alpha,i}(s)|=1$ for all $s$ and also 
$\big|L_{\varepsilon,i}-L_{\varepsilon,\alpha,i}\big|\leq M_1\,\varepsilon$ by (\ref{eqn 2.5.20}).
In view of this, it follows from 
(\ref{eqn 2.5.31}) that 
\begin{equation}
\label{eqn 2.5.33}
\left|\frac{\partial \tilde{f}_{\varepsilon,i}(\alpha,s)}{\partial\alpha}\right|\leq \frac{\big|q_{\varepsilon,\alpha,i}^{\prime}(\pi_{\varepsilon,\alpha,i}(s))-\hat{q}_{\varepsilon,i}^\prime(\hat{\pi}_{\varepsilon,i}(s))\big|+O(\varepsilon)}{|{\nabla} U_{\varepsilon,i}(q_{\varepsilon,\alpha,i}(\pi_{\varepsilon,\alpha,i}(s))|}
\end{equation}
for all $\alpha\in (0,\varepsilon]$, $\varepsilon\in S$, and $s\in\Re$, where the reciprocal of $|{\nabla} U_{\varepsilon,i}(p)|$ is bounded above by a uniform constant $M_5$ (see Lem. \ref{lem 2.2.4}). By again using the fact that $|\pi_{\varepsilon,\alpha,i}(s)-s|,|\hat{\pi}_\varepsilon(s)-s|\leq M_1\varepsilon$ (see (\ref{eqn 2.5.22*})), where $M_1\varepsilon\leq\sqrt{\varepsilon}$ for $\varepsilon\in S$ sufficiently small, we conclude that
\begin{equation}
\label{eqn 2.5.35}
\big|q_{\varepsilon,\alpha,i}^{\prime}(\pi_{\varepsilon,\alpha,i}(s))-\hat{q}_{\varepsilon,i}^\prime(\hat{\pi}_{\varepsilon,i}(s))\big|\leq
\big|q_{\varepsilon,\alpha,i}^\prime(\pi_{\varepsilon,\alpha,i}(s))-q_{\varepsilon,\alpha,i}^\prime(s)\big|
\end{equation}
$$+\big|q_{\varepsilon,\alpha,i}^\prime(s)-\hat{q}_{\varepsilon,i}^\prime(s))\big|+\big|\hat{q}_{\varepsilon,i}^\prime(s)-\hat{q}_{\varepsilon,i}^\prime(\hat{\pi}_{\varepsilon,i}(s))\big|$$
$$\leq e_{\varepsilon,\alpha,i}(s)+h_{\varepsilon,\alpha,i}(s)+\hat{e}_{\varepsilon,i}(s),$$

\noindent
for all $s\in\Re$, $\alpha\in(0,\varepsilon]$ and sufficiently small $\varepsilon\in S$, where 
$$e_{\varepsilon,\alpha,i}(s):=\sup\,\big\{\big|q_{\varepsilon,\alpha,i}^\prime(s+\lambda)-q_{\varepsilon,\alpha,i}^\prime(s)\big|:|\lambda|\leq\sqrt{\varepsilon}\,\big\},$$
$$\hat{e}_{\varepsilon,i}(s):=\sup\,\big\{\big|\hat{q}_{\varepsilon,i}^\prime(s+\lambda)-\hat{q}_{\varepsilon,i}^\prime(s)\big|:|\lambda|\leq\sqrt{\varepsilon}\,\big\},$$
$$h_{\varepsilon,\alpha,i}(s):=\big|q_{\varepsilon,\alpha,i}^\prime(s)-\hat{q}_{\varepsilon,i}^\prime(s)\big|.$$

\noindent
By substituting each of the $L_{\varepsilon,i}$-periodic functions: $q'_{\varepsilon,\alpha,i}(s)$ and $\hat{q}'_{\varepsilon,i}(s)$ for $f(s)$ in Lem. \ref{lem 2.5.4}, Eq. (\ref{eqn 2.5.17}) and interpreting the integrals $\int_0^{L_{\varepsilon,i}}\big|q_{\varepsilon,\alpha,i}''(s)\big|\,ds$ and $\int_0^{L_{\varepsilon,i}}\big|\hat{q}_{\varepsilon,i}''(s)\big|\,ds$ in terms of the total curvatures per $P$-period (in $x$) of the respective arcs $\Gamma_{\varepsilon,\alpha,i}$ and $\hat{\Gamma}_{\varepsilon,i}$, we conclude that
\begin{equation}
\label{eqn 2.5.36}
\int_0^{L_{\varepsilon,i}}e_{\varepsilon,\alpha,i}(s)\,ds\leq M_5\sqrt{\varepsilon};\,\,\int_0^{L_{\varepsilon,i}}\hat{e}_{\varepsilon,i}(s)\,ds\leq M_6\sqrt{\varepsilon}.
\end{equation}
for all $\alpha\in(0,\varepsilon]$ and sufficiently small $\varepsilon\in S$. Similarly, by substituting the $L_{\varepsilon,i}$-periodic function $f(s):=q_{\varepsilon,\alpha,i}(s)-\hat{q}_{\varepsilon,i}(s)$ into (\ref{eqn 2.5.16}), we get
\begin{equation}
\label{eqn 2.5.37}
\Big(\int_0^{L_{\varepsilon,i}}\big|q_{\varepsilon,\alpha,i}^\prime(s)-\hat{q}_{\varepsilon,i}^\prime(s)\big|\,ds\,\Big)^2\leq
\end{equation}
$$
8\int_0^{L_{\varepsilon,i}}\big|q_{\varepsilon,\alpha,i}(s)-\hat{q}_{\varepsilon,i}(s)\big|\,ds\int_0^{L_{\varepsilon,i}}\big(\big|q^{\prime\prime}_{\varepsilon,\alpha,i}(s)\big|+\big|\hat{q}^{\prime\prime}_{\varepsilon,i}(s)\big|\big)\,ds,$$
from which it easily follows by estimating the first integral on the second line of (\ref{eqn 2.5.37}) and estimating the second integral by the total curvatures of  $\Gamma_{\varepsilon,\alpha,i}$ and $\hat{\Gamma}_{\varepsilon,i}$ that
\begin{equation}
\label{eqn 2.5.38}
\int_0^{L_{\varepsilon,i}}h_{\varepsilon,\alpha,i}(s)\,ds\leq M_7\sqrt{\varepsilon}
\end{equation}
for all $\alpha\in(0,\varepsilon]$ and sufficiently small $\varepsilon\in S$. By substituting the inequalities (\ref{eqn 2.5.36}(a)), (\ref{eqn 2.5.36}(b)), and (\ref{eqn 2.5.38}) into (\ref{eqn 2.5.33}) and (\ref{eqn 2.5.35}), we conclude that
\begin{equation}
\label{eqn 2.5.39}
\int_0^{L_{\varepsilon,i}}\big|(\partial/\partial \alpha)\,\tilde{f}_{\varepsilon,i}(\alpha,s)\big|\,ds\leq M_8\sqrt{\varepsilon}.
\end{equation}
for all $\alpha\in(0,\varepsilon]$ and sufficiently small $\varepsilon\in S$. Finally, it follows from (\ref{eqn 2.5.30a}), (\ref{eqn  2.5.39}), and Fubini's theorem that

\begin{equation}
\label{eqn 2.5.40}
\int_0^{L_{\varepsilon,i}}\,\tilde{f}_{\varepsilon,i}(s)\,ds\leq\int_0^{L_{\varepsilon,i}}\int_0^\varepsilon \big|(\partial/\partial\alpha)\,\tilde{f}_{\varepsilon,i}(\alpha,s)\big|\,d\alpha\,ds
\end{equation}
$$=\int_0^\varepsilon\int_0^{L_{\varepsilon,i}}\big|(\partial\big/\partial\alpha)\tilde{f}_{\varepsilon,i}(\alpha,s)\big|\,ds\,d\alpha\leq\int_0^\varepsilon M_8\sqrt{\varepsilon}\,\,d\alpha\leq M_3\varepsilon^{3/2},$$

\noindent
for all sufficiently small $\varepsilon\in S$, as was required in (\ref{eqn 2.5.30}).

\begin{corollary}
\label{cor 2.5.1}
{(Generalization of Prop. \ref{prop 2.5.1})}
Given $\varepsilon\in S$, ${\bm\Gamma}_\varepsilon\in{\bm{\cal F}}_\varepsilon$, $\delta\in[0,1)$, and the capacitary potentials $U_{\varepsilon,i}(p):=U_i({\bm\Gamma}_\varepsilon;p)$ (defined in the closures of the periodic domains $\Omega_\varepsilon:=\Omega({\bm\Gamma}_\varepsilon)$), let the related capacitary potentials $U_{\varepsilon,\delta,i}(p):{\rm Cl}(\Omega_{\varepsilon,\delta,i})\rightarrow\Re$ be defined such that $U_{\varepsilon,\delta,i}(p):=\big((U_{\varepsilon,i}(p)-\varepsilon\delta)\big/(1-\varepsilon\delta)\big)$
in the closure of $\Omega_{\varepsilon,\delta,i}:=\{0<U_{\varepsilon,\delta,i}(p)<1\}$, we use ${p}_{\varepsilon,\delta,i}(s):\Re\rightarrow{\Gamma}_{\varepsilon,\delta,i}$ to denote positively-oriented arc-length parametrizations of the arcs ${\Gamma}_{\varepsilon,\delta,i}:=\{U_{\varepsilon,\delta,i}(p)=0\}$, $i=1,2$, and we use $\tilde{p}_{\varepsilon,\delta,i}(s):\Re\rightarrow\tilde{\Gamma}_{\varepsilon,\delta,i}$ to denote the mapping of each $s\in\Re$ into the unique point of intersection $\tilde{p}_{\varepsilon,\delta,i}(s)$ of the arc $\tilde{\Gamma}_{\varepsilon,\delta,i}:=\{U_{\varepsilon,\delta,i}(p)=\varepsilon\}$ with the unique arc $\gamma_{\varepsilon,\delta,i}(s)$ of steepest ascent of the function $U_{\varepsilon,\delta,i}(p)$ originating at the point $p_{\varepsilon,\delta,i}(s)$. Then there exist constants $C$ and $M$, independent of small $\delta\in(0,1)$ such that 
\begin{equation}
\label{eqn 2.5.1a'}
\int_0^{L_{\varepsilon,\delta,i}}f_{\varepsilon,\delta,i}(s)\,ds\leq M\varepsilon^{3/2},
\end{equation}
uniformly for $i=1,2$, and for all sufficiently small $\varepsilon\in S$, where 
we define $L_{\varepsilon,\delta,i}:=||\Gamma_{\varepsilon,\delta,i}||$ and 
\begin{equation}
\label{eqn 2.5.2a'}
0\leq f_{\varepsilon,\delta,i}(s):=\big|\tilde{p}_{\varepsilon,\delta,i}(s)-{p}_{\varepsilon,\delta,i}(s)\big|-{\rm dist}\big({p}_{\varepsilon,\delta,i}(s),\tilde{\Gamma}_{\varepsilon,\delta,i}\big)\leq C\varepsilon.
\end{equation}
\end{corollary}
\vspace{.1in}

\noindent
{\bf Proof.} For any given values $\varepsilon\in S$ and $\delta\in[0,1)$, we redefine all the notation in Prop. \ref{prop 2.5.1} and its proof by replacing all the definitions previously expressed in terms of (or based on) the capacitary potentials $U_{\varepsilon,i}(p):{\rm Cl}(\Omega_{\varepsilon,i})\rightarrow\Re$ by the corresponding definitions based on the corresponding capacitary potentials $U_{\varepsilon,\delta,i}(p):{\rm Cl}(\Omega_{\varepsilon,\delta,i})\rightarrow\Re$. Thus Cor. \ref{cor 2.5.1} reduces in the case $\delta=0$ to Prop. \ref{prop 2.5.1}, and the proof valid for the case $\delta=0$ easily extends to establish uniform estimates valid for small $\varepsilon\in S$ and small $\delta\in(0,1)$.

\subsection{Sharp positive lower bounds for $|\nabla U(p)|$ on $\Gamma$}
\label{subsection 3.3}
\begin{theorem} 
{(Sharp uniform lower bound for boundary gradient of the capacitary potential)}
\label{thm 2.6.1}
Let $U(p):=U({\bf\Gamma};p)$ in the closure of $\Omega:=\Omega({\bf\Gamma})$, where ${\bf\Gamma}=(\Gamma_1,\Gamma_2)\in{\boldsymbol{\cal F}}$ denotes a weak solution of Prob. \ref{prob2.1.1} (see Def. \ref{def 2.3.1} and Thm. \ref{thm 2.3.1}), obtained as the limit of a convergent sequence of operator "fixed points" ${\bf \Gamma}_n\in{\bm{\cal F}}_{\varepsilon_n}$ corresponding to a null-sequence of values $\varepsilon_n$ in the interval $(0,\varepsilon_0)$.
Assume the component-wise arc-length and total curvature per $P$-period of the "fixed points" ${\bf\Gamma}_n, n\in N$, are uniformly bounded (as they will be under the assumptions of Thm. \ref{thm 2.1.3}, where we set ${\bf X}:={\bf Y}$). Then: 
\begin{equation}
\label{eqn 2.6.1}
\liminf_{p\rightarrow p_i}|{\nabla} U(p)|\geq a_i(p_i),
\end{equation}
$i=1,2$, where $p\in\Omega$ and $p_i$ denotes any point in $\Gamma_i$.
In fact we have
\begin{equation}
\label{eqn 2.6.2}
|{\nabla} U_i(p)|\geq E_i(p):=\hat{a}_i(p)\,{\rm exp}\big(-2\,C_i\,U_i(p)\big),
\end{equation}
in $\hat{\Omega}_i:=\{p\in\Omega:U_i(p)<1/2\}$ for $i=1,2$, where 
$\hat{a}_i(p):{\rm Cl}(\hat{\Omega}_i)\rightarrow\Re$ denotes the continuous function such that $\hat{a}_i(p)=a_i(p)$ on $\partial\hat{\Omega}_i$ and ${\rm ln}\big(\hat{a}_i(p)\big)$ is harmonic in $\hat{\Omega}_i$, and where   
$C_i=\sup\big\{{\rm ln}\big({a}_i(q)\big/|{\nabla} U(q)|\big):q\in\hat{\Gamma}\big\}$, and where $\hat{\Gamma}:=\{U(q)=1/2\}$. Observe that in (\ref{eqn 2.6.2}), $E_i(p)$, $i=1,2$, denotes a continuous function in ${\rm Cl}(\hat{\Omega}_i)
$, such that $E_i(p)=a_i(p)$ on $\Gamma_i$. 
\end{theorem}
\begin{remark}
\label{rem 2.6.1}
For $\phi\in L^1(\Re)$ and $g\in L^p(\Re)$ with $p\geq 1$, we have: 
\begin{equation}
\label{eqn 2.6.3}
\int_{-\infty}^\infty \big|(g*\phi)(x)\big|^p dx\leq\Big(\int_{-\infty}^\infty |\phi(x)|dx\Big)^p\Big(\int_{-\infty}^\infty|g(x)|^p dx\Big),
\end{equation}
where "$g*\phi$" denotes the convolution of $\phi$ and $g$ (see [SL], Thm. 6p, p. 374). If $g\in L^p(\Re)$ and $\phi(x)$ denotes an integrable, $\kappa$-periodic function, then 
\begin{equation}
\label{eqn 2.6.4}
\int_{\tau}^{\tau+\kappa}\big|(g*\phi)(x)\big|^p dx\leq\Big(\int_{\tau}^{\tau+\kappa}|\phi(x)|dx\Big)^p\Big(\int_{-\infty}^\infty|g(x)|^p dx\Big)
\end{equation}
for any $\tau\in\Re$; in fact (\ref{eqn 2.6.4}) follows by setting $\phi(x):=\phi_n(x)$ in ( \ref{eqn 2.6.3}), where we define $\phi_n(x)=\phi(x)$ in $[-n\kappa,n\kappa]$, $\phi_n(x)=0$ in $\Re\setminus[-n\kappa,n\kappa]$, and using the fact that $g*\phi_n\rightarrow g*\phi$ as $n\rightarrow\infty$, where $g*\phi$ is $\kappa$-periodic.
\end{remark}
\begin{lemma}
\label{lem 2.6.1}
Let be given a family of fixed points ${\bf\Gamma}_\varepsilon=(\Gamma_{\varepsilon,1},\Gamma_{\varepsilon,2})\in {\boldsymbol{\cal F}}_\varepsilon$, $\varepsilon\in S\subset(0,\varepsilon_0)$, such that the total curvature (in one $P$-period) of the curves $\Gamma_{\varepsilon,i}$, $i=1,2$, is bounded above by a constant $M$ (independent of $\varepsilon\in S$). Then for any $\delta\in(0,1)$, there exists a constant $C(\delta)$ such that
\begin{equation}
\label{eqn 2.6.5}
\int_{\Gamma_{\varepsilon,\delta,i}}\big|{\nabla} U_{\varepsilon,i}(p)\big|^2\,ds_{\varepsilon,\delta,i}\leq C(\delta),
\end{equation}
for $i=1,2$, uniformly for all $\varepsilon\in S$, where 
$U_{\varepsilon,i}(p):=U_i({\bf\Gamma}_\varepsilon;p)$ and 
$\Gamma_{\varepsilon,\delta,i}:=$ $\Phi_{\varepsilon\delta,i}({\bf\Gamma}_\varepsilon)$ $=\{U_{\varepsilon,i}(p)=\varepsilon\delta\}$.
\end{lemma}
\noindent
{\bf Proof.} We fix the value $\varepsilon\in S$, and suppress the subscript $\varepsilon$. We also assume without loss of generality that the curve-pair ${\bf\Gamma}=(\Gamma_1,\Gamma_2)\in{\bf X}$ is analytic. We let $\kappa=\kappa({\bf\Gamma})$ denote the capacity of one $P$-period (in $x$) of  the domain $\Omega:=\Omega({\bf\Gamma})$, while $K_i=K(\Gamma_{i}):={\rm total\,\ curvature\,\, of\,\,one}\,P\,{\rm -period\,\,of\,\,}\Gamma_i$. Let $w=F(z):\Re\times[0,1]\rightarrow{\rm Cl}(\Omega)$ be an analytic, $\kappa$-periodic, onto function, whose restriction to $[\tau,\tau+\kappa)\times[0,1]$ is one-to-one for each $\tau\in\Re$. Let $G(x,y)={\rm ln}(|F'(z)|)$ in $\Re\times[0,1]$. Then the function  
$$G_y(x,y)=(\partial/\partial y)\,{\rm ln}\,(|F'(z)|)=-(\partial/\partial x)\,{\rm arg}\,(F'(z))$$
is harmonic in $\Re\times[0,1]$, and we have $\int_\tau^{\tau+\kappa}\phi(x)\,dx=0$ and $\int_\tau^{\tau+\kappa} |\phi(x)|\,dx\leq M$  for any $\tau\in\Re$, where we define $\phi(x):=G_y(x,0)$ in $\Re$. We define the function $H(x,y)$ in $\Re\times[0,\infty)$ such that $H_y(x,0):=\phi(x)$ in $\Re\times\{0\}$ and such that $H_y$ is defined by the convolution integral: 
\begin{equation}
\label{eqn 2.6.6}
H_y(x,y):=\frac{1}{\pi}\int_{-\infty}^\infty\frac{y}{(x-t)^2+y^2} \,\,\phi(t)\,dt
\end{equation}
in $\Re\times(0,\infty)$. Observe that $H_y(x,y)$ is bounded, $\kappa$-periodic (in $x$), continuous in $\Re\times[0,\infty)$, and harmonic in $\Re\times(0,\infty)$. Also $H_y=G_y$ in $\Re\times\{0\}$ and $|H_y-G_y|\leq C$ on $\Re\times\{\frac{1}{2}\}$ for some constant $C$, from which it follows that $|H_y(x,y)-G_y(x,y)|\leq 2Cy$ in $\Re\times[0,\frac{1}{2}]$ by the maximum principle. Therefore, 
\begin{equation}
\label{eqn 2.6.7}
|(G(x,b)-G(x,a))-(H(x,b)-H(x,a))|\leq C (b^2-a^2)
\end{equation}
for any $x\in\Re$ and $0<a<b\leq \frac{1}{2}$. By integrating (\ref{eqn 2.6.6}) over $y\in[a,b]\subset[0,1/2]$, we get
$$H(x,b)-H(x,a)=\frac{1}{2\pi}\int_{-\infty}^\infty {\rm ln}\left (\frac{(x-t)^2+b^2}{(x-t)^2+a^2}\right )\phi(t)\,dt,$$
\noindent
from which it follows (see Remark \ref{rem 2.6.1}) that
\begin{equation}
\label{eqn 2.6.8}
\int_\tau^{\tau+\kappa}\left |H(x,b)-H(x,a)\right |^N dx\leq \left({\frac{M}{2\pi}}\right)^N\int_{-\infty}^\infty\left(
{\rm ln}\,\frac{x^2+b^2}{x^2+a^2}\right)^N dx
\end{equation}
for any $N\geq 1$. By substituting (\ref{eqn 2.6.8}) into the series expansion for ${\rm exp}\,\big(H(x,b)-H(x,a)\big)$, we see that
\begin{equation}
\label{eqn 2.6.9}
\int_\tau^{\tau+\kappa}\left (e^{(H(x,b)-H(x,a))}-1\right )dx\leq I(M,a,b),
\end{equation}
where
$$I(M,a,b):=\int_{-\infty}^\infty \left( \left(\frac{x^2+b^2}{x^2+a^2}\right)^{M/2\pi}-1\right)dx<\infty.$$
\noindent
By substituting (\ref{eqn 2.6.7}) into (\ref{eqn 2.6.9}), we find that
\begin{equation}
\label{eqn 2.6.10}
\int_{\Gamma_a}\left (\frac{|{\nabla} U(p_a)|^2}{|{\nabla} U(p_b)|}\right)ds_a=
\int_\tau^{\tau+\kappa}\left (\frac{|F'(x+jb)|}{|F'(x+ja)|}\right )dx
\leq (I(M,a,b)+\kappa)\,e^{C(b^2-a^2)}
\end{equation}
for any $0<a<b\leq 1/2$, where $j=\sqrt{-1}$, and where, for any fixed $\alpha\in (0,1/2]$, $p_\alpha:=F(x+j\alpha)$ is the $\kappa$-periodic mapping of the line $\Re\times\{\alpha\}$ onto the level curve of $U_i(p)=U_i({\bf\Gamma};p)$ at altitude $\alpha$, and the arc $\Gamma_\alpha$ is one $P$-period (in $x$) of the image of that mapping. At this point, one can see that if ${\bf\Gamma}$ is not analytic, then ${\bf\Gamma}$ can be approximated by sequence of analytic curve-pairs in ${\bf X}$ having the same total-curvature bounds, and that the estimate (\ref{eqn 2.6.10}) follows in the limit. We now reintroduce the parameter $\varepsilon\in S$. The assertion (\ref{eqn 2.6.5}) follows from (\ref{eqn 2.6.10}) by setting $U_i=U_{\varepsilon,i}$, $a=\delta\varepsilon$, and $b=\varepsilon$, and by observing that (a) the constant $C=C_\varepsilon$ can be chosen independent of small $\varepsilon\in S$, and (b) $|{\nabla} U_{\varepsilon,i}(p)|$ is bounded above on $\Gamma_\varepsilon:=\{U_{\varepsilon,i}(p)=\varepsilon\}$, with an upper bound independent of small $\varepsilon\in S$ (see Lem. \ref{lem 2.2.1}(a)), and (c) $I(M,\delta\varepsilon,\varepsilon)$ is bounded from above as $\varepsilon\rightarrow 0+$ for fixed $\delta\in(0,1)$ and $M\in\Re_+$.
\vspace{.1in}

\noindent
{\bf Proof of Thm. \ref{thm 2.6.1}.} In addition to the notation in the statement of Thm. \ref{thm 2.6.1}, we also assume that ${\bf\Gamma}_\varepsilon\rightarrow{\bf\Gamma}=(\Gamma_1,\Gamma_2)\in{\boldsymbol{\cal F}}$ as $\varepsilon\rightarrow 0+$ in $S$, where $S$ denotes a positive, countable subset of $(0,1/2]$, with sole accumulation point at $0$, and where ${\bf\Gamma}_\varepsilon=(\Gamma_{\varepsilon,1},\Gamma_{\varepsilon,2})\in{\bm{\cal F}}_\varepsilon$ is a "fixed point" of one of the operators ${\bf T}_\varepsilon^\pm$ for each $\varepsilon\in S$. In addition to the notation in the statement of Thm. \ref{thm 2.6.1}, for any $\varepsilon\in S$, we define the "full domains" $\Omega_\varepsilon:=\Omega({\bm\Gamma}_\varepsilon)$, the capacitary potentials $U_{\varepsilon,i}(p):=U_i({\bm \Gamma}_{\varepsilon};p):{\rm Cl}(\Omega_\varepsilon)\rightarrow\Re$, and the "half-domains" $\hat{\Omega}_{\varepsilon,i}:=\{p\in\Omega_\varepsilon:U_{\varepsilon,i}(p)<1/2\}$ and "center arcs" $\hat{\Gamma}_\varepsilon:=\{U_{\varepsilon,i}(p)=1/2\}$. For $\varepsilon\in S$ and $\delta\in(0,1)$ we define the arcs $\Gamma_{\varepsilon,\delta,i}:=\Phi_{\varepsilon\delta,i}({\bm
\Gamma}_\varepsilon)=\{U_{\varepsilon,i}(p)=\varepsilon\delta\}$, 
the full domains $\Omega_{\varepsilon,\delta,i}:=\{p\in\Omega_\varepsilon:U_{\varepsilon,i}(p)>\varepsilon\delta\}$, the harmonic mappings $U_{\varepsilon,\delta,i}(p):{\rm Cl}(\Omega_{\varepsilon,\delta,i})\rightarrow\Re$ such that $U_{\varepsilon,\delta,i}(p):=\big((U_{\varepsilon,i}(p)-\varepsilon\delta)\big/(1-2\varepsilon\delta)\big)$, the half domains $\hat{\Omega}_{\varepsilon,\delta,i}:=\{\varepsilon\delta<U_{\varepsilon,i}(p)<1/2\}$, and center arcs $\hat{\Gamma}_{\varepsilon}$. 
\vspace{.1in}

\noindent
In the notation of the proof of Lem. \ref{lem 2.4.1}, 
the function $E_{\varepsilon,0}(z):={\rm ln}\big(\varepsilon\big/\big|g_\varepsilon(f_\varepsilon(z)+\varepsilon)-z\big|\big)$ is harmonic for sufficiently small $\varepsilon\in S$. It follows that the continuous function
\begin{equation}
\label{eqn 2.6.11'}
\phi_{\varepsilon,i}(p):={\rm ln}\big(\varepsilon\big/\big|\Pi_{\varepsilon,i}(p)-p\,\big|\big):{\rm Cl}(\dot{\Omega}_{\varepsilon,i})\rightarrow\Re
\end{equation}
is harmonic in the domain ${\dot{\Omega}}_{\varepsilon,i}:=\{p\in\Omega_\varepsilon:U_{\varepsilon,i}(p)\leq 1-\varepsilon\}$ for any sufficiently small $\varepsilon\in S$, where for each $p\in\dot{\Omega}_{\varepsilon,i}$, $\Pi_{\varepsilon,i}(p)$ denotes the point of intersection of the arc of steepest ascent of $U_{\varepsilon,i}$ through $p$ with the level curve of $U_{\varepsilon,i}$ at the altitude $U_{\varepsilon,i}(p)+\varepsilon$. In view of this definition, it 
follows by Prop. \ref{prop 2.5.1} and Cor. \ref{cor 2.5.1} (especially Eqs. (\ref{eqn 2.5.1a'}) and (\ref{eqn 2.5.2a'}), where $\big|\tilde{p}_{\varepsilon,\delta,i}(s)-p_{\varepsilon,\delta,i}(s)\big|=\big|\Pi_{\varepsilon,i}(p_{\varepsilon,\delta,i}(s))-p_{\varepsilon,\delta,i}(s)\big|$ for all $s\in\Re$) that 
\begin{equation}
\label{eqn 2.6.15}
0\leq\big|\Pi_{\varepsilon,i}(p)-p\big|-{\rm dist}\big(p,\tilde{\Gamma}_{\varepsilon,\delta,i}\big)\leq f_{\varepsilon,\delta,i}(p)=\varepsilon\,\, z_{\varepsilon,\delta,i}(p),
\end{equation}
\noindent
uniformly for all $p\in\Gamma_{\varepsilon,\delta,i}$ and all sufficiently small $\varepsilon\in S$ and $\delta\in(0,1)$, where the measurable functions $z_{\varepsilon,\delta,i}(p):\Gamma_{\varepsilon,\delta,i}\rightarrow\Re$ satisfy 
\begin{equation}
\label{eqn 2.6.15c}
0\leq z_{\varepsilon,\delta,i}(p)\leq C_1\,\,{\rm in}\,\,\Gamma_{\varepsilon,\delta,i}\,\,{\rm and}\,\,\int_{\Gamma_{\varepsilon,\delta,i}}\,z_{\varepsilon,\delta,i}(p)\,ds\leq M\sqrt{\varepsilon},
\end{equation} 
and where $C_1$ and $M$ denote uniform constants independent of sufficiently small $\varepsilon\in S$ and $\delta\in(0,1)$. (Here the arc integral refers to one $P$-period (in $x$).) 
\vspace{.1in}

\noindent
We have $({\bm{\rm i}})$ $a_i(p)\,{\rm dist}\big(p,\tilde{\Gamma}_{\varepsilon,i}\big)=\varepsilon$ for all sufficiently small $\varepsilon\in S$ and all points $p\in\Gamma_{\varepsilon,i}$ (see Thm. \ref{thm 2.1.1}, Eq. (\ref{eqn 2.1.13})). Therefore, for any fixed value $\delta\in(0,1)$, there exists a value $\eta=\eta(\delta)>0$ so small that 
\begin{equation}
\label{eqn 2.6.12}
a_{\delta,i}(p)\,\,{\rm dist}(p,\tilde{\Gamma}_{\varepsilon,i})\leq\varepsilon
\end{equation}
uniformly for all sufficiently small $\varepsilon\in S$ and all points $p\in\Gamma_{\varepsilon,\delta,i}$, where we define $a_{\delta,i}(p):=\eta(\delta)\,a_i(p)$ throughout $\Re^2$. For the remainder of the proof, we choose the value $\eta=\eta(\delta)$ (corresponding to any given value $\delta\in(0,1)$) to be as large as possible subject to the requirement that (\ref{eqn 2.6.12}) holds for all sufficiently small $\varepsilon\in S$ and all $p\in\Gamma_{\varepsilon,\delta,i}$. Toward an upper bound for $\eta(\delta)$, given any point $p\in\Gamma_{\varepsilon,\delta,i}$, we choose the point $q=q_{\varepsilon,\delta,i}(p)\in\Gamma_{\varepsilon,i}$ such that $p$ and $q$ both lie on the same arc of steepest ascent of the function $U_{\varepsilon,i}$. By Cor. \ref{cor 2.2.1}, there exists a positive constant $C$, independent of small $\varepsilon\in S$, such that $({\bm{\rm ii}})$ $|q_{\varepsilon,\delta,i}(p)-p|\leq (\varepsilon\delta/C)$. By using (\ref{eqn 2.6.12}), $({\bm{\rm i}})$, and $({\bm{\rm ii}})$, we see that
\begin{equation}
\label{eqn 2.6.13}
\varepsilon\geq a_{\delta,i}(p)\,{\rm dist}(p,\tilde{\Gamma}_{\varepsilon,i})\geq a_{\delta,i}(p)\big({\rm dist}(q,\tilde{\Gamma}_{\varepsilon,i})-|p-q|\big)
\end{equation}
$$\geq a_{\delta,i}(p)\big([\varepsilon/a_i(q)]-(\varepsilon\delta/C)\big)$$
for all $p\in\Gamma_{\varepsilon,\delta,i}$, where $q=q_{\varepsilon,\delta,i}(p)$. It is easily seen that if Eq. (\ref{eqn 2.6.12}), and therefore (\ref{eqn 2.6.13}), is satisfied for all sufficiently small $\varepsilon\in S$ at a particular value $\delta\in(0,1)$ such that $a_i(p)\,\delta\leq\overline{A}\delta<C$ for all $p\in\Gamma_{\varepsilon,i}$, then we have that 
\begin{equation}
\label{eqn 2.6.14}
\eta(\delta)<\liminf_{\varepsilon\rightarrow 0+}\Big(\frac{C\,a_i\big(q_{\varepsilon,\delta,i}(p)\big)}{a_i(p)\big[C-a_i\big(q_{\varepsilon,\delta,i}(p)\big)\delta\big]}\Big)\leq\frac{C}{C-a_i(p)\delta}\leq\frac{C}{C-\overline{A}\delta},
\end{equation}
for the same constant $C>0$.
\vspace{.1in}

\noindent
Due to the fact that the function $\phi_{\varepsilon,i}(p):{\rm Cl}(\dot{\Omega}_{\varepsilon,i})\rightarrow\Re$, which was defined in Cor. \ref{cor 2.4.1aaa}, Eq. (\ref{eqn 2.6.11}), is harmonic in $\dot{\Omega}_{\varepsilon,i}$, the related continuous function  
\begin{equation}
\label{eqn 2.6.16}
\phi_{\varepsilon,\delta,i}(p):={\rm ln}\left(\frac{\varepsilon}{\hat{a}_{\varepsilon,\delta,i}(p)\big|\Pi_{\varepsilon,i}(p)-p\big|}\right):{\rm Cl}(\hat{\Omega}_{\varepsilon,\delta,i})\rightarrow\Re
\end{equation}
\noindent
is harmonic at least in the domain $\hat{\Omega}_{\varepsilon,\delta,i}$,
where we define the continuous function $\hat{a}_{\varepsilon,\delta,i}(p):{\rm Cl}\big(\hat{\Omega}_{\varepsilon,\delta,i}\big)\rightarrow\Re$ such that $\hat{a}_{\varepsilon,\delta,i}(p)=a_{\delta,i}(p):=\eta(\delta)a_i(p)$ on $\partial\hat{\Omega}_{\varepsilon,\delta,i}$ and ${\rm ln}\big(\hat{a}_{\varepsilon,\delta,i}(p)\big)$ is harmonic in $\hat{\Omega}_{\varepsilon,\delta,i}$. (Here, the value $\delta\in(0,1)$ is fixed and sufficiently small, and $\eta(\delta)$ satisfies (\ref{eqn 2.6.14}).) It follows from (\ref{eqn 2.6.15}), (\ref{eqn 2.6.15c}), (\ref{eqn 2.6.12}), and (\ref{eqn 2.6.16}) that $\hat{a}_{\varepsilon,\delta,i}(p)\big|\Pi_{\varepsilon,i}(p)-p\big|\leq \big(1+\hat{a}_{\varepsilon,\delta,i}(p)\,z_{\varepsilon,\delta,i}(p)\big)\,\varepsilon,$
and therefore that
\begin{equation}
\label{eqn 2.6.17}
\phi_{\varepsilon,\delta,i}(p)\geq {\rm ln}\big(1/(1\,+\,a_{\delta,i}(p)\,z_{\varepsilon,\delta,i}(p)\big)\geq-a_{\delta,i}(p)\,z_{\varepsilon,\delta,i}(p)
\end{equation}
\noindent
both for all points $p\in\Gamma_{\varepsilon,\delta,i}$. We also consider the continuous function
\begin{equation}
\label{eqn 2.6.18}
\psi_{\varepsilon,\delta,i}(p):=\phi_{\varepsilon,\delta,i}(p)+2\,C_{\varepsilon,\delta,i}\,{U}_{\varepsilon,\delta,i}(p):{\rm Cl}(\hat{\Omega}_{\varepsilon,\delta,i})\rightarrow\Re
\end{equation}
\noindent
(harmonic in $\hat{\Omega}_{\varepsilon,\delta,i}$), where we define 
\begin{equation}
\label{eqn 2.6.19}
{U}_{\varepsilon,\delta,i}(p):=\left(2(U_{\varepsilon,i}(p)-\varepsilon\delta)/(1-2\varepsilon\delta)\right),
\end{equation}
$${\nabla} {U}_{\varepsilon,\delta,i}(p)=(2{\nabla} U_{\varepsilon,i}(p)/(1-2\varepsilon\delta)),
$$
and where we use the constant $C_{\varepsilon,\delta,i}$ to denote the least value such that
\begin{equation}
\label{eqn 2.6.20}
C_{\varepsilon,\delta,i}\geq -\phi_{\varepsilon,\delta,i}(p)={\rm ln}\Big(\big(a_{\delta,i}(p)|\Pi_{\varepsilon,i}(p)-p|\big)\big/\varepsilon\Big)
\end{equation}
uniformly for all $p\in\hat{\Gamma}_{\varepsilon}$. Then we have 
\begin{equation}
\label{eqn 2.6.21}
\psi_{\varepsilon,\delta,i}(p)\geq 0 \,\,{\rm on}\,\,\hat{\Gamma}_{\varepsilon};\,\,\psi_{\varepsilon,\delta,i}(p)\geq-a_{\delta,i}(p)\,z_{\varepsilon,\delta,i}(p)\,\,{\rm on}\,\,\Gamma_{\varepsilon,\delta,i},
\end{equation}
\noindent
due to (\ref{eqn 2.6.17}), (\ref{eqn 2.6.18}), (\ref{eqn 2.6.19}), and (\ref{eqn 2.6.20}).
\vspace{.1in}

\noindent
We also introduce the $P$-periodic (in $x$) Green's function $G_{\varepsilon,\delta,i}(p_{0,i},q)\geq 0$ for the Laplace operator in the $P$-periodic domain $\hat{\Omega}_{\varepsilon,\delta,i}$. For any specified point $p_{0,i}\in\hat{\Omega}_i$ and value $\rho>0$ such that ${\rm Cl}\big(B_\rho(p_{0,i})\big)\subset\hat{\Omega}_i$, and for any sufficiently small values $\varepsilon\in S$ and $\delta\in(0,1)$ (so small that ${\rm Cl}\big(B_\rho(p_{0,i})\big)\subset\hat{\Omega}_{\varepsilon,\delta,i}$), we use $G_{\varepsilon,\delta,i}(p_{0,i},q)$ to denote the unique $P$-periodic (in $x$) continuous function of $q$ in ${\rm Cl}(\hat{\Omega}_{\varepsilon,\delta,i})$ such that 

$$\Delta_q {G}_{\varepsilon,\delta,i}(p_{0,i},q)+\sum_{n=-\infty}^\infty\delta(q-p_{0,i}-n(P,0))=0$$
in $\hat{\Omega}_{\varepsilon,\delta,i}$ and ${G}_{\varepsilon,\delta,i}(p_{0,i},q)=0 \,\,{\rm\,\,for\,\, all\,\,}\,\,q\in\partial\hat{\Omega}_{\varepsilon,\delta,i}$.
where $\delta(\cdot)$ denotes the Dirac delta function.
We also use $C_\rho(p_{0,i})$ to denote a positive constant such that 

\begin{equation}
\label{eqn 2.6.22}
{G}_{\varepsilon,\delta,i}(p_{0,i},q)\leq C_\rho(p_{0,i})\,{U}_{\varepsilon,\delta,i}(q)
\end{equation}
for all $q\in \partial Q_\rho(p_{0,i})$, and uniformly for all sufficiently small values $\varepsilon\in S$ and $\delta\in(0,1)$, where we define the $P$-periodic set $Q_\rho(p_{0,i}):= \big(\bigcup_{n=-\infty}^\infty B_\rho(p_{0,i}+n(P,0))\big)$. Since both sides of (\ref{eqn 2.6.22}) are harmonic functions in $\hat{\Omega}_{\varepsilon,\delta,i}\setminus Q_\rho(p_{0,i})$ which vanish on $\Gamma_{\varepsilon,\delta,i}$, and since ${G}_{\varepsilon,\delta,i}(p_{0,i},q)=0\leq {U}_{\varepsilon,\delta,i}(q)=1/2$ for all $q\in\hat{\Gamma}_{\varepsilon}$, it follows that the inequality (\ref{eqn 2.6.22}) holds throughout $q\in\hat{\Omega}_{\varepsilon,\delta,i}\setminus Q_\rho(p_{0,i})$ by the maximum principle. Therefore, we also have
 
\begin{equation}
\label{eqn 2.6.23}
|{\nabla}_q {G}_{\varepsilon,\delta,i}(p_{0,i},q)|\leq C_\rho(p_{0,i})\,|{\nabla} {U}_{\varepsilon,\delta,i}(q)|
\end{equation}
$$=(2C_\rho(p_{0,i})/(1-2\varepsilon\delta))\,|{\nabla}{U}_{\varepsilon,i}(q)|$$
\noindent
for $q\in\Gamma_{\varepsilon,\delta,i}$.
\noindent
Since $\psi_{\varepsilon,\delta,i}(q)$ is harmonic throughout $\hat{\Omega}_{\varepsilon,\delta,i}$, it follows from Green's second identity that

\begin{equation}
\label{eqn 2.6.24}
\psi_{\varepsilon,\delta,i}(p_{0,i})=\int_{\partial\hat{\Omega}_{\varepsilon,\delta,i}}\big(\partial{G}_{\varepsilon,\delta,i}(p_{0,i},q)\big/\partial{\nu}_q\big)
\,\psi_{\varepsilon,\delta,i}(q)\,d\hat{s}_{\varepsilon,\delta,i},
\end{equation}
where $\hat{\Omega}_{\varepsilon,\delta,i}$ refers to one $P$-period (containing the point $p_{0,i}$) of the corresponding $P$-periodic region, and $\partial\big/\partial\boldsymbol{\nu}_q$ refers to the normal derivative at $q\in\partial\hat{\Omega}_{\varepsilon,\delta,i}$ in the direction of the interior normal vector. Since $(\partial\big/\partial\boldsymbol{\nu}_q)\,G_{\varepsilon,\delta,i}(p_{0,i},q)=|{\nabla}_q\,G_{\varepsilon,\delta,i}(p_{0,i},q)|\geq 0$ for all $q\in\Gamma_{\varepsilon,\delta,i}\cup\hat{\Gamma}_{\varepsilon}$, it follows from (\ref{eqn 2.6.24}) via the inequalities (\ref{eqn 2.6.21} a,b) and (\ref{eqn 2.6.23}) that

\begin{equation}
\label{eqn 2.6.25}
\psi_{\varepsilon,\delta,i}(p_{0,i})\geq-\int_{\Gamma_{\varepsilon,\delta,i}}\big|{\nabla}_q\,{G}_{\varepsilon,\delta,i}(p_{0,i},q)\big|\,a_{\delta,i}(q)\,z_{\varepsilon,\delta,i}(q)\,ds_{\varepsilon,\delta,i}
\end{equation}

$$\geq-\big(2C_{\rho}(p_{0,i})\big/{(1-2\varepsilon\delta)\big)}\int_{\Gamma_{\varepsilon,\delta,i}}\big|{\nabla}_q\,
{U}_{\varepsilon}(q)\,\big|\,\,a_{\delta,i}(q)\,z_{\varepsilon,\delta,i}(q)\,ds_{\varepsilon,\delta,i}$$

$$\geq-\frac{2C_\rho(p_{0,i})}{1-2\varepsilon\delta}\sqrt{\int_{\Gamma_{\varepsilon,\delta,i}}\big|{\nabla}_q\, {U}_{\varepsilon,i}(q)\big|^2\,ds_{\varepsilon,\delta,i}}\sqrt{\int_{\Gamma_{\varepsilon,\delta,i}}a_{\delta,i}^2(q)\,z_{\varepsilon,\delta,i}^2(q)\,ds_{\varepsilon,\delta,i}},$$
where all the arc integrals are restricted to one $P$-period (in $x$) and the final step in (\ref{eqn 2.6.25}) is based on the Schwartz inequality.
\vspace{.1in}

\noindent
Finally, we observe that the first integral in the final line of (\ref{eqn 2.6.25}) remains uniformly bounded as $\varepsilon\rightarrow 0+$ provided that the constant $\delta\in(0,1)$ is sufficiently small (see Lem. \ref{lem 2.6.1}, Eq. (\ref{eqn 2.6.5})), while, concerning the second integral, we have
 
$$a_{\delta,i}(p)=\eta(\delta)\,a_i(p)\leq \big(C\overline{A}\,\big/(C-\overline{A}\,\delta)\big)$$
in $\Re^2$ uniformly for all $0<\delta<\min\{1,(C/\,\overline{A}\,)\}$ (see (\ref{eqn 2.6.14})), while it follows from (\ref{eqn 2.6.15c}) that
$$\int_{\Gamma_{\varepsilon,\delta,i}}\,z_{\varepsilon,\delta,i}^2(p)\,ds_{\varepsilon,\delta,i}\leq C_1\int_{\Gamma_{\varepsilon,\delta,i}}\,z_{\varepsilon,\delta,i}(p)\,ds_{\varepsilon,\delta,i}\leq C_1\,M\sqrt{\varepsilon},$$
uniformly for sufficiently small $\delta\in(0,1)$ and $\varepsilon\in S$. Therefore, for any sufficiently small $\delta\in(0,1)$, the second integral approaches $0$ as $\varepsilon\rightarrow 0+$ relative to $S$. For any specified point $p_{0,i}\in\hat{\Omega}_i$, we abbreviate Eq. (\ref{eqn 2.6.24}) to state that 
$$\psi_{\varepsilon,\delta,i}(p_{0,i})\geq -\zeta(\varepsilon,\delta,i),$$
where $\zeta(\varepsilon,\delta,i)\rightarrow 0$ as $\varepsilon\rightarrow 0+$ relative to $S$ for any fixed, sufficiently small value $\delta\in(0,1)$. By substituting the definition of the function $\psi_{\varepsilon,\delta,i}(p)$ (see (\ref{eqn 2.6.17}), (\ref{eqn 2.6.18}), (\ref{eqn 2.6.19}), and (\ref{eqn 2.6.20})), we obtain the equivalent inequality:

\begin{equation}
\label{eqn 2.6.26}
\frac{\varepsilon}{\big|\Pi_{\varepsilon,i}(p_{0,i})-p_{0,i}\big|}\geq \hat{a}_{\varepsilon,\delta,i}(p_{0,i})\,{\rm exp}\big(-2\,C_{\varepsilon,\delta,i}\,{U}_{\varepsilon,\delta,i}(p_{0,i})\big)\,{\rm exp}\big(-\zeta(\varepsilon,\delta,i)\big),
\end{equation}
valid for any point $p_{0,i}\in\hat{\Omega}_i$ if $\varepsilon\in S$ is sufficiently small (depending on $p_{0,i})$. We have that $U_{\varepsilon,\delta,i}(p)\rightarrow U_i(p)$ and ${\nabla} {U}_{\varepsilon,\delta,i}(p)\rightarrow{\nabla}{U}_i(p)$, both in any compact subset of $\hat{\Omega}_i$ as $\varepsilon\rightarrow 0+$ relative to $S$ (see Def. \ref{def 2.3.1} and Thm \ref{thm 2.3.1}), while, by the Theorem of the mean, we also have for any point $p_{0,i}\in\hat{\Omega}_i$ that

$$\big(\varepsilon\big/\big|\Pi_{\varepsilon,i}(p_{0,i})-p_{0,i}\big|\big)\leq|{\nabla} U_{\varepsilon,i}(p_{0,i}^*)|$$
for some point $p_{0,i}^*$ in the line-segment joining $p_{0,i}$ to $\Pi_{\varepsilon,i}(p_{0,i})$.
In view of this, it follows from (\ref{eqn 2.6.26}) in the limit as 
$\varepsilon\rightarrow 0+$ relative to $S$ that if $\delta\in(0,1)$ is sufficiently small, then
\begin{equation}
\label{eqn 2.6.27}
\big|{\nabla} U_i(p_{0,i})\big|\geq \hat{a}_{\delta,i}(p_{0,i})\,\,{\rm exp}\big(-2\,C_{\delta,i}\,U_i(p_{0,i})\big)
\end{equation}
throughout $\hat{\Omega}_i$, where, in terms of notation introduced in the statement of Thm. \ref{thm 2.6.1}, we have set $\hat{a}_{\delta,i}(p_{0,i})=\eta(\delta)\,\hat{a}_i(p_{0,i})$ in $\hat{\Omega}_i$ (in terms of notation introduced in the statement of Thm. \ref{thm 2.6.1}), and where ${\rm exp}\big(C_{\varepsilon,\delta,i}\big)$ $\rightarrow{\rm exp}\big(C_{\delta,i}\big)$ $=\eta(\delta)\,{\rm exp}\big(C_i\big)$ as $\varepsilon\rightarrow 0+$, where we define $C_i$ such that ${\rm exp}\big(C_i\big):=\sup\big\{\big({a}_{i}(q)\big/\big|\nabla U_{i}(q)\big|\big):$ $q\in\hat{\Gamma}\big\}$. For $p$ near $\Gamma_i$ in $\hat{\Omega}_i$, the inequality (\ref{eqn 2.6.27}) reduces to 

\begin{equation}
\label{eqn 2.6.28}
|{\nabla} U_i(p)|\geq \hat{a}_i(p)-z\big({\rm dist}\big(p,\Gamma_i\big)\big),
\end{equation}
where $z(t)\rightarrow 0$ as $t\rightarrow 0+$, while, relative to the arc $\hat{\Gamma}$, (\ref{eqn 2.6.27}) reduces to the inequality: 
\begin{equation}
\label{eqn 2.6.29}
|{\nabla} U_i(p)|\geq \hat{a}_i(p)\,{\rm exp}\big(-C_i\big)=\hat{a}_i(p)\,\inf\big\{|{\nabla} U_i(p)|/a_i(p):p\in\hat{\Gamma}\big\}
\end{equation}
\vspace{.1cm}
\noindent
Since ${\rm ln}\big(|{\nabla} U_i(p)|\big/\hat{a}_i(p)\big)$ is harmonic in $\hat{\Omega}_i$, it follows from (\ref{eqn 2.6.28}) and (\ref{eqn 2.6.29}) by the maximum principle that Eq. (\ref{eqn 2.6.2}) holds, as was asserted.

\subsection{Existence of classical solutions for Probs. \ref{prob 1.1} and \ref{prob2.1.1}}
\label{subsection 3.4}
\begin{theorem}
\label{thm 2.6.2}{(Extension of Thms. \ref{thm1.2} and \ref{thm2.1.0}; existence of classical solutions between strict lower and strict upper solutions)}
In the context of Prob. \ref{prob2.1.1}, let be given a strict lower solution (or sub-solution) $\tilde{\bf\Gamma}\in {\bf X}\cap C^2$ and a strict upper solution (or super-solution) $\hat{{\bf\Gamma}}\in{\bf X}\cap C^2$ such that $\tilde{\bf\Gamma}<\hat{{\bf\Gamma}}$. Then:
\vspace{.1in}

\noindent
(a) There exists at least one weak solution ${\bm \Gamma}=(\Gamma_1,\Gamma_2)\in{\boldsymbol{\cal F}}$ such that 
$\tilde{{\bm\Gamma}}<{\bm\Gamma}<\hat{{\bm\Gamma}}$, and such that, for $i=1,2$, the arc-length per $P$-period and total curvature per $P$-period of $\Gamma_i$ are both bounded by $M$ (i.e. $||\Gamma_i||\leq M$ and $K(\Gamma_i)\leq M$), where $M$ denotes a constant depending only on $\underline{A},\overline{A},A_1, A_2$. 
\vspace{.1in}

\noindent
(b) Any weak solution ${\bm\Gamma}$ with the properties in Part (a) is a classical solution in the sense that 
${\bf\Gamma}=(\Gamma_1,\Gamma_2)\in{\bf X}\cap C^1$ and that, for the function $U(p):=U({\bf\Gamma};p)$, the derivative ${\nabla} U(p)$ has a continuous extension to the closure of $\Omega:=\Omega({\bf \Gamma})$ such that $|{\nabla} U(p)|=a_i(p)$ point-wise in $\Gamma_i, i=1,2$. 
\end{theorem}

\noindent
{\bf Proof.} In the context of Def. \ref{def 2.1.5}, it follows from Thm. \ref{thm 2.1.3} (see also Thm. \ref{thm 2.1.2}, Thm. \ref{thm 3.1.2}(c)), Eq. (\ref{eqn 3.1.6a}) and Lem. \ref{lem 2.2.5}, Eq. (\ref{eqn 2.2.X})) that there exists a constant $M$ such that for any sufficiently small value $\varepsilon\in(0,\varepsilon_1]$, the operator ${\bf T}_\varepsilon$ (denoting either ${\bf T}_{\varepsilon}^+$ or ${\bf T}_{\varepsilon}^-)$ has a fixed point ${\bf\Gamma}_{\varepsilon}=(\Gamma_{\varepsilon,1},\Gamma_{\varepsilon,2})\in{\bf Y}$ such that $||\Gamma_{\varepsilon,1}||, ||\Gamma_{\varepsilon,2}||, K(\Gamma_{\varepsilon,1}), K(\Gamma_{\varepsilon,2})\leq M$. In view of this, it follows from Thm. \ref{thm 2.3.1} that, given any null-sequence $\big(\varepsilon(n)\big)_{n=1}^\infty$ of values in $(0,\varepsilon_1]$, there exists a subsequence $\big(\varepsilon(n(k))\big)_{k=1}^\infty$ of the natural numbers such that the pairs ${\bf \Gamma}_k:={\bf\Gamma}_{n(k)}\in{\bf Y}$ converge as $k\rightarrow\infty$ to a weak solution ${\bf\Gamma}\in{\bf Y}$. In view of Thm. \ref{thm 2.1.3}, Thm. \ref{thm 2.4.1} implies that $\limsup_{p\rightarrow p_i}|{\nabla} U(p)|\leq a_i(p_i)$ for $p\in\Omega$ and $i=1,2$, and Thm. \ref{thm 2.6.1} implies that $\liminf_{p\rightarrow p_i}|{\nabla} U(p)|\geq a_i(p_i)$, for $p\in\Omega$ and $i=1,2$, where $U(p):=U({\bf\Gamma};p)$ in ${\rm Cl}(\Omega)$. The smoothness of ${\bf\Gamma}$ follows from Thm. \ref{thm 2.7.1} below. It also follows from Thm. \ref{thm 2.1.3}(b) that $||\Gamma_1||, ||\Gamma_2||, K(\Gamma_1), K(\Gamma_2)\leq M$.

\begin{theorem} 
{(Uniform curvature bounds for solutions)}
\label{thm 2.7.1}
For any given classical solution ${\bf \Gamma}=(\Gamma_1,\Gamma_2)$ of Prob. \ref{prob2.1.1} (which exists by Thm. \ref{thm 2.6.2}), let $U(p)=U({\bf\Gamma};p)$ denote the capacitary potential in the closure of the domain $\Omega:=\Omega({\bm\Gamma})$. Then $\big|{\nabla}\big({\rm ln}(|{\nabla} U(p)|)\big)\big|\leq L_0$ and $\big|{\nabla}\big({\rm arg}({\nabla} U(p))\big)\big|\leq L_0$, both in $\Omega:=\Omega({\bf\Gamma})$, where the bound $L_0$ depends only on $\overline{A}, \underline{A}, A_1$, and $A_2$, in fact we can choose: 
\begin{equation}
\label{eqn 2.7.5aa}
L_0={\rm ln}\big(\,\overline{A}\,\big/\,\underline{A}\,\big)\,+\,\big(\,A_1\,\big/\underline{A}\,\big)+\big(2\,\overline{A}\,\big(\,\overline{A}\,A_2\,+\,A_1^2)\big)\big/\,\underline{A}^4\big).
\end{equation}
It follows that the conjugate harmonic functions: ${\rm ln}(|{\nabla} U(p)|)$ and ${\rm arg}({\nabla} U(p))$, are uniformly Lipschitz-continuous in ${\rm Cl}(\Omega)$, with the same Lipschitz constant $L_0$ 
(Here $|{\nabla} U|=a_i(p)$ on $\Gamma_i$.) In particular, the curvatures of the level curves $\Gamma_\alpha:=\{U(x)=\alpha\}$, $\alpha\in(0,1)$ are uniformly bounded independent of $\alpha\in (0,1)$, and the same curvature bound (dependent on $\overline{A}, \underline{A}, A_1$, and $A_2$ only) carries over to the directed boundary arcs $\Gamma_1,\Gamma_2$. 
\end{theorem}

\begin{lemma} 
{(Boundary derivative estimates)}
\label{lem 2.7.1}
(a) For any classical solution ${\bf \Gamma}\in{\bf X}$ of Prob. \ref{prob2.1.1} and any continuous function $f(p):{\rm Cl}(\Omega)\rightarrow\Re$ (where $\Omega:=\Omega({\bf\Gamma})$), the solution $\phi(p)$ of the Dirichlet problem: $\Delta \phi(p)=f(p)$ in $\Omega$, $\phi(\partial\Omega)=0$ is such that $|\phi(p)|\leq\big({\boldsymbol |}f{\boldsymbol |}\big/2\,\underline{A}^2\big)$ and $|\nabla \phi|\leq \big(2\overline{A}{\boldsymbol |}f{\boldsymbol |}\big/\underline{A}^2\big)$, both in $\Omega$, where ${\boldsymbol |}f{\boldsymbol |}:=\sup\{|f(p)|:p\in\Omega\}$. 
\vspace{.1in}

\noindent
(b) Let $\phi(p)$ solve the above Dirichlet problem, where $f(p)\geq \varepsilon$ for some $\varepsilon>0$. Then $|\nabla\phi(p)|\geq (\underline{A}\varepsilon/4\overline{A}^2)$ for all $p\in\partial\Omega$.
\end{lemma}

\noindent
{\bf Proof.} We map $\Omega$ conformally onto the strip $(0,1)\times\Re$ under the analytic mapping: $w=u+jv=F(z)=U(z)+jV(z)$, where $V(z)$ is the harmonic conjugate of $U:=U({\bf\Gamma};z)$ in $\Omega$. We define the function $\psi(u,v):(0,1)\times\Re\rightarrow\Re$ such that $\phi(z)=\psi(F(z))$ in $\Omega$. Then 
\begin{equation}
\label{eqn 2.7.6}
|{\nabla}_w \psi|=(|{\nabla}_z \phi|/|F'(z)|),
\end{equation}
\begin{equation}
\label{eqn 2.7.7}
f(x,y)=\big(\partial_x^2+\partial_y^2\big)\,\phi(x,y)=|F'(z)|^2\big(\partial_u^2+\partial_v^2\big)\,\psi(u,v),
\end{equation}
\noindent
where
\begin{equation}
\label{eqn 2.7.8}
0<\underline{A}\leq|F'(z)|\leq \overline{A}
\end{equation}
\noindent
in $\Omega$ (by the maximum principle, since ${\rm ln}(|F'|)$ is harmonic in $\Omega$ and (\ref{eqn 2.7.8}) holds on $\partial\Omega$). By (\ref{eqn 2.7.6}), (\ref{eqn 2.7.7}), and maximum principles for the Poisson equation, we have 
\begin{equation}
\label{eqn 2.7.9}
|\psi(u,v)|\leq \Psi(u,v):=\big({\boldsymbol |}f{\boldsymbol |}/2\underline{A}^2\big)u(1-u)
\end{equation}
\noindent
in $(0,1)\times\Re$. The first estimate in Part (a) follows from this. Also, it follows from (\ref{eqn 2.7.9}) and a gradient estimate for solutions of the Poisson Equation: $\Delta\psi=g$ (namely $|\nabla \psi(z_0)|\leq(2\sqrt{2}/r)\sup\{|\psi(z)|:|z-z_0|\leq r\}+(r/2\sqrt{2})|g|$; see [GT], Sect. 3.4, p. 37) that $|\nabla_w\psi(u,v)|\leq\big(\boldsymbol{|}f\boldsymbol{|}\big/2\underline{A}^2\big)$ uniformly in $(0,1)\times\Re$, from which the second estimate in Part (a) follows by (\ref{eqn 2.7.6}) and (\ref{eqn 2.7.8})).

\begin{remark}
\label{rem 7.4.3+}
{(Strengthening of Thm. \ref{thm 2.4.1}.)} 
After the estimates (\ref{eqn 2.4.1}) and (\ref{eqn 2.4.2}) have been proved, it follows from them that for any weak solution ${\bf\Gamma}\in{\boldsymbol{\cal F}}$ of Prob. \ref{prob2.1.1}, we have $|\nabla U(p)|\leq\overline{A}$ in $\Omega:=\Omega({\bf\Gamma})$, since $|\nabla U(p)|=|\nabla U({\bf\Gamma};p)|$ is sub-harmonic in $\Omega$. Also, $\hat{a}_i(p)\geq \underline{A}$ in $\hat{\Omega}_i$, since ${\rm ln}\big(\hat{a}_i(p)\big)$ is harmonic in $\hat{\Omega}_i$ and the same estimate holds on $\partial\hat{\Omega}_i$. In view of this, we can actually choose $C_i:={\rm ln}\big(\,\overline{A}\big/\underline{A}\,)$ in (\ref{eqn 2.4.2}). Thus, the estimate (\ref{eqn 2.4.2}) is improved to state that 
\begin{equation}
\label{eqn 2.7.1}
{\rm ln}\big(|{\nabla}U(p)|\big/\hat{a}_i(p)\big)\leq 2\,{\rm ln}\big(\,\overline{A}\big/\underline{A}\,\big)\,U_i(p)\leq 2\overline{A}\,{\rm ln}\big(\,\overline{A}\big/\underline{A}\,\big)\,{\rm dist}(p,\Gamma_i)
\end{equation}

\noindent
in $\hat{\Omega}_i=\{U_i(p)<1/2\}$ for any weak solution $\bf\Gamma\in{\boldsymbol{\cal F}}$.
\label{rem 2.4.1}
\end{remark}
\vspace{.1in}

\noindent
{\bf Proof of Thm. \ref{thm 2.7.1}.} To obtain uniform upper bounds for $\big|\nabla\big({\rm ln}\big(|{\nabla}U_i(p)|)\big)\big|$, $i=1,2$, we estimate the terms in the obvious inequality:
\begin{equation}
\label{eqn 2.7.3}
\big|\nabla\big({\rm ln}\big(|{\nabla} U_i(p)|)\big)\big|\leq \big|\nabla\,{\rm ln}\big(|{\nabla}U_i(p)|\big/\hat{a}_i(p)\big)\big|+\big|{\nabla}\big({\rm ln}\big(a_i(p)\big)\big|+\big|{\nabla}\hat{\phi}_i(p)\big|,
\end{equation}
\noindent
relative to $\hat{\Omega}_i:=\{p\in\Omega: U_i(p)<1/2\}$, $i=1,2$. Here, we define the function: $\hat{\phi}_i(p):={\rm ln}\big(a_i(p)\big/\hat{a}_i(p)\big):{\rm Cl}(\hat{\Omega}_i)\rightarrow\Re$, where $\hat{a}_i(p):{\rm Cl}(\hat{\Omega}_i)\rightarrow\Re$ denotes a logarithmically harmonic function which coincides with ${a_i(p)}$ on $\partial\hat{\Omega}_i$. Concerning the first term on the right side of (\ref{eqn 2.7.3}), it follows from (\ref{eqn 2.7.1}) by a standard derivative estimate that $\big|{\nabla}\big({\rm ln}\big(|{\nabla}U(p)|\big/\hat{a}_i(p))\big)\big|\leq {\rm ln}\big(\,\overline{A}\big/\underline{A}\,\big)$ near $\hat{\Gamma}_i$ in $\Omega$, $i=1,2$. For the second term, we have $\big|{\nabla}\big({\rm ln}\big(a_i(p)\big)\big|\leq\big(A_1\big/\underline{A}\,\big)$ in $\Re^2$ by assumption. For the third term, we observe that $\hat{\phi}_i(p)$ solves the Dirichlet problem: $\Delta\,\hat{\phi}_i(p)=f_i(p):=\Delta\,\big({\rm ln}\big(a_i(p)\big)$ in $\Omega$, $\hat{\phi}_i(\partial\hat{\Omega}_i)=0$. Therefore, it follows from Lem. \ref{lem 2.7.1} that 
$$\big|{\nabla}\hat{\phi}_i(p)\big|\leq 2(\,\overline{A}\big/\underline{A}^2\,){\boldsymbol|}f_i{\boldsymbol|}\leq 2\big(\,\overline{A}\big/\underline{A}^4\big)\big(\,\overline{A}A_2+A_1^2\big)$$
\noindent
in $\hat{\Omega}_i$. The assertion (\ref{eqn 2.7.5aa}) follows by substituting these three estimates into (\ref{eqn 2.7.3}).
\vspace{.1in}
\subsection{Uniform $C^{\,3,\tilde{\varrho}}$-continuations}
\label{subsection 3.5}
\begin{theorem} 
{(Smooth uniform continuations of capacitary potentials of solutions of Prob.\ref{prob2.1.1})}
\label{thm 2.7.2}
Given $\varrho\in(0,1]$, let ${\bm{\mathfrak{M}}}(\varrho)$ denote the family of all classical solutions ${\bf\Gamma}=(\Gamma_1,\Gamma_2)\in{\bf X}$ of Prob. \ref{prob2.1.1} corresponding to a family of functions $a_1(p),a_2(p)\in \boldsymbol{{\cal A}}\cap C^{\,3,\varrho}$, whose norms relative to $\boldsymbol{{\cal A}}$ and $C^{\,3,\varrho}$ are uniformly bounded from above by a single value $\overline{B}$. Then:
\vspace{.1in}

\noindent
(a) The set of all periodic conformal mappings $w=F(z)$ of the strip $\overline{\omega}:=[0,1]\times\Re$ onto the closed regions ${\rm Cl}(\Omega({\bf\Gamma}))$, ${\bf\Gamma}\in{\bm{\mathfrak{M}}}(\varrho)$ is uniformly bounded in the $C^{\,3,\tilde{\varrho}}(\overline{\omega})$-norm (for some value $\tilde{\varrho}\in(0,\varrho]$), and they can, therefore, all be continued as $C^{\,3,\tilde{\varrho}}$-functions to a single wider strip ${\omega}_\delta:=(-\delta,1+\delta)\times\Re$ (where $\delta>0$), in such a way as to remain uniformly bounded in the $C^{\,3,\tilde{\varrho}}$-norms on $\omega_\delta$.
\vspace{.1in}

\noindent
(b) Let ${\bm{\mathfrak{M}}}^*(\varrho)$ denote the sub-family of pairs ${\bf\Gamma}\in{\bm{\mathfrak{M}}}(\varrho)$ such that $|p_i(t)-p_i(\tau)|\geq 
\mathfrak{E}_0\,\eta(|t-\tau|)$ for $i=1,2$ and all $t,\tau\in\Re$, where $\mathfrak{E}_0$ is any fixed positive constant, $\eta(t):=\min\{1,t\}$, and $p_i(t):\Re\rightarrow\Gamma_i$, $i=1,2$, denote arc-length parametrizations of the corresponding arcs $\Gamma_i$. Then there exists a value $\delta>0$ so small that the continued functions $w=F(z)$ map the strip $\omega_\delta$ invertably onto the extended regions $\Omega_\delta({\bf\Gamma})$ in such a way that the inverse mappings $z=F^{-1}(w):\Omega_\delta({\bf\Gamma})\rightarrow\omega_\delta$, ${\bf \Gamma}\in{\boldsymbol{\mathfrak {M}}}(\varrho)$ are also uniformly bounded in the $C^{\,3,\tilde{\varrho}}$-norm. The continued capacitary potentials $U(p)=U({\bf\Gamma};p):\Omega_\delta({\bf\Gamma})\rightarrow(-\delta,1+\delta)$, ${\bf\Gamma}\in{\bm{\mathfrak{M}}}^*(\varrho)$, defined such that $U(p):={\rm Re}(F^{-1}(z))$, are also uniformly bounded in their respective $C^{\,3,\tilde{\varrho}}(\Omega_\delta({\bf\Gamma}))$-norms.
\end{theorem}
\begin{lemma}
\label{lem 2.7.3}
Given any functions $f_i\in C^{\,1,\varrho}(\Re)$, $i=1,2$, such that ${\boldsymbol |}f_i{\boldsymbol |}, {\boldsymbol |}f_i^{\,\prime}{\boldsymbol |}\leq M$ and $|f_i^{\,\prime}(x_1)-f_i^{\,\prime}(x_2)|\leq L\,|x_1-x_2|^{\varrho}$ for all $x_1,x_2\in\Re$ (for some constants $M,L\in\Re_+$ and $\varrho\in(0,1]$), let $U(x,y):{\rm Cl}(\omega)\rightarrow\Re$ denote the bounded solution of the Dirichlet problem: $\Delta U=0$ in $\omega:=\Re\times(0,1)$, $U=f_i$ on $\gamma_i$, $i=1,2$ (where $\gamma_1=\Re\times\{0\}$, $\gamma_2:=\Re\times\{1\}$). 
Then there exist new constants $M^*,L^*\in\Re_+$ and $\varrho^*\in(0,\varrho)$ such that 
$$|\nabla U(p)|\leq M^*;\,\,|\nabla U(p)-\nabla U(q)|\leq L^*|p-q|^{\varrho^*},$$
both for all points $p,q\in{\rm Cl}(\omega)$. (In fact one can first choose $\varrho^*$ arbitrarily close to $\varrho$ in $(0,\varrho)$, and then choose $0<M^*=C_0\big(M+(L/\varrho)\big)$ and $0<L^*= C_1\big((M+(L/(\varrho-\varrho^*))\big)$, where $C_0,C_1$ denote universal constants.
\end{lemma}

\noindent
{\bf Proof outline.} 
Under the assumptions, we have $|U_x(x,y)|\leq M$ and $|U_x(x+h,y)-U_x(x,y)|\leq L|h|^\varrho$, both for all $(x,y),(x+h,y)\in\omega$, since the left sides are continuous in ${\rm Cl}(\omega)$ and sub-harmonic in $\omega$, and since the same inequalities hold on $\partial\omega$. Therefore $|U_x(x+h,y)-U_x(x,y)|\leq (L+2M)|h|^{\varrho^*}$ for $(x,y),(x+h,y)\in\Re\times[0,1]$ and $\varrho^*\in(0,\varrho)$, as follows directly if $|h|\leq 1$, and from: $|U_x|\leq M$ otherwise. For the proof of the remaining estimates, we observe that $U(p)=U_1(p)+U_2(p)$, where $U_1$ (resp. $U_2$) solves the Dirichlet problem in the case where $f_2(x)=0$ ($f_1(x)=0$). It suffices to consider only the first case, for which we set $f=f_1$ and $U:=U_1$. For this case, the solution of the Dirichlet problem is given (for $0<y\leq 1$) by the convolution-integral formula: 
\begin{equation}
\label{eqn 2.7.10}
U(x,y)=\int_{-\infty}^\infty K(t,y)\,f(x-t)\,dt,
\end{equation}
where the kernal is the harmonic function: $K(x,y)=\frac{1}{4}\big({\rm sin}(2\alpha y)\big)\big/\big({\rm sin}^2(\alpha y)+{\rm sinh}^2(\alpha x)\big)=(1/\pi)(\partial/\partial y)({\rm Re}({\rm \ln}({\rm sinh}(\alpha z)))$, in which $\alpha=\pi/2$. One can show by applying parameter differentiation and integration by parts to (\ref{eqn 2.7.10}), that: 
\begin{equation}
\label{eqn 2.7.11}
U_x(x,y)-f^{\,\prime}(x)=\int_{-\infty}^\infty K(t,y)\,(f^{\,\prime}(x-t)-f^{\,\prime}(x))\,dt,
\end{equation}
\begin{equation}
\label{eqn 2.7.11a}
U_y(x,y)=-2f(x)+\int_{-\infty}^\infty \big(\tilde{K}(t,y)+{\rm sign}(t)\big)(f^{\,\prime}(x-t)-f^{\,\prime}(x))\,dt,
\end{equation}

\begin{equation}
\label{eqn 2.7.11aa}
U_y(x,y)-U_y(x,0)=\int_{-\infty}^\infty \big(\tilde{K}(t,y)-\tilde{K}(t,0)\big)(f^{\,\prime}(x-t)-f^{\,\prime}(x))\,dt,
\end{equation}
all for $0<y\leq 1$, where (\ref{eqn 2.7.11}), (\ref{eqn 2.7.11aa}) follow from (\ref{eqn 2.7.10}), (\ref{eqn 2.7.11a}), respectively, and where (\ref{eqn 2.7.11a}), with $y=0$, defines the function $U_y(x,0):\Re\rightarrow\Re$. Here we define $\tilde{K}(x,y):=(\partial/\partial y)\int_{0}^x K(s,y)ds$  $=-\big(\big({\rm sinh}(\alpha x){\rm cosh}(\alpha x)\big/\big({\rm sinh}^2(\alpha x)+{\rm sin}^2(\alpha y)\big)$, $\tilde{K}(x,0)=-{\rm coth}(\alpha x)$, and $\tilde{K}(x,y)-\tilde{K}(x,0)={\rm coth}(\alpha x)\big({\rm sin}^2(\alpha y)\big/({\rm sinh}^2(\alpha x)+{\rm sin}^2(\alpha y))\big)$. 
Under the assumptions, one can show by estimating the singular integrals (\ref{eqn 2.7.11}), (\ref{eqn 2.7.11a}), and (\ref{eqn 2.7.11aa}) that $|U_y(x,y)|\leq C_0\,\big(\boldsymbol{|}f\boldsymbol{|}+\boldsymbol{|}f^{\,\prime}\boldsymbol{|}+(L/\rho)\big)$ for $y=0$, and that $|U_x(x,y+h)-U_x(x,y)|, |U_y(x,y+h)-U_y(x,y)|\leq C_1\,\big((L/\varrho)+\boldsymbol{|}f^{\,\prime}\boldsymbol{|}\big)|h|^\varrho$, both for $y=0$ and $|h|\leq 1$, where $C_0,C_1$ are constants independent of $L,M,\varrho$. The first inequality holds for all $(x,y)\in\Re\times[0,1]$, as one sees by odd continuation of $U$ across $\gamma_2$ (where $U=0$) and the maximum principle (since the left side is bounded and sub-harmonic). For any fixed $0<h<1$, the second and third inequalities also hold for all $(x,y)\in\Re\times[0,1-h]$, by the same odd continuation argument. Similarly, a double application of (\ref{eqn 2.7.11a}) yields a formula for $U_y(x+h,0)-U_y(x,0)$, which one can estimate to show (after increasing $C_1$ if necessary) that, for any given $h\in [-1,1]\setminus\{0\}$, we have: $|U_y(x+h,y)-U_y(x,y)|\leq C_1\big(\boldsymbol{|}f^{\,\prime}\boldsymbol{|}+(L/\varrho)+L\,{\rm ln}(1/|h|)\big)|h|^\varrho$ for all $x\in\Re$, provided that $y=0$. This inequality extends to $(x,y)\in\Re\times[0,1]$ by odd continuation (of $U$ across $\gamma_1$), since the left side is bounded and subharmonic. Finally, for $h\in[-1,1]\setminus\{0\}$ and $0<\varrho^*<\varrho^*+\delta=\varrho$, we have $|U_y(x+h,y)-U_y(x,y)|\leq C_1\,\big(\boldsymbol{|}f^{\,\prime}\boldsymbol{|}+(L/\varrho)+(L/e(\varrho-\varrho^*))\big)|h|^{\varrho^*}$ for all $(x,y)\in\Re\times[0,1]$, since $|h|^\delta{\rm ln}(1/|h|)\leq (1/e\delta)$. Finally, the same inequality holds for $|h|\geq 1$ after suitably increasing $C_1$ while keeping $C_0$ fixed.
\vspace{.1in}

\noindent
{\bf Proof of Thm. \ref{thm 2.7.2}, Part(a).} 
As in the proof of Lem. \ref{lem 2.6.1}, for any $\varrho\in(0,1]$ and any classical solution ${\bf\Gamma}\in{\bm{\mathfrak{M}}}(\varrho)$ of Prob. \ref{prob2.1.1}, let $w=F(z)$ denote the corresponding periodic conformal mapping of the strip $\overline{\omega}:={\rm Cl}(\omega):=\Re\times[0,1]$ onto ${\rm Cl}(\Omega({\bf\Gamma}))$. (The term "admissible" expresses a correspondence of the mapping $w=F(z)$ to some arc-pair ${\bf\Gamma}\in\bm{\mathfrak{M}}(\varrho)$.) 
For any admissible function $w=F(z)$ and any $z\in\overline{\omega}$, we have that 
\begin{equation}
\label{eqn 2.7.1-8p}
\big(1\big/\,\overline{A}\,)\leq |F'(z)
|=\big(1\big/|\nabla U(F(z)|\big)\leq\big(1\big/\underline{A}\,).
\end{equation}
We have that (${\bf i}$): $\phi(x,y)=A\big(F(x,y)\big)$ in $\Re\times\{0,1\}$, where we define $\phi(x,y):={\rm ln}\big(|F'(x,y)|\big)$ in $\Re\times[0,1]$ and $A(p)=\big(1\big/a(p)\big)=\big(1\big/a_i(p)\big)$ on $\Gamma_i$ for $i=1,2$. By differentiating (${\bf i}$), we determine that (${\bf ii}$): $\phi_x(x,y)=A'\big(F(x,y)\big)F_x(x,y)$ in $\Re\times\{0,1\}$, from which it follows that if $A'$ is uniformly bounded and the function $F_x(x,y)$ is in $C^{\,\varrho}\big(\Re\times\{0,1\}\big)$, then the function $\phi_x(x,y)$ is in $C^{\,\varrho}\big(\Re\times\{0,1\}\big)$. It follows by Lem. \ref{lem 2.7.3} that the function $\nabla \phi(x,y)$ is in $C^{\,\hat{\varrho}}\big(\Re\times[0,1]\big)$ for a smaller value $\hat{\varrho} \in(0,1]$, and therefore, by the Cauchy-Riemann equations, that the function $\nabla\big({\rm log}\big(F'(x,y)\big)\big)$ is also in $C^{\,\hat{\varrho}}\big(\Re\times[0,1]\big)$. Finally, one shows using (\ref{eqn 2.7.1-8p}), that $\nabla(F'(x,y))$ is in $C^{\,\hat{\varrho}}\big(\Re\times[0,1]\big)$, and therefore that all second order partial derivatives of $F(x,y)$ are in $C^{\,\hat{\varrho}}\big(\Re\times[0,1]\big)$, completing the first cycle of the argument.
\vspace{.1in}

\noindent
At this point, it follows from (${\bf ii}$) by differentiation that (${\bf iii}$): $\phi_{xx}(x,y)=D_x \big(A'(F)F_x\big)=A'(F)F_{xx}(x,y)+A''(F)F_x^2$ in $\Re\times\{0,1\}$. Since $A'$ and $A''$ are uniformly bounded and the functions $F_x(x,y)$ and $F_{xx}(x,y)$ are both in $C^{\,\hat{\varrho}}\big(\Re\times\{0,1\}\big)$, it follows from (${\bf iii}$) that the function $\phi_{xx}(x,y)$ is in $C^{\,\hat{\varrho}}\big(\Re\times\{0,1\}\big)$, from which it further follows by Lem. \ref{lem 2.7.3} that $(\partial/\partial x)\nabla\phi(x,y)$ is in $C^{\,\tilde{\varrho}}\big(\Re\times[0,1]\big)$ for some value $\tilde{\varrho}\in(0,\hat{\varrho}]$, and then by the Cauchy-Riemann equations that first the function $(\partial/\partial x)\nabla\big({\rm log}\big(F'(x,y)\big)$ is in $C^{\,\tilde{\varrho}}\big(\Re\times[0,1]\big)$. Finally, it follows by applying (\ref{eqn 2.7.1-8p}) that the function $(\partial/\partial x)\nabla\big(F'(x,y)\big)$ is in $C^{\,\tilde{\varrho}}\big(\Re\times[0,1]\big)$. Therefore, all third-order partial derivatives of the function $F(x,y)$ are in $C^{\,\tilde{\varrho}}\big(\Re\times[0,1]\big)$, completing the second cycle.
\vspace{.1in}

\noindent
A closer look at the same steps (especially the applications of Lem. \ref{lem 2.7.3}) shows that the family of all admissible functions $w=F(z)$ is uniformly bounded in the $C^{\,3,\tilde{\varrho}}$-norm in $\Re\times[0,1]$. 
The remaining claim about uniform $C^{\,3,\tilde{\varrho}}$-continuations of admissible functions $F$ to $\omega_\delta$ now follows by a well-known result (see [GT], Lemma 6.37), which states that if $U\in C^{\,3,\tilde{\varrho}}\big(\overline{\omega}\big)$, then $U$ has a continuation to $C^{\,3,\tilde{\varrho}}\big(\omega_\delta\big)$ such that the $C^{\,3,\tilde{\varrho}}\big(\omega_\delta\big)$-norm of $U$ cannot exceed the $C^{\,3,\tilde{\varrho}}\big(\overline{\omega}\big)$-norm of $U$ multiplied by some value depending only on $\delta$.
\vspace{.1in}

\noindent
{\bf Proof of Thm. \ref{thm 2.7.2}, Part (b).} 
Given an admissible transformation $w=u+jv=F(z)=\phi(z)+j\psi(z):\overline{\omega}\rightarrow{\rm Cl}\big(\Omega({\bf\Gamma})\big)$, (corresponding to ${\bf\Gamma}\in{\bm{\mathfrak{M}}}^*(\varrho)$), we also use $w=F(z)$ to denote the continuation of $F$ to $\omega_\delta$. In view of the uniform bounds for the $C^{\,3,\tilde{\varrho}}$-norms of the continuations to $\omega_{\delta_0}$ of the admissible function $F$, it follows that there exists a constant $\delta_1\in(0,\delta_0)$ such that
\begin{equation} 
\label{eqn 2.7.1-8pq}\big(F(z)-F(z_0)\big)=(z-z_0)\big(F'(z_0)+E(z_0,z)\big)
\end{equation} for any $z_0\in\Re\times\{0, 1\}$, where $E(z_0,z)$ is a remainder term such that $E(z_0,z)\rightarrow 0$ as $z\rightarrow z_0$. By double-application of the remainder formula (\ref{eqn 2.7.1-8pq}) in the cases where $z=z_1$ and $z=z_2$, one sees that there exist positive constants $\delta_2\in(0,\delta
_1)$ and $0<C_1<C_2<\infty$ such that
\begin{equation}
\label{eqn 2.7.1-8q}
C_1\sqrt{(t-\tau)^2+(\alpha-\beta)^2}\leq\big|F(t,\alpha)-F(\tau,\beta)\big|
\end{equation}
$$\leq C_2\sqrt{(t-\tau)^2+(\alpha-\beta)^2}$$
for all $t,\tau, \alpha,\beta\in\Re$ such that $|t-\tau|\leq\delta_2$ and $|\alpha|,|\beta|\leq\delta_2$. By assumption, there is a value $\mathfrak{E}_0>0$ such that $|F(t,0)-F(\tau,0)|\geq\mathfrak{E}_0\,\eta\,(|t-\tau|)$ for $t,\tau\in\Re$, where we define $\eta(t):={\rm min}\{t,1\}$. Therefore
\begin{equation}
\label{eqn 2.7.1-8r}  
|F(t,\alpha)-F(\tau,\beta)|\geq|F(t,0)-F(\tau,0)|
\end{equation}
$$-|F(t,\alpha)-F(t,0)|-|F(\tau,\beta)-F(\tau,0)|$$
$$\geq\mathfrak{E}_0\,\eta\,(|t-\tau|)-2C_2\overline{\delta}\geq(\mathfrak{E}_0/2)\sqrt{(\eta(|t-\tau|))^2+(\alpha-\beta)^2}$$
for all $t,\tau,\alpha,\beta\in\Re$ such that $|t-\tau|\geq\delta_2$ and $|\alpha|,|\beta|\leq\overline{\delta}:=\min\big\{\delta_2, \big(\mathfrak{E}_0\big/(4C_2+\mathfrak{E}_0)\big)\big\}$. 
By (\ref{eqn 2.7.1-8q}) and (\ref{eqn 2.7.1-8r}), there exists a value $M:=(1/{\rm min}\{C_1,\mathfrak{E}_0\})$ such that for all admissible $F$, the restrictions of the mappings $F(z)$ to $\Re\times(-\overline{\delta},\overline{\delta})$ are all one-to-one mappings of $\Re\times(-\overline{\delta}, \overline{\delta})$ onto the $M\overline{\delta}$-neighborhoods of the corresponding arcs $\Gamma_i$. It follows from this that if $\Omega_\delta({\bf\Gamma})$ denotes the union of $\Omega({\bf\Gamma})$ with the $M\overline{\delta}$-neighborhoods of the arcs $\Gamma_i$, $i=1,2$, then the functions $w=F(z):\omega_\delta\rightarrow\Omega_\delta({\bf\Gamma})$ are one-to-one, and therefore globally invertible, 
since their restrictions $w=F(z):\overline{\omega}\rightarrow{\rm Cl}(\Omega({\bf\Gamma}))$ were already known to be globally invertible conformal mappings. It follows from all this that the inverse mappings $z=G(w):=F^{-1}(w):\Omega_\delta({\bf\Gamma})\rightarrow\omega_\delta$ are all Lipschitz-continuous functions with Lipschitz constant $(1/M)$, whose domains contain $M\overline{\delta}$-neighborhoods of the arcs $\Gamma_i$, $i=1,2$, and which are uniformly bounded in the $C^{\,3,\tilde{\varrho}}$-norm relative to the $(\overline{\delta}/M)$-neighborhood of $\Omega({\bf\Gamma})$. 
In view of this, the assertions of Thm. \ref{thm 2.7.2}(b) hold, and the assertions regarding the uniform continuations of the capacitary potentials $U({\bf\Gamma};p)$ follows from this by defining $U(z)={\rm Re}(F^{-1})(z)$ for $F$ corresponding to ${\bf\Gamma}$. 
\section{Qualitative properties of solutions}
\label{section 4}
\subsection{Main qualitative results}
\label{subsection 4.1}
Chapter \ref{section 4}, for which this section serves as the summary of main results, focuses primarily on results in the context of the double-free-boundary problem relative to a "valley" $G$ of a single logarithmically-subharmonic flow-speed function $a(p):G\rightarrow\Re_+$ (see Prob. \ref{prob 4I.+}). The main topics in this context include existence and non-existence of solutions (see Thms. \ref{thm 4I.*} and \ref{thm 4I.0}), as well as two successfull approaches to the resolution of the uniqueness question under suitable assumptions The second (less direct) approach proves the existence of continuously and monotonically (in a positive sense) varying one-parameter local solution families, as defined in Def. \ref{def 4.1.1}. In fact these solution families can be made the basis for a local uniqueness proof (see Thm. \ref{thm 4.1.2}).  Our local uniqueness assertions follow from Thm. \ref{thm 4.1.2} via Thm. \ref{thm 4.1.1} and Cor. \ref{cor 4.1.2}). 

\begin{problem} {(Double-free-boundary flow problem with two flow-speed functions in a periodic strip-like region)}. 
\label{prob 4I} 
In the context of Prob. \ref{prob2.1.1}, let be given a simply-connected, $P$-periodic (in $x$), strip-like domain $G$ having interior and exterior tangent balls of uniform radius at all boundary points. We use $\partial^+G$ and $\partial^-G$) to denote the upper and lower boundary components of $G$, respectively, both of which are in  ${\rm X}$. We also use ${\bf X}(G)$ (resp. ${\bf X}(\overline{G}))$ to denote the set of all pairs ${\bf\Gamma}=(\Gamma_1,\Gamma_2)\in{\bf X}$ such that $\partial^-G<\,(\leq)\,\Gamma_1<\Gamma_2<\,(\leq)\,\partial^+G$. Given the $P$-periodic (in $x$) flow-speed functions $a_1(p), a_2(p):\Re^2\rightarrow\Re_+$, both in ${\boldsymbol{\cal A}}$, we seek a pair ${\bf\Gamma}=(\Gamma_1,\Gamma_2)\in{\bf X}(\overline{G})$ such that $|{\nabla}U({\bf\Gamma};p)|=a_i(p)$ for all $p\in\Gamma_i$, $i=1,2$.
\end{problem}
\begin{problem}
\label{prob 4I.} 
This refers to Prob. \ref{prob 4I} in the case where the functions $a_1(p)$ and $a_2(p)$ both coincide with a single $P$-periodic (in $x$) function $a(p):\Re^2\rightarrow\Re_+$ in the class  $\boldsymbol{{\cal{A}}}\cap C^{\,3,\varrho}$ for some $\varrho\in(0,1]$.
\end{problem}

\begin{problem}
\label{prob 4I.+}
This is Prob. \ref{prob 4I.} in the case where the boundary flow-speed function $a(p):G\rightarrow\Re_+$ is (weakly) logarithmically subharmonic.
\end{problem}
\begin{theorem} 
{(Non-existence of solutions of Prob. \ref{prob 4I})}
\label{thm 4I.0} 
(a) In Prob. \ref{prob 4I}, let $A_i(G):=\inf\{a_i(p):p\in G\}$ and $B_i(G):=\sup\{a_i(p):p\in G\}$, $i=1,2$. If $A_1(G)>B_2(G)$ or $A_2(G)>B_1(G)$, then there does not exist any classical solution ${\bf\Gamma}\in{\bf X}(G)$. In general, Prob. \ref{prob 4I} has no classical solution ${\bf\Gamma}\in{\bf X}(G)$ at the vector ${\bm\lambda}=(\lambda_1,\lambda_2)\in\Re^2_+$ if $A_1(G)\lambda_1>B_2(G)\lambda_2$ or $A_2(G)\lambda_2>B_1(G)\lambda_1$. (b) In Prob. \ref{prob 4I.+}, no classical solution exists at ${\bm\lambda}\in\Re_+^2$ if $\big(\lambda_{3-i}\big/\lambda_i\big)<\big(\,\underline{A}\,\big/\,\overline{A}\,\big)$ for either $i=1$ or $i=2$.
\end{theorem}
\vspace{.1in}

\noindent
{\bf Proof.} Let ${\bf\Gamma}=(\Gamma_1,\Gamma_2)\in{\bf X}(G)$ denote a classical solution at ${\bm\lambda}$ such that $\lambda_1\, A_1(G)>\lambda_2\,B_2(G)$. Let $\phi(p):={\rm ln}\big({\nabla}U(p)|\big)$ in ${\rm Cl}(\Omega)$, where $U(p):=U({\bf\Gamma};p)$. Then 
$\phi(p)$ is a harmonic function in $\Omega:=\Omega({\bf\Gamma})$ such that $\inf\{\phi(p):p\in\Gamma_2\}=\inf\big\{{\rm ln}\big(\lambda_2\,a_2(p)\big):p\in\Gamma_2\big\}>\sup\big\{{\rm ln}\big(\lambda_1\,a_1(p)\big):p\in\Gamma_1\big\}=\sup\{\phi(p):p\in\Gamma_1\}.$ Therefore $\int_{\gamma_i}K_i(p)\,ds=\int_{\gamma_i}\phi_{\boldsymbol{\nu}}(p)\,ds>0$ for $i=1,2$, where $K_i(p)$ denotes the signed curvature of $\Gamma_i$ at $p\in\Gamma_i$, and where $\gamma_i$ denotes one $P$-period of $\Gamma_i$. But this contradicts the $P$-periodicity of the curves $\Gamma_i$. The contradiction in the alternate case is similar.
\begin{theorem} 
{(Existence of classical solutions of Prob. \ref{prob 4I} between weakly-lower and weakly-upper classical solutions)} 
\label{thm 4I.*} 
In the context of Prob. \ref{prob 4I}, let be given a (weakly) lower (classical) solution $\tilde{\bf\Gamma}=(\tilde{\Gamma}_1,\tilde{\Gamma}_2)\in{\bf X}(\overline{G})$ and a (weakly) upper (classical) solution  $\hat{\bf\Gamma}=(\hat{\Gamma}_1,\hat{\Gamma}_2)\in{\bf X}(\overline{G})$ such that $\tilde{\bf\Gamma}<\hat{\bf\Gamma}$. Then there exists a classical solution ${\bf\Gamma}\in{\bf{X}}(\overline{G})$ such that $\tilde{\bf\Gamma}\leq{\bf\Gamma}\leq\hat{\bf\Gamma}$. 
\end{theorem}
\vspace{.1in}

\noindent
{\bf Proof.} We define the function-sequences $\big(a_{n,1}\big)_{n=1}^\infty$ and $\big(a_{n,2}\big)_{n=1}^\infty$, in which all the functions are in the class ${\boldsymbol{\cal A}}$; also such that $a_{n,i}(p)\rightarrow a_i(p)$ uniformly in any compact subset of $\Re^2$ as $n\rightarrow\infty$, $i=1,2$, and, finally, such that for each $n\in N$, we have $(-1)^i\big(a_{n,i}(p)-a_i(p)\big)>0$ in $\tilde{\Gamma}_i$ and $(-1)^i\big(a_{n,i}(p)-a_i(p)\big)<0$ in $\hat{\Gamma}_i$, $i=1,2$, 
so that $\tilde{\bf\Gamma}$ (resp. $\hat{\bf\Gamma}$) is a strict lower (upper) solution of Prob. \ref{prob 4I} relative to the functions $a_{n,1}$ and $a_{n,2}$. By Thm. \ref{thm 2.7.1}, for each $n\in N$, there exists a classical solution ${\bf\Gamma}_n=(\Gamma_{n,1},\Gamma_{n,2})\in{\bf X}(\overline{G})$ of Prob. \ref{prob 4I} corresponding to the new functions $a_{n,1}$ and $a_{n,2}$, such that $\tilde{\bf\Gamma}_n\leq{\bf\Gamma}_n\leq\hat{\bf\Gamma}_n$. Moreover, by Thm. \ref{thm 2.7.1}, the curvatures of the arcs $\Gamma_{n,1}$ and $\Gamma_{n,2}$ are uniformly bounded over $n\in N$, and we have $\underline{A}\leq|{\nabla}U_n(p)|\leq\overline{A}$ and $|{\nabla}\phi_n(p)|\leq L_0$, both uniformly in $\Omega_n:=\Omega({\bf\Gamma}_n)$ and independent of $n\in N$, where $U_n(p):=U({\bf\Gamma}_n;p)$ and $\phi_n(p):={\rm ln}\big(|{\nabla}U_n(p)|\big)$. In view of this, one can conclude by passing to a subsequence (expressed in the original notation) that ${\bf\Gamma}_n\rightarrow{\bf\Gamma}$ for some solution ${\bf\Gamma}\in{\bf X}(\overline{G})$. Clearly $\tilde{\bf\Gamma}\leq{\bf\Gamma}\leq\hat{\bf\Gamma}$, and ${\bf\Gamma}$ has the properties listed in Thm. \ref{thm 2.7.1}.

\begin{lemma} 
{(Main curvature estimate for solutions of Prob. \ref{prob 4I.+}}
\label{lem Z.1}
For any pair ${\bm\lambda}=(\lambda_1,\lambda_2)\in\Re_+^2$ and any classical-solution-pair ${\bf\Gamma}=(\Gamma_1,\Gamma_2)\in{\bf X}(G)\cap C^{2}$ of Prob. \ref{prob 4I.+} at ${\bm\lambda}$, we have
\begin{equation}
\label{eqn Z.1}
\big|K(p)-(\partial/\partial\boldsymbol{\nu})\,{\rm ln}\big(a(p)\big)\big|\leq\lambda_i\,a(p)|\mu|\,+\,\big(2H\overline{A}\,\overline{\lambda}\big/\underline{A}^2\underline{\lambda}^2\,\big),
\end{equation}
pointwise on $\Gamma_{i}$, $i=1,2$, where for $p\in\Gamma_i$, we set $\boldsymbol{\nu}(p):=\big({\nabla}U(p)\big/|{\nabla}U(p)|\big)$, and we use $K(p)$ to denote the counter-clockwise-oriented curvature of $\Gamma_{i}$ at the point $p\in\Gamma_{i}$. Also, we set $\mu:={\rm ln}\big(\lambda_2/\lambda_1\big)$ and $\overline{\mu}:={\rm ln}\big(\,\overline{\lambda}\big/\underline{\lambda}\,\big)$, 
and we assume for some constant $H>0$ that $0\leq\Delta{\rm ln}\big(a(p)\big)\leq H$ throughout $G$. 
Actually, the estimate \ref{eqn Z.1} holds for all points $p\in{\rm Cl}\big(\Omega({\bm \Gamma})\big)$, where $K(p)$ (resp. $\boldsymbol{\nu}$) denotes the left curvature of (resp. the left normal to) the level curve of $U$ through the point $p\in\,{\rm Cl}(\Omega)$.
\end{lemma}

\noindent
{\bf Proof.} We define the function $\phi(p):={\rm ln}\big(|{\nabla}U(p)|\big/\lambda_1 a(p)\big)-\mu\, U(p)$ in the closure of $\Omega:=\Omega({\bf\Gamma})$. Then $0\leq-\Delta \phi(p)=\Delta\,{\rm ln}\big(a(p)\big)\leq H$ in $\Omega$, and it follows from the fact that $|{\nabla}U(p)|=\lambda_ia(p)$ on $\Gamma_i$, $i=1,2$, that $\phi(p)=0$ on 
$\partial\Omega$. In view of this, it follows from Lem. \ref{lem 2.7.1}(a) that 
\begin{equation}
\label{eqn Z.1aa}
\big|{\nabla}\phi(p)\big|\leq\big(2H\overline{A}\,\overline{\lambda}\,\big/\underline{A}^2\underline{\lambda}^2\big)=\big(2H\overline{A}\,{\rm exp}(\overline{\mu})\big/\underline{A}^2\underline{\lambda}\big).
\end{equation} 
Now the inequality (\ref{eqn Z.1}) follows from (\ref{eqn Z.1aa}) in view of the facts that 
\begin{equation}
\label{eqn Z.1aaa}
\partial\phi\big/\partial\boldsymbol{\nu}=K(p)-(\partial\big/\partial\boldsymbol{\nu})\,{\rm ln}\big(a(p)\big)-\mu\,\partial U\big/\partial\boldsymbol{\nu}
\end{equation}
in ${\rm Cl}(\Omega)$, and where $\partial U(p)/\partial\boldsymbol{\nu}=\lambda_ia(p)$ on $\Gamma_i$, $i=1,2$. Similarly, the estimate (\ref{eqn Z.1}) follows from (\ref{eqn Z.1aa}) and (\ref{eqn Z.1aaa}), in view of the fact that $|{\nabla} U(p)|\leq\overline{\lambda}\,\,\overline{A}$ for all $p\in{\rm Cl}(\Omega)$.

\begin{theorem} {(Non-existence of classical solutions of Prob. \ref{prob 4I.+} at {\bm\lambda})}
\label{thm 4I.2}
There does not exist any classical solution of Prob. \ref{prob 4I.+} at a pair ${\bm\lambda}\in\Re_+^2$ such that either $\mu\lambda_1>E$ or $-\mu\lambda_2>E$, where we define $E:=\sup\big\{\big(|{\nabla}a(p)|\big/a^2(p)\big):p\in G\big\}>0$ and $\mu={\rm ln}(\lambda_2/\lambda_1)$. In other words, if Prob. \ref{prob 4I.+} does have a classical solution at ${\bm\lambda}$, then ${\rm min}\{(\lambda_1,\lambda_2\}\,|\mu|\,\leq\, E$. 
\end{theorem}
\noindent
{\bf Proof.} Let ${\bf\Gamma}=(\Gamma_1,\Gamma_2)\in{\bf X}(G)$ denote a $C^{2}$-solution of Prob. \ref{prob 4I.+} at ${\bm\lambda}\in\Re_+^2$. In the context of the proof of Lem. \ref{lem Z.1} (where $\Delta\,{\rm ln}\big(a(p)\big)\geq 0$ in $\Omega$ by assumption), we define $\phi(p):={\rm ln}\big(|{\nabla}U(p)|\big/\lambda_1a(p)\big)-\mu U(p)$. Then $\phi(p)=0$ on $\partial\Omega$ and therefore $\phi(p)\geq 0$ in $\Omega$ (since $\Delta\phi\leq 0$ there), from which it follows by (\ref{eqn Z.1aa}) that
$(-1)^i(\partial\phi\big/\partial\boldsymbol{\nu})=(-1)^i\big(K_i(p)-(\partial/\partial\boldsymbol{\nu}){\rm ln}\big(a(p)\big)-\mu\lambda_i a(p)\big)\leq 0
$ on $\Gamma_i$ for $i=1,2$, where ${\boldsymbol{\nu}}$ is the left normal to $\Gamma_i$. Therefore $K_1(p)\geq a(p)\big(\mu\lambda_1-\big(|{\nabla} a(p)|/a^2(p)\big)\big)$ on $\Gamma_1$ and $K_2(p)\leq a(p)\big(\mu\lambda_2+\big(|{\nabla} a(p)|/a^2(p)\big)\big)$ on $\Gamma_2$. Since for any ($P$-periodic) classical solution ${\bf\Gamma}\in{\bf X}(G)$, the arc-length integral of the signed curvature $K_i(p)$ over any one $P$-period (in $x$) of either of the component arcs $\Gamma_i$, $i=1,2,$ must vanish, it is impossible for a classical solution to exist if $\mu\,\lambda_1\,a^2(p)>|{\nabla}a(p)|$ for all $p\in\Gamma_1$ or if $-\mu\,\lambda_2\,a^2(p)>|{\nabla} a(p)|$ for all $p\in\Gamma_2$. The assertion follows.

\begin{definition} {(Continuously-varying and positively-ordered families of solutions of Prob. \ref{prob 4I})}
\label{def 4.1.1}
Given the positive continuous functions $a_1(p),a_2(p):\Re^2\rightarrow\Re_+$ and the continuous vector-valued function $\boldsymbol{\lambda}(t)=(\lambda_1(t),\lambda_2(t)):\Re\rightarrow\Re_+^2$ such that 
\begin{equation}
\label{eqn 4.1.2}
\lambda_1(t)\,({\rm resp.}\,\lambda_2(t))\,\,{\rm is}\,\,{\rm strictly}\,\,{\rm increasing}\,({\rm decreasing})\,\,{\rm in}\,\,\Re,
\end{equation}
we use the phrase "parametrized solution family" to refer to an open interval $I$ and a mapping ${\bf\Gamma}(t)=(\Gamma_{1}(t),\Gamma_{2}(t)):I\rightarrow{\bf X}(G)$ such that for each $t\in I$, the arc-pair ${\bf\Gamma}(t)\in{\bf X}(G)\cap C^{\,3,\tilde{\varrho}}$ solves Prob. \ref{prob 4I} at $\boldsymbol{\lambda}(t)\in\Re_+^2$, which means that the capacitary potential $U(t;p):=U({\bf\Gamma}(t);p)$ defined in the closure of $\Omega(t):=\Omega({\bf\Gamma}(t))$ is a $C^1$-function which satisfies the condition
\begin{equation}
\label{eqn 4.1.1}
\big|{\nabla}_p\,U(t;p)\big|=\lambda_i(t)\,a_i(p)\,\,{\rm on}\,\,\Gamma_{i}(t)\,
\end{equation}
for $i=1,2$ and $t\in I$. The parametrized solution-family $\{{\bf \Gamma}(t)=(\Gamma_{1}(t),\Gamma_{2}(t)):t\in I\}$ of Prob. \ref{prob 4I} is called "positively-ordered" (or "elliptically-ordered" (see Beurling [AB1]) if 
\begin{equation}
\label{eqn 4.1.3}
{\bf\Gamma}(\alpha)<{\bf\Gamma}(\beta)\,\,{\rm whenever}\,\,\alpha<\beta\,\,{\rm in}\,\,I,
\end{equation}
and it is called "continuously-varying" if
\begin{equation}
\label{eqn 4.1.4}
{\bf\Gamma}(t)\rightarrow{\bf\Gamma}(\tau)\,\,{\rm as}\,\, t\rightarrow \tau\,\,{\rm for}\,\,{\rm  any}\,\,\tau\in\,\,I,
\end{equation}
in the sense that that $H\big({\bf\Gamma}(t),{\bf\Gamma}(\tau)\big)\rightarrow 0\,\,{\rm as}\,\,t\rightarrow\tau$, where $H\big({\bf\Gamma}(t),{\bf\Gamma}(\tau)\big)$ denotes the Hadamard distance between the arcs.
\end{definition}
\begin{theorem}
\label{thm 4.1.1}
{(Uniqueness of solutions of Probs. \ref{prob 1.1}, \ref{prob2.1.1}, or \ref{prob 4I}, viewed as members of positively-ordered, continuously-varying parametrized solution-families)} In the context of Def. \ref{def 4.1.1}, given a real open interval ${I}$ and the continuous function ${\bm\lambda}(t)=(\lambda_1(t),\lambda_2(t)):{\rm Cl}(I)\rightarrow\Re_+^2$ satisfying (\ref{eqn 4.1.2}) relative to ${I}$, let be given a continuously-varying, positively-ordered, parametrized solution-family ${\bf \Gamma}(t)=(\Gamma_{1}(t),\,$ $\Gamma_{2}(t)):t\in I\}$ for Prob. \ref{prob 1.1}, \ref{prob2.1.1}, or \ref{prob 4I}.  
For a given, fixed value $\tau\in I$, let $\tilde{\bf\Gamma}$ denote ${\bf any}$ classical solution of (\ref{eqn 4.1.1}) at $t:=\tau$ (therefore not necessarily a member of the above continuous, monotone solution-family). Then in fact $\tilde{{\bf\Gamma}}={\bf\Gamma}(\tau)$, provided that there exist values $\alpha,\beta\in I$ such that $\alpha<\tau<\beta$ and
\begin{equation}
\label{eqn 4.1.5}
{\bf\Gamma}(\alpha)\leq \tilde{\bf\Gamma}\leq{\bf\Gamma}(\beta).
\end{equation}
\end{theorem}

\noindent
{\bf Proof.} For $\alpha,\beta\in I$, we set $U_\alpha(p):=U({\bf\Gamma}_\alpha;p)$, $U_\beta(p):=U({\bf\Gamma}_{\beta};p)$, and $\tilde{U}(p):=U(\tilde{{\bf\Gamma}};p)$ in the closures of the (either always annular or always $P$-periodic and strip-like) domains $\Omega_\alpha:=\Omega({\bf\Gamma}_{\alpha})$, $\Omega_\beta:=\Omega({\bf\Gamma}_\beta)$ and $\tilde{\Omega}:=\Omega(\tilde{{\bf\Gamma}})$, respectively (where ${\bf\Gamma}_t:={\bf\Gamma}(t)$ for all $t\in\Re$). By (\ref{eqn 4.1.3}), (\ref{eqn 4.1.4}), and (\ref{eqn 4.1.5}), there exists a maximum (resp. minimum) value $\alpha\in I$ (resp. $\beta\in I$) such that the first (resp. second) inequality in (\ref{eqn 4.1.5}) holds. For maximum $\alpha\in I$ such that ${\bf\Gamma}_\alpha\leq\tilde{{\bf\Gamma}}$, there exists a point $p_{\alpha,1}\in\Gamma_{\alpha,1}\cap\tilde{\Gamma}_1$, or else a point $p_{\alpha,2}\in\Gamma_{\alpha,2}\cap\tilde{\Gamma}_2$, or both. In the first case, we have $U_\alpha(p_{\alpha,1})=\tilde{U}(p_{\alpha,1})=0$, where $U_\alpha(p)\geq \tilde{U}(p)\geq 0$ in ${\rm Cl}(\Omega_\alpha\cap\tilde{\Omega})$ by comparison principles. In view of (\ref{eqn 4.1.4}), it follows from this that 
$$\lambda_1(\tau)\,a_1(p_{\alpha,1})=\big|{\nabla} \tilde{U}(p_{\alpha,1})\big|\leq\big|{\nabla}U_\alpha(p_{\alpha,1})\big|=\lambda_1(\alpha)\,a_1(p_{\alpha,1}).$$
Therefore, $\lambda_1(\alpha)\geq \lambda_1(\tau)$, which contradicts (\ref{eqn 4.1.2}) if $\alpha<\tau$ in $I$. In the second case, we have $U_\alpha(p_{\alpha,2})=\tilde{U}(p_{\alpha,2})=1$, where $\tilde{U}(p)\leq U_\alpha(p)\leq 1$ in ${\rm Cl}(\Omega_\alpha\cap\tilde{\Omega})$ by comparison principles. In view of (\ref{eqn 4.1.1}), it follows from this that $$\lambda_2(\tau)\,a_2(p_{\alpha,2})=\big|{\nabla} \tilde{U}(p_{\alpha,2})\big|\geq\big|{\nabla} U_\alpha(p_{\alpha,2})\big|=\lambda_2(\alpha)\,a_2(p_{\alpha,2}).$$
Therefore $\lambda_2(\alpha)\leq \lambda_2(\tau)$, which again contradicts (\ref{eqn 4.1.2}) if $\alpha<\tau$ in $I$. We conclude in either case that if $\alpha$ is maximum in $I$ subject to (\ref{eqn 4.1.5}), then $\alpha\geq \tau$. A similar argument shows that if $\beta\in I$ is minimum subject to (\ref{eqn 4.1.5}), then $\beta\leq \tau$. In view of (\ref{eqn 4.1.2}) and (\ref{eqn 4.1.5}), which together imply that $\alpha\leq\beta$, it follows that $\alpha=\beta=\tau$, and therefore that $\tilde{\bf\Gamma}={\bf\Gamma}(\tau)$.
\vspace{.1in}

\begin{lemma} 
{(Continuation of positively-ordered, continuously-varying, para-metrized solution-families)}
\label{lem 4.1.1}
In the context of Def. \ref{def 4.1.1}, let be given a real open interval ${I}$ and a continuous function ${\bm\lambda}(t)=(\lambda_1(t),\lambda_2(t)):{\rm Cl}(I)\rightarrow\Re_+^2$ satisfying (\ref{eqn 4.1.2}) relative to ${I}$. Given two  intersecting open intervals $\tilde{I}$ and $\hat{I}$, both contained in $I$, let $\tilde{\bf\Gamma}(t):\tilde{I}\rightarrow{\bf X}$ and $\hat{\bf\Gamma}(t):\hat{I}\rightarrow{\bf X}$ denote two parametrized families of ${\bf\lambda}(t)$-solutions of Prob. \ref{prob 4I.}, each satisfying (\ref{eqn 4.1.3}) and (\ref{eqn 4.1.4}) relative to it's respective (open) $t$-interval. Assume there exist values $\alpha\in\tilde{I}$, $\beta\in\hat{I}$, and $\tau\in\tilde{I}\cap\hat{I}$ such that $\alpha<\tau<\beta$ and either ($\boldsymbol{i}$)
$\tilde{\bf\Gamma}(\alpha)<\hat{\bf\Gamma}(\tau)<\tilde{\bf\Gamma}(\beta)$ or ($\boldsymbol{i}\boldsymbol{i}$) $\hat{\bf\Gamma}(\alpha)<\tilde{\bf\Gamma}(\tau)<\hat{\bf\Gamma}(\beta)$ (in terms of the partial ordering in ${\bf X}$). Then $\tilde{\bf\Gamma}(t)=\hat{\bf\Gamma}(t)$ for all $t\in \tilde{I}\cap\hat{I}$. Therefore we can define a parametrized family $\dot{\bf\Gamma}(t):\dot{I}\rightarrow{\bf X}$ of solutions of Prob. \ref{prob 4I.}, satisfying (\ref{eqn 4.1.3}) and (\ref{eqn 4.1.4}) relative to the larger interval $\dot{I}:=\tilde{I}\cup\hat{I}$, such that $\dot{\bf\Gamma}(t)=\tilde{\bf\Gamma}(t)$ in $\tilde{I}$ and $\dot{\bf\Gamma}(t)=\hat{\bf\Gamma}(t)$ in $\hat{I}$. 
\end{lemma}
\vspace{.1in}

\noindent
{\bf Proof} If the assumptions are satisfied by the value $\tau\in\tilde{I}\cap\hat{I}$ (for some $\alpha$, $\beta$), then it follows from Thm. \ref{thm 4.1.1} that $\tilde{\bf\Gamma}(\tau)=\hat{\bf\Gamma}(\tau)$, from which it follows by continuity that the assumptions are satisfied by any value $t\in\tilde{I}\cap\hat{I}$ such that $|t-\tau|$ is sufficiently small. Thus if $\tilde{\bf\Gamma}(\tau)=\hat{\bf\Gamma}(\tau)$, then by Thm. \ref{thm 4.1.1} we have $\tilde{\bf\Gamma}(t)=\hat{\bf\Gamma}(t)$ for all $t\in \tilde{I}\cap\hat{I}$ such that $|t-\tau|$ is sufficiently small. The assertion follows. 
\vspace{.1in}

\begin{definition} (with application to Thm. \ref{thm 4.1.2})
\label{def 2.7.2}
In the context of Prob. \ref{prob 4I.+}, given a vector ${\bm \lambda}=(\lambda_1,\lambda_2)\in\Re_+^2$ and a positive value $\rho_0>0$, we use the notation $\boldsymbol{\cal R}({\bm\lambda};\rho_0)$ to denote the set of all solution pairs ${\bf\Gamma}\in{\bf X}(\overline{G})$ at $\boldsymbol{\lambda}$ such that 
\begin{equation}
\label{eqn 2.7.1-1} {\rm capacity}({\Omega})\,|\mu| \leq\pi-\rho_0,
\end{equation}

\noindent
where
$${\rm capacity}\big({\Omega}):=
\lambda_1\int_{\Gamma_1}\,a(p)\,ds=\lambda_2\,
\int_{\Gamma_2}\,a(p)\,ds$$
\noindent 
denotes the capacity of one $P$-period (in $x$) of the domain $\Omega:=\Omega({\bf\Gamma})$, 
and where we set $\mu={\rm  ln}\big(\lambda_2\big/\lambda_1\big)$. Given a positive interval $\Lambda_+:=\big[\,\underline{\lambda},\overline{\lambda}\,\big]$ (where $0<\underline{\lambda}<\overline{\lambda}<\infty$), we use $\boldsymbol{\cal R}(\Lambda_+;\rho_0)$ to denote the union of the sets $\boldsymbol{\cal R}({\bm \lambda};\rho_0)$ over all ${\bm\lambda}$ such that $\lambda_1, \lambda_2\in\Lambda_+$
(compare to Thms. \ref{thm 4I.0}, \ref{thm 4I.2}, and Lem. \ref{lem 4B.2}).
\end{definition}

\begin{definition} (with application to Thm. \ref{thm 4.1.2})
\label{def 4A}
In the context of Prob. \ref{prob 4I.+}, for any given vector ${\bm\lambda}=(\lambda_1,\lambda_2)\in\Re_+^2$ and positive scalar $r>0$, we use the notation: ${\boldsymbol{\cal S}}({\bm\lambda};r)$ to denote the set of all classical solutions ${\bf\Gamma}\in{\bf X}(\overline{G})$ at ${\bm\lambda}$ such that the positive flow-speed function $a(p):\Re^2\rightarrow\Re_+$ in $C^{2,1}$ is weakly logarithmically subharmonic in $\Omega:=\Omega({\bf\Gamma})$, and is such that
\begin{equation}
\label{eqn 4A.1} 
E({\bf\Gamma})\,>\,r\,+\,\frac{\mu^2}{1+\sqrt{\mu^2+1}}\,=\,r\,+\,\sqrt{\mu^2+1}\,-\,1,
\end{equation}
where $\mu:={\rm ln}\big(\lambda_2\big/\lambda_1\big)$, and 
\begin{equation}
\label{eqn 4A.3}
E({\bf\Gamma}):=\inf\big\{\big(|{\nabla} W(p)|\big/|{\nabla} U(p)|\big);p\in \partial\Omega\big\}.
\end{equation} 
Here, we define the functions $U(p):=U({\bf\Gamma};p)$ and $W(p):=W({\bf\Gamma};p)$, both in the closure of $\Omega:=\Omega({\bf\Gamma})$, such that $W({\bf\Gamma};p)$ solves the Dirichlet problem:
\begin{equation}
\label{eqn 4A.4}
\Delta W=\Delta\,{\rm ln}(a(p))\geq 0\,\,{\rm in}\,\, \Omega;\,\, W(\partial\Omega)=0.
\end{equation}
\end{definition}
\begin{remark}
\label{rem 4I.Lo}
(a) For $E=E({\bm\Gamma})$, the inequality: (\ref{eqn 4A.1}) is equivalent to the following inequalities (related to (\ref{eqn 4E.31}), (\ref{eqn 4E.32}), (\ref{eqn 4E.33})):
\begin{equation}
\label{eqn 4A.2}
r<1+E-|\mu|\,\,{\rm and}\,\,(1+E+\mu-r)(1+E-\mu-r)>1.
\end{equation}
(b) Let $S$ denote any subset of $\Re_+^2$. Then, in either of Defs. \ref{def 2.7.2} and Def. \ref{def 4A}, the defining condition (namely condition (\ref{eqn 2.7.1-1}) or (\ref{eqn 4A.1}), resp.), is satisfied by all classical solutions ${\bf\Gamma}\in{\bf X}(G)\cap C^{\,3,\tilde{\varrho}}$ of Prob. \ref{prob 4I.+} at ${\bm\lambda}\in S$, provided that it is satisfied by all classical solutions at ${\bm\lambda}\in S$ such that $|\mu|\leq \mu_0:={\rm ln}\big(\,\overline{A}\,\big/\underline{A}\,\big)$ (because there are no classical solutions at pairs ${\bm\lambda}\in S$ such that $|\mu|>\mu_0$; see Thm. \ref{thm 4I.0}). 
\end{remark}

\noindent
\begin{lemma}
\label{lem Lie} 
Assume in Prob. \ref{prob 4I.+} that there exist constants $0<\delta\leq H$ such that $\delta\leq\Delta\,{\rm ln}\big(a(p)\big)\leq H$ throughout $G$. Then for any fixed constant $\rho_0\in(0,\pi)$ and any sufficiently small constant $r_0>0$, we have that (a) ${\bf\Gamma}\in\bm{\cal{R}}\big({\bm\lambda};\rho_0\big)$ and (b) 
${\bf\Gamma}\in\bm{\cal{S}}\big({\bm\lambda};r_0\big)$, both uniformly among all vectors 
${\bm\lambda}\in\Re_+^2$ such that $|{\bm\lambda}|$ (resp. $|{\bm\lambda}|\,|\mu|$) is sufficiently large (sufficiently small), and among all classical solutions ${\bf\Gamma}=(\Gamma_1,\Gamma_2)\in{\bf X}(G)\cap C^{\,3,\tilde{\varrho}}$ of Prob. \ref{prob 4I.+} corresponding to these vectors.
\end{lemma}

\noindent
{\bf Proof.} Regarding assertion (a), by Thm. \ref{thm Z.2} there exists a constant $L$ such that $||\Gamma_1||,||\Gamma_2||\leq L$, uniformly for all classical solutions corresponding to vectors ${\bm\lambda}$ such that $|{\bm\lambda}|$ (resp. $|{\bm\lambda}||\mu|$) is sufficiently large (small), where $||\Gamma_i||$ denotes the arc-length of one $P$-period of a periodic arc-component $\Gamma_i$ of ${\bf\Gamma}$ and $\mu:={\rm ln}\big(\lambda_2\big/\lambda_1\big)$). Therefore, we have
\begin{equation}
\label{eqn 2.7.1-2}
{\rm capacity}(\Omega)\,|\mu|\leq {\rm max}\{\lambda_1,\lambda_2\}\,\overline{A}\,L\,|\mu|\leq \pi-\rho_0,
\end{equation}
uniformly for sufficiently large $|{\bm\lambda}|$ and small $|{\bm\lambda}|\,|\mu|$, where the first inequality in (\ref{eqn 2.7.1-2}) is an estimate of the capacity of $\Omega=\Omega({\bf\Gamma})$.
\vspace{.1in}

\noindent
Turning to assertion (b), it follows from the curvature estimate (\ref{eqn Z.1}) in Lem. \ref{lem Z.1} that the absolute curvatures of the boundary arcs $\Gamma_1, \Gamma_2$ of classical solutions ${\bf\Gamma}$ of Prob. \ref{prob 4I.+} at the vectors ${\bm\lambda}\in\Re_+^2$ are uniformly bounded from above relative to all ${\bm\lambda}\in\Re_+^2$ such that $|{\bm\lambda}|$ (resp. $|{\bm\lambda}|\,|\mu|$) is sufficiently large (sufficiently small). On the other hand, we have ${\rm dist}(\Gamma_1,\Gamma_2)\geq\big(1\big/\,\overline{A}\,\,|{\bm\lambda}|\big)$ for any solution at ${\bm\lambda}\in\Re_+^2$. Therefore, if $|{\bm\lambda}|$ (resp. $|{\bm\lambda}|\,|\mu|$) is sufficiently large (sufficiently small), then for any point $p_0\in\partial\Omega$ (where $\Omega=\Omega({\bf\Gamma}))$, there exists a ball $B_R(q_0)$ with center-point $q_0\in\Omega$ and radius $R=\big(1\big/2\overline{A}\,\,|{\bm\lambda}|\big)$, such that $B_R(q_0)\subset\Omega$ and $p_0\in\partial B_R(q_0)$. For any $\delta>0$, we define the functions $W_{R,\delta}(p):{\rm Cl}\big({B}_R(q_0)\big)\rightarrow\Re$, such that $W_{R,\delta}(p):=(\delta/4)\big(|p-q_0|^2-R^2\big)$. Observe that $W_{R,\delta}(p)<0$, $\Delta W_{R,\delta}(p)=\delta$, and $\nabla W_{R,\delta}(p)=(\delta/2)(p-q_0)$, all at any point $p\in B_{R}(q_0)$, 
while obviously $W_{R,\delta}(p)=0$ for all $p\in\partial B_R(q_0)$. For a given function $W(p):{\rm Cl}\big(\Omega\big)\rightarrow\Re$ such that $W(\partial\Omega)=0$ and $\delta\leq\Delta\,W(p)\leq H$ in $\Omega$ for some small value $\delta\in(0,H]$, one sees that $W(p)<0$ in $\Omega$ and $W(p)\leq 0$ in ${\rm Cl}(\Omega)$, both by the maximum principle. By an application of the comparison principle, one sees that $W(p)\leq W_{R,\delta}(p)\leq 0$ for all $p\in B_R(q_0)$. It follows by the Hoepf boundary-point lemma that $|\nabla W(p_0)|\geq|\nabla W_{R,\delta}(p_0)|=(R\delta/4)$, and therefore that $\big(|\nabla W(p_0)|\big/|\nabla U(p_0)|\big)\geq \big(\delta\big/8\,\overline{A}^2\,|{\bm\lambda}|^2\big)$, where the point $p_0$ is arbitrary in $\partial\Omega$. In view of (\ref{eqn 2.7.1-2}), the assertion follows from this.
\vspace{.1in}

\begin{theorem} 
{(Existence of a locally Lipschitz-continuously-varying, strictly-positively-ordered, global family of solutions of Prob. \ref{prob 4I.+} containing a specified solution)}
\label{thm 4.1.2}
In the context of Prob. \ref{prob 4I.+} (in which the function $a\in\boldsymbol{\cal A}\cap C^{\,3,\tilde{\varrho}}$ satisfies $\Delta\,{\rm ln}\big(a(p)\big)\geq 0$ in $G$), let be given a classical solution $\hat{\bf\Gamma}=(\hat{\Gamma}_1,\hat{\Gamma}_2)\in{\bf X}(G)\cap C^{3, \tilde{\varrho}}$ at the vector value $\hat{\bm\lambda}=(\hat{\lambda}_1,\hat{\lambda}_2)\in\Re_+^2$ such that $\hat{\bm\Gamma}\in\boldsymbol{\cal R}\big(\hat{\bm\lambda};\rho_0\big)\cap\boldsymbol{\cal S}\big(\hat{\bm\lambda}; r_0\big)$ for some values $\rho_0\in(0,\pi)$ and $r_0>0$ (see defs \ref{def 2.7.2} and \ref{def 4A}). 
Also let be given a value $t_0\in\Re$ and a smooth vector-valued function
$${\bm\lambda}(t)=(\lambda_1(t),\lambda_2(t)): \Re\rightarrow\Re_+^2$$
such that ${\bm\lambda}(t_0)=\hat{\bm\lambda}$ and such that the first (resp. second) component $\lambda_1(t)$ (resp. $\lambda_2(t)$) strictly decreases (increases) from $+\infty$ (resp. $0$) to $0$ (resp, $+\infty$) as $t$ increases in $\Re$. Then:
\vspace{.1in}

\noindent
(a) There exist a real open interval $I$ containing $t_0$ and a local Lipschitz-continu-ously varying, positively ordered parametrized solution family
$${\bf\Gamma}(t)=(\Gamma_1(t),\Gamma_2(t)): I\rightarrow{\bf X}(G)\cap C^{3;\varrho}$$
such that ${\bf\Gamma}(t_0)=\hat{\bf\Gamma}$ and also such that for each $t\in I$, 
${\bf\Gamma}(t)$ is a classical solution of Problem \ref{prob 4I.+} at ${\bm\lambda}(t)$. 
\vspace{.1in}

\noindent
(b) The above parametrization can be continued beyond either endpoint $\tau$ of the parameter interval $I$ (which is therefore not maximal) if ($\boldsymbol{i}$) there exists a value $\eta>0$ such that $N_\eta\big(\Omega({\bf\Gamma}(t)\big)\subset G$ for all $t\in I$, (where $N_\eta(S)$ denotes the $\eta$-neighborhood of $S$), and if ($\boldsymbol{i}\boldsymbol{i}$): there exist constants $r_0>0$ and $\rho_0\in(0,\pi)$ such that ${\bf\Gamma}(t)\in\boldsymbol{\cal R}\big({\bm\lambda}(t); \rho_0\big)\cap\boldsymbol{\cal S}\big({\bm\lambda}(t);\, r_0\big)$ for all $t\in I$.
\end{theorem}

\noindent
The proof as given at the conclusion of Section \ref{subsection 4.3}.
 
\begin{corollary}
\label{cor 4.1.2}
In the context of Thm. \ref{thm 4.1.2}, it follows from Thm. \ref{thm 4.1.1} that the arc-pair $\hat{{\bm\Gamma}}\in{\bf X}(G)\cap C^{\,3,\tilde{\varrho}}$ is the unique solution of Prob. \ref{prob 4I.+} at $\hat{{\bm\lambda}}$ among all curve-pairs $\tilde{{\bm \Gamma}}\in{\bf X}(G)$ such that ${\bm\Gamma}(t^-)<\tilde{{\bm\Gamma}}<{\bm\Gamma}(t^+)$ for some values $t^\pm\in I$ such that $\pm t^\pm>0$.
\end{corollary}
\begin{remark}
{(Local uniqueness result for Prob. \ref{prob 4I} in a more general context)}
\label{rem 4I.Q}
In the context of Thm. \ref{thm 4.1.1}, let $\big(\tilde{\Omega},\tilde{{\bf\Gamma}}\big)$ denote any $P$-periodic, multiply-connected, smooth open set $\tilde{\Omega}$ whose boundary $\partial\tilde{\Omega}$ is $P$-periodically partitioned into a pair of sets $\tilde{{\bf\Gamma}}=(\tilde{\Gamma}_1,\tilde{\Gamma}_2)$, each of which is a disjoint union of one infinite $P$-periodic arc and an infinite collection of simple closed curves. We let $|\nabla \tilde{U}(p)|=\lambda_i(\tau)\,a_i(p)$ on $\tilde{\Gamma}_i$, $i=1,2$, where $\tilde{U}(p):{\rm Cl}(\tilde{\Omega})\rightarrow\Re_+$ is a harmonic function in $\tilde{\Omega}$ satisfying the boundary conditions $\tilde{U}(\tilde{\Gamma}_i)=i-1$. Then $\tilde{{\bf\Gamma}}={\bf\Gamma}_\tau$. The proof is the same as before.
\end{remark}

\noindent
\begin{remark}
\label{rem }
{\it (Applicability of Thm. \ref{thm 4.1.2} and Cor. \ref{cor 4.1.2} to one-dimensional flows)} In the one-dimensional flow-model, the solution pairs $(x_1,x_2)$,  in which $x_1<x_2$, satisfy the equation: $\lambda_1 a(x_1)=(1/(x_2-x_1))=\lambda_2 a(x_2)$ corresponding to a given pair $(\lambda_1,\lambda_2)\in\Re_+^2$, where for specificity, we choose the flow-speed function: $a(x):={\rm exp}\big((\delta/2)\,x^2\big)$ (so that $\big({\rm ln}(a(x))\big)''=\delta$). In the context of Defs. \ref{def 2.7.2} and \ref{def 4A}, we thus have ${\rm capacity}(x_1,x_2)=\big(1/(x_2-x_1)\big)$, $E=(\delta/2)(x_2-x_1)^2$, and $$\mu:={\rm ln}\big(\lambda_2/\lambda_1\big)={\rm ln}\big(a(x_1)/a(x_2)\big)=(\delta/2)(x_1^2-x_2^2).$$
Therefore the condition: ${\rm capacity}(x_1,x_2)\,|\mu|<(\pi-\rho_0)$ is satisfied if 
$$(\delta/2)(x_2-x_1)^2<(\pi-\rho_0),$$
and the condition: $\mu^2<E$ is satisfied if
$$\delta\,(x_2+x_1)^2<2.$$
In other words, in the present example, these results are applicable to streams of width less than $\sqrt{2\,(\pi-\rho_0)/\delta}$ whose mid-lines are at a distance less than
$1/\sqrt{2\delta}$ from the center of the valley of the flow-speed function.
\end{remark}
\subsection{Uniqueness of classical solutions}
\label{subsection 4.2}
\begin{theorem}
\label{thm 4I.bib}
{(Proximity of two classical solutions)} 
In Prob. \ref{prob 4I.} (where the flow-speed function $a(p):G\rightarrow\Re_+$ is (weakly) logarithmically subharmonic), we let $\tilde{{\bf\Gamma}}=(\tilde{\Gamma}_1,
\tilde{\Gamma}_2)$ and $\hat{{\bf\Gamma}}=(\hat{\Gamma}_1,\hat{\Gamma}_2)$, both in ${\bf X}(G)$, denote a classical inner solution and a classical outer solution, respectively, such that $|{\nabla}\tilde{U}(p)|=a(p)$ on $\tilde{\Gamma}_2$ and $|{\nabla}\hat{U}(p)|=a(p)$ on $\hat{\Gamma}_1$. Then we have that $\tilde{\Gamma}_2\geq\hat{\Gamma}_1$ unless $a(p)$ is logarithmically harmonic throughout $\tilde{\Omega}\cup\hat{\Omega}$, where $\tilde{\Omega}:=\Omega(\tilde{\bf\Gamma})$ and $\hat{\Omega}:=\Omega(\hat{\bf\Gamma})$.
Similarly, if $\hat{\bf\Gamma}=(\hat{\Gamma}_1,\hat{\Gamma}_2)$ is the classical inner solution such that $|\nabla\hat{U}(p)|=a(p)$ on $\hat{\Gamma}_2$ and $\tilde{\bf\Gamma}=(\tilde{\Gamma}_1,\tilde{\Gamma}_2)$ is the classical outer solution such that $|\nabla\tilde{U}(p)|=a(p)$ on $\tilde{\Gamma}_1$, then $\hat{\Gamma}_2\geq\tilde{\Gamma}_1$ unless $a(p)$ is logarithmically harmonic throughout $\tilde{\Omega}\cup\hat{\Omega}$. Therefore, if $\tilde{\bf\Gamma}$ and $\hat{\bf\Gamma}$ are both classical solutions and $a(p)$ is not logarithmically harmonic in $\tilde{\Omega}\cup\hat{\Omega}$, then we have both $\tilde{\Gamma}_2\geq\hat{\Gamma}_1$ and $\hat{\Gamma}_2\geq\tilde{\Gamma}_1$, from which it follows that the two configurations are close enough so that the intersection of the closures of the two (doubly connected) stream beds contains a curve $\gamma\in{\rm X}$. 
\end{theorem}
\noindent
{\bf Proof.} For the proof of the first assertion, we let $\tilde{U}(p):=U(\tilde{{\bf\Gamma}};p)$ and $\tilde{\phi}(p):={\rm ln}\big(|{\nabla}\tilde{U}(p)|\big/a(p)\big)$, both in the closure of $\tilde{\Omega}:=\Omega(\tilde{{\bf\Gamma}})$. Then $\Delta\tilde{\phi}(p)\leq 0$ in $\tilde{\Omega}$, $\tilde{\phi}(p)=0$ on $\tilde{\Gamma}_2$, and $\tilde{\phi}(p)\geq 0$ on $\tilde{\Gamma}_1$. Therefore $\tilde{\phi}\geq 0$ in $\tilde{\Omega}$, from which it follows that 
\begin{equation}
\label{eqn 4B.1}
\tilde{\phi}_{\boldsymbol{\nu}}(p)=\tilde{K}_2(p)-\big(a_{\boldsymbol{\nu}}(p)\big/a(p)\big)\leq 0
\end{equation}
on $\tilde{\Gamma}_2$, and also, by a similar argument, that
\begin{equation}
\label{eqn 4B.1a}
\hat{\phi}_{\boldsymbol{\nu}}(p)=\hat{K}_1(p)-\big(a_{\boldsymbol{\nu}}(p)\big/a(p)\big)\geq 0
\end{equation}
on $\hat{\Gamma}_1$, where $\hat{\phi}(p):={\rm ln}\big(|{\nabla}\hat{U}(p)|\big/a(p)\big)$ in the closure of $\hat{\Omega}$, and where $\tilde{K}_2(p)$ and $\hat{K}_1(p)$ denote respectively the signed curvatures of the arcs $\tilde{\Gamma}_2$ and $\hat{\Gamma}_1$ at their respective points. If the first claim is false, then the open set $E(\tilde{\Gamma}_2)\cap D(\hat{\Gamma}_1)$ is nonempty. For specificity, we assume that ${\Gamma}_1$ intersects $\tilde{\Gamma}_2$ at least once. Then there exists a non-empty, connected region $\omega$ contained in the region $E(\tilde{\Gamma}_2)\cap D(\hat{\Gamma}_1)$, whose boundary $\partial\omega$ is the disjoint union $\partial\omega=\tilde{\gamma}_2\cup\hat{\gamma}_1$, where $\tilde{\gamma}_2:={\rm Cl}(\omega)\cap\tilde{\Gamma}_2$ and $\hat{\gamma}_1:={\rm Cl}(\omega)\cap\hat{\Gamma}_1$. Also the arcs $\tilde{\gamma}_2$ and $\hat{\gamma}_1$ have common initial and terminal endpoints $p_1$ and $p_2$. We have 
\begin{equation}
\label{eqn 4B.2}
\int_{\tilde{\gamma}_2}\tilde{K}_2(p)\,ds-\int_{\hat{\gamma}_1}\hat{K}_1(p)\,ds\leq\int_{\tilde{\gamma}_2}\big(a_{\boldsymbol{\nu}}(p)\big/a(p)\big)\,ds-\int_{\hat{\gamma}_1}\big(a_{\boldsymbol{\nu}}(p)\big/a(p)\big)\,ds,
\end{equation}
by integrating (\ref{eqn 4B.1}) and (\ref{eqn 4B.1a}). We also have
\begin{equation}
\label{eqn 4B.3}
\int_{\tilde{\gamma}_2} \tilde{K}_2(p)\,ds-\int_{\hat{\gamma}_1} \hat{K}_1p)\,ds=\big(\tilde{\theta}_2(p_2)-\tilde{\theta}_2(p_1)\big)-\big(\hat{\theta}_1(p_2)-\hat{\theta}_1(p_2)\big)\geq 0,
\end{equation}
where $\tilde{\theta}_2(p)$ and $\hat{\theta}_1(p)$ denote continuous arguments of the forward tangent vectors to the c*rves $\tilde{\gamma}_2$ and $\hat{\gamma}_1$ and where $p_1$ and $p_2$ are their joint initial and terminal points. We also have that 
\begin{equation}
\label{eqn 4B.4}
0\geq\int\int_{\omega}\Delta\,{\rm ln}(a(p))\,dA=\int_{\hat{\gamma}_2}\big(a_{\boldsymbol{\nu}}(p)\big/a(p)\big)\,ds-\int_{\tilde{\gamma}_1}\big(a_{\boldsymbol{\nu}}(p)\big/a(p)\big)\,ds,
\end{equation}
as follows by applying the divergence Theorem to the sub-harmonic function ${\rm ln}\big(a(p)\big)$ in $\omega$. By stringing the inequalities in (\ref{eqn 4B.2}), (\ref{eqn 4B.3}), and (\ref{eqn 4B.4}) into a chain of inequalities which are all in the same direction and begin and end with zero, we determine that all the inequalities in the chain are in fact equalities. Therefore, the double-integral over $\omega$ in (\ref{eqn 4B.4}) vanishes, and the two arc integrals (in (\ref{eqn 4B.4})) are equal; also the two signed curvature integrals in (\ref{eqn 4B.3})) are equal; in fact $\tilde{\theta}_2(p_1)-\hat{\theta}_1(p_1)=0$ and $\hat{\theta}_1(p_2)-\tilde{\theta}_2(p_2)=0$, since the latter two expressions are non-negative and sum to zero. Since ${\rm ln}\big(a(p)\big)$ is sub-harmonic in $\omega$, it follows that it is harmonic in $\omega$. Finally, we conclude by integrating the equations  
in (\ref{eqn 4B.1}), (\ref{eqn 4B.1a}) that
\begin{equation}
\label{eqn 4B.5}
\int_{\tilde{\gamma}_2}\tilde{\phi}_{\boldsymbol{\nu}}(p)\,ds=0=\int_{\hat{\gamma}_1}\hat{\phi}_{\boldsymbol{\nu}}(p)\,ds.
\end{equation}
Since $\tilde{\phi}_{\boldsymbol{\nu}}(p)\leq 0$ on $\tilde{\gamma}_2$ and $\hat{\phi}_{\boldsymbol{\nu}}(p)\geq 0$ on $\hat{\gamma}_1$ (by (\ref{eqn 4B.1}), (\ref{eqn 4B.1a})), it follows that $\tilde{\phi}_{\boldsymbol{\nu}}(p)=0$ on $\tilde{\gamma}_2$ and $\hat{\phi}_{\boldsymbol{\nu}}(p)=0$ on $\hat{\gamma}_1$.
Since the functions $\tilde{\phi}$ and $\hat{\phi}$ are already known to be non-negative $  $and (weakly) super-harmonic in the respective regions $\tilde{\Omega}$ and $\hat{\Omega}$, and to vanish on $\tilde{\Gamma}_2$ and $\hat{\Gamma}_1$, resp., that $\tilde{\phi}(p)=0$ in $\tilde{\Omega}$ and $\hat{\phi}(p)=0$ in $\hat{\Omega}$. Therefore, ${\rm ln}\big(a(p)\big)$ coincides with the harmonic functions ${\rm ln}\big(|{\nabla}\tilde{U}(p)|\big)$ (resp. ${\rm ln}\big(|{\nabla}\hat{U}(p)|\big)$) in $\tilde{\Omega}$ (resp. $\hat{\Omega}$). Therefore ${\rm ln}\big(a(p)\big)$ is a $C^2$-function which is harmonic in $\tilde{\Omega}\cup\,\omega\,\cup\hat{\Omega}$, proving the first claim in the case where $\hat{\Gamma}_1$ and $\tilde{\Gamma}_2$ intersect. In the alternate case where they do not intersect, one sets $\omega=E(\tilde{\Gamma}_2)\cap D(\hat{\Gamma}_1)$ and procedes with  the same argument. Finally, the second assertion is equivalent to the first and the third assertion is equivalent to the first two. 
\begin{theorem}
\label{thm 4B.bjb}
{(Uniqueness of classical solutions)} Assume in the context of Prob. \ref{prob 4I.} (where we assume the given flow-speed function $a(p):G\rightarrow\Re_+$ to be weakly logarithmically subharmonic in the annular region $G$) that there exist two (not necessarily distinct) classical solutions ${\bf\Gamma}_1:=(\Gamma_{1,1},\Gamma_{1,2})\in{\bf X}(G)\cap C^{\,3,\tilde{\varrho}}$ and ${\bf\Gamma}_2:=(\Gamma_{2,1},\Gamma_{2,2})\in{\bf X}(G)\cap C^{\,3,\tilde{\varrho}}$. For $i=1,2$, let $\Omega_i:=\Omega({\bf\Gamma}_i)$ and $U(p):=U({\bf\Gamma}_i;p)$ denote the corresponding flow channels and stream functions. Then:
\vspace{.1in}

\noindent
(a) If $\Gamma_{1,1}=\Gamma_{2,1}$ or $\Gamma_{1,2}=\Gamma_{2,2}$, then in fact 
$\bm\Gamma_1=\bm\Gamma_2$.
\vspace{.1in}

\noindent
(b)  If the inequality $\Gamma_{1,1}\leq\Gamma_{2,2}$ holds, but not the  inequality  $\Gamma_{1,2}\leq\Gamma_{2,2}$, and  if the inequality $\Gamma_{2,1}\leq\Gamma_{1,2}$ holds, but  not the inequality $\Gamma_{2,1}\geq\Gamma_{1,1}$, then the function $a(p):G\rightarrow\Re_+$ is identically logarithmically harmonic in $\Omega_1\cup\Omega_2$.
The same argument, but with the roles of ${\bf\Gamma}_1$ and ${\bf\Gamma}_2$ reversed, shows that if $\Gamma_{2,1}\leq\Gamma_{1,2}$ but not $\Gamma_{2,2}\leq\Gamma_{1,2}$ and if $\Gamma_{1,1}\leq\Gamma_{2,2}$ but not $\Gamma_{1,1}\geq\Gamma_{2,1},$, then again the function $a(p):G\rightarrow\Re_+$ is identically logarithmically harmonic in $\Omega_1\cup\Omega_2$.

\noindent
(c) If the flow-speed function $a(p)$ is not identically logarithmically harmonic throughout $\Omega_1\cup\Omega_2$, then one of the stream beds $\Omega_1$ and $\Omega_2$ must be a subset of the other.
\vspace{.1in}

\noindent
(d) If we have either $\Omega_1\subset\Omega_2$ or $\Omega_2\subset\Omega_1$, then $\Omega_1=\Omega_2$, except possibly if the function $a(p):G\rightarrow\Re_+$ is identically harmonic in $\Omega_1\cup\Omega_2$
                           
\noindent
(e) If $\Omega_1\subset\Omega_2$ and ${\rm ln}(a(p))$ is not identically harmonic throughout $\Omega_1\cup\Omega_2$, then either $\Omega_1=\Omega_2$ or else ${\rm Cl}(\Omega_1)\subset\Omega_2$.
\vspace{.1in}

\noindent
(f) Obviously if the function $a(p):G\rightarrow\Re_+$ is strictly logarithmically subharmonic, then $a(p)$ is not logarithmically harmonic throughout 
$\Omega_1\cup\Omega_2$. Therefore, it suffices to assume in Parts (a)-(e) above, and also in Thm. \ref{thm 4I.bib}, that $a(p):G\rightarrow\Re_+$ is strictly logarithmically subharmonic throughout $G$. 
\end{theorem}
\noindent
{\bf Proof of Part (a)}
Assume that $\Gamma_{1,2}=\Gamma_{2,2}$, but that it is not true that $\Gamma_{1,1}\geq\Gamma_{2,1}$. Then there exists a level curve $\dot{\Gamma}_{2,1}$ of the capacitary potential $U_2(p):{\rm Cl}(\Omega_2)\rightarrow[0,1]$ at minimum altitude $\delta>0$ subject to the requirement that $\dot{\Gamma}_{2,1}\geq\Gamma_{1,1}$. We have that $U_1(p)\geq\tilde{U}_2(p):=((U_2(p)-\delta)/(1-\delta))$ throughout $\tilde{\Omega}_2:=\{p\in\Omega_2:U_2(p)>\delta\}$ by maximum principles for harmonic functions, from which it follows that $|\nabla U_1(\tilde{p})|\geq|\nabla\tilde{U}_2(\tilde{p})|$ at any point $\tilde{p}\in\Gamma_{1,1}\cap\dot{\Gamma}_{2,1}$. Finally, we have that: $|\nabla U_2(p)|\geq a(p)$ throughout ${\rm Cl}(\Omega_2)$, since the superharmonic function $\phi_2(p)={\rm ln}(|\nabla U_2(p)|/a(p)):{\rm Cl}(\Omega_2)\rightarrow\Re$ vanishes on 
$\partial\Omega_2$. As a consequence of the above, we have that

$$(1-\delta)\,a(\tilde{p})=(1-\delta)|\nabla U_1(\tilde{p})|\geq(1-\delta)|\nabla \tilde{U}_2(\tilde{p})|=|\nabla U_2(\tilde{p})|=a(\tilde{p}).$$
This contradiction proves that in fact $\Gamma_{1,1}\geq\Gamma_{2,1}$. By reversing the 
roles of ${\bf\Gamma}_1$ and ${\bf\Gamma}_2$, we also conclude that $\Gamma_{2,1}\geq\Gamma_{1,1}$. Therefore if $\Gamma_{1,2}=\Gamma_{2,2}$, then $\Gamma_{1,1}=\Gamma_{2,1}$. 
\vspace{.1in}

\noindent
{\bf Proof of Part (b)} Assuming that the inequality $\Gamma_{1,1}\leq\Gamma_{2,2}$ is satisfied, but the inequality  $\Gamma_{1,2}\leq\Gamma_{2,2}$ is not satisfied, there exists a level curve $\dot{\Gamma}_{1,2}$ of the capacitary potential $U_1(p):{\rm Cl}(\Omega_1)\rightarrow[0,1]$ at maximum altitude $C<1$ subject to the requirement that $\dot{\Gamma}_{1,2}\leq\Gamma_{2,2}$. Also, assuming that the inequality: $\Gamma_{2,1}\leq\Gamma_{1,2}$ is satisfied, but the inequality: $\Gamma_{2,1}\geq\Gamma_{1,1}$ is not satisfied, there exists a level curve $\dot{\Gamma}_{2,1}$ of the capacitary potential $U_2(p):{\rm Cl}(\Omega_2)\rightarrow[0,1]$ at the minimum altitude $\delta>0$ subject to the requirement 
that $\dot{\Gamma}_{2,1}\geq\Gamma_{1,1}$. Thus, we have that
$$\dot{\Gamma}_{1,2}=\{p\in\Omega_1:U_1(p)=C\}\,\,{\rm and}\,\,\dot{\Gamma}_{2,1}=\{p\in\Omega_2:U_2(p)=\delta\},$$ where 
$$C:=\min\{U_1(p):p\in\Gamma_{2,2}\cap\Omega_1\}\,{\rm and}\,\,
\delta:=\max\{U_2(p):p\in\Gamma_{1,1}\cap\Omega_2\}.$$
We define the arc-pairs $\tilde{\bf\Gamma}_i:=(\tilde{\Gamma}_{i,1},\tilde{\Gamma}_{i,2})$, $i=1,2$, such that
-$$\tilde{\bf\Gamma}_1:=(\Gamma_{1,1},\dot{\Gamma}_{1,2})\,\,{\rm and}\,\,
\tilde{\bf\Gamma}_2:=(\dot{\Gamma}_{2,1},\Gamma_{2,2}).$$
It directly follows from the above that 
$$\tilde{U}_1(p):=U(\tilde{\bf\Gamma}_1;p)=U_1(p)/C$$
in the closure of $\tilde{\Omega}_1:=\Omega(\tilde{\bf\Gamma}_1)$ and 
$$\tilde{U}_2(p):=U(\tilde{\bf\Gamma}_2;p)=((U_2(p)-\delta)/(1-\delta))$$
the closure of $\tilde{\Omega}_2:=\Omega(\tilde{\bf\Gamma}_2)$. Therefore, we have 
$$|\nabla\tilde{U}_1(p)|=(|\nabla U_1(p)|/C)\,\,{\rm in}\,\,\tilde{\Omega}_1\,\,{\rm and}
\,|\nabla\tilde{U}_2(p)|=(|\nabla U_2(p)|/(1-\delta))\,\,{\rm in}\,\, {\rm Cl}(\tilde{\Omega}_2)$$. We also have from the above that
$$\tilde{\Gamma}_{1,1}=\Gamma_{1,1}\leq\dot{\Gamma}_{1,1}=\tilde{\Gamma}_{2,1}\,\,{\rm and}\,\,\tilde{\Gamma}_{1,2}=\dot{\Gamma}_{1,2}\leq\Gamma_{2,2}=\tilde{\Gamma}_{2,2},$$ which is the proof by components that:
$$\tilde{\bf\Gamma}_1\leq\tilde{\bf\Gamma}_2.$$
In view of this inequality, it follows by maximum and comparison principles for harmonic functions that
$$\tilde{U}_1(p)\geq\tilde{U}_2(p)$$
throughout the closure of the region $\tilde{\omega}:=\tilde{\Omega}_1\cap\tilde{\Omega}_2$. Moreover, one has that $\tilde{U}_1(p)=\tilde{U}_2(p))=0$ for all points 
$p\in\tilde{\sigma}_1:=\tilde{\Gamma}_{1,1}\cap\tilde{\Gamma}_{2,1}\subset\partial\tilde{\omega}$ and 
$\tilde{U}_1(p)=\tilde{U}_2(p)=1$ for all points
$p\in\tilde{\sigma}_2=\tilde{\Gamma}_{1,2}\cap\tilde{\Gamma}_{2,2}\subset\partial\tilde{\omega}$, where it is clear that $\tilde{\sigma}_i\not=\emptyset$ for $i=1,2$. It  follows from this by the Hopf boundary point lemma that
$$|\nabla\tilde{U}_1(\tilde{p}_1)|\geq|\nabla\tilde{U}_2(\tilde{p}_1)|\, {\rm and}\,\,|\nabla\tilde{U}_1(\tilde{p}_2)|\leq|\nabla\tilde{U}_2(\tilde{p}_2)|$$ 
at any point $\tilde{p}_i\in\tilde{\sigma}_i$ for $i=1,2$. Finally, we have that 
$$|\nabla U_i(p)|\geq a(p)$$
throughout ${\rm Cl}(\Omega_i),$
$i=1,2$, because the superharmonic functions $\phi_i(p)={\rm ln}(|\nabla U_i(p)|/a(p)):{\rm Cl}(\Omega_i)\rightarrow\Re$, $i=1,2$, vanish on their respective domain boundaries $\partial\Omega_i$, $i=1,2$. As an application of all of the above, we conclude that
\begin{equation}
\label{eqn EZ1}
\frac{a(\tilde{p}_1)}{1-\delta}=\frac{|\nabla U_1(\tilde{p}_1)|}{1-\delta}=|\nabla\tilde{U}_1(\tilde{p}_1)|\geq|\nabla\tilde{U}_2(\tilde{p}_1)|=\frac{|\nabla U_2(\tilde{p}_1)|}{C}\geq\frac{a(\tilde{p}_1)}{C}
\end{equation}
and
\begin{equation}
\label{eqn EZ2}
\frac{a(\tilde{p}_2)}{1-\delta}\leq\frac{|\nabla U_1(\tilde{p}_2)|}{1-\delta}=|\nabla\tilde{U}_1(\tilde{p}_2)|\leq|\nabla\tilde{U}_2(\tilde{p}_2)|=\frac{|\nabla U_2(\tilde{p}_2)|}{C}=\frac{a(\tilde{p}_2)}{C},
\end{equation}
for $\tilde{p}_i\in\tilde{\sigma}_i$, $i=1,2$. It follows from Equations (\ref{eqn EZ1}) and (\ref{eqn EZ2}) together that $C=1-\delta$ and that all inequalities in (\ref{eqn EZ1}) and (\ref{eqn EZ2}) reduce to equations. In particular, we have that $|\nabla U_2(\tilde{p}_1)|=a(\tilde{p}_1)$ and $|\nabla U_1(\tilde{p}_2)|=a(\tilde{p}_2)$. Since $\tilde{p}_i\in\Omega_i$ for $i=1,2$, it follows by application of the strict maximum principle to the non-negative superharmonic functions $\phi_i(p):{\rm Cl}(\Omega_i)\rightarrow\Re$, $i=1,2$, that $|\nabla U_i(p)|=a(p)$ throughout $\Omega_i$ for $i=1,2$, and therefore that the logarithmically subharmonic function $a(p):G\rightarrow\Re_+$ is in fact logarithmically harmonic throughout 
$\Omega_1\cup\Omega_2$. 
\vspace{.1in}

\noindent
{\bf  Proof  of Part (c))}. In Part (b), let $A, B, C$, and $D$ denote the inequalities $\Gamma_{1,2}\leq\Gamma_{2,2}$, $\Gamma_{2,1}\geq\Gamma_{1,1}$, 
$\Gamma_{2,2}\leq\Gamma_{1,2}$, and $\Gamma_{1,1}\geq\Gamma_{2,1}$ in the given order. 
Under the assumptions of Thm. \ref{thm 4I.bib}, the inequalities $\Gamma_{1,1}\leq\Gamma_{2,2}$ and $\Gamma_{2,1}\leq\Gamma_{1,2}$ both hold automatically. Therefore we can slightly revise Part (b) to state that if both $A$ and $B$ do not hold or if both $C$ and $D$ do not hold, then the function ${\rm ln}(a(p))$ is identically harmonic in $\Omega_1\cup\Omega_2$. Therefore if the function ${\rm ln}(a(p))$ is not identically harmonic in $\Omega_1\cup\Omega_2$, then either $A$ or $B$ holds and either $C$ or $D$ holds. In other words, either $A$ and $C$ both hold or $A$ and $D$ both hold, or $B$ and $C$ both hold or else $B$ and $D$ both hold. But $A$ and $C$ (resp. $B$ and $D$) imply that $\Gamma_{1,2}=\Gamma_{2,2}$ (resp. $\Gamma_{1,1}=\Gamma_{2,1}$, from which it  follows by Part (a) that $\Omega_1=\Omega_2$. The inequalities $A$ and $D$ are satisfied 
only by $\Gamma_{2,1}\leq\Gamma_{1,1}<\Gamma_{1,2}\leq\Gamma_{2,2}$, i.e. $\Omega_1\subset\Omega_2$, while $B$ and $C$ are satisfied only by $\Gamma_{1,1}\leq\Gamma_{2,1}<\Gamma_{2,2}\leq\Gamma_{1,2}$, i.e. $\Omega_2\subset\Omega_1$
\vspace{.1in}

\noindent
{\bf Proof of Part (d)}We assume for the purpose of obtaining a contradiction that the classical solutions ${\bf\Gamma}_1=(\Gamma_{1,1},\Gamma_{1,2}))\in{\bf X}$ and  ${\bf\Gamma}_2=(\Gamma_{2,1}, \Gamma_{2,2})\in{\bf X}$ are distinct, the corresponding stream beds are such that $\Omega_1\subset\Omega_2$, and the given weakly logarithmically sub-harmonic flow-speed function $a(p):G\rightarrow\Re_+$ is not identically logarithmically harmonic throughout $\Omega_2$. For each $i=1,2$, we define the continuous function $\phi_i(p):{\rm Cl}\big(\Omega_i\big)\rightarrow\Re$ such that
\begin{equation}
\label{eqn 3ob}
\phi_i(p):={\rm ln}\big(|\nabla U_i(p)|\big/a(p)\big).
\end{equation}
We have that $|\nabla U_i(p)|=a(p)$ on $\partial\Omega_i$, from which it follows that $\phi_i(\partial\Omega_i)=0$ for $i=1,2$. Also for $i=1,2$ the function: ${\rm ln}(\big|\nabla U_i(p)\big|):\Omega_i\rightarrow\Re$ is harmonic. It follows from this by the assumed properties of $a(p):G\rightarrow\Re_+$ that the functions $\phi_i(p):{\rm Cl}(\Omega_i)\rightarrow\Re_+$, $i=1,2$, are weakly super-harmonic in their respective domain interiors and are therefore non-negative by the maximum principle. Also the non-negative, weakly super-harmonic function $\phi_2(p):{\rm Cl}(\Omega_2)\rightarrow\Re_+$ is not identically harmonic in $\Omega_2$ and is therefore strictly positive in $\Omega_2$ by the strict maximum principle. For the function $\psi(p):=\phi_2(p)-\phi_1(p):{\rm Cl}(\Omega_1)\rightarrow\Re$, we have $\Delta\psi(p)=\Delta\phi_2-\Delta\phi_1=0$ in $\Omega_1$ and $\psi(p)=\phi_2(p)-\phi_1(p)\geq 0$ on $\partial\Omega_1$, so that we have $\psi(p)\geq 0$ in ${\rm Cl}(\Omega_1)$ by the maximum principle. Also $\psi(p)=\phi_2(p)-\phi_1(p)=\phi_2(p)>0$ on $\Omega_2\cap\partial\Omega_1$. Therefore, we 	can't have $\psi(p)=0$ identically in ${\rm Cl}(\Omega_1)$, from which it follows that $\psi(p)>0$ throughout $\Omega_1$ by the strict maximum principle. To sum up, we have
\begin{equation}
\label{eqn 4ob}
\phi_1(p)\,\,<\,\,\phi_2(p)
\end{equation}
in $\Omega_1$, from which it follows that
\begin{equation}
\label{eqn 5ob}
|\nabla U_1(p)|\,\,<\,\,|\nabla U_2(p)|
\end{equation}
throughout $\Omega_1$, since $|\nabla U_i(p)|\,\,=\,\,a(p)\,{\rm exp}\big(\phi_i(p)\big)$ in $\Omega_i$ for $i=1,2$. Therefore, we have
\begin{equation}
\label{eqn 6ob}
\int_{\gamma_1}\,|\nabla U_1(p)|\,ds\,\,<\,\,\int_{\gamma_1}\,|\nabla U_2(p)|\,ds\,\,\leq\,\,\int_{\gamma_2}\,|\nabla U_2(p)|\,ds, 
\end{equation}
where $\gamma_2$ denotes any arc of steepest ascent of the  function $U_2(p)$ joining $\Gamma_{2,1}$ to $\Gamma_{2,2}$ through $\Omega_2,$ and $\gamma_1$ denotes the shortest sub-arc of $\gamma_2$ which connects $\Gamma_{1,1}$ to $\Gamma_{1,2}$ through $\Omega_1$.
On the other hand, by (\ref{eqn 1.2}) and the fundamental theorem of calculus we have that 
\begin{equation}
\label{eqn 7ob}
\int_{\gamma_i}\,\nabla U_i(p)\cdot {\bf n}_i(p)\,ds\,=\,U_i(p_{i,2})-U(p_{i,1})\,=\,1
\end{equation}
for $i=1,2$, where $\gamma_i$ denotes any smooth arc joining the initial point $p_{i,1}\in\Gamma_{i,1}$ to the terminal point $p_{i,2}\in\Gamma_{i,2}$ through $\Omega_i$, and where ${\bf n}_i(p)$ denotes the forward unit tangent vector to $\gamma_i$ at the point $p\in\gamma_i$. It follows directly from (\ref{eqn 7ob})in the case where $i=1$ by absolute value estimates that 
\begin{equation}
\label{eqn 8ob}
\int_{\gamma_1}\,|\nabla U_1(p)|)\,ds\,\,\geq\,\,1,
\end{equation}
where $\gamma_1$ is any smooth arc joining $\Gamma_{1,1}$ to $\Gamma_{1,2}$ through $\Omega_1$. 
In the context of (\ref{eqn 6ob}), we have 
\begin{equation}
\label{eqn 9ob}
\int_{\gamma_2}\,|\nabla U_2(p)|\,ds\,\,=\,\,1,
\end{equation}
as follows from (\ref{eqn 7ob}) in the case where $i=2$ by making the the substitution: ${\bf n}_2(p)=\big(\nabla U_2(p)\big/|\nabla U_2(p)|\big)$. Clearly the equations (\ref{eqn 6ob}), (\ref{eqn 8ob}), and (\ref{eqn 9ob}) are inconsistent and cannot all hold.  This contradiction proves the assertion.
\vspace{.1in}

\noindent 
{\bf Proof of Part (e)} Let $\phi_i(p):={\rm ln}\big(|\nabla U_i(p)|\big/a(p)\big)$ in ${\rm Cl}(\Omega_i)$ for $i=1,2$. Then $\phi_i(\partial\Omega_i)=0$ and $\Delta\phi_i(p)=\,-\Delta{\rm ln}\big(a(p)\big)\leq 0$. Since $\Delta{\rm ln}(a(p))$ does not vanish throughout $\Omega_i$, we have that $\phi_i(p)>0$ throughout $\Omega_i$ by the strict maximum principle. Therefore, assuming that $\Omega_1\not=\Omega_2$ and $\Omega_1\subset\Omega_2$, we have that $\phi_1(p)\leq\phi_2(p)$ first on $\partial\Omega_1$ and then throughout  ${\rm Cl}(\Omega_1)$, by maximum and comparison principles, and the fact that $\Delta(\phi_2-\phi_1)=0$ in $\Omega_1$. In fact  we have $\phi_1(p)<\phi_2(p)$ in $\Omega_1$ by the strict maximum principle, since we do not have $\phi_1(p)=\phi_2(p)$ everywhere on $\partial\Omega_1$. For the purpose of obtaining a contradiction, let $q$ denote a point in $\partial\Omega_1\cap\partial\Omega_2$. Observe that $\phi_1(q)=0=\phi_2(q)$. In the above context, the Hopf boundary-point Lemma implies that ($\boldsymbol{i}$): $0<\partial\phi_1(q)/\partial\nu<\partial\phi_2(q)/\partial\nu$, where $\nu$ denotes the interior normal to both $\partial\Omega_1$ and $\partial\Omega_2$ at $q$. But we must also have that ($\boldsymbol{i}\boldsymbol{i}$): $K_2(q)\leq K_1(q)$, where $K_i(q)$ denotes the counter-clockwise-oriented curvature of $\partial\Omega_i$ at $q$, $i=1,2$. But ($\boldsymbol{i}$) and ($\boldsymbol{i}\boldsymbol{i}$) together contradict the equations: $\partial\phi_i(q)/\partial\nu=K_i(q)-a_\nu(q)/a(q)$, $i=1,2$.
\vspace{.1in}

\begin{remark}
\label{rem 3.2} {(Uniqueness example)}
Given a positive function $a(r):\Re_+\rightarrow\Re_+$, the anulus $\Omega(r_1,r_2):=\{p\in\Re^2:r_1<|p|<r_2\}$ solves Problem \ref{prob 1.1} corresponding to the flow-speed function $A(p):=a(|p|):\Re^2\setminus\{0\}\rightarrow\Re_+$ if and only if the values $0<r_1<r_2<\infty$ together satisfy the equation: 
\begin{equation}
\label{eqn ex}
r_1\,a(r_1)=\big(1/{\rm ln}(r_2/r_1)\big)=r_2\,a(r_2).
\end{equation}
Now assume for some value $r_0>0$ that the related function $\phi(r):=r\,a(r):\Re_+\rightarrow\Re_+$ strictly decreases (resp. increases) from $\infty$
to $\phi_0:=\phi(r_0)>0$ (resp. from $\phi_0:=\phi(r_0)$ to $\infty$) as $r$ increases from $0$ to $r_0$ (resp. from $r_0$ to $\infty$). Then for any value $t>\phi_0$, there exist unique values $r_1(t)\in (0,r_0)$ and $r_2(t)\in(r_0,\infty)$ such that $\phi(r_1(t))=t=\phi(r_2(t))$. Moreover, the 
continuous function $r_1(t):(\phi_0,\infty)\rightarrow(0,r_0)$ (resp. $r_2(t):(\phi_0,\infty)\rightarrow(r_0,\infty)$) is decreasing (increasing), and the annulus $\Omega(r_1(t),r_2(t))$ satisfies (\ref{eqn ex}) if and only if $(r_2(t)/r_1(t))={\rm exp}(1/t)$, But as $t$ increases from $\phi_0$ to $\infty$, $(r_2(t)/r_1(t))$ increases from $1$ to $\infty$ while ${\rm exp}(1/t)$ decreases from ${\rm exp}(1/\phi_0)>1$ to $0$ Therefore, the two graphs intersect at at a unique value $t_0\in(\phi_0,\infty)$.
Therefore under the present assumptions, Problem \ref{prob 1.1} has one and only one annular solution centered at the origin, namely the annulus $\Omega(r_1(t_0),r_2(t_0))$. \end{remark}
\subsection{Continuously and monotonically 
varying solution families}
\label{subsection 4.3}

\begin{lemma}
\label{lem 4B.2} 
{(Uniform self-separation of solutions of Prob. \ref{prob 4I.+}. Uniform continuation of capacitary potentials)} 
Assume in the context of Prob. \ref{prob 4I.+} that the fixed positive $P$-periodic flow-speed function $a(p):\Re^2\rightarrow\Re_+$ is in the class $\boldsymbol{\cal {A}}\cap C^{\,3,\varrho}$ for some $\varrho\in(0,1]$, and is weakly-logarithmically-subharmonic relative to $G$. Then: (a) for any vector ${\bm\lambda}=(\lambda_1,\lambda_2)\in\Re_+^2$ and any classical solution ${\bm\Gamma}_{\boldsymbol{\lambda}}\in {\bf X}(G)\cap C^{\,3,\tilde{\rho}}$ of Prob. \ref{prob 4I.+} at $\boldsymbol{\lambda}$ such that ${\bf\Gamma}_{\boldsymbol{\lambda}}\in{\boldsymbol{\cal R}}(\Lambda_+;\rho_0)$ (for some given value $\rho_0\in(0,\pi$); see Def. \ref{def 2.7.2}, Lem. \ref{lem Lie} and Thm. \ref{thm Z.2}), we have:
\begin{equation}
\label{eqn 2.7.1-22}
|p_{\boldsymbol{\lambda},i}(t)-p_{\boldsymbol{\lambda},i}(\tau)|\geq\min\Big\{\frac{2}{\pi}\big|t-\tau\big|,\frac{2}{K_0}, \Big(\frac{\rho_0}{B_2}\Big)^2, R_0, R_1, 1\Big\}
\end{equation}
$i=1,2$, uniformly for all $t,\tau\in\Re$, where $K_0$ denotes a uniform upper bound for the absolute curvatures of the components of these classical solutions (the bound exists by Thm. \ref{thm 2.7.1}), also  $p_{\boldsymbol{\lambda},i}(t):\Re\rightarrow\Gamma_{\boldsymbol{\lambda},i}$ denotes any arc-length parametrization of $\Gamma_{\boldsymbol{\lambda},i}$, and, finally, the positive constants $B_2, R_0, R_1>0$, which remain to be specified, are such that $B_2, (1/R_0), (1/R_1)$ have upper bounds which depend only on $\underline{A}\,, \overline{A}, A_1$, $A_2$,  and ${\bm \lambda}$.
\vspace{.1in}

\noindent
(b) Assume in Prob. \ref{prob 4I.+} that the flow-speed function $a(p):\Re^2\rightarrow \Re_+$ in $\boldsymbol{{\cal A}}\cap C^{\,3,\varrho}$ is logarithmically subharmonic in $G$. Then, given the positive closed interval $\Lambda_+:=[\,\underline{\lambda}\,,\overline{\lambda}\,]\subset\Re_+$, the main assertions of Thm. \ref{thm 2.7.2} and Part (a) above both apply uniformly with respect to all pairs ${\bm\lambda}\in\Lambda_+^2$ and all solutions ${\bf\Gamma}_{\boldsymbol{\lambda}}\in{\bf X}(G)$ of Prob. \ref{prob 4I.+} at ${\bm\lambda}$. Therefore, the capacitary potentials $U_{\boldsymbol{\lambda}}(p):=U({\bf\Gamma}_{\boldsymbol{\lambda}};p)$ of all of these solutions can be continued as single-valued $C^{\,3,\tilde{\varrho}}$-functions in simply-connected domains $\Omega_{\boldsymbol{\lambda}}^*$ satisfying the same uniform bound on their $C^{\,3,\tilde{\varrho}}$-norms relative to their respective domains $\Omega_{{\boldsymbol\lambda},\delta}$ uniformly containing the $\delta$-neighborhoods $N_\delta(\Omega_{{\boldsymbol{\lambda}}})$ of the respective regions  $\Omega_{{\boldsymbol\lambda}}:=\Omega({\bf\Gamma}_{\boldsymbol{\lambda}})$ (with the same value $\delta>0$ in all cases). 
\end{lemma}

\noindent
{\bf Proof of Part (a).} Let be given a classical solution ${\bf\Gamma}=(\Gamma_1,\Gamma_2)\in{\bf X}(G)$ of Prob. \ref{prob 4I.+} at a fixed pair ${\bm \lambda}\in\Lambda_+^2$, and let $p_i(t):\Re\rightarrow\Gamma_i$ be the arc-length parametrizations of the components. We choose the positive direction on $\Gamma_1$ and $\Gamma_2$ such that $\boldsymbol{\nu}$ always points locally to the left, where we use $\boldsymbol{\nu}(p)$ to denote the unit vector pointing in the direction of $\nabla U_i(p)$ on $\partial\Omega$, and where $U_i(p):=U_i({\bf\Gamma};p)$ in (the closure of)  $\Omega:=\Omega({\bf\Gamma})$. (Thus the positive direction on $\Gamma_1$ (resp. $\Gamma_2$) corresponds to increasing (decreasing) $t$.) For any fixed $t_0\in\Re$, it follows from the absolute curvature bound that 
\begin{equation}
\label{eqn ZW}
|p_i(t)-p_i(t_0)|\geq(2/\pi)|t-t_0|\,\,{\rm whenever}\,\,|t-t_0|\leq(\pi/K_0)
\end{equation}
for $i=1,2$. Since $|p_i(t)-p_i(t_0)|\rightarrow +\infty$ as $t\rightarrow\pm\infty$, it suffices to show that if ${\bf\Gamma}\in {\boldsymbol{\cal R}}(\Lambda_+;\rho_0)$ for some fixed value $\rho_0\in(0,\pi)$ (see (\ref{eqn 2.7.1-1})), then for any $t_0\in\Re$, any local minimum of the function $r(p_i(t)):=|p_i(t)-p_i(t_0)|:\Re\rightarrow\Re_+$ must exceed a certain uniform positive lower bound. Given the points $p_{i,0}\in\Gamma_i$, $i=1,2$, it suffices to consider the point-sets ${\cal M}_i^\pm(p_{i,0})$, $i=1,2$, whose respective elements $p_i^\pm\in\Gamma_i$ are such that $p_i^+$ is a local minimizer of the function  $r_i(p):=|p-p_{i,0}|$ relative to $p\in\Gamma_i\cup\Omega$, whereas $p_i^-$ is a local minimizer of $r_i(p):=|p-p_{i,0}|$ relative to $p\in \Gamma_i\cup D_i(\Gamma_i)$, $i=1,2$. Then for $p_{i,0}\in\Gamma_i$ and $p_i^\pm\in{\cal M}_i^\pm(p_{i,0})$, we have
\begin{equation}
\label{eqn 2.7.1-3}
{\boldsymbol{\nu}}(p_i^\pm)|p_i^\pm-p_{i,0}|=\mp(p_i^\pm-p_{i,0}).
\end{equation}
It also follows from the uniform bound on the curvature of $\Gamma_1, \Gamma_2$ that 
\begin{equation}
\label{eqn 2.7.1-4}
\big|\theta\big({\boldsymbol{\nu}}(p_{i,0}),{\boldsymbol{\nu}}(p_i^\pm)\big)\big|\leq B_0\sqrt{\big|p_{i,0}-p_i^\pm\big|},
\end{equation}
where $B_0:=E_0\sqrt{K_0}$ for a dimensionless constant $E_0>0$, and where $\theta\big({\boldsymbol{\nu}}(p_{i,0}),$ ${\boldsymbol{\nu}}(p_i^\pm)\big)$ denotes the angle between ${\boldsymbol{\nu}}(p_{i,0})$ and ${\boldsymbol{\nu}}(p_i^\pm)$. We assert that, given points $p_{i,0}\in\Gamma_i$ and $p_i^+\in{\cal M}_i^+(p_{i,0})$ such that $r_i:=|p_i^+-p_{i,0}|\leq 1$, we must have 
\begin{equation}
\label{eqn 2.7.1-5}
B_2\sqrt{r_i}\geq B_1 r_i+B_0 \sqrt{r_i}\geq\pi-\mu_i\,{\rm capacity}(\Omega)\geq\pi-|\mu_i|\,{\rm capacity}(\Omega)\geq\rho_0
\end{equation}
(for $\mu_i:={\rm ln}\big(\lambda_{3-i}\big/\lambda_i\big)$), where $B_1=(A_1/\underline{A})$ and $B_2=B_0+B_1$. Here, the final inequality follows from the assumption that  ${\bf\Gamma}\in{\boldsymbol{\cal R}}(\Lambda_+;\rho_0)$ (for some given value $\rho_0\in(0,\pi$)). For the proof of the second (and main) inequality in (\ref{eqn 2.7.1-5}), we assume that $p_i^+\in{\cal M}_i^+(p_{i,0})$, where $p_{i,0}$ precedes $p_i^+$ in terms of the ordering on $\Gamma_i$. We let $\gamma_i$ denote an arc-segment of $\Gamma_i$ having $p_{i,0}$ (resp. $p_i^+$) as its initial (resp. terminal) end-point, and we use $L_i$ to denote the straight line-segment connecting the same points in the same order. Obviously $r_i:=||L_i||=|p_{i,0}-p_i^+|$. We can assume that $\gamma_i$ does not cross $L_i$, since otherwise one could choose a new relative minimum point $p_i^+\in\Gamma_i$ which is closer to $p_{i,0}$. One also sees by using (\ref{eqn 2.7.1-4}) and a uniform upper bound on the curvature of the fixed boundary $\partial G$ that there exists a value $R_0=R_0(K_0)>0$ such that ${\rm Cl}(\omega_i)\subset G$ whenever $r_i\in(0,R_0]$. We also have that $\phi_i(\Gamma_i)=0$, $\phi_i(p)=\mu_i:={\rm ln}(\lambda_{3-i}/\lambda_i)$ on $\Gamma_{3-i}$, and $\Delta\phi_i(p)=\,-\Delta\,{\rm ln}\big(a(p)\big)\leq 0$ in $\Omega$, where we define the function $\phi_i(p):={\rm ln}\big(|\nabla U_i(p)|\big/\lambda_i\,a(p)\big)$ in terms of the capacitary potential $U_i(p):=U_i({\bf \Gamma};p)$. Therefore, we have $\phi_i(p)\geq\mu_i\,U_i(p)$ in ${\rm Cl}(\Omega)$ by the comparison principle, since the same inequality holds on $\partial\Omega$. It follows that
\begin{equation}
\label{eqn 2.7.1-6}
\phi_{i,\boldsymbol{\nu}}=K_i(p)-({\rm ln}(a(p))_{\boldsymbol{\nu}}\geq \mu_i\,U_{i,\boldsymbol{\nu}}(p)=\mu_i\lambda_i\,a(p)
\end{equation}
on $\Gamma_i$, where $K_i(p)$ denotes the counter-clockwise-oriented curvature of $\Gamma_i$ at $p\in\Gamma_i$. In fact we have $\big|\int_{\gamma_i} K_i(p)\,|dt|-\pi\big|\leq B_0\sqrt{r_i}$, due to (\ref{eqn 2.7.1-3}) and (\ref{eqn 2.7.1-4}), both in the "$+$" case. Again assuming that $r_i\leq 1$, it follows by integrating (\ref{eqn 2.7.1-6}) on $\gamma_i$, rearranging the terms, and applying the divergence theorem to the function $\psi(p):={\rm ln}\big(a(p)\big)$ in the bounded connected region $\omega_i$ such that $\partial\omega_i=\gamma_i\cup L_i$ (where $\Delta\psi\geq 0$ in $\omega_i$) that
\begin{equation}
\label{eqn 2.7.1-7}
\pi-\mu_i\lambda_i||\gamma_i||_a-B_0\sqrt{r_i}\leq\int_{\gamma_i} \psi_{\boldsymbol{\nu}}\,|dt|\leq\int_{\partial\omega_i}\psi_{\boldsymbol{\nu}}\,|dt|+B_1\,r_i
\end{equation}
$$=\int\int_{\omega_i}(-\Delta\psi)\,dA+B_1\,r_i\leq\,B_1\,r_i,$$
where $||\gamma_i||_a$ is the weighted arc-length of $\gamma_i$ (with weight-function $a(p)$) and ${\boldsymbol{\nu}}$ denotes the interior normal to $\omega_i$ on $\partial\omega_i$. It follows that
\begin{equation}
\label{eqn 2.7.1-5'}
B_1 r_i+B_0 \sqrt{r_i}\geq\pi-\lambda_i\mu_i\,||\gamma_i||_a\geq
\pi-\mu_i\,{\rm capacity}(\Omega_{\boldsymbol{\lambda}}),
\end{equation}
from which Eq. (\ref{eqn 2.7.1-5}) follows.
\vspace{.1in}

\noindent
Similarly, for $p_{i,0}\in\Gamma_i$ and $p_i^-\in{\cal M}_i^-(p_{i,0})$, $i=1,2$, we will show that \begin{equation}
\label{eqn 2.7.1-8}
K_0\overline{\lambda}\,\,\overline{A}\,r_i\geq\underline{\lambda}\,\,\underline{A}\,\big(\pi-B_0\sqrt{r_i}\,\big),
\end{equation}
$i=1,2$, where $K_0$ denotes a bound on the absolute curvature of $\Gamma_1$ and $\Gamma_2$. In fact the flow across the segment $L_i$ joining $p_{i,0}$ to $p_i^-$ is bounded from above by $\overline{\lambda}\,\,\overline{A}\,r_i$, since $|\nabla U_i|\leq\overline{\lambda}\,\,\overline{A}$ in $\Omega$. On the other hand, the flow across the arc $\gamma_i\subset \Gamma_i$ joining $p_{i,0}$ to $p_i^-$ exceeds $\underline{\lambda}\,\underline{A}\,||\gamma_i||$, where, in view of the absolute curvature bound $K_0$, it follows from (\ref{eqn 2.7.1-3}) and (\ref{eqn 2.7.1-4}) in the "$-$"-cases, that the length of the arc $\gamma_i$ is at least $\big(\big[\pi-B_0\sqrt{r_i}\,\,\big]\big/K_0\big)$. Eq. (\ref{eqn 2.7.1-7}) now follows from the fact that the flows across $L_i$ and $\gamma_i$ are equal, where the flow across $\gamma_i$ exceeds $\underline{\lambda}\,\,\underline{A}\,||\gamma_i||$. It follows from this that 
\begin{equation}
\label{eqn 2.7.1-8a}
{\rm if}\,\,r_i\leq 1,\,\,\,{\rm then}\,\,\,r_i\geq R_1:=\big(\pi\,\underline{\lambda}\,\,\underline{A}\,/(K_0\,\overline{\lambda}\,\,\overline{A}+C_0\underline{\lambda}\,\,\underline{A})\big)^2.
\end{equation}

\noindent
By (\ref{eqn 2.7.1-5}) and (\ref{eqn 2.7.1-8a}), we conclude that either $r_i>1$ or else $r_i\geq(\rho_0/B_2)^2$ and $r_i\geq R_1$ (both). Finally, in view of Def. \ref{def 2.7.2} (see (\ref{eqn 2.7.1-1})), the assertion (\ref{eqn 2.7.1-2}) follows by substituting these inequalities into (\ref{eqn ZW}).
\vspace{.1in}

\noindent
{\bf Proof of Part (b).} In view of Thms. \ref{thm 2.6.2}, \ref{thm 2.7.1} and \ref{thm 2.7.2}, this follows from Part (a). 

\begin{definition} 
{(Notation)}
\label{def 4D}
Given a pair $\hat{{\bm\lambda}}\in{\bm\Lambda}_+$ and a classical solution $\hat{{\bf\Gamma}}\in{\bf X}\cap C^{\,3,\tilde{\varrho}}$ of Prob. \ref{prob 4I.+} at ${\bm\lambda}$, let $\hat{U}(p):=U(\hat{{\bf\Gamma}};\,p): {\rm Cl}(\Omega)\rightarrow\Re$ denote the capacitary potential in the domain $\Omega:=\Omega({\bf\Gamma})$. Assume that the mapping $\hat{U}(p): {\rm Cl}(\hat{\Omega})\rightarrow\Re$ has a real-valued $C^{\,3,\tilde{\varrho}}$-continuation $\hat{U}(p)$ to a neighborhood $\Omega$ of ${\rm Cl}(\hat{\Omega})$, and that there exists an interval $A=(-\delta_0^*, 1+\delta_0^*)$ (where $\delta_0^*>0$) such that the family of sets $\Gamma_\alpha:=\{p\in \Omega: \hat{U}(p)=\alpha\}$, $\alpha\in A$, constitutes a continuously and monotonically-varying family of double-point free periodic arcs. 
Then for any vector $\boldsymbol{\delta}=(\delta_1,\delta_2)\in\Re^2$ such that $\delta:=|\boldsymbol{\delta}|<\delta_0^*$, we define the curve-pair 
$${\bf M}_{\boldsymbol{\delta}}({\bf\Gamma}):=\big(M_{1,\,\delta_1}({\bf\Gamma}),M_{2,\,\delta_2}({\bf\Gamma})\big)\in{\bf X}(G)$$
\noindent
component-wise such that 
\begin{equation}																			\label{eqn 4D.1}
M_{i,\,\delta_i}({\bf\Gamma}):=\big\{p\in\Omega: \hat{U}(p)=i+\delta_i-1\big\}\in X,
\end{equation}
for each $i=1,2$.  
\end{definition}
\begin{lemma} 
{(A family of solutions of Prob. \ref{prob 4I.} in a neighborhood of a given solution of Prob. \ref{prob 4I.})}
\label{lem 4E.1}
In Prob. \ref{prob 4I.}, let $\hat{{\bf\Gamma}}:=(\hat{\Gamma}_1,\hat{\Gamma}_2)\in{\bf X}\cap C^{\,3,\tilde{\varrho}}$ denote a fixed classical solution
at $\hat{{\bm\lambda}}:=(\hat{\lambda}_1,\hat{\lambda}_2)\in\Re_+^2$ such that:
\begin{equation}
\label{eqn 4E.1}
\Delta\,\,{\rm ln}\big(a(p)\big)\geq 0\,\,{\rm in}\,\,{\rm Cl}(\hat{\Omega}), 
\end{equation}
where $\hat{\Omega}:=\Omega(\hat{{\bm\Gamma}})$. 
Assume that the capacitary potential $\hat{U}(p):=U(\hat{\bf\Gamma};p)$, originally defined in the closure of the domain $\hat{\Omega}:=\Omega(\hat{{\bf\Gamma}})$, can be continued as a real-valued $C^{\,3,\tilde{\varrho}}$-function (still called $\hat{U}(p)$) defined in some neighborhood $\Omega$ of ${\rm Cl}(\hat{\Omega})$ (conditions guaranteeing this continuation are given in Thm. \ref{thm 2.7.2} and Lem. \ref{lem 4B.2}(a)(b)). Finally, assume that $\hat{{\bf\Gamma}}\in{\boldsymbol{\cal S}}
(\hat{{\bm\lambda}};r_0)$ for some value $r_0>0$ (see Def. \ref{def 4A}). Then for any given $r\in(0,r_0)$, there exist values $\delta_1^*\in(0,\delta_0^*]$ and $C_0>0$, and a unit vector
\begin{equation}
\label{eqn 4E.2}
{\bm v}:=(v_1,v_2)\in(C_0,\infty)^2,
\end{equation}
such that for any $\delta\in(0,\delta_1^*]$ and any pair ${\bm\lambda}=(\lambda_1,\lambda_2)\in\Re_+^2$ satisfying 
\begin{equation}
\label{eqn 4E.3}
|\lambda_i-\hat{\lambda}_i|\,<\,r\,\hat{\lambda}_i\,v_i\,\delta\,\,\,\,{\rm for}\,\,\,i=1,2,
\end{equation}
where the double-arc-pair $\hat{{\bf\Gamma}}_{{\bm v}\delta}:={\bf M}_{{\bm v}\delta}(\hat{{\bf\Gamma}})$ (resp. $\hat{{\bf\Gamma}}_{-{\bm v}\delta}:={\bf M}_{-{\bm v}\delta}(\hat{{\bf\Gamma}})$) (see Def.
\ref{def 4D}) is a strict upper (resp. lower) classical solution of Prob. \ref{prob 4I.} at ${\bm\lambda}$ 
such that $\hat{{\bm\Gamma}}_{{\bm v}\delta}>\hat{{\bm\Gamma}}$ (resp. $\hat{{\bm\Gamma}}>\hat{{\bm\Gamma}}_{-{\bm v}\delta}$). Therefore, for any ${\bm\lambda}\in\Re^2_+$ satisfying (\ref{eqn 4E.3}), there exists (by Thms. \ref{thm 2.6.2} and \ref{thm 2.7.1}) a classical solution ${\bf\Gamma}({\bm\lambda})\in{\bf X}\cap C^{\,3,\tilde{\varrho}}$ of Prob. \ref{prob 4I.} at ${\bm\lambda}$ such that
\begin{equation}
\label{eqn 4E.4}
\hat{\bf\Gamma}_{-{\bm v}\delta}<{\bf\Gamma}({\bm\lambda})<\hat{{\bf\Gamma}}_{{\bm v}\delta}.
\end{equation}
\end{lemma}
\noindent
{\bf Proof.} 
Given the $C^{\,3,\tilde{\varrho}}$-solution $\hat{{\bf\Gamma}}\in{\boldsymbol{\cal S}}
(\hat{{\bm\lambda}};r_0)$ of Prob. \ref{prob 4I.}, for any vector $\boldsymbol{\delta}=(\delta_1,\delta_2)\in\Re^2$ such that $\delta:=|\boldsymbol{\delta}|=\sqrt{\delta_1^2+\delta_2^2}$ is sufficiently small, we define the capacitary potential $\hat{U}_{\boldsymbol{\delta}}(p):=\,U(\hat{{\bf\Gamma}}_{\boldsymbol{\delta}};\,p)$ in the closure of the periodic strip-like domain $\hat{\Omega}_{\boldsymbol{\delta}}:=\,\Omega(\hat{{\bf\Gamma}}_{\boldsymbol{\delta}}),$ where $\hat{{\bf\Gamma}}_{\boldsymbol{\delta}}:=\,(\hat{\Gamma}_{1,\,\delta_1},\hat{\Gamma}_{2,\,\delta_2}):={\bf M}_{\boldsymbol{\delta}}(\hat{{\bf\Gamma}})\in{\bf X}$ (see Def. \ref{def 4D}). 
We begin, in this context, by developing estimates for the functions $\hat{F}_{i,\,\boldsymbol{\delta}}(p):\hat{\Gamma}_{i,\,\delta_i}\rightarrow\Re$, $i=1,2$, defined such that
\begin{equation}
\label{eqn 4E.5}
\hat{F}_{i,\boldsymbol{\delta}}(p):=\,|\nabla{\hat{U}}_{\boldsymbol{\delta}}(p)|-\hat{\lambda}_i\,a(p)\,=\,\hat{U}_{\boldsymbol{\delta},\boldsymbol{\nu}}(p)-\hat{\lambda}_i\,a(p),
\end{equation}
(where $\boldsymbol{\nu}\,=\,\boldsymbol{\nu}(p)$ refers to the unit normal vector to $\partial\,\hat{\Omega}_{\boldsymbol{\delta}}$ at any point 
$p\in\partial\,\hat{\Omega}_{\boldsymbol{\delta}}$ and in the direction $\nabla \hat{U}_{\boldsymbol{\delta}}(p)$). For any sufficiently small ${\bm\delta}=(\delta_1,\delta_2)\in\Re^2$, the arc-pair $\hat{{\bm\Gamma}}_{\boldsymbol{\delta}}:=\,(\hat{\Gamma}_{1,\delta_1},\,\hat{\Gamma}_{2,\,\delta_2})\in{\bf X}$ is an upper (resp. lower) classical solution of Prob. \ref{prob 4I.} at $\hat{\bm\lambda}$ (see Defs. \ref{def 2.1.2} and \ref{def 2.1.3}) provided that, in the notation of (\ref{eqn 4E.5}), we have
\begin{equation}
\label{eqn 4E.25}
(-1)^i\hat{F}_{i,\boldsymbol{\delta}}(p)\,<\,(\,>\,)\,\,0\,\,{\rm on}\,\,\hat{\Gamma}_{i,\delta_i}\
\end{equation}
for $i=1,2$. More generally, given a vector ${\bm\lambda}=(\lambda_1,\lambda_2)\in\Re_+^2$ such that
\begin{equation}
\label{eqn 4E.4+}
|\lambda_i-\hat{\lambda}_i|\leq r\,\hat{\lambda}_i\,|\delta_i|
\end{equation}
for some value $r\in(0,\,r_0]$, and for $i=1,2,$ the same arc-pair $\hat{{\bm\Gamma}}_{\boldsymbol{\delta}}:={\bf M}_{\boldsymbol{\delta}}(\hat{{\bm\Gamma}})$ will also be an upper (resp. lower) solution of Prob. \ref{prob 4I.} at ${\bm\lambda}$ provided that
\begin{equation}
\label{eqn 4E.25a}
(-1)^i\hat{F}_{i,\boldsymbol{\delta}}(p)\,<\,-r\,\hat{\lambda}_i|\delta_i|a(p)\,\,\,\big(\,>\,r\,\hat{\lambda}_i|\delta_i|a(p)\,\,\big)
\,\,{\rm on}\,\,\hat{\Gamma}_{i,\delta_i}
\end{equation}
for sufficiently small ${\bm\delta}\in\Re^2$ and for $i=1,2$ (see (\ref{eqn 4E.17a})). We also define the related functions
\begin{equation}
\label{eqn 4E.6}
\hat{F}_i(p)\,:=\,|\nabla\hat{U}(p)|-\hat{\lambda}_i\,a(p)\,=\,\hat{U}_{\boldsymbol{\nu}}(p)-\hat{\lambda}_i\,a(p),
\end{equation}
$i=1,2$, in neighborhoods of the corresponding arcs $\hat{\Gamma}_i$ (where $\boldsymbol{\nu}=\boldsymbol{\nu}(p)$ again denotes the unit normal vector to $\partial\,\hat{\Omega}$ at the point $p\in\partial\,\hat{\Omega}$ in the direction of the vector $\nabla\,\hat{U}(p)$). Observe that for each $i=1,2$ and for all points $p\in\hat{\Gamma}_i$, we have that $\hat{F}_i(p)=0$, and
it follows from (\ref{eqn 4E.6}) by differentiation that
\begin{equation}
\label{eqn 4E.7}
\hat{F}_{i,\,\boldsymbol{\nu}}(p)\,=\,\big(\hat{U}_{\boldsymbol{\nu}\boldsymbol{\nu}}(p)-\hat{\lambda}_i\,a_{\boldsymbol{\nu}}(p)\big)
\end{equation}
$$=\,\hat{\lambda}_i\,a(p)\,\big((\hat{U}_{\boldsymbol{\nu}\boldsymbol{\nu}}(p)\big/\hat{U}_{\boldsymbol{\nu}}(p))\,-\,(a_{\boldsymbol{\nu}}(p)\big/a(p))\big)$$
$$=\,\hat{\lambda}_i\,a(p)\big(\hat{K}_i(p)\,-\,(a_{\boldsymbol \nu}(p)\big/a(p))\big)$$
$$=\,\hat{\lambda}_i\,a(p)\,\hat{\phi}_{i,\,\boldsymbol{\nu}}(p)\,=\,\hat{\lambda}_i\,a(p)\,\mathfrak{C}_i(p),$$
where, in terms of the constant $\hat{\mu}:=\,{\rm ln}\big(\hat{\lambda}_2\big/\hat{\lambda}_1\big)$, the counter-clockwise-oriented curvature $\hat{K}_i(p)$ of $\hat{\Gamma}_i$ at $p\in\hat{\Gamma}_i$, the weakly-subharmonic function $\hat{W}(p):\,{\rm Cl}\,(\hat{\Omega})\rightarrow\Re$, and the boundary function $\hat{E}(p):\,\partial\,\hat{\Omega}\rightarrow\Re_+$, we define the continuous and weakly-superharmonic functions $\hat{\phi}_i(p):\,{\rm Cl}\,(\hat{\Omega})\rightarrow\Re$, $i=1,2$, such that
\begin{equation}
\label{eqn 4E.8}
\hat{\phi}_i(p)\,:=\,{\rm ln}\,\big(|\nabla \hat{U}(p)|\,\big/\,\hat{\lambda}_i\,a(p)\big)\,=\,(-1)^{i+1}\,\hat{\mu}\,\hat{U}_i(p)-\hat{W}(p),
\end{equation}
and the corresponding boundary-derivative functions $\mathfrak{C}_i(p):\,\hat{\Gamma}_i\rightarrow\Re$, $i=1,2$, such that
\begin{equation}
\label{eqn 4E.9}
\mathfrak{C}_i(p)\,:=\,\hat{\phi}_{i,\,\boldsymbol{\nu}}(p)\,=\,\hat{\mu}\,\hat{U}_{\boldsymbol{\nu}}(p)-\hat{W}_{\boldsymbol{\nu}}(p)
\,=\,\hat{\lambda}_i\,a(p)\big(\hat{\mu}\,-\,(-1)^i\,\hat{E}(p)\big).
\end{equation}
In the above context, and for small ${\bm\delta}=(\delta_1,\delta_2)\in\Re^2,$ we will prove for $i=1,2$ that 
\begin{equation}
\label{eqn 4E.10}
\big|\hat{F}_i(\hat{p}_{i,\,\delta_i})\,-\,\mathfrak{C}_i(\hat{p}_i)\,\delta_i\,\big|\,\leq\,\delta\,z(\delta),
\end{equation}
\begin{equation}
\label{eqn 4E.11}
\big|\hat{U}_{\boldsymbol{\delta},\,\boldsymbol{\nu}}(\hat{p}_{i,\,\delta_i})-\hat{U}_{\boldsymbol{\nu}}(\hat{p}_{i,\,\delta_i})-\hat{\lambda}_i\,a(\hat{p}_{i,\,\delta_i})(\delta_1-\delta_2)\big|\,\leq\,\delta\,z(\delta),
\end{equation}
\begin{equation}
\label{eqn 4E.12}
\big|\hat{F}_{i,\,\boldsymbol{\delta}}(\hat{p}_{i,\,\delta_i})\,-\,\hat{\lambda}_i\,a(\hat{p}_i)\,\big([\,\hat{\mu}\,-\,(-1)^i\,\hat{E}(\hat{p}_i)]\,\delta_i\,+\,[\delta_1\,-\,\delta_2]\,\big)\,\big|\,\leq\,\delta\,z(\delta),
\end{equation}
where $\delta=|{\bm\delta}|$ and $z(\delta):\Re_+\rightarrow\Re_+$ denotes some positive function (distinct in each application) such that $z(\delta)\rightarrow 0+$ as $\delta\rightarrow 0+$. Here, for any points $\hat{p}_i\in\hat{\Gamma}_i$, $i=1,2$, and any vector ${\bm\delta}=\,(\delta_1,\,\delta_2)$ such that $\delta=\,|{\bm \delta}|$ is sufficiently small, we set
\begin{equation}
\label{eqn 4E.13}
\hat{p}_{i,\,\delta_i}:=\,\hat{p}_i+h_i\,{\bm\nu}\in\hat{\Gamma}_{i,\,\delta_i}:=M_{i,\,\delta_i}(\hat{{\bf\Gamma}}),
\end{equation}
where the magnitude of $h_i=\,h_i(\delta_i)>0$ is minimum subject to (\ref{eqn 4E.13}). Toward the proof of the estimate (\ref{eqn 4E.10}), we observe that 
\begin{equation}
\label{eqn 4E.14}
\delta_i:=\hat{U}(\hat{p}_{i,\,\delta_i})-\hat{U}(\hat{p}_i)=\hat{U}_{\boldsymbol{\nu}}(\hat{p}_i)\,h_i\,+\,R_i\,h_i,
\end{equation}
\begin{equation}
\label{eqn 4E.15} 
\hat{F}_i(\hat{p}_{i,\,\delta_i})-\hat{F}_i(\hat{p}_i)=\hat{F}_{i,\,\boldsymbol{\nu}}(\hat{p}_i)\,h_i+\,R_i\,h_i
\end{equation}
for $i=1,2$, from which it follows that 
\begin{equation}
\label{eqn 4E.16}
\hat{F}_i(\hat{p}_{i,\,\delta_i})-\hat{F}_i(\hat{p}_i)\,=\,
\big(\hat{F}_{i,\,\boldsymbol{\nu}}(\hat{p}_i)\big/\hat{U}_{\boldsymbol{\nu}}(\hat{p}_i)\big)\,\delta_i+R_i\,\delta_i.
\end{equation}
Here (in (\ref{eqn 4E.14}), (\ref{eqn 4E.15}), and (\ref{eqn 4E.16})), we have $\hat{F}_i(\hat{\Gamma}_i)\,=\,0$, 
and we use $R_i$ to denote the coefficient of $h_i$ in the second-order Taylor remainder. Notice that $R_i$ can be estimated in terms of $\hat{\boldsymbol{\lambda}}$ and the secondd derivatives, respectively, of the functions $a(p)$ and $\hat{U}(p)$, which are both uniformly bounded by the $C^{\,3,\tilde{\varrho}}$-norms of the $C^{\,3,\tilde{\varrho}}$-continuation $\hat{U}(p):\Omega\rightarrow \Re$. 
This completes the proof of the estimate (\ref{eqn 4E.10}), which now follows fromw turn to the proof of the estimate (\ref{eqn 4E.11}), which is based on the following simple identity: We have
\begin{equation}
\label{eqn 4E.18}
\big(1+\delta_2-\delta_1\big)\,\hat{U}_{\boldsymbol{\delta}}(p)\,=\,\hat{U}(p)-\delta_1
\end{equation}
for all points $p\in\partial\,\hat{\Omega}_{\boldsymbol{\delta}}$, provided that the value $\delta:=\,|\boldsymbol{\delta}|$ is small enough to guarantee that ${\rm Cl}\big(\hat{\Omega}_{\boldsymbol{\delta}}\big)\subset\Omega$. Also, it follows from Thm. \ref{thm 2.7.2}(a) that if $\Delta\,\hat{U}(p)\,=\,0$ in ${\rm Cl}\big(\hat{\Omega}\big)$, and therefore that $\big|\Delta\,\hat{U}(p)\big|\leq\delta z(\delta)$ in $\hat{\Omega}_{\boldsymbol{\delta}}\setminus\hat{\Omega}$, where $z(\delta):\Re_+\rightarrow\Re_+$ denotes a positive null-function depending on the $C^{\,3,\,\tilde{\varrho}}$-norm of the $C^{\,3,\,\tilde{\varrho}}$-continuation $\hat{U}(p):\Omega\rightarrow\Re$.
Therefore, the function $\hat{\psi}_{\boldsymbol{\delta}}(p):=\,\hat{U}(p)-\hat{U}_{\boldsymbol{\delta}}(p):{\rm Cl}\big(\hat{\Omega}_{\boldsymbol{\delta}}\big)\rightarrow\Re$ must be such that
\begin{equation}
\label{eqn 4E.19}
\hat{\psi}_{\boldsymbol{\delta}}\big(\partial\hat{\Omega}_{\boldsymbol{\delta}}\big)\,=\,0,\,\, \Delta\,\hat{\psi}_{\boldsymbol{\delta}}\big(\hat{\Omega}\big)=0,\,\,{\rm and}\,\,\,|\Delta \,\hat{\psi}_{\boldsymbol{\delta}}(p)|\leq\delta\,z(\delta)\,\,{\rm for}\,\,{\rm all}\,\,p\in \hat{\Omega}_{\boldsymbol{\delta}},
\end{equation}
\noindent
and therefore such that 
\begin{equation}
\label{eqn 4E.20}
\big|\hat{\psi}_{\boldsymbol{\delta}}(p)\big|\leq \delta\,z(\delta)
\end{equation}
for any point $p\in\hat{\Omega}_{\boldsymbol{\delta}}$. For any point $p_0\in\Re^2$, let $u(p_0;\,p):\omega(p_0)\rightarrow\Re$ denote a harmonic "barrier" function such that $u(p_0;p)\,:=\,\big(\delta\,z(\delta)\,{\rm ln}\big(|p-p_0|\big/R_0\big)\big/{\rm ln}\big(R_1/R_0)\big),$ in the closed annular domain $\omega(p_0)\,:=\,\big\{p\in\Re^2:\,R_0\leq|p-p_0|\leq R_1\big\}$. Then, we have $u(p_0;p)\,=\,0$ (resp. $u(p_0;p)\,=\,\delta\,z(\delta))$ on the interior (resp. exterior) boundary component of $\omega(p_0)$. Also, assuming that in the definition of $\omega(p_0)$ we have $0<R_0\,:=\,{\rm dist}\big(p_0,\,\hat{\Omega}_{\boldsymbol{\delta}}\big)<R_1$ and $R_1<R_0+{\rm dist}\big(\hat{\Gamma}_1,\hat{\Gamma}_2\big)$, it follows by the comparison principle that $|\hat{\psi}_\delta(p)|\leq u(p_0;\,p)$, first for all $p\in\partial\,(\hat{\Omega}_{\boldsymbol{\delta}}\cap \omega(p_0))$, then for all $p\in\hat{\Omega}_{\boldsymbol{\delta}}\cap\omega(p_0)$, from which one finally concludes that 
\begin{equation}
\label{eqn 4E.21}
|\nabla\hat{\psi}_{\delta}(p)|\leq|\nabla u(p_0;p)|\,=\,\big(\delta\,z(\delta)/R_0{\rm ln}(R_1/R_0)\big).
\end{equation}
for all points $p\in\partial\hat{\Omega}_{\boldsymbol{\delta}}$ such that $p\in\partial\omega(p_0)$ for some ball $\omega(p_0)$ such that $\hat{\Omega}_{\boldsymbol{\delta}}\cap\omega(p_0)=\emptyset$. Also, in view of (\ref{eqn 4E.18}), it follows from (\ref{eqn 4E.20}) that
\begin{equation}
\label{eqn 4E.22}
\big|(1+\delta_2-\delta_1)\,\hat{U}_{\boldsymbol{\delta}}(p)\,-\,(\hat{U}(p)-\delta_1)\big|\leq\delta\, z(\delta)
\end{equation} 
for all points $p\in{\rm Cl}\big(\hat{\Omega}_{\boldsymbol{\delta}}\big),$
and, in view of (\ref{eqn 4E.21}), it follows from (\ref{eqn 4E.22}) that 
\begin{equation}
\label{eqn 4E.23}
\big|(1+\delta_2-\delta_1)\nabla\,\hat{U}_{\boldsymbol{\delta}}(p)\,-\,\nabla\,\hat{U}(p)\big|\leq \delta\,z(\delta)
\end{equation}
for sufficiently small $\delta>0$ and all points $p\in\partial\,\hat{\Omega}_{\boldsymbol{\delta}}$, or, alternatively, such that for some other null-function $z(\delta):\Re_+\rightarrow\Re_+$, we have 
\begin{equation}
\label{eqn 4E.24}
\big|\nabla\,\hat{U}_{\boldsymbol{\delta}}(p)\,-(1+\delta_1-\delta_2)\,\nabla\,\hat{U}(p)\big|\leq \delta\,z(\delta),
\end{equation}
for any sufficiently small $\delta>0$ and for all $p\in\partial\,\hat{\Omega}_{\boldsymbol{\delta}}$, from which the estimate (\ref{eqn 4E.11}) follows. Finally, the estimate (\ref{eqn 4E.12}) follows directly from (\ref{eqn 4E.10}) and (\ref{eqn 4E.11}), in view of (\ref{eqn 4E.9}) and the following simple identity:
\begin{equation}
\label{eqn 4E.17}
\hat{F}_{i,\,\boldsymbol{\delta}}(p)\,=\,
\hat{F}_i(p)\,+\,\big(\hat{F}_{i,\,\boldsymbol{\delta}}(p))-\hat{F}_i(p)\big)
\end{equation}
$$=\,\hat{F}_i(p)\,+\,\big(\hat{U}_{\boldsymbol{\delta},\,\boldsymbol{\nu}}(p)
-\hat{U}_{\boldsymbol{\nu}}(p)\big),$$
for $i=1,2$ and for all points $p\in \hat{\Gamma}_{i,\,\delta_i}:=\,M_{i,\,\delta_i}(\hat{{\bm\Gamma}})$. 
At this point, a direct comparison of the Eqs. (\ref{eqn 4E.25a}) and (\ref{eqn 4E.12}) shows that for sufficiently small ${\bm\delta}\in\Re^2$, the arc-pair $\hat{{\bf\Gamma}}_{\boldsymbol{\delta}}=(\hat{\Gamma}_{1,\delta_1},
\hat{\Gamma}_{2,\delta_2})$ is an upper (resp. lower) classical solution of Prob. \ref{prob 4I.} at any
 vector ${\bm\lambda}=(\lambda_1,\lambda_2)\in\Re_+^2$ satisfying (\ref{eqn 4E.4+}) provided that
\begin{equation}
\label{eqn 4E.17a}
\big(\,[\,(-1)^i\,\hat{\mu}-\,\hat{E}(\hat{p}_i)\,]\,\delta_i\,+\,(-1)^i\,(\delta_1-\delta_2)\big)<-r|\delta_i|-\delta\,z(\delta)\,\,\big(\,>r|\delta_i|+\delta\,z(\delta)\big),
\end{equation}
both for $i=1,2$, where $z(\delta):[0,\infty)\rightarrow[0,\infty)$ denotes some continuous, strictly-increasing function such that $z(0)=0$. Let $\delta>0$ and ${\bm v}=(v_1,v_2)$ denote respectively any (positive) scalar and any unit vector in the first quadrant. It is easy to see that if ${\bm\delta}=\delta{\bm v}$ (resp. ${\bm\delta}=-\delta{\bm v}$), then $\hat{{\bf\Gamma}}_{\boldsymbol{\delta}}\,>\,\hat{{\bf\Gamma}}$ (resp. $\hat{{\bf\Gamma}}_{\boldsymbol{\delta}}\,<\,\hat{{\bf\Gamma}}$). Moreover, it follows by studying the four possible cases of (\ref{eqn 4E.17a}) that if a positive scaler $\delta>0$ and a corresponding unit vector ${\bm v}=(v_1,v_2)\in\Re_+^2$ together satisfy both of the following conditions:
\begin{equation}
\label{eqn 4E.31}
v_1\,+\,z(\delta)\,<\,\big(1+\hat{E}-\hat{\mu}-r\big)\,v_2;\,\,\,v_2\,+\,z(\delta)\,<\,\big(1+\hat{E}+\hat{\mu}-r\big)\,v_1,
\end{equation}
where $\hat{E}=\,E(\hat{\bm\Gamma})=\,\min\big\{\big(|\nabla\hat{W}(p)|\big/|\nabla\hat{U}(p)|\big):p\in\hat{{\bm\Gamma}}\big\}$ and $z(\delta)$ denotes the same fixed positive function as in (\ref{eqn 4E.17}), then the arc-pair $\hat{{\bf\Gamma}}_{\boldsymbol{\delta}}:={\bf M}_{\boldsymbol{\delta}}(\hat{{\bf\Gamma}})$
is an upper (resp. lower) classical solution of Prob. \ref{prob 4I.} at any vector ${\bm\lambda}\in\Re_+^2$ satisfying (\ref{eqn 4E.4+}) provided that ${\bm\delta}=\delta{\bm v}$ (resp. ${\bm\delta}=-\delta{\bm v}$).
\vspace{.1in}

\noindent
To discuss the existence of vectors ${\bm v}$ satisfying (\ref{eqn 4E.31}), we begin with the related problem of finding all the unit vectors ${\bm v}=(v_1,v_2)\in\Re_+^2$ such that
\begin{equation}
\label{eqn 4E.32}
\big(1\big/(1+\hat{E}-\hat{\mu}-r_0)\big)\,<\,(v_2/v_1)\,<\,\big(1+\hat{E}+\hat{\mu}-r_0\big).
\end{equation}
Clearly (\ref{eqn 4E.32}) has infinitely-many solutions, since
\begin{equation}
\label{eqn 4E.33} 
\big(1+\hat{E}-\hat{\mu}-r_0\big)\big(1+\hat{E}+\hat{\mu}-r_0\big)=\big(1+\hat{E}-r_0)^2-\hat{\mu}^2>1
\end{equation}
by assumption (see (\ref{eqn 4A.2})). Furthermore, any pair ${\bm v}\in\Re_+^2$ solves (\ref{eqn 4E.32}) if and only if it also satisfies both inequalities in (\ref{eqn 4E.31}) in the special case where $r:=r_0$ and $\delta=0$ (and $z(0)=0$). In view of this, one sees that any solution ${\bm v}$ of (\ref{eqn 4E.32}) also solves both inequalities in (\ref{eqn 4E.31}) provided that $\delta\in(0,\delta_1^*]$, where $\delta_1^*\in(0,\delta_0^*]$ is small enough so that $z(\delta_1^*)<\,C_0\,(r_0-r)$, where $C_0:=\min\{M_1, M_2\}$, and $M_i$, $i=1,2,$ denotes the greatest lower bound of $v_i$ among all unit vectors ${\bm v}=(v_1,v_2)$ satisfying (\ref{eqn 4E.32}). One can also choose $C_0=\min\big\{(1\big/\sqrt{A_1^2+1}\big),\big(1\big/\sqrt{A_2^2+1}\big)\big\}$ in (\ref{eqn 4E.2}), where $A_i=(1+\hat{E}+(-1)^i\hat{\mu}-r_0)\leq M_i$ for $i=1,2$. Finally, for the purpose of choosing one single solution ${\bm v}$ of (\ref{eqn 4E.32}) which depends continuously on ${\bm\lambda}$, (perhaps named the "optimal" solution), it is natural to define: 
\begin{equation}
\label{eqn 4E.34}
{\bm v}=({\rm cos}(\theta),{\rm sin}(\theta)),\,\,{\rm where}\,\,\,\theta\,=\,(\theta_1+\theta_2)/2,\,\,\, \theta_1\,=\,{\rm arctan}(1/A_1),
\end{equation}
$${\rm and}\,\,\,\theta_2\,=\,{\rm arctan}(A_2).$$
\begin{lemma} {(Strict monotonicity and Lipschitz continuity of parametrized solutions of Prob. \ref{prob 4I.} in a neighborhood of a given solution)}
\label{lem 4F}
Let $\hat{{\bf\Gamma}}=(\hat{\Gamma}_1,\hat{\Gamma}_2)\in{\bf X}$ denote a fixed solution of Prob. \ref{prob 4I.} at $\hat{{\bm\lambda}}$ $=(\hat{\lambda}_1,\hat{\lambda}_2)\in\Re_+^2$, where $a(p):\Re^2\rightarrow\Re_+$ and $\hat{\bf \Gamma}$ have the additional properties assumed in Lemms. \ref{lem 4B.2} and \ref{lem 4E.1} (in particular, (\ref{eqn 4E.2}) holds, and $\hat{{\bf\Gamma}}\in{\boldsymbol{\cal S}}(\hat{{\bm\lambda}};r_0)$ for some $r_0>0$). Let the unit vector ${\bm v}$ satisfy the condition (\ref{eqn 4E.32}) or (\ref{eqn 4E.34}) for some fixed $r\in(0,r_0)$ and for all $\delta\in(0,\delta^*_1]$ (where $\delta^*_1\in(0,\delta_0^*]$ is chosen such that $|R_i(\delta)|\leq C_0(r_0-r)$ for $\delta\in(0,\delta^*_1]$). Then: 
\vspace{.1in}

\noindent
(a) For any classical solution ${\bf\Gamma}={\bf\Gamma}({\bm\lambda})$ of Prob. \ref{prob 4I.} at ${\bm\lambda}\in\Re^2_+$ such that
\begin{equation}
\label{eqn 4F.1}
{\bf M}_{-\delta_0^*{\bm v}}(\hat{{\bf\Gamma}})<{\bf\Gamma}<{\bf M}_{\delta_0^*{\bm v}}(\hat{{\bf\Gamma}})
\end{equation}
(in terms of Def. \ref{def 4D}), there exist values $\alpha=\alpha({\bm\lambda})$,  
$\beta=\beta({\bm\lambda})\in I_{\delta_0^*}:=(-\delta_0^*,\,\delta_0^*)$ such that $0\leq\alpha\leq\beta$ and 
\begin{equation}
\label{eqn 4F.2}
\hat{{\bm \Gamma}}\,\leq\,\hat{{\bf\Gamma}}_{\alpha{\bm v}}:={\bf M}_{\alpha{\bm v}}(\hat{{\bf\Gamma}})\leq{\bf\Gamma}\leq\hat{{\bf\Gamma}}_{\beta{\bm v}}:={\bf M}_{\beta{\bm v}}(\hat{{\bf\Gamma}}),
\end{equation}
\noindent
and where $\alpha$ (resp. $\beta$) is maximum (resp. minimum) subject to (\ref{eqn 4F.2}).
\vspace{.1in}

\noindent
(b) There exist uniform positive constants $\mathfrak{A},\mathfrak{B}>0$ and a value $\delta_2^*=\,\delta_2^*(\hat{\bm{\lambda}})\in(0,\delta^*_1]$ such that for any classical solution ${\bf\Gamma}={\bf\Gamma}({\bm\lambda}):=\,(\Gamma_1,\Gamma_2)\in{\bf X}\cap C^{2}$ of Prob. \ref{prob 4I.} at a pair ${\bm\lambda}=\,(\lambda_1,\lambda_2)\in\Lambda_+^2$ (such that $\lambda_1>\hat{\lambda}_1,\,\,\lambda_2<\hat{\lambda}_2$), and for any values $\alpha,\beta\in I_{\delta_2^*}:=(-\delta_2^*,\delta_2^*)$ such that the condition (\ref{eqn 4F.2}) holds, we have:
\begin{equation}
\label{eqn 4F.3}
\mathfrak{A}\,\,{\rm min}\Big\{\frac{(\lambda_1-\hat{\lambda}_1)}{\hat{\lambda}_1},\frac{(\hat{\lambda}_2-\lambda_2)}{\hat{\lambda}_2}\Big\}\leq\alpha\leq\beta\leq \mathfrak{B}\,\,{\rm max}\Big\{\frac{(\lambda_1-\hat{\lambda}_1)}{\hat{\lambda}_1}, \frac{(\hat{\lambda}_2-\lambda_2)}{\hat{\lambda}_2}\Big\}.
\end{equation}
In fact one can choose $\mathfrak{A}$ such that $\big(2\,\overline{A}\,\,\overline{E}\,\,\big)\,\mathfrak{A}=\underline{A}$.
\end{lemma}
\noindent
{\bf Proof.} Assuming  that (\ref{eqn 4F.1}) holds for a classical solution ${\bf\Gamma}:={\bf\Gamma}({\bm\lambda})$ of Prob. \ref{prob 4I.+} at ${\bm\lambda}$, it follows from (\ref{eqn 4F.2}) by the comparison principle that  
\begin{equation}
\label{eqn 4F.5}
\hat{U}_{\alpha{\bm v}}(p)\geq U(p)\,\,{\rm in}\,\,\hat{\Omega}_{\alpha{\bm v}}\cap\Omega\,\,\,{\rm and}\,\,\,\hat{U}_{\beta{\bm v}}(p)\leq U(p)\,\,{\rm in}\,\, \hat{\Omega}_{\beta{\bm v}}\cap\Omega,
\end{equation} 
where $U(p):=U({\bf\Gamma};p)$ in the closure of $\Omega:=\Omega({\bf\Gamma})$, and where, for either $\kappa=\alpha$ or $\kappa=\beta$, we set $\hat{{\bf\Gamma}}_{\kappa{\bm v}}={\bf M}_{\kappa{\bm v}}(\hat{\bf\Gamma})=(\hat{\Gamma}_{\kappa v_1,1},\hat{\Gamma}_{\kappa v_2,2})$, and $\hat{U}_{\kappa{\bm v}}(p):=U(\hat{\bf \Gamma}_{\kappa{\bm v}};p)$ in the closure of $\hat{\Omega}_{\kappa{\bm v}}:=\Omega(\hat{{\bf\Gamma}}_{\kappa{\bm v}})$. For $\alpha$ maximum subject to (\ref{eqn 4F.2}), there exists a point $\hat{p}_{\alpha,1}\in\Gamma_1\cap\hat{\Gamma}_{\alpha v_1,1}$ such that $\hat{U}_{\alpha{\bm v}}(\hat{p}_{\alpha,1})=U(\hat{p}_{\alpha,1})=0$, and therefore 
\begin{equation}
\label{eqn 4F.6}
|{\nabla}\hat{U}_{\alpha{\bm v}}(\hat{p}_{\alpha,1})|\geq|{\nabla} U(\hat{p}_{\alpha,1})|
\end{equation}
(due to (\ref{eqn 4F.5})), or else there exists a point $\hat{p}_{\alpha,2}\in\Gamma_2\cap\hat{\Gamma}_{\alpha,2}$ at which $\hat{U}_{\alpha{\bm v}}(\hat{p}_{\alpha,2})=U(\hat{p}_{\alpha,2})=1$, from which it follows via (\ref{eqn 4F.5}) that 
\begin{equation}
\label{eqn 4F.7}
|{\nabla}\hat{U}_{\alpha{\bm v}}(\hat{p}_{\alpha,2})|\leq|{\nabla} U(\hat{p}_{\alpha,2})|.
\end{equation}
In the first case, one sees by combining (\ref{eqn 4E.12}) with (\ref{eqn 4F.6}) that
\begin{equation}
\label{eqn 4F.8}
\hat{\lambda}_1\,a(\hat{p}_{\alpha,1})\,+\,\hat{\lambda}_1\,a(\tilde{p}_{\alpha,1})\,{\cal E}_1\alpha\,+\,R_1\,\alpha\,
\end{equation}
$$=|{\nabla}\hat{U}_{\alpha{\bm v}}(\hat{p}_{\alpha,1})|\geq|{\nabla} U(\hat{p}_{\alpha,1})|=\lambda_1\,a(\hat{p}_{\alpha,1}),
$$
at points $\hat{p}_{\alpha,1}\in\hat{\Gamma}_{\alpha,1}\cap\Gamma_1$ and $\tilde{p}_{\alpha,1}\in\hat{\Gamma}_1$ such that $|\tilde{p}_{\alpha,1}-\hat{p}_{\alpha,1}|={\rm dist}(\hat{p}_{\alpha,1},\hat{\Gamma}_1)$, where ${\cal E}_1:=\big((E(\tilde{p}_{\alpha,1})+1+\hat{\mu})v_1-v_2\big) >rv_1>C_0 r>0$ (see (\ref{eqn 4E.31}) and (\ref{eqn 4E.2})). For $\lambda_1>\hat{\lambda}_1$, the assumptions that $\alpha\leq 0$ and $|R_1(\alpha)|\leq\hat{\lambda}_1a(\tilde{p}_{\alpha,1}){\cal E}_{1}$ lead to a contradiction in (\ref{eqn 4F.8}). Therefore, we have $\alpha>0$ if $|R_1(\alpha)|\leq\hat{\lambda}_1a(\tilde{p}_{\alpha,1}){\cal E}_{1}$. But for $\alpha>0$, (\ref{eqn 4F.8}) implies that
\begin{equation}
\label{eqn 4F.9}
\alpha({\bm\lambda})\geq \frac{a(\hat{p}_{\alpha,1})(\lambda_1-\hat{\lambda}_1)}{\hat{\lambda}_1a(\tilde{p}_{\alpha,1})\,{\cal E}_{1}+|R_1(\alpha)|}\geq\frac{\underline{A}(\lambda_1-\hat{\lambda}_1)}{2\,\overline{A}\,\overline{{\cal E}}_1\,\hat{\lambda}_1},
\end{equation}
if we assume that $|R_1(\alpha)|<\hat{\lambda}_1\,\underline{A}C_0r<\hat{\lambda}_1a(\tilde{p}_{\alpha,1})\,{\cal E}_{1}$ and $\lambda_1>\hat{\lambda}_1$, where $\overline{{\cal E}}_1:=\sup\{{\cal E}_1(p):p\in\hat{\Gamma}_1\}$. In the second case, it follows from (\ref{eqn 4E.12}) and (\ref{eqn 4F.7}), that
\begin{equation}
\label{eqn 4F.10}
\hat{\lambda}_2\,a(\hat{p}_{\alpha,2})-\hat{\lambda}_2\,a(\tilde{p}_{\alpha,2})\,{\cal E}_2\,\alpha-R_2\,\alpha\,
\end{equation}
$$=\,|{\nabla}\hat{U}_{\alpha{\bm v}}(\hat{p}_{\alpha,2})|\leq|{\nabla} U(\hat{p}_{\alpha,2})|\,=\,\lambda_2\,a(\hat{p}_{\alpha,2}),$$ 
at points $\hat{p}_{\alpha,2}\in\hat{\Gamma}_{\alpha,2}\cap\Gamma_2$ and $\tilde{p}_{\alpha,2}\in\hat{\Gamma}_2$ such that 
$|\tilde{p}_{\alpha,2}-\hat{p}_{\alpha,2}|={\rm dist}(\hat{p}_{\alpha,2},\hat{\Gamma}_2)$, and ${\cal E}_2=\big((E(\tilde{p}_{\alpha,2})+1-\hat{\mu})v_2-v_1\big)>r v_2> C_0r>0$ (see (\ref{eqn 4E.31}) and (\ref{eqn 4E.2})). For $\hat{\lambda}_2>\lambda_2$, this equation is contradicted by the assumptions that $\alpha\leq 0$ and $|R_1(\alpha)|<\hat{\lambda}_2a(\tilde{p}_{\alpha,2})\,{\cal E}_{2}$. For $\alpha>0$, it follows from (\ref{eqn 4F.10}) that
\begin{equation}
\label{eqn 4F.11}
\alpha({\bm\lambda})\geq\frac{a(\hat{p}_{\alpha,2})(\hat{\lambda}_2-\lambda_2)}{\hat{\lambda}_2 a(\tilde{p}_{\alpha,2})\,{\cal E}_{2}+|R_2(\alpha)|}\geq\frac{\underline{A}(\hat{\lambda}_2-\lambda_2)}{2\,\overline{A}\,\overline{{\cal E}}_2\,\hat{\lambda}_2}
\end{equation}
if we assume that $|R_2(\alpha)|\leq C_0\hat{\lambda}_2\,\underline{A}r\leq\hat{\lambda}_2a(\tilde{p}_{\alpha,2})\,{\cal E}_{2}$ and $\hat{\lambda}_2>\lambda_2$, where $\overline{{\cal E}}_2:=\sup\{{\cal E}_2(p):p\in\hat{\Gamma}_2\}$. Therefore, by (\ref{eqn 4F.9}) and (\ref{eqn 4F.11}), we have
\begin{equation}
\label{eqn 4F.12a}
\alpha({\bm\lambda})\geq\frac{\underline{A}}{2\overline{A}\,\max\{\overline{{\cal E}}_1,\overline{{\cal E}}_2\}}\min\Big\{\frac{(\hat{\lambda}_2-\lambda_2)}{\hat{\lambda}_2},\frac{(\lambda_1-\hat{\lambda}_1)}{\hat{\lambda}_1}\Big\}
\end{equation}
for $|R_2(\alpha)|\leq C_0\hat{\lambda}_2\,\underline{A}r$, completing the proof of the existence of the positive uniform lower bound for $\alpha$ in (\ref{eqn 4F.3}).
\vspace{.1in}

\noindent
We turn now to the existence of the uniform upper bound for $\beta$ in (\ref{eqn 4F.3}). By (\ref{eqn 4E.12}) and (\ref{eqn 4F.5}), we have that
\begin{equation}
\label{eqn 4F.12}
\hat{\lambda}_1\,a(\hat{p}_{\beta,1})\,+\,\hat{\lambda}_1\,a(\tilde{p}_{\beta,1})\,{\cal E}_1\,\beta\,+\,R_1\,\beta
\end{equation}
$$=|{\nabla}\hat{U}_{\beta{\bf v}}(\hat{p}_{\beta,1})|\leq|{\nabla} U(\hat{p}_{\beta,1})|=\lambda_1a(\hat{p}_{\beta,1})$$
at points $\hat{p}_{\beta,1}\in\Gamma_1\cap\hat{\Gamma}_{\beta,1}$ and $\tilde{p}_{\beta,1}\in\hat{\Gamma}_1$ such that $|\tilde{p}_{\beta,1}-\hat{p}_{\beta,1}|={\rm dist}(\hat{p}_{\beta,1},\hat{\Gamma}_1)$, where ${\cal E}_1:=\big((E(\tilde{p}_{\beta,1})+1+\hat{\mu})v_1-v_2\big)\geq v_1 r\geq C_0 r>0$, or else, alternatively, we have that
\begin{equation}
\label{eqn 4F.13}
\hat{\lambda}_2a(\hat{p}_{\beta,2})-\hat{\lambda}_2a(\tilde{p}_{\beta,2})\,{\cal E}_2\beta+R_2\,\beta\,
\end{equation}
$$
=|{\nabla}\hat{U}_{\beta{\bm v}}(\hat{p}_{\beta,2})|\geq|{\nabla} U(\hat{p}_{\beta,2})|=\lambda_2a(\hat{p}_{\beta,2})$$
at points $\hat{p}_{\beta,2}\in\Gamma_2\cap\hat{\Gamma}_{\beta,2}$ and $\tilde{p}_{\beta,2}\in\hat{\Gamma}_2$ such that $|\tilde{p}_{\beta,2}-\hat{p}_{\beta,2}|={\rm dist}(\hat{p}_{\beta,2},\hat{\Gamma}_2)$, where ${\cal E}_2:=\big((E(\tilde{p}_{\beta,2})+1-\hat{\mu})v_2-v_1\big)\geq v_2 r\geq C_0 r>0$. Of course we have $\beta\geq\alpha>0$ in both cases, by the already established part of (\ref{eqn 4F.3}).) Therefore, if either (\ref{eqn 4F.12}) or (\ref{eqn 4F.13}) holds, and if we assume that $\lambda_1-\hat{\lambda}_1, \hat{\lambda}_2-\lambda_2>0$, and that $2\,|R_i(\beta)|< \min\{\hat{\lambda}_1,\hat{\lambda}_2\}\,\underline{A}C_0 r\leq\min\{\hat{\lambda}_1a(\tilde{p}_{\beta,1})\,{\cal E}_1,\hat{\lambda}_2a(\tilde{p}_{\beta,2})\,{\cal E}_2\}$, $i=1,2$, then we have:
\begin{equation}
\label{eqn 4F.14}
\beta({\bm\lambda})\leq\max_{i=1,2}\frac{a(\hat{p}_{\beta,i})\,|\lambda_i-\hat{\lambda}_i|}{\hat{\lambda}_i a(\tilde{p}_{\beta,i})\,{\cal E}_i-|R_i(\beta)|}
\end{equation}
$$\leq\Big\{\frac{2\,\overline{A}}{\underline{A}\,\min\{\underline{{\cal E}}_1,\underline{{\cal E}}_2\}}\Big\}\max\Big\{{\frac{\lambda_1-\hat{\lambda}_1}{\hat{\lambda}_1},\frac{\hat{\lambda}_2-\lambda_2}{\hat{\lambda}_2}}\Big\},$$
completing the proof.

\begin{lemma}
\label{lem 4G}
(Existence of strictly positively-ordered and locally Lipschitz-con-tinuously varying, local parametrized solution-families for Prob. \ref{prob 4I.+}) In the context of Prob. \ref{prob 4I.+} and Thm. \ref{thm 4.1.2}, let be given a solution $\hat{\bf\Gamma}\in{\bf X}(G)\cap C^{3,\tilde{\varrho}}$ of Prob. \ref{prob 4I.+} at $\boldsymbol{\lambda}(t_0)$ such that $({\bf i})$:$\,\hat{\bf\Gamma}\in\boldsymbol{\cal R}\big({\bm\lambda}(t_0);\rho_0\big)$ and $({\bf ii})$: $\hat{\bf\Gamma}\in\boldsymbol{\cal S}\big({\bm\lambda}(t_0); r_0\big)$ for some given values $r_0>0$ and $\rho_0\in(0,\pi)$. Then there exist a real open interval $I$ and a Lipschitz-continuously and  strictly monotonically varying parametrized solution family ${\bf\Gamma}(t)\,:\,I\rightarrow{\bf X}(G)\cap C^{3,\,\tilde{\varrho}}$ such that $t_0\in I$ and ${\bf\Gamma}(t_0)=\hat{\bf\Gamma}$, also such that for every $t\in I$, the pair ${\bf\Gamma}(t)\in{\bf X}(G)\cap C^{3,\,\tilde{\varrho}}$ solves Prob. \ref{prob 4I.+} at $\boldsymbol{\lambda}(t)$, and finally such that ${\bf\Gamma}(t)\in\boldsymbol{\cal R}\big({\bm\lambda}(t);\rho_1\big)$ and $\,{\bf\Gamma}(t)\in\boldsymbol{\cal S}\big({\bm\lambda}(t); r_1\big)$. both for all $t\in I$, where $r_1$ and $\rho_1$ denote suitable positive constants.
\end{lemma}

\noindent
{\bf Proof;} We define $U(t_0;p):=\,U({\bf\Gamma}(t_0);p)$ for all points $p$ in the closure of the region $\Omega(t_0):=\Omega\big({\bf\Gamma}(t_0)\big)$ (see Defs. \ref{def 2.1.1} and \ref{def 4A}). In view of assumption (${\bf i}$) above, it follows from Lem. \ref{lem 4B.2}(a),(b) that 
\vspace{.1in}

\noindent
(${\bf iii}$): there exists a value $\eta=\eta(t_0)>0$ such that the capacitary potential $U(t_0;p):{\rm Cl}\big(\Omega(t_0)\big)\rightarrow[0,1]$ has a single-valued $C^{\,3,\tilde{\varrho}}$-continuation by the same name $U(t_0;p)$ to a simply-connected domain containing the $\eta$-neighborhood of ${\rm Cl}\big(\Omega(t_0\big)$. In view of Def. \ref{def 4D} and the property (${\bf iii}$), it follows from Lem. \ref{lem 4E.1} and assumption (${\bf ii}$) that there exists a positive constant $\delta_1^*\in(0,\delta_0^*]$ small enough so that for any $\delta\in(0,\delta_1^*]$, and a corresponding unit vector ${\bm v}={\bm v}(t_0)$ (uniquely defined by (\ref{eqn 4E.34})), the two arc-pairs 
\begin{equation}
\label{eqn G.1}
{\bf \Gamma}_{-{\bm v}\delta}(t_0):={\bf M}_{-{\bm v}\delta}\big({\bf\Gamma}(t_0)\big)\,\,\,{\rm and}\,\,\,{\bf \Gamma}_{{\bm v}\delta}(t_0):={\bf M}_{{\bm v}\delta}\big({\bf\Gamma}(t_0)\big)
\end{equation}
constitute respective strict lower and strict upper classical solutions of Prob. \ref{prob 4I.} relative to any vector ${\bm\lambda}=(\lambda_1,\lambda_2)\in\Re^2_+$ such that (${\bf iv}$): $|\lambda_i-\lambda_i(t_0)|\leq\,r\lambda_i(t_0)v_i\delta$ for $i=1,2$ (see (\ref{eqn 4E.3})). By substituting the given locally Lipschitz-continuous mapping ${\bm\lambda}={\bm\lambda}(t):\Re\rightarrow\Re_+^2$ into the condition (${\bf iv}$), one sees that for any $\delta\in(0,\delta_1^*]$, the condition (${\bf iv}$) is satisfied by any vector ${\bm\lambda}:={\bm\lambda}(t)$ such that (${\bf v}$): $|\lambda_i(t)-\lambda_i(t_0)\big|\leq L\,|t-t_0|\leq r\lambda_i(t_0)\,v_i\,\delta$ for $i=1,2$, where $L$ denotes a uniform Lipschitz constant for both of the functions $\lambda_1(t),\lambda_2(t):I\rightarrow\Re_+$. Therefore, (${\bf vi}$): there exists a constant $C:=(r/L)\min\{\lambda_i(t_0)\,v_i:i=1,2\}>0$ such that for any $\delta\in(0, \delta_1^*]$, any and any $t\in I_\delta(t_0):=(t_0-C\delta, t_0+C\delta)$, the vector ${\bm\lambda}={\bm\lambda}(t)$ satisfies (${\bf iv}$). In view of Thms. \ref{thm 2.6.2} and \ref{thm 2.7.1}, it follows from (${\bf vi}$) that (${\bf vii}$): for any $\delta\in(0,\delta_1^*]$, $t_0\in I$, and $t\in I_{\delta}(t_0)$, there exists at least one classical solution ${\bf\Gamma}(t)\in{\bf X}\cap C^{\,1,1}$ of Prob. \ref{prob 4I.} at ${\bm\lambda}(t)$ such that 
\begin{equation}
\label{eqn G.2}
{\bf\Gamma}_{-{\bm v}\delta}(t_0)\,<\,{\bf\Gamma}(t)\,<\,{\bf \Gamma}_{{\bm v}\delta}(t_0)
\end{equation}

\noindent
and ${\bf\Gamma}(t)\in\boldsymbol{\cal R}\big({\bm\lambda}(t);\rho_0\big)\,\cap\,\boldsymbol{\cal S}\big({\bm\lambda}(t);r_0\big)$ (see Defs. \ref{def 2.7.2} and \ref{def 4A}). In this context, it follows from Lem. \ref{lem 4F}, Eq. (\ref{eqn 4F.2}) (by the substitutions: $\hat{\bm\lambda}={\bm\lambda}(t)$ and ${\bm\lambda}={\bm\lambda}(t+h)$), that
\begin{equation}
\label{eqn chris.1a}
{\bf\Gamma}(t)\leq{\bf M}_{\alpha h{\bm v}}\big({\bf\Gamma}(t)\big)\leq{\bf\Gamma}(t+h)\leq{\bf M}_{\beta h{\bm v}}\big({\bf\Gamma}(t)\big)
\end{equation}
for any $\delta\in(0,\delta_2^*]$, any $t\in I_{\delta}(t_0)$ such that $(t+h)\in I_{\delta}(t_0)$ for sufficiently small $h>0$, depending on $\delta\in(0,\delta_2^*]$, and for a suitable unit vector ${\bm v}={\bm v}(t)$ (see (\ref{eqn 4E.31}) and (\ref{eqn 4E.32})). In view of (\ref{eqn G.2}) and (\ref{eqn  chris.1a}), it follows from Lem. \ref{lem 4F} that (${\bf viii}$): for any $t_0\in I$ and $\delta\in(0,\delta_1^*]$, and for any classical solution $\hat{{\bf\Gamma}}\in{\bf X}\cap C^{\,3,\tilde{\varrho}}$ of Prob. \ref{prob 4I.} at $\hat{{\bm\lambda}}:={\bm \lambda}(t_0)$, there exists a strictly-
positively-ordered and locally Lipschitz-continuously-varying parametrized solution family ${\bf\Gamma}(t):I_{\delta}(t_0)\rightarrow{\bf X}\cap C^{\,3,\tilde{\varrho}}$ such that ${\bf\Gamma}(\tau)=\hat{{\bf\Gamma}}$, and also such that for each $t$ in the parameter-interval $I_{\delta}(t_0)$ with center-point 8$t_0\in I$ and constant length $2C\delta$, the ordered arc-pair ${\bf\Gamma}(t)\in{\bf X}\cap C^{\,3,\tilde{\varrho}}$ denotes a classical solution of Prob. \ref{prob 4I.+} at ${\bm\lambda}(t)$ (see Def. \ref{def 4.1.1}). 
\vspace{.1in} 

\noindent
{\bf Proof of Thm. \ref{thm 4.1.2}} In the context of Prob. \ref{prob 4I.+} and Thm \ref{thm 4.1.2}, we let $\tau_1$ (resp. $\tau_2$) denote the initial (terminal) endpoint of the interval $I$, and we define the classical solutions ${\bf\Gamma}(\tau_i)$, $i=1,2$, at the corresponding vector values ${\bf\lambda}(\tau_i)$, $i=1,2$, by continuity, so that the mapping ${\bf\Gamma}(t):[\tau_1,\tau_2]\rightarrow{\bf X}(G)\cap\,C^{\,3,\,\tilde{\varrho}}$ is continuous. By Lem. \ref{lem 4G}, for any sufficiently small values $\delta_1$, $\delta_2>0$, there exist the positively-ordered and locally Lipschitz-continuously varying local solution families ${\bf\Gamma}^*_i(t):\big(\tau_i-\delta_i,\tau_i+\delta_i\big)\rightarrow {\bf X}(G)\cap\,C^{\,3,\,\tilde{\varrho}}$,  $i=1,2$, such that ${\bf \Gamma}_i(\tau_i)={\bf\Gamma}(\tau_i)$ for $i=1,2$. In terms of these solution families, we define the positively-ordered, locally Lipschitz-continuously varying solution family $\dot{\bf\Gamma}(t):\big(\tau_1-\delta_1,\tau_2+\delta_2\big)\rightarrow {\bf X}(G)\cap\,C^{\,3,\,\tilde{\varrho}}$ such that $\dot{\bf\Gamma}(t)={\bf\Gamma}(t)$ for $\tau_1\leq t\leq \tau_2$, $\dot{\bf\Gamma}(t)={\bf\Gamma}^*_1(t)$ for $\tau_1-\delta_1<t\leq\tau_1$, and $\dot{{\bf\Gamma}}(t)={\bf\Gamma}^*_2(t)$ for $\tau_2\leq t<\tau_2+\delta_2$. 

\section{Arc-length and total curvature estimates} 
\label{section 5}
\subsection{Main operator and fixed-point estimates}
\label{subsection 5.1}
\begin{theorem} 
{(Main arc-length and total curvature estimates for the operators ${\bf T}_\varepsilon$)}

\label{thm 3.1.1}
\noindent
In the context of Defs. \ref{def 2.1.4} and \ref{def 2.1.5}, for any ${\bf\Gamma}=(\Gamma_1,\Gamma_2)\in{\bf X}$ and $\varepsilon\in(0,\varepsilon_0)$ (where $\varepsilon_0:={\rm min}\{1/2,$ $(\underline{A}^2/2A_1)\}$ throughout this paper), we let ${\bf T}_\varepsilon({\bf\Gamma})=$ $({T}_{\varepsilon,1}({\bf\Gamma}),{T}_{\varepsilon,2}({\bf\Gamma}))$ $={\bf\Psi}_\varepsilon\circ{\bf\Phi}_\varepsilon({\bf\Gamma})$, where ${\bf T}_\varepsilon$ (resp. ${\bf\Psi}_\varepsilon$) denotes either ${\bf }_\varepsilon^+$ (resp. ${\bf\Psi}_\varepsilon^+$) or ${\bf T}_\varepsilon^-$ (resp. ${\bf\Psi}_\varepsilon^-$). Then:
\vspace{.1in}

\noindent
(a) We have ${T}_{\varepsilon,i}({\bf\Gamma})\in{\rm X}$ for either $i=1,2$, any ${\bf\Gamma}:=(\Gamma_1,\Gamma_2)\in{\bf X}$, and any sufficiently small $\varepsilon\in(0,\varepsilon_0)$. 
\vspace{.1in}

\noindent
(b) There exist constants $A,B,C,H\geq 1$ sufficiently large (depending only on $\underline{A},\overline{A},A_1$, and $A_2$) such that for any $h\in(0,1/2)$, we have 
\begin{equation}
\label{eqn 3.1.1}
K(T_{\varepsilon,i}({\bf\Gamma}))\leq K(\Gamma_i)+\Big[K_{h,i}+(Ah-1)K(\Gamma_i)+Bh||\Gamma_i||\Big](\varepsilon/h)
\end{equation}
$$+HR(\Gamma_i,\Gamma_{h,i})(\varepsilon^2/h),$$
\begin{equation}
\label{eqn 3.1.2}
||{T}_{\varepsilon,i}({\bf\Gamma})||\leq ||\Gamma_i||+\Big[||\Gamma_{h,i}||+Ch\,K(\Gamma_i)-||\Gamma_i||\Big](\varepsilon/h)
\end{equation}
$$+H\,R(\Gamma_i,\Gamma_{h,i})(\varepsilon^2/h),$$
$i=1,2$, both uniformly for all curve-pairs ${\bf\Gamma}\in{\bf X}$ such that $||\Gamma_i||, K(\Gamma_i)<\infty$ and for all $0<\varepsilon\leq\varepsilon_0(h):={\rm min}\{\varepsilon_0,h\}$, where $\Gamma_{h,i}:=\Phi_{h,i}({\bf\Gamma})$, $K_{h,i}:=K(\Gamma_{h,i})$, and $R(\Gamma_i,\Gamma_{h,i}):={\rm max}\{||\Gamma_i||,||\Gamma_{h,i}||,K_i,K_{h,i}\}$. 
\vspace{.1in}

\noindent
(c) For $i=1,2$, and for any functional $F:{\rm X}\rightarrow \Re$ of the form: $F(\Gamma_i):=\lambda K(\Gamma_i)+\mu ||\Gamma_i||$, with $\lambda,\mu>0$ and $\lambda+\mu=1$, it follows from Eqs. (\ref{eqn 3.1.1}) and (\ref{eqn 3.1.2}) that
\begin{equation}
\label{eqn 3.1.3}
F(T_{\varepsilon,i}({\bf\Gamma}))\leq F(\Gamma_i)+F(\Gamma_{h,i})(\varepsilon/h)
\end{equation}
$$-\Big[\big((1-Ah)-(\mu/\lambda) Ch\big)\lambda K(\Gamma_i)+\big(1-(\lambda/\mu)\,Bh\big)\mu||\Gamma_i||\Big](\varepsilon/h)$$
$$+HR(\Gamma_i,\Gamma_{h,i})(\varepsilon^2/h)$$
for $i=1,2$, $\varepsilon\in(0,\varepsilon_0(h)]$, and ${\bf\Gamma}\in {\bf X}$ such that $||\Gamma_i||, K(\Gamma_i)<\infty$.
\vspace{.1in}

\noindent
(d) In (\ref{eqn 3.1.3}), we define:
\begin{equation}
\label{eqn 3.1.3a}
\lambda:=\frac{2C}{\sqrt{A^2+4BC}+2C-A}\,\,{\rm and}\,\,\mu:=\frac{\sqrt{A^2+4BC}-A}{\sqrt{A^2+4BC}+2C-A}.
\end{equation}
Then the values $\lambda,\mu>0$ are such that $\lambda+\mu=1$,  
\begin{equation}
\label{eqn 3.1.3aa}
A+(\mu/\lambda)C=(\lambda/\mu) B=P_0:=\big(A+\sqrt{A^2+4BC}\,\big)/2>0,
\end{equation}
and $\lambda,\mu\geq(2\big/LP_0)$, where $L$ denotes the denominator in (\ref{eqn 3.1.3a}). It follows from (\ref{eqn 3.1.3}) and (\ref{eqn 3.1.3aa}) that for any $h\in(0,1/2)$, we have
\begin{equation}
\label{eqn 3.1.4}
F(T_{\varepsilon,i}({\bf\Gamma}))\leq F(\Gamma_i)+\big[F(\Gamma_{h,i})+(P_0 h-1) F(\Gamma_i)\big](\varepsilon/h)+HR(\Gamma_i,\Gamma_{h,i})(\varepsilon^2/h),
\end{equation}
for $i=1,2$, all $\varepsilon\in(0,\varepsilon_0(h)]$, and all ${\bf\Gamma}\in {\bf X}$ such that $||\Gamma_i||, K(\Gamma_i)<\infty$.
\end{theorem}
\begin{theorem} 
{(Arc-length and total curvature estimate for fixed points of the operators ${\bf T}_\varepsilon$)}
\label{thm 3.1.2}
(a) In the context of Thm. \ref{thm 3.1.1}, choose a value $h\in(0,1/2)$ such that $2P_0 h<1$. Then there exists a value $\hat{\varepsilon}_0(h)\in(0,\varepsilon_0(h)]$, depending only on $\underline{A},\overline{A},A_1, A_2$, and $h$, such that 
\begin{equation}
\label{eqn 3.1.5}
F(T_{\varepsilon,i}({\bf\Gamma})
)\leq F(\Gamma_i)+\Big(2F(\Phi_{h,i}({\bf\Gamma}))+(2hP_0-1)F(\Gamma_i)\Big)(\varepsilon/h)
\end{equation}
for $i=1,2$ and any $\varepsilon\in(0,\hat{\varepsilon}_0(h)]$ and ${\bf\Gamma}\in {\bf X}$. 
\vspace{.1in}

\noindent
(b) Let ${\bf Y}$ be the invariant set in Def. \ref{def 2.1.8} such that ${\bf T}_\varepsilon:{\bf Y}\rightarrow{\bf Y}$ for $\varepsilon\in(0,\varepsilon_1]$, where $\varepsilon_1\in(0,\varepsilon_0)$ and ${\bf T}_\varepsilon$ denotes either ${\bf T}_\varepsilon^+$ or ${\bf T}_\varepsilon^-$. Let be given $h\in(0,1/2)$ such that $2P_0 h<1$ and $0<\hat{\varepsilon}_0(h)\leq\varepsilon_1$. Then for any $\varepsilon\in(0,\hat{\varepsilon}_0(h)]$, ${\bf T}_\varepsilon$ has a "fixed point" ${\bf \Gamma}_\varepsilon\in{\bf Y}$ (see  Thms. \ref{thm 2.1.2} and \ref{thm 2.1.3}) such that $F\big(\Gamma_{\varepsilon,\,i}\big)<\infty$ for $i=1,2$.
\vspace{.1in}

\noindent
(c) For any $\varepsilon\in(0,\hat{\varepsilon}_0(h)]$ and any fixed point ${\bf\Gamma}_\varepsilon:=(\Gamma_{\varepsilon,1},\Gamma_{\varepsilon,2})\in{\bf Y}$ of ${\bf T}_\varepsilon$ such that $F\big(\Gamma_{\varepsilon,\,i}\big)<\infty$, $i=1,2$, it follows from  (\ref{eqn 3.1.5}) that
\begin{equation}
\label{eqn 3.1.6}
\big(1-2hP_0\big)F\big(\Gamma_{\varepsilon,i})\leq 2F(\Phi_{h,i}({\bf\Gamma}_\varepsilon)\big)
\end{equation}
for $i=1,2$. In view of Lem. \ref{lem 2.2.5}, there is a constant $M(h)$ such that $F(\Phi_{h,i}({\bf\Gamma}_\varepsilon))$ $\leq M(h)$ for $i=1,2$, $\varepsilon\in(0,\hat{\varepsilon}_0(h)]$, and all fixed points ${\bf\Gamma_\varepsilon}\in{\bf Y}$ of the operators ${\bf T}_\varepsilon$. 
In view of this, and the estimate: $\lambda,\mu\geq\big(2\big/L\,P_0\big)$ (see Thm. \ref{thm 3.1.1}(d)), it follows from (\ref{eqn 3.1.6}) that 
\begin{equation}
\label{eqn 3.1.6a}
\big(2\big/L\,P_0\big)\max\big\{||\Gamma_{\varepsilon,\,i}||,K(\Gamma_{\varepsilon,\,i})\big\}\leq F(\Gamma_{\varepsilon,\,i})\leq\big(M(h)\big/(1-2hP_0)\big)
\end{equation}
for $i=1,2$, all $\varepsilon\in(0,\hat{\varepsilon}_0(h)]$, and all fixed points ${\bf\Gamma}_{\varepsilon,\,i}\in{\bf Y}$ of ${\bf T}_{\varepsilon,\,i}$
such that $F({\Gamma}_{\varepsilon,\,i})<\infty$, $i=1,2$.
\end{theorem}
\begin{corollary}
\label{cor 3.1.1a} 
{(Further estimates)} 
In Prob. \ref{prob2.1.1} and Def. \ref{def 2.1.3}, let the functions $a_1(p), a_2(p)$ be replaced by the related functions $\hat{a}_i(p):=\kappa a_i(p)$, $i=1,2$. Then in Thms. \ref{thm 3.1.1} and \ref{thm 3.1.2}, the estimates (\ref{eqn 3.1.1})-(\ref{eqn 3.1.6a}) continue to hold, where the constants $A, B, C$ are replaced by the new constants $(A/\kappa)$, $(B/\kappa)$, $(C/\kappa)$.
\end{corollary}

\subsection{Length and turning-angle estimates for arc-partitions}
\label{subsection 5.2}
\begin{definition}
\label{def 3.2.1} {(Operator and polygon notation)}
(a) In slightly revised notation from Def. \ref{def 2.1.5}, for any $\varepsilon\in(0,\varepsilon_0)$ and any arc $\Gamma\in\tilde{{\rm X}}_\varepsilon$, we define the arcs $\Gamma_{\varepsilon,i}=\Psi_{\varepsilon,i}(\Gamma):=\partial G_{\varepsilon,i}(\Gamma)$, $i=1,2$, where the sets $G_{\varepsilon,i}(\Gamma)$, $i=1,2,$ are defined in Def. \ref{def 2.1.4}. We observe that $G_{\varepsilon,i}(\Gamma)=G_{\varepsilon,i}\big(\hat{\Gamma}_{\varepsilon,i}\big)$, $i=1,2$, where we define the arcs $\hat{\Gamma}_{\varepsilon,i}:=\partial \hat{G}_{\varepsilon,i}(\Gamma))$ (see Def. \ref{def 2.1.4} and Lem. \ref{lem 2.1.4*}). Observe that $||\hat{\Gamma}_{\varepsilon,i}||\leq||\Gamma||$ and $K\big(\hat{\Gamma}_{\varepsilon,i}\big)\leq K(\Gamma)$ for $i=1,2$ (see ). In view of this, it suffices to Prove the assertions of this section in the special case where $\Gamma$ has the properties of the arc $\hat{\Gamma}_{\varepsilon,i}$.
\vspace{.1in}

\noindent
(b) Given a fixed, $P$-periodic (in $x$), double-point free polygonal arc $\Gamma\in{\rm X}$, we use $V(\Gamma)$ to denote the (finite) set of all vertices $q$ of $\Gamma$, while $S(\Gamma)$ denotes the (finite) set of all sides $L$ of $\Gamma$ (these being straight line-segments of $\Gamma$ which, by convention, do not contain their endpoints).
\end{definition}
\begin{definition}
\label{def 3.2.2} {(Ordering of points of an arc $\Gamma\in{\rm X})$} Let $\Gamma\in{\rm X}$ be a continuous, double-point-free arc. Then for any points $p_1,p_2\in\Gamma$, we use the notation: "$p_1\leq p_2$" to mean that a point traveling in the positive direction (from $x=-\infty$ to $x=+\infty$) along $\Gamma$ will not pass $p_2$ before passing $p_1$. Of course the notation "$p_1<p_2$" means that $p_1\leq p_2$ and that $p_1\not=p_2$. We remark that any open segment $A_{\varepsilon,i}$ of $\Gamma$ (i.e. relatively-open, connected proper subset $A$ of ${\Gamma}$ has the form: $A=\{p\in{\Gamma}_{\varepsilon,i}:p_1<p<p_2\}$ for suitable points $p_1,p_2\in{\Gamma}$.    
\end{definition}

\begin{definition}
\label{def 3.2.3} {(A partition of ${\Gamma}_{\varepsilon,i}:=\Psi_{\varepsilon,i}(\Gamma)$, where $\Gamma\in{\rm X}$ is a polygonal arc)} For each $\varepsilon\in(0,\varepsilon_0)$ and $i\in\{1,2\}$, each polygonal arc $\Gamma\in{\rm X}$, and any point $p\in{\Gamma}_{\varepsilon,i}:=\Psi_{\varepsilon,i}(\Gamma)$ (see Def. \ref{def 3.2.1}(a)), we set
\begin{equation}
\label{eqn 3.2.1}
\Pi_{\varepsilon,i}(p)=\big\{q\in\Gamma:|p-q|={\rm dist}(p,\Gamma)\big\},
\end{equation}
which is clearly a closed, non-empty subset of $\Gamma$. For each side $L\in S(\Gamma)$ and vertex $q\in V(\Gamma)$, we define 
\begin{equation}
\label{eqn 3.2.1r} 
A_{\varepsilon,i}(L)=\big\{p\in{\Gamma}_{\varepsilon,i}:\Pi_{\varepsilon,i}(p)\subset L\big\},
\end{equation}
\begin{equation}
\label{eqn 3.2.1r2}
B_{\varepsilon,i}(q)=\big\{p\in{\Gamma}_{\varepsilon,i}:\Pi_{\varepsilon,i}(p)=\{q\}\,\big\},
\end{equation}
\begin{equation}
\label{eqn 3.2.1r1}
C_{\varepsilon,i}=\big\{p\in{\Gamma}_{\varepsilon,i}:\Pi_{\varepsilon,i}(p)\,\,\, {\rm contains\,\, at\,\, least\,\, two\,\, points}\big\}.
\end{equation} 
\end{definition}
\begin{theorem}
\label{thm 3.2.1}
(a) For either $i=1$ or $i=2$, any fixed value $\varepsilon\in(0,\varepsilon_0)$, and any fixed $P$-periodic (in $x$) polygonal arc $\Gamma\in{\rm X}$, the sets
\begin{equation}
\label{eqn 3.2.2}
A_{\varepsilon,i}(L),L\in S(\Gamma);\,\, B_{\varepsilon,i}(q),q\in V(\Gamma);\,\big\{p\big\},\,\,p\in C_{\varepsilon,i},
\end{equation}
some of which may be empty, constitute a partition of the arc ${\Gamma}_{\varepsilon,i}:=\Psi_{\varepsilon,i}(\Gamma)$.
\vspace{.1in}

\noindent
(b) For each side $L\in S(\Gamma)$, $A_{\varepsilon,i}(L)$ is an open set relative to ${\Gamma}_{\varepsilon,i}$.
\vspace{.1in}

\noindent
(c) For any $L\in S(\Gamma)$ (resp. $q\in V(\Gamma)$), the set $A_{\varepsilon,i}(L)$ (resp. $B_{\varepsilon,i}(q)$) is a segment (i.e. connected subset) of ${\Gamma}_{\varepsilon,i}$.
\vspace{.1in}

\noindent
(d) The partition (\ref{eqn 3.2.2}) of ${\Gamma}_{\varepsilon,i}$ consists of at most finitely many connected sets relative to any $P$-period of ${\Gamma}_{\varepsilon,i}$.
\end{theorem}
\begin{lemma}
\label{lem 3.2.1}
(a) If $\emptyset\not=\Pi_{\varepsilon,i}(p)\subset L$ (for given $\varepsilon\in(0,\varepsilon_0)$, $i\in\{1,2\}$, polygonal arc $\Gamma\in{\rm X}$, $p\in{\Gamma}_{\varepsilon,i}$, and $L\in S(\Gamma)$), then $\Pi_{\varepsilon,i}(p)\cap L=\{q\}$, where the point $q\in L$ is the unique orthogonal projection of $p$ onto $L$. 
\vspace{.1in}

\noindent
(b) For fixed $\varepsilon\in(0,\varepsilon_0)$ and $i\in\{1,2\}$, let $p_1,p_2\in\Gamma_{\varepsilon,i}$ denote points such that $p_1<p_2$ relative to $\Gamma_{\varepsilon,i}$. Then $q_1\leq q_2$ relative to $\Gamma$, where $q_1,q_2\in\Gamma$ denote any points such that $q_1\in\Pi_{\varepsilon,i}(p_1)$ and $q_2\in\Pi_{\varepsilon,i}(p_2)$
\end{lemma}
\noindent
{\bf Proof.} Part (a) is self-explanitory. Regarding Part (b), let $\gamma_1$ (resp. $\gamma_2$) denote the straight-line-segment having the initial endpoint $p_1$ (resp. $p_2$) in $\Gamma_{\varepsilon,i}$ and the terminal endpoint $q_1$ (resp. $q_2$) in $\Gamma$. We remark that apart from their initial endpoints, the line-segments $\gamma_1,\gamma_2$ lie entirely in the open set $D_{3-i}\big(\Gamma_{\varepsilon,i}\big)$ (due to Lem. \ref{lem 3.3.2}(a)), and that apart from their terminal endpoints, they both lie in the open set $D_i(\Gamma)$ (since $|q_j-p_j|=
{\rm dist}\big(p_j,\Gamma\big)$). Using these facts, and assuming that the assertion is false, we have that $q_2<q_1$ in $\Gamma$. Therefore, the line-segments $\gamma_1$ and $\gamma_2$ intersect at a point $r$. For $\alpha_j=|p_j-r|$ and $\beta_j=|q_j-r|$, $j=1,2$, we have $\alpha_1+\beta_1\leq \alpha_1+\beta_2$, since no arc joining $p_1$ to $\Gamma$ is shorter than $\gamma_1$. Thus $\beta_1\leq\beta_2$. It follows that $|p_2-q_1|<\alpha_2+\beta_1\leq\alpha_2+\beta_2={\rm dist}(p_2,\Gamma)$, a contradiction.
\vspace{.1in}

\noindent
{\bf Proof of Thm. \ref{thm 3.2.1} (b).} {(The sets $A_{\varepsilon,i}(L)$ are open (relative to $\Gamma_{\varepsilon,i}$)).} For fixed $\varepsilon\in(0,\varepsilon_0)$ and $i\in\{1,2\}$, let $p_0$ denote a point in $A_{\varepsilon,i}(L)$ for some fixed $L\in S(\Gamma)$. Then $\emptyset\not=\Pi_{\varepsilon,i}(p_0)\subset L$. Therefore $\Pi_{\varepsilon,i}(p_0)=\{q_0\}$, where $q_0$ denotes the unique orthogonal projection of $p_0$ onto $L$ and therefore the unique point closest to $p_0$ in $L$ (see Lem. \ref{lem 3.2.1}(a)). Clearly, we have $|q-p_0|>|q_0-p_0|$ for all points $q\in\Gamma\setminus L$, since $\Pi_{\varepsilon,i}(p_0)\subset L$ (see (\ref{eqn 3.2.1r}, \ref{eqn 3.2.1r2}, and (\ref{eqn 3.2.1r1})). Since the set $\Gamma\setminus L$ is closed, it follows that there exists a positive constant $\eta_0>0$ (depending on $p_0$) such that $|q-p_0|\geq|q_0-p_0|+\eta_0$ uniformly for all points $q\in\Gamma\setminus L$. Therefore, for any point $p\in{\Gamma}_{\varepsilon,i}$ such that $|p-p_0|<(\eta_0/3)$, we have that $|p-q_0|\leq\big(|p_0-q_0|+|p-p_0|\big)\leq \big(|p_0-q_0|+(\eta_0/3)\big)$ and $|p-q|\geq\big(|p_0-q|-|p-p_0|\big)\geq\big(|p_0-q_0|+\eta_0\big)-|p-p_0|\geq\big(|p_0-q_0|+(2\eta_0/3)\big)$ for all $q\in\Gamma\setminus L$, from which it follows that $p\in A_{\varepsilon,i}(L)$. This completes the proof that the set $A_{\varepsilon,i}(L)$ is open relative to $\Gamma_{\varepsilon,i}$, since $p_0$ is arbitrary in $A_{\varepsilon,i}(L)$.
\vspace{.1in}

\noindent
{\bf Proof of Thm. \ref{thm 3.2.1}(c).} (The sets $A_{\varepsilon,i}(L)$, $B_{\varepsilon,i}(q)$ are connected)) To prove that the set  $A_{\varepsilon,i}:=A_{\varepsilon,i}(L)$ (which is open relative to $\Gamma_{\varepsilon,i}$) is connected (for any fixed $\varepsilon\in(0,\varepsilon_0)$ and $i\in\{1,2\}$), we will show that $A_{\varepsilon,i}=I_{\varepsilon,i}$, where $I_{\varepsilon,i}$ denotes the smallest open arc-segment of ${\Gamma}_{\varepsilon,i}$ containing $A_{\varepsilon,i}$. For any point $p\in I_{\varepsilon,i}$, there exist points $p_1,p_2\in A_{\varepsilon,i}$ such that and $p_1<p<p_2$. Thus, for any $q\in\Pi_{\varepsilon,i}(p)$ we have $q_1\leq q\leq q_2$ (by double application of Lem. \ref{lem 3.2.1}(b)), where $\Pi_{\varepsilon,i}(p_j)=\{q_j\}\subset L$, $j=1,2$. Thus $\Pi_{\varepsilon,i}(p)\subset L$, whence $p\in A_{\varepsilon,i}$. Thus $A_{\varepsilon,i}$ is connected. Finally, regarding the connectedness of the set $B_{\varepsilon,i}(q)$, it follows from Lem. \ref{lem 3.2.1}(b) that for any points $p_1, p, p_2\in\Gamma_{\varepsilon,i}$ such that $p_1,p_2\in B_{\varepsilon,i}(q)$ and $p_1<p<p_2$ in terms of the natural ordering in $\Gamma_{\varepsilon,i}$, we have that $p\in B_{\varepsilon,i}(q)$.
\begin{lemma}
\label{lem 3.2.2}
For fixed $\varepsilon\in(0,\varepsilon_0)$ and $i\in\{1,2\}$ and any fixed double-point-free polygonal arc $\Gamma\in{\rm X}$, there exists a value $\rho>0$ such that for any point $p\in{\Gamma}_{\varepsilon,i}:=\Psi_{\varepsilon,i}(\Gamma)$ and any distinct points $q_1,q_2\in\Pi_{\varepsilon,i}(p)\subset\Gamma$, we have $|q_1-q_2|\geq\rho$.
\end{lemma}
\noindent
{\bf Proof.} If $q_1,q_2\in V(\Gamma)$, then $|q_1-q_2|\geq\rho_1>0$, where $\rho_1>0$ is the minimum distance between distinct vertices of $\Gamma$. We assume for the remainder of the proof that $q_1\in L_1$ for some side $L_1\in S(\Gamma)$. Then $q_2\notin{\rm Cl}(L_1)$, by Lem. \ref{lem 3.2.1}(a). If $q_2$ does not belong to either of the sides of $\Gamma$ which are adjacent to $L_1$, then $|q_1-q_2|\geq\rho_2$, where $\rho_2>0$ is the minimum distance between non-adjacent sides of $\Gamma$. Finally, assume $q_2$ belongs to a side $L_2$ of $\Gamma$ which is adjacent to $L_1$. Then it is easily seen that $|q_2-q_1|\geq 2|p-q_1|{\rm sin}(|\theta|/2)$, where $|p-q_1|\geq(\varepsilon/a(p))\geq(\varepsilon\big/\,\overline{A}\,\,)$ and $\theta\in(-\pi,0)$ is the turning angle of the forward tangent to $\Gamma$ at the common vertex of $L_1$ and $L_2$. Therefore $|q_2-q_1|\geq\rho_3:=(2\,b\,\varepsilon/\,\overline{A}\,\,)>0$, where $b$ is the minimum of ${\rm sin}(|\theta|/2)$ over the turning angles of $\Gamma$ at all vertices of $\Gamma$. We remark that in all three above cases, the (strict) positivity of the constants $\rho_1, \rho_2, \rho_3$ depends on the $P$-periodicity of the fixed polygonal arc $\Gamma\in{\rm X}$. This completes the proof, where we set $\rho={\rm min}\{\rho_1,\rho_2,\rho_3\}$.
\begin{lemma}
\label{lem 3.2.3}
For fixed $\varepsilon\in(0,\varepsilon_0)$, $i\in\{1,2\}$, and any fixed polygonal arc $\Gamma\in{\rm X}$, the set $C_{\varepsilon,i}$ contains at most finitely-many points relative to any $P$-period (in $x$) of the arc ${\Gamma}_{\varepsilon,i}$.
\end{lemma}
\noindent
{\bf Proof.} Let $p_1,p_2,\cdots,p_n$ denote $n\geq 3$ distinct points of $C_{\varepsilon,i}$, restricted to one $P$-period of $\Gamma_{\varepsilon,i}$. We can assume that $p_1<p_2<\cdots<p_n$ in terms of the natural ordering on ${\Gamma}_{\varepsilon,i}$, where $p_{n}=p_1+(P,0)$. By Lem. \ref{lem 3.2.1}(b), there exist points $q_j^\pm\in\Pi_{\varepsilon,i}(p_j),j=1,\cdots,n,$ such that
\begin{equation}
\label{eqn 3.2.3}
q_1^-<q_1^+\leq q_2^-<q_2^+\leq\cdots\leq q_{n-1}^-<q_{n-1}^+\leq q_n^-<q_n^+.
\end{equation}
For each $j=1,\cdots,n$, the arc-segment $\gamma_j$ of $\Gamma$ which joins $q_j^-$ to $q_j^+$ (without intersecting the other points) has length $||\gamma_j||\geq\rho$, due to Lem. \ref{lem 3.2.2}. This is a contradiction unless $n\leq(||\Gamma||/\,\rho)$.
\vspace{.1in}

\noindent
{\bf  Proof of Thm. \ref{thm 3.2.1}, Parts (a) and (d).} \,\,Taking Lem. \ref{lem 3.2.1}(a) into account, Thm. \ref{thm 3.2.1}(a) simply re-expresses the fact that for any fixed $p\in{\Gamma}_{\varepsilon,i}$, exactly one of the following three alternatives is true: $({\bf i})$: $\Pi_{\varepsilon,i}(p)$ is a one-point set contained in one of the sides $L\in S(\Gamma)$, or $({\bf ii})$: $\Pi_{\varepsilon,i}(p)$ is a one-point set containing one of the vertices $q\in V(\Gamma)$, or $({\bf iii})$: $\Pi_{\varepsilon,i}(p)$ is not a one-point set. Concerning Part (d), the partition (\ref{eqn 3.2.2}) of $\Gamma_{\varepsilon,i}$ consists of finitely many connected components per $P$-period of $\Gamma_{\varepsilon,i}$, since each $L\in S(\Gamma)$ (resp. $q\in V(\Gamma)$) corresponds to at most one connected set $A_{\varepsilon,i}(L)$ (resp. $B_{\varepsilon,i}(q)$) of ${\Gamma}_{\varepsilon,i}$, and since the set $C_{\varepsilon,i}$ contains at most finitely-many points per $P$-period. 

\begin{lemma} {(Arc-length and total curvature estimates for the image of a line-segment)}
\label{lem 3.3.1}
In the context of Def. \ref{def 3.2.1}, for given $i\in\{1,2\}$, $\varepsilon\in(0,\varepsilon_0)$, and any given $P$-periodic polygonal arc $\Gamma\in {\rm X}$, let $\gamma_{\varepsilon,i}$ denote any (connected) arc-segment of the arc $\Gamma_{\varepsilon,i}:=\Psi_{\varepsilon,i}(\Gamma)$ which projects orthogonally onto a connected segment $I$ of a straight-line-segment $L$ of $\Gamma$. Then there exist constants $C_1:=(2\,A_1\,/\,\underline{A}^2)$ and $C_2:=\big((2\,\underline{A}A_2+4A_1^2)\big/\,\underline{A}^3\big)\leq\big((2\underline{A}+8\delta)A_2\big/\,\underline{A}^3\big)$ such that 
$$\gamma_{\varepsilon,i}=\{(x,y_{\varepsilon,i}(x)):x\in I\},$$
in suitable Cartesian coordinates (such that $I=(\alpha,\beta)\times\{0\}\subset\Gamma$ and $\omega_{\varepsilon,i}:=\{(x,y):0<y<y_{\varepsilon,i}(x),x\in I\}\subset\Omega(\Gamma_{\varepsilon,i})$), where $y=y_{\varepsilon,i}(x):I\rightarrow \Re$ is a $C^2$ function such that
\begin{equation}
\label{eqn 3.3.1}
|y'_{\varepsilon,i}(x)|\leq C_1\varepsilon\,\,{\rm and}\,\, |\kappa_{\varepsilon}(x,y_{\varepsilon,i}(x))|\leq C_2\varepsilon
\end{equation}
for $x\in I$. Thus, we have
\begin{equation}
\label{eqn 3.3.2}
||\gamma_{\varepsilon,i}||\leq \sqrt{1+C_1^2\,\varepsilon^2}||I||\leq (1+(C_1^2/2)\,\varepsilon^2)\,||I||
\end{equation}
and
\begin{equation}
\label{eqn 3.3.3}
K(\gamma_{\varepsilon,i})\leq C_2\,||\gamma_{\varepsilon,i}\,||\leq \overline{C}_2||I||\,\varepsilon,
\end{equation}
where $\overline{C}_2:=C_2\,\sqrt{1+C_1^2\,\varepsilon_0^2}\leq\sqrt{2}\,C_2$.
\end{lemma}
\noindent
{\bf Proof.}
For convenience, we fix the index $i\in\{1,2\}$ and suppress the subscript "$i$" in the following notation: We consider the equation:
\begin{equation}
\label{eqn 3.3.4}
\phi(x,y):=y\,\,a(x,y)=\varepsilon,
\end{equation}
where we assume $a\in\boldsymbol{{\cal A}}, x\in \Re$ and $y>0$. By the intermediate value theorem, (\ref{eqn 3.3.4}) has at least one solution $y=y_\varepsilon(x)$ for each $x$. Any solution must be such that
\begin{equation}
\label{eqn 3.3.5}
(\varepsilon/\,\overline{A}\,)\leq y_\varepsilon(x)\leq(\varepsilon\,/\underline{A}\,).
\end{equation}
We have
\begin{equation}
\label{eqn 3.3.6}
\phi_y(x,y)=a+y\,a_y\geq \underline{A}-(A_1\varepsilon/\underline{A}\,)\geq (\underline{A}\,/2)
\end{equation}
for $0<y\leq(\varepsilon\,/\underline{A}\,)$ and $0<\varepsilon<\varepsilon_0$.
If $\varepsilon\in(0,\varepsilon_0)$, then the equation:
$$y_\varepsilon(x)\,a(x,y_\varepsilon(x))=\varepsilon$$
is solved by a unique $C^2$-function $y=y_\varepsilon(x):\Re\rightarrow \Re_+$ satisfying the condition (\ref{eqn 3.3.5}). By differentiation, we have $y'_{\varepsilon}(x)=\left((-y_{\varepsilon}\,a_x)/(a+y_{\varepsilon}\,a_y)\right)$. In view of (\ref{eqn 3.3.6}), it follows from this that
$$|y'_\varepsilon(x)|\leq \left((A_1\,y_{\varepsilon}(x))/(\underline{A}/2)\right)\leq C_1\,\varepsilon,$$
for $0<\varepsilon\leq\varepsilon_0$, where $C_1\varepsilon_0\leq 1$.
\vspace{.1in}

\noindent
To estimate the curvature $\kappa_{\varepsilon}(p)$ of $\gamma_{\varepsilon}$ at $p\in\gamma_{\varepsilon}$, we use the formula: $\kappa(p)=-(\phi_{\boldsymbol{\tau}\boldsymbol{\tau}}/\phi_{\boldsymbol{\nu}})$, where $\boldsymbol{\tau}$ points in the tangent direction on $\gamma_{\varepsilon}$, and $\boldsymbol{\nu}$ is the corresponding exterior normal (to $\omega_{\varepsilon}$). We have  
\begin{equation} 
\label{eqn 3.3.8}
\phi_{\boldsymbol{\tau}}=y_{\boldsymbol{\tau}}\,a+y\, a_{\boldsymbol{\tau}}=0;\quad \phi_{\boldsymbol{\nu}}=|\nabla\phi|=y_{\boldsymbol{\nu}}\,a+y\,a_{\boldsymbol{\nu}};\quad \phi_{\boldsymbol{\tau}\boldsymbol{\tau}}=y\,  a_{\boldsymbol{\tau}\boldsymbol{\tau}}+2\,a_{\boldsymbol{\tau}} \,y_{\boldsymbol{\tau}}. 
\end{equation}
It follows from the first two equations in (\ref{eqn 3.3.8}) that
\begin{equation}
\label{eqn 3.3.9}
|\nabla\phi|\geq\sqrt{1-y_{\boldsymbol{\tau}}^2}-y|a_{\boldsymbol{\nu}}|=a\,\sqrt{1-(y\,a_{\boldsymbol{\tau}}/a)^2}-y|a_{\boldsymbol{\nu}}|\geq\,(a/2)
\end{equation}
for $0<2y|\nabla a|\leq a$. It follows from (\ref{eqn 3.3.9}) and the final equation in (\ref{eqn 3.3.8}) that
$$|\kappa_{\varepsilon}(p)|=\frac{|y\,a_{\boldsymbol{\tau}\boldsymbol{\tau}}+2\,y_{\boldsymbol{\tau}}\,a_{\boldsymbol{\tau}|}}{|\nabla\phi|}=
\frac{y\,\left|a \,a_{\boldsymbol{\tau}\boldsymbol{\tau}}-2\,a_{\boldsymbol{\tau}}^2\right|}{a|\nabla\phi|}\leq
\frac{2y|a a_{\boldsymbol{\tau}\boldsymbol{\tau}}-2\,a_{\boldsymbol{\tau}}^2|}{a^2}
\leq C_2\,\varepsilon,
$$
for all $0<\varepsilon<\varepsilon_0$ and $p\in\gamma_\varepsilon$.

\begin{lemma}
\label{lem 3.3.2} {(Polar-coordinate arc-length and total curvature estimates)}
(a) For any $\varepsilon\in(0,\varepsilon_0)$, any function $a\in\boldsymbol{{\cal A}}$, and any point $p_\varepsilon\in\Re^2$ such that $|p_\varepsilon|\,a(p_\varepsilon)=\varepsilon$ (of which there exists exactly one on every radial), we have $|q|\,a(q)<\varepsilon$ for all $q\in\gamma_\varepsilon$, where $\gamma_\varepsilon$ denotes the straight line joining $p_\varepsilon$ to the origin.
\vspace{.1in}

\noindent
(b) Let the closed arc $\gamma_\varepsilon$ be the solution set for the equation $|p|\,a(p)=\varepsilon$, and, for any $\alpha\leq\beta\leq\alpha+2\pi$, let $\gamma_\varepsilon(J)$ denote the intersection of $\gamma_\varepsilon$ with the polar-coordinate angular sector $J:=[\alpha,\beta]$. Then there exist constants $C_3:=(A_1/\underline{A}^3)+(4A_1^2\varepsilon_0/\underline{A}^5)$ and $C_4$ (depending only on $\underline{A}, \overline{A}, A_1,A_2$; to be defined in (\ref{eqn 3.3.20a}) and (\ref{eqn 3.3.21a}).) such that for any $\varepsilon\in(0,\varepsilon_0)$, we have 
\begin{equation}
\label{eqn 3.3.10}
|\,||\gamma_\varepsilon(J)||-(\varepsilon/a(0))(\beta-\alpha)|\leq C_3 (\beta-\alpha)\,\varepsilon^2,
\end{equation}
\begin{equation}
\label{eqn 3.3.11}
|K(\gamma_\varepsilon(J))-(\beta-\alpha)|\leq C_4(\beta-\alpha)\,\varepsilon.
\end{equation}
\end{lemma}
\noindent
{\bf Proof.}
In polar coordinates, the equation $\phi(r,\theta):=r\,a(r,\theta)=\varepsilon$ has at least one solution $r=r_\varepsilon(\theta):\Re\rightarrow\Re_+$ at $\varepsilon>0$, where we have $(\varepsilon\big/\,\overline{A}\,)\leq r_\varepsilon(\theta)\leq(\varepsilon\big/\underline{A}\,)$ for any solution at $\varepsilon>0$. By differentiation of the function: $\phi(r,\theta)=r\,a(r,\theta)$, we see that $$\phi_r(r,\theta)=a(r,\theta)+r\,a_r(r,\theta)\geq\underline{A}-A_1 r\geq \underline{A}-(A_1\big/\underline{A}\,)\,\varepsilon\geq(\,\underline{A}\,/2)$$ for $0<r<(\varepsilon\big/\underline{A}\,)$ and $0<\varepsilon<\varepsilon_0:=\min\{1/2, (\,\underline{A}^2\,)/2A_1)\}$, which shows that for each $\theta\in\Re$, there exists exactly one value $r_{\varepsilon}(\theta)>0$ such that $(r_{\varepsilon}(\theta),\theta)\in\gamma_\varepsilon$, and that, for this value, we have $(r,\theta)\in\omega_\varepsilon$ for all values $r\in[0,r_\varepsilon(\theta))$, implying that the region $\omega_\varepsilon$ is simply connected. By the theorem of the mean, we have
\begin{equation}
\label{eqn 3.3.12}
\big|a(r_{\varepsilon}(\theta),\theta)-a(0)\big|\leq A_1 r_{\varepsilon}(\theta)\leq\big(A_1/\underline{A}\,\big)\,\varepsilon,
\end{equation} 
and it follows by substituting the equation $r_{\varepsilon}(\theta)\,a(r_{\varepsilon}(\theta),\theta)=\varepsilon$ into (\ref{eqn 3.3.12}) that

$$\left|\frac{r_{\varepsilon}(\theta)}{\varepsilon}-\frac{1}{a(0)}\right|=\left|\frac{1}{a(r_{\varepsilon},\theta)}-\frac{1}{a(0)}\right|\leq\frac{\left|a(0)-a(r_{\varepsilon}(\theta),\theta)\right|}{a(r_{\varepsilon}(\theta),\theta)\,a(0)}\leq\frac{A_1\,varepsilon}{a(r_\varepsilon(\theta),\theta)\,a(0)\underline{A}},$$
from which it follows by multiplication by $\varepsilon$ that 
\begin{equation}
\label{eqn 3.3.13}
\big|r_{\varepsilon}(\theta)-(\varepsilon\big/a(0)\big)\big|\leq \big(A_1\big/\underline{A}\,\big)\,\big(\varepsilon^2\big/a(r_\varepsilon(\theta),\theta)\,a(0)\big)
\end{equation}
for all $\theta$ and $0<\varepsilon<\varepsilon_0$. The exterior normal $\boldsymbol{\nu}_\varepsilon$ to the arc $\gamma_\varepsilon:=\{r(p)a(p)=\varepsilon\}$ at any point $p_\varepsilon\in\gamma_\varepsilon$ is given by
$$\boldsymbol{\nu}_\varepsilon=\frac{\nabla(r\,a)}{|\nabla(ra)|}=\frac{\boldsymbol{\nu}_0+(r/a)\nabla\,a}{|\boldsymbol{\nu}_0+(r/a)\nabla a|},$$
where $\boldsymbol{\nu}_0=(p/r)$ and $\nabla(r\,a)$ is evaluated at $p_\varepsilon$. It is easily follows that ${\rm tan}(\phi)\leq 2(A_1/\underline{A})\,r$ for $2 r\,A_1<a(p)$, where $\phi$ is the angle between $\boldsymbol{\nu}_0$ and $\boldsymbol{\nu}_\varepsilon$. Therefore, we have
\begin{equation}
\label{eqn 3.3.14}
|r'_\varepsilon(\theta)|\leq 2\,(A_1/\underline{A}\,)\,r_\varepsilon^2(\theta)\leq 2\,(A_1/\underline{A}^3)\,\varepsilon^2
\end{equation}
\begin{equation}
\label{eqn 3.3.15}
\left|r_\varepsilon(\theta)\,\sqrt{1+\big(r_\varepsilon'(\theta)\big/r_\varepsilon(\theta))^2}-(\varepsilon/a(0))\right|
\end{equation}
$$\leq r_\varepsilon(\theta)\left(\sqrt{1+(r_\varepsilon'(\theta)/r_\varepsilon(\theta))^2}-1\right)+\left|r_\varepsilon(\theta)-(\varepsilon/a(0)\right)|$$
$$\leq r_\varepsilon(\theta)(r_\varepsilon'(\theta)/r_\varepsilon(\theta))^2+\left|r_\varepsilon(\theta)-(\varepsilon/a(0)\right)|$$
$$\leq (\varepsilon/\underline{A}\,)(2A_1\,\varepsilon/\underline{A}^2)^2+(A_1\,\varepsilon^2/\underline{A}^3)\leq C_3\,\varepsilon^2.$$
The first assertion (\ref{eqn 3.3.10}) now follows from (\ref{eqn 3.3.15}), in view of the identity:
$$
||\gamma_\varepsilon(J)||=\int_\alpha^\beta r_\varepsilon(\theta)\sqrt{1+(r_\varepsilon'(\theta)/r_\varepsilon(\theta))^2}\,d\theta
$$
\vspace{.1in}

\noindent
We now turn to the total curvature estimate (\ref{eqn 3.3.11}). For any specified point $p_\varepsilon\in\gamma_\varepsilon$, one can choose new Cartesian coordinates $(x,y)$, with the origin unchanged, such that $p_\varepsilon=(x_\varepsilon,y_\varepsilon)$ and $\gamma_\varepsilon$ is locally (near $p_\varepsilon$) the graph of the function $y=y_\varepsilon(x)$ such that $p_\varepsilon=(x_\varepsilon,y_\varepsilon(x_\varepsilon))$ and $y'(x_\varepsilon)=0$. By twice differentiating the equation: $r(x,y_\varepsilon(x))\,a(x,y_\varepsilon(x))=\varepsilon$, and using the identities: $r_x=x/r$, $r_y=y/r$, $r_{xx}=y^2/r^3$, and $y'(x_\varepsilon)=0$, one sees that
\begin{equation}
\label{eqn 3.3.16}
y'_\varepsilon(x)=r_x\,a+r\,a_x=(x/r)\,a+a_x\,r=0,
\end{equation}
\begin{equation}
\label{eqn 3.3.17}
-y''_\varepsilon(x)=\frac{r_{xx}\,a+2\,r_x\,a_x+a_{xx}\,r}{r_y\,a+a_y r}=\frac{(y^2/r^3)\,a+2\,(x/r)\,a_x+r\,a_{xx}}{(y/r)\,a+a_y\,r}
\end{equation}
$$=\frac{1}{r}\,\frac{(1-(x/r)^2)\,a+2\,r\,(x/r)\,a_x+r^2\,a_{xx}}{a\,\sqrt{1-(x/r)^2}+a_y\,r},$$
\begin{equation}
\label{eqn 3.3.18}
(|x|/r)\leq(|\nabla a|/a)\,r\leq(A_1/\underline{A}\,)\,r\leq(A_1/\underline{A}^2)\,\varepsilon.
\end{equation}
for $x=x_\varepsilon$, $y=y_\varepsilon$. By using (\ref{eqn 3.3.18}) (which follows from (\ref{eqn 3.3.16})), one can show that 
\begin{equation}
\label{eqn 3.3.19}
a\,\sqrt{1-(x/r)^2}+a_y\,r\geq ((\sqrt{3}-1)/2)\,a\geq (1/3)\,a
\end{equation}
for $0<\varepsilon\leq\varepsilon_0$. It follows from (\ref{eqn 3.3.17}), (\ref{eqn 3.3.18}), and (\ref{eqn 3.3.19}) that 
\begin{equation}
\label{eqn 3.3.19a}
\Big|r\,y_\varepsilon''(x)+\frac{a}{a\,\sqrt{1-(x/r)^2}+a_y\,r}\Big|\leq \frac{3}{a}\,\,\Big(a(x/r)^2+2\,|a_x|\,|x|\,|x/r|+|a_{xx}|\,r^2\Big),
\end{equation} for $x=x_\varepsilon$, $y=y_\varepsilon$, and $0<\varepsilon<\varepsilon_0$. We also have, due to (\ref{eqn 3.3.19}), that 
\begin{equation}\label{eqn 3.3.19b}
\Big|\frac{a}{\sqrt{1-(x/r)^2}+a_y\,r}-1\Big|\leq\frac{3}{a}\,\Big((x/r)^2+(|a_y|/a)\,r\Big)
\end{equation}
for $x=x_\varepsilon$, $y=y_\varepsilon$, and $0<\varepsilon< \varepsilon_0$, and it follows by combining (\ref{eqn 3.3.19a}) and (\ref{eqn 3.3.19b}) and estimating the terms that,
for the signed curvature $\kappa_\varepsilon(p)$ of $\gamma_\varepsilon$, we have:
\begin{equation}
\label{eqn 3.3.20}
|r\,\kappa_\varepsilon(p)-1|\leq C_5\,\varepsilon,
\end{equation}
uniformly for $p\in\gamma_\varepsilon$ and $0<\varepsilon<\varepsilon_0$, where 
\begin{equation}
\label{eqn 3.3.20a} C_5:=\big(3/\underline{A}^6\big)\Big(\,\underline{A}^3\,A_1+\Big[\,\underline{A}\,A_1^2+\underline{A}^2 A_1^2+\underline{A}^3A_2\Big]
\varepsilon_0+2\,A_1^3\,\varepsilon_0^2\Big).
\end{equation}
It follows that 
\begin{equation}
 \label{eqn 3.3.21}
\left|r_\varepsilon\kappa_\varepsilon\sqrt{1+(r_\varepsilon'/r_\varepsilon)^2}-1\right|
\leq r_\varepsilon|\,\kappa_\varepsilon|\left(\sqrt{1+(r_\varepsilon'/r_\varepsilon)^2}-1\right)+\left|r_\varepsilon\kappa_\varepsilon-1\right|
\end{equation}
$$\leq r_\varepsilon\,|\kappa_\varepsilon|\,(r'_\varepsilon/r_\varepsilon)^2+\left|r_\varepsilon\kappa_\varepsilon-1\right| 
\leq (1+C_5\,\varepsilon)(r_\varepsilon'/r_\varepsilon)^2+C_5\,\varepsilon\leq C_4\,\varepsilon$$
for $0<\varepsilon<\varepsilon_0$, where we set
\begin{equation}
\label{eqn 3.3.21a}
C_4:=\left((1+C_5\,\varepsilon_0)(2A_1/\underline{A}^2)^2\,\varepsilon_0+C_5\right).
\end{equation}
The assertion (\ref{eqn 3.3.11}) follows from (\ref{eqn 3.3.21}), in view of the fact that
$$K(\gamma_\varepsilon(J))=\int_\alpha^\beta r_\varepsilon(\theta)\kappa_\varepsilon\,(r_\varepsilon(\theta),\theta)\,\sqrt{1+(r'_\varepsilon(\theta)/r_\varepsilon(\theta))^2}\,d\theta.$$

\begin{lemma} {(Turning-angle estimate)}
\label{lem 3.4.1} 
For given $\varepsilon\in(0,\varepsilon_0)$ and for any given $P$-periodic polygonal arc $\Gamma\in{\rm X}$, let $\Gamma_{\varepsilon,i}:=\Psi_{\varepsilon,i}(\Gamma),$ $i=1,2,$ (see Def. \ref{def 3.2.1}). Then the following assertions hold:  
\vspace{.1in}

\noindent
(a) The arcs ${\Gamma}_{\varepsilon,i}$ are piecewise-smooth, $P$-periodic, simple arcs. In fact for $i=1,2$,  ${\Gamma}_{\varepsilon,i}$ is a $C^2$-arc except at a collection of "vertices"  located in the set $C_{\varepsilon,i}\subset{\Gamma}_{\varepsilon,i}$ (see Defs. \ref{def 3.2.1} and \ref{def 3.2.3}, and Thm. \ref{thm 3.2.1}), where the set $C_{\varepsilon,i}$ contains finitely-many points relative to any $P$-period (in $x$) of the arc $\Gamma_{\varepsilon,i}$ 
\vspace{.1in}

\noindent
(b) Given a particular point ${p}_{\varepsilon,i}\in C_{\varepsilon,i}\subset{\Gamma}_{\varepsilon,i}=\Psi_{\varepsilon,i}(\Gamma)$, in terms of the natural ordering in $\Gamma$, let $q_{\varepsilon,i}^-$ (resp. $q_{\varepsilon,i}^+$) denote the maximal (resp. minimal) point in $\Gamma$ such that $q_{\varepsilon,i}^-\leq q\leq q_{\varepsilon,i}^+$ for all points $q\in\Pi_{\varepsilon,i}({p}_{\varepsilon,i})\subset\Gamma$. Then the turning angle ${\mathfrak{A}}_{\varepsilon,i}({p}_{\varepsilon,i})$ of ${\Gamma}_{\varepsilon,i}$ at the point ${p}_{\varepsilon,i}\in C_{\varepsilon,i}\subset\Gamma_{\varepsilon,i}$ is such that: 
\begin{equation}
\big|\mathfrak{A}_{\varepsilon,i}(p_{\varepsilon,i})\big|\,\leq\big(1+C_6\varepsilon\big)\,\phi\big(\hat{\bm\nu}^+_{\varepsilon,i},\hat{\bm\nu}^-_{\varepsilon,i}\big),
\end{equation}
where $C_6:=\big(6A_1\big/\underline{A}^2\big)$, $\hat{\bm\nu}_{\varepsilon,i}^\pm:=\big(p_{\varepsilon,i}-q^\pm_{\varepsilon,i})\big/|p_{\varepsilon,i}-q_{\varepsilon,i}^\pm|\big)$, and $\phi(\boldsymbol{\nu}_1,\boldsymbol{\nu}_2)$ denotes the counter-clockwise turning angle from the unit vector $\boldsymbol{\nu}_1$ to the unit vector $\boldsymbol{\nu}_2$. 
\end{lemma}
\vspace{.1in}

\noindent
{\bf Proof.} Concerning Part (a), we observe that ${\Gamma}_{\varepsilon,i}$ is the double-point-free boundary of a  closed, simply-connected, $P$-periodic (in $x$) domain (see Lems. \ref{lem 2.1.4} and \ref{lem 2.1.4*}). The smoothness of the arc ${\Gamma}_{\varepsilon,i}$ relative to $\Re^2\setminus C_{\varepsilon,i}$ follows from Lems. \ref{lem 3.3.1} and \ref{lem 3.3.2}. At this point, the assertion follows by Thm. \ref{thm 3.2.1}.
\vspace{.1in}

\noindent
We turn to the proof of Part (b): In terms of the natural ordering of points in the arc ${\Gamma}_{\varepsilon,i}$ (for fixed $\varepsilon\in(0,\varepsilon_0)$ and fixed $i\in\{1,2\}$; see Def. \ref{def 3.2.2}), the point ${p}_{\varepsilon,i}$ is the terminal (resp. initial) end-point of an arc-segment ${\gamma}_{\varepsilon,i}^-$ (resp. ${\gamma}_{\varepsilon,i}^+$) in ${\Gamma}_{\varepsilon,i}$. Here, we can assume ${\gamma}_{\varepsilon,i}^+=B(q_{\varepsilon,i}^+)\subset\Gamma_{\varepsilon,i}$ if $q_{\varepsilon,i}^+\in V(\Gamma_i)$ and that ${\gamma}_{\varepsilon,i}^+=A(L_{\varepsilon,i}^+)\subset\Gamma_{\varepsilon,i}$ if $q_{\varepsilon,i}^+\in L_{\varepsilon,i}^+\in S(\Gamma_i)$. Similarly, we can assume ${\gamma}_{\varepsilon,i}^-=B(q_{\varepsilon,i}^-)\subset\Gamma_{\varepsilon,i}$ if $q_{\varepsilon,i}^-\in V(\Gamma_i)$, while ${\gamma}_{\varepsilon,i}^-=A(L_{\varepsilon,i}^-)\subset\Gamma_{\varepsilon,i}$ if $q_{\varepsilon,i}^-\in L_{\varepsilon,i}^-\in S(\Gamma_i)$. Then ${\gamma}_{\varepsilon,i}^\pm$ is a portion of the level curve at altitude $\varepsilon$ of the function $r_{\varepsilon,i}^\pm(p)\,a_i(p)$, where $r_{\varepsilon,i}^\pm(p):=|p-q_{\varepsilon,i}^\pm|$ if $q_{\varepsilon,i}^\pm\in V(\Gamma_i)$ and $r_
{\varepsilon,i}^\pm(p):={\rm dist}(p,L_{\varepsilon,i}^\pm)$ if $q_{\varepsilon,i}^\pm
\in L_{\varepsilon,i}^\pm\in S(\Gamma_i)$. The upper normal $\boldsymbol{\nu}_{\varepsilon,i}^\pm$ to the arc ${\gamma}_{\varepsilon,i}^\pm\subset\Gamma$ at the point $p_{\varepsilon,i}^\pm\in\gamma_{\varepsilon,i}^\pm$ is given by
\begin{equation}
\label{eqn 3.4.3}
\boldsymbol{\nu}_{\varepsilon,i}^\pm
:=\frac{\nabla(r_{\varepsilon,i}^\pm\, a_i)}{|\nabla(r_{\varepsilon,i}^\pm\,a_i)|}
=\frac{\hat{\boldsymbol{\nu}}_{\varepsilon,i}^\pm\,+\,\boldsymbol{h}_{\varepsilon,i}}{|\hat{\boldsymbol{\nu}}_{\varepsilon,i}^\pm\,+\,\boldsymbol{h}_{\varepsilon,i}|},
\end{equation}
where
\begin{equation}
\label{eqn 3.4.4}
\hat{\boldsymbol{\nu}}_{\varepsilon,i}^\pm:=\big(({p}_{\varepsilon,i}-q_{\varepsilon,i}^\pm)/|{p}_{\varepsilon,i}-q_{\varepsilon,i}^\pm|\big)=\nabla r_{\varepsilon,i}^\pm({p}_{\varepsilon,i}),
\end{equation}
\begin{equation}
\label{eqn 3.4.5}
\boldsymbol{h}_{\varepsilon,i}:=\big(r_{\varepsilon,i}^\pm({p}_{\varepsilon,i})\nabla\,a_i({p}_{\varepsilon,i})/a_i({p}_{\varepsilon,i})\big)=\big(\varepsilon\,\nabla a_i({p}_{\varepsilon,i})\big/a_i^2({p}_{\varepsilon,i})\big).
\end{equation}
One can verify that $|\hat{\boldsymbol{\nu}}_{\varepsilon,i}^\pm|=1$ and $|\boldsymbol{h}_{\varepsilon,i}|\leq \big(\varepsilon\,|\nabla\,a_i({p}_{\varepsilon,i})|/a_i^2({p}_{\varepsilon,i})\big)\leq \big(A_1\big/\underline{A}^2\big)\,\varepsilon$. It follows from (\ref{eqn 3.4.3}), (\ref{eqn 3.4.4}), and (\ref{eqn 3.4.5}) that the equation:
\begin{equation}
\label{eqn 3.4.6}
{\rm tan}\big(\theta_{\varepsilon,i}^\pm-\hat{\theta}_{\varepsilon,i}^\pm\big)=\frac{|\boldsymbol{h}_{\varepsilon,i}|\,{\rm sin}\big(\psi_{\varepsilon,i}-\hat{\theta}_{\varepsilon,i}^\pm\big)}{1+|\boldsymbol{h}_{\varepsilon,i}|\,{\rm cos}\big(\psi_{\varepsilon,i}-\hat{\theta}_{\varepsilon,i}^\pm\big)}
\end{equation}
is satisfied by setting $\theta_{\varepsilon,i}^\pm:={\rm arg}(\boldsymbol{\nu}_{\varepsilon,i}^\pm)$, corresponding to $\hat{\theta}_{\varepsilon,i}^\pm:={\rm arg}(\hat{\boldsymbol{\nu}}_{\varepsilon,i}^\pm)$, and $\psi_{\varepsilon,i}:={\rm arg}(\boldsymbol{h}_{\varepsilon,i})$. It follows from (\ref{eqn 3.4.6}) by differentiation (holding $\psi_{\varepsilon,i}$ fixed) that
\begin{equation}
\label{eqn 3.4.7}
\frac{\partial\theta_{\varepsilon,i}^{\pm}}{\partial\hat{\theta}_{\varepsilon,i}^{\pm}}=1-\frac{{\rm cos}^2\big(\theta_{\varepsilon,i}^\pm-\hat{\theta}_{\varepsilon,i}^\pm\big)\big({\rm cos}\big(\psi_{\varepsilon,i}-\hat{\theta}_{\varepsilon,i}^\pm\big)+|\boldsymbol{h}_{\varepsilon,i}|\big)\,|\boldsymbol{h}_{\varepsilon,i}|}{\big(1+|\boldsymbol{h}_{\varepsilon,i}|{\rm cos}\big(\psi_{\varepsilon,i}-\hat{\theta}_{\varepsilon,i}^\pm\big)\big)^2}.
\end{equation}
It follows from (\ref{eqn 3.4.7}) that
$$\left|\frac{\partial\theta_{\varepsilon,i}^\pm}{\partial\hat{\theta}_{\varepsilon,i}^\pm}\right|\leq 1+\frac{\big(1+|\boldsymbol{h}_{\varepsilon,i}|\big)|\,\boldsymbol{h}_{\varepsilon,i}|}{\big(1-|\boldsymbol{h}_{\varepsilon,i}|\big)^2}\leq \big(1+6\,|\boldsymbol{h}_{\varepsilon,i}|\big)$$
for $|\boldsymbol{h}_{\varepsilon,i}|\leq 1/2$, independent of $\theta_{\varepsilon,i}^\pm,\hat{\theta}_{\varepsilon,i}^\pm$, from which it follows by the theorem of the mean that 
\begin{equation}
\label{eqn 3.4.8}
|\theta_{\varepsilon,i}^+-\theta_{\varepsilon,i}^-|\leq \big(1+6\,|\boldsymbol{h}_{\varepsilon,i}|\big)|\hat{\theta}_{\varepsilon,i}^+-\hat{\theta}_{\varepsilon,i}^-|\leq \big(1+(6\,A_1\big/\underline{A}^2)\,\varepsilon\,\big)|\hat{\theta}_{\varepsilon,i}^+-\hat{\theta}_{\varepsilon,i}^-|
\end{equation}
for $|\boldsymbol{h}_{\varepsilon,i}|\leq(A_1/\underline{A}^2)\,\varepsilon\leq(1/2)$, completing the proof of Part (b). We omit the straight-forward proof of Part (c). 

\subsection{Proofs of main estimates}
\label{subsection 5.3}
\begin{lemma} {(Main estimates for the operators $\Phi_{\varepsilon,i}$)}
\label{lem 3.5.1}
Let the periodic arc-pair ${\bf\Gamma}=(\Gamma_1,\Gamma_2)\in {\bf X}$ be such that $||\Gamma_i||, K(\Gamma_i)<\infty$, $i=1,2$, where $||\cdot||$ (resp. $K(\cdot)$) refers to arc-length (total curvature) relative to one $P$-period (in $x$). Then:
\begin{equation}
\label{eqn 3.5.1}
||\Gamma_{\varepsilon,i}||\leq(1-(\varepsilon/h))||\Gamma_i||+(\varepsilon/h)||\Gamma_{h,i}||=||\Gamma_i||+\big(||\Gamma_{h,i}||-||\Gamma_i||\big)(\varepsilon/h),
\end{equation}
\begin{equation}
\label{eqn 3.5.2}
K(\Gamma_{\varepsilon,i})\leq(1-(\varepsilon/h))K(\Gamma_i)+(\varepsilon/h) K(\Gamma_{h,i})
\end{equation}
$$=K(\Gamma_i)+\big(K(\Gamma_{h,i})-K(\Gamma_i )\big)(\varepsilon/h),$$
for $i=1,2$, and for all $0<\varepsilon\leq h\leq 1$, where we set $\Gamma_{\varepsilon,i}:=\Phi_{\varepsilon,i}({\bf\Gamma})$ and $\Gamma_{h,i}:=\Phi_{h,i}({\bf\Gamma})$.
\end{lemma}
-
\noindent
{\bf Proof.}
Given ${\bf\Gamma}\in {\bf X}$, let $\Omega=\Omega({\bf\Gamma})$. For $j=\sqrt{-1}$, let $w=F(x+jy)$ be a continuous, periodic mapping of the strip $\Re\times[0,1]$ onto ${\rm Cl}(\Omega),$ whose restriction to $\Re\times(0,1)$ is a conformal mapping onto $\Omega$. For the proof of (\ref{eqn 3.5.2}), we assume w.l.o.g. that $\Omega$ is sufficiently regular, so that the harmonic function $f(x,y):={\rm arg}(F'(x+jy)):\Re\times(0,1)\rightarrow\Re$ extends continuously to $\Re\times[0,1]$. The total curvature of the level curve of $U({\bf\Gamma};p)$ corresponding to the line $L_y:=[0,1]\times\{y\}$ is given by $k(y):=K(\Gamma_y)=\int_0^L g(x,y)\,dx$ for $y\in[0,1]$ where $g(x,y):=|f_x(x,y)|$ is sub-harmonic, and where $L$ denotes the period of $F$. Thus
$$k''(y)=\int_0^L g_{yy}(x,y)\,dx=\int_0^L (\Delta\,g-g_{xx})\,dx\quad$$
$$\quad \geq-\int_0^L g_{xx}\,dx=g_x(0,y)-g_x(L,y)=0$$
for all $y\in(0,1)$, where $\Delta$ denotes the Laplace operator. The assertion (\ref{eqn 3.5.2}) follows directly from this. The assertion (\ref{eqn 3.5.1}) follows by the same argument, but with $g(x,y):=|F'(x+jy)|$ (which is again a sub-harmonic function).
\begin{lemma} {(Main estimates for the operators $\Psi_{\varepsilon,i}$)}
\label{lem 3.5.2}
There exist uniform constants $A,B,C,D$, and $E$, depending only on the uniform constants $\underline{A}$, $\overline{A}$, $A_1$, and $A_2$, such that
\begin{equation}
\label{eqn 3.5.3}
K\big(\Psi_\varepsilon(\Gamma)\big)\leq K(\Gamma)+\Big[A\,\,K(\Gamma)+B\,||\Gamma||\Big]\,\varepsilon,
\end{equation}
\begin{equation}
\label{eqn 3.5.4}
||\Psi_\varepsilon(\Gamma)||\leq ||\Gamma||+C\,\,K(\Gamma)\,\varepsilon+\Big[D\,||\Gamma||+E\,\,K(\Gamma)\Big]\,\varepsilon^2,
\end{equation}
both uniformly for all $\varepsilon\in(0,\varepsilon_0)$ and all $P$-periodic (in $x$) arcs $\Gamma\in\tilde{{\rm X}}_\varepsilon$. Here, for any $\varepsilon\in(0,\varepsilon_0)$, the mapping: $\Psi_{\varepsilon}(\Gamma):{\rm X}\rightarrow{\rm X}$ is defined by any one of the following rules: Either $\Psi_{\varepsilon}^-(\Gamma)=\partial\,G_{\varepsilon,1}(\Gamma)$, or $\Psi_{\varepsilon}^+(\Gamma)=\partial\,G_{\varepsilon,2}(\Gamma)$, or $\hat{\Psi}_{\varepsilon}^-(\Gamma)=\partial\,H_{\varepsilon,1}(\Gamma)$, or else $\hat{\Psi}_{\varepsilon}^+(\Gamma)=\partial\,H_{\varepsilon,2}(\Gamma)$ for all $\Gamma\in{\rm X}$ (see Def. \ref{def 3.2.1}).
\end{lemma}

\noindent
{\bf Proof.} For any fixed $\varepsilon\in(0,\varepsilon_0)$, we use $\Psi_{\varepsilon}(\Gamma):{\rm X}\rightarrow{\rm X}$ to denote a mapping defined by either one of the following rules: Either $\Psi_{\varepsilon}^-(\Gamma)=\partial\,G_{\varepsilon,1}(\Gamma)$ or $\Psi_{\varepsilon}^+(\Gamma)=\partial\,G_{\varepsilon,2}(\Gamma)$. For any given $P$-periodic (in $x$) polygonal arc $\Gamma\in{\rm X}$, we define the arc $\Gamma_\varepsilon:=\Psi_{\varepsilon}(\Gamma)\in{\rm X}$. Then $\Gamma_{\varepsilon}$ can be partitioned according to Thm. \ref{thm 3.2.1} (Eq. (\ref{eqn 3.2.2})). Concerning the arc-length estimate, Eq. (\ref{eqn 3.5.3}), the length of one $P$-period (in $x$) of the arc ${\Gamma}_{\varepsilon}$ is the sum of the lengths of the disjoint connected components: $A_{\varepsilon}(L)$ and $B_{\varepsilon}(q)$ corresponding distinct sides $L\in S(\Gamma)$ and vertices $q\in V(\Gamma)$, all restricted to one $P$-period (in $x$) of the arc $\Gamma$. In view of Lem. \ref{lem 3.3.1}, Eq. (\ref{eqn 3.3.2}) and Lem. \ref{lem 3.3.2}, Eq. (\ref{eqn 3.3.10}), it follows that 
\begin{equation}
\label{eqn 3.5.5}
||\Gamma_\varepsilon||=||\Psi_\varepsilon(\Gamma)||\leq\big(1+(C_1^{\,2}/2)\,\varepsilon^2\big)\,\sum ||L||+\big((\varepsilon/\underline{A}\,)+C_3\,\varepsilon^2\big)\,\sum\,\mathfrak{A}(q)
\end{equation}
$$\leq \big(1+(C_1^{\,2}/2)\,\varepsilon^2\big)||\Gamma||+\big((\varepsilon/\underline{A}\,)+C_3\,\varepsilon^2\big)\,K(\Gamma)$$
$$\leq||\Gamma||+\,C\,K(\Gamma)\,\varepsilon+\big(D\,||\Gamma||+\,E\,K(\Gamma)\big)\,\varepsilon^2,$$
where, in (\ref{eqn 3.5.5}), the first sum is over all sides $L\in S(\Gamma)$ and the second sum is over all vertices $q\in V(\Gamma)$, both subject to the restriction to one $P$-period (in $x$) of $\Gamma$. Here, we use $\mathfrak{A}(q)$ to denote the absolute turning angle of $\Gamma$ at a vertex $q\in\Gamma$. Thus, the estimate (\ref{eqn 3.5.4}) holds for any $\varepsilon\in(0,\varepsilon_0)$ and any $P$-periodic polygonal arc $\Gamma\in{\rm X}$, where we set $C:=(1/\underline{A})$, $D:=C_3$, and $E:=C_1^{\,2}/2$.
\vspace{.1in}

\noindent
Turning now to the proof of total-curvature estimate (\ref{eqn 3.5.4}), we begin again with a fixed, but arbitrary, value $\varepsilon\in(0,\varepsilon_0)$ and a fixed, but arbitrary, $P$-periodic (in $x$), double-point-free, polygonal arc $\Gamma\in{\rm X}$, and we partition the corresponding $P$-periodic arc $\Gamma_{\varepsilon}:=\Psi_{\varepsilon}(\Gamma)\in{\rm X}$ as in (\ref{eqn 3.2.2}). Then the total curvature of one $P$-period (in $x$) of the arc ${\Gamma}_{\varepsilon}$ equals the sum of the total curvatures of all the arcs $A_{\varepsilon}(L), B_{\varepsilon}(q)\subset\Gamma_{\varepsilon}$ corresponding to all the sides $L\in S(\Gamma)$ and vertices $q\in V(\Gamma)$ located in any one $P$-period (in $x$) of $\Gamma$, together with the sum of the absolute turning angles ${\mathfrak A}_{\varepsilon}(p_{\varepsilon})$ of ${\Gamma}_{\varepsilon}$ at all the points $p_{\varepsilon}$ located in any one $P$-period (in $x$) of $\Gamma_{\varepsilon}$ such that $p_{\varepsilon,i}\in C_{\varepsilon}\subset{\Gamma}_{\varepsilon}$ (see Thm. \ref{thm 3.2.1}).
\vspace{.1in}

\noindent
We let $\big(p_{\varepsilon,i}\big)_{i\in Z}$ denote the sequence of all points $p_{\varepsilon,i}\in C_
{\varepsilon}$, indexed such that $p_{\varepsilon,i}<p_{\varepsilon,i+1}$ for all $i\in i\in Z:=\{0,\pm 1,\pm 2, \pm 3,\cdots\}$ (in terms of the natural ordering in $\Gamma_{\varepsilon}$). For each $i\in Z$, we define the three arc-segments \begin{equation}
\label{3.5.6}
\gamma_{\varepsilon,i}=\big\{p\in\Gamma_{\varepsilon}: p_{\varepsilon,i-1}<p<p_{\varepsilon,i}\big\};\,\,\, \ddot{\gamma}_{\varepsilon,i}=\big\{q\in\Gamma:q_{\varepsilon,i}^-\leq q\leq q_{\varepsilon,i}^+\big\},
\end{equation} 
$$\dot{\gamma}_{\varepsilon,i}=\big\{q\in\Gamma: q_{\varepsilon,i-1}^+\,<\,q\,<\,q_{\varepsilon,i}^-\big\},$$ where the points $q_{\varepsilon,i}^\pm\in\Gamma$ are chosen such that $q_{\varepsilon,i}^-$ is maximal and $q_{\varepsilon,i}^+$ is minimal, both in terms of the natural ordering in $\Gamma$ and subject to the requirements that $q_{\varepsilon,i}^-\leq \Pi_{\varepsilon}\big(p_{\varepsilon,i}\big)\leq q_{\varepsilon,i}^+$. Obviously the points and arc-segments $p_{\varepsilon,i}$ and $\gamma_{\varepsilon,i}$, $i\in Z$ (resp. the arc-segments $\dot{\gamma}_{\varepsilon,i}$ and $\ddot{\gamma}_{\varepsilon,i}$, $i\in Z$) constitute a partition of the arc $\Gamma_{\varepsilon}$ (resp. $\Gamma$) such that, in terms of the natural ordering in $\Gamma_{\varepsilon}$ (resp. $\Gamma$), we have  $\gamma_{\varepsilon,i}<p_{\varepsilon,i}<\gamma_{\varepsilon,i+1}$ for all $i\in Z$ (resp. $\dot{\gamma}_{\varepsilon,i}<\ddot{\gamma}_{\varepsilon,i}<\dot{\gamma}_{\varepsilon,i+1}$ for all $i\in Z$). 
Due to the $P$-periodicity (in $x$) of the arc $\Gamma_{\varepsilon}\in{\rm X}$, there exists a natural number $n=n(\varepsilon)\in N$ such that $p_{\varepsilon,i+n}=p_{\varepsilon,i}+(P,0)$ for all $i\in Z$, where the same sequence of equations also holds with $p_{\varepsilon,i}$ replaced by the points
$q_{\varepsilon,i}^\pm$ or the arc-segments ${\gamma}_{\varepsilon,i}$, $\dot{\gamma}_{\varepsilon,i}$, or $\ddot{\gamma}_{\varepsilon,i}$, etc. Finally, we have
\begin{equation}
\label{eqn 3.5.7}
\Pi_{\varepsilon}\big(p_{\varepsilon,i}\big)\subset\ddot{\gamma}_{\varepsilon,i}\,\,\,{\rm and}\,\,\,\Pi_{\varepsilon}\big(\gamma_{\varepsilon,i}\big)=\dot{\gamma}_{\varepsilon,i},
\end{equation}
due to Def. \ref{def 3.2.3} and Lem. \ref{lem 3.2.1}(b). In view of (\ref{eqn 3.5.7}), it follows from  Thm. \ref{thm 3.2.1} (Eq. (\ref{eqn 3.2.2})) that for each $i\in Z$ the arc-segment $\gamma_{\varepsilon,i}$ is the disjoint union of all the arc-segments $A_{\varepsilon}(L)$ and $B_{\varepsilon}(q)$ corresponding to $L\in S(\Gamma)$ and $q\in V(\Gamma)$ such that $L, \{q\}\subset\dot{\gamma}_{\varepsilon,i}$ (which constitute a finite partition of $\dot{\gamma}_{\varepsilon,i}$). Since each arc-segment $\gamma_{\varepsilon,i}$ has no vertices, the total curvature of $\gamma_{\varepsilon,i}$ must be equal to the sum of the total curvatures of the arcs in this partition of $\gamma_{\varepsilon,i}$. It follows by applying the estimates Lem. \ref{lem 3.3.1}, Eq. (\ref{eqn 3.3.3}) and Lem. \ref{lem 3.3.2}, Eq. (\ref{eqn 3.3.11}) to the individual arc-segments in this partition and $\gamma_{\varepsilon,i}$ and summing the terms that
\begin{equation}
\label{eqn 3.5.8}
K\big(\gamma_{\varepsilon,i}\big)\,\leq\,\big(1+C_4\,\varepsilon\big)\,K\big(\dot{\gamma}_{\varepsilon,i}\big)\,+\,\overline{C}_2\,||\dot{\gamma}_{\varepsilon,i}||\,\varepsilon.
\end{equation}

By Lem. \ref{lem 3.4.1}, the absolute turning angle ${\mathfrak{A}}_{\varepsilon}(p_{\varepsilon,i})$ of the arc $\Gamma_{\varepsilon}$ at the vertex $p_{\varepsilon,i}$ is estimated by:
\begin{equation}
\label{eqn 3.5.9}
{\mathfrak{A}}_{\varepsilon}(p_{\varepsilon,i})\leq (1+C_6\varepsilon)\,\phi(\hat{\boldsymbol{\nu}}_{\varepsilon,i}^-,\hat{\boldsymbol{\nu}}_{\varepsilon,i}^+)
\end{equation}
for each $i\in Z$, where $C_6=6(A_1/\underline{A}^2)$, $\hat{\boldsymbol\nu}^\pm_{\varepsilon,i}:=\big((p_{\varepsilon,i}-q_{\varepsilon,i}^\pm)\big/|p_{\varepsilon,i}-q_{\varepsilon,i}^\pm|\big)$, and where 
$\phi(\boldsymbol{\nu}_1,\boldsymbol{\nu}_2)$ denotes the magnitude of the angle between the unit vectors $\boldsymbol{\nu}_1$ and $\boldsymbol{\nu}_2$. We also use $\tilde{\boldsymbol{\nu}}_{\varepsilon,i}^\pm$, $i\in Z$, to denote the upper normal to $\Gamma$ on the straight-line-segment $l_{\varepsilon,i}^\pm$ of $\ddot{\gamma}_{\varepsilon,i}$, chosen such that $q_{\varepsilon,i}^\pm$ is an endpoint of $l_{\varepsilon,i}^\pm$ (Thus $\tilde{\boldsymbol{\nu}}_{\varepsilon,i}^\pm=\hat{\boldsymbol{\nu}}_{\varepsilon,i}^\pm$ if $q_{\varepsilon,i}^\pm$ is not a vertex of $\Gamma$). Then it is easily seen that
\begin{equation}
\label{eqn 3.5.10}
\phi(\hat{\boldsymbol{\nu}}_{\varepsilon,i}^+,\hat{\boldsymbol{\nu}}_{\varepsilon,i}^-)\leq \phi(\hat{\boldsymbol{\nu}}_{\varepsilon,i}^-,\tilde{\boldsymbol{\nu}}_{\varepsilon,i}^-)+K(\ddot{\gamma}_{\varepsilon,i})+\phi(\tilde{\boldsymbol{\nu}}_{\varepsilon,i}^+,\hat{\boldsymbol{\nu}}_{\varepsilon,i}^+).
\end{equation}
It is also easily seen that $\phi(\hat{\boldsymbol{\nu}}_{\varepsilon,i}^+,\tilde{\boldsymbol{\nu}}_{\varepsilon,i}^+)=0$ (resp. $\phi(\tilde{\boldsymbol{\nu}}_{\varepsilon,i}^-,\hat{\boldsymbol{\nu}}_{\varepsilon,i}^-)=0$) if the point $q_{\varepsilon,i}^+$ (resp. $q_{\varepsilon,i}^-$) is not a vertex of the arc $\Gamma$. On the other hand, if one or both of the points $q_{\varepsilon,i}^\pm$ are vertices of the arc $\Gamma$, then \begin{equation}
\label{eqn 3.5.11}
0\leq\phi(\hat{\boldsymbol{\nu}}_{\varepsilon,i}^+,\tilde{\boldsymbol{\nu}}_{\varepsilon,i}^+)+\phi(\hat{\boldsymbol{\nu}}_{\varepsilon,i}^-, \tilde{\boldsymbol{\nu}}_{\varepsilon,i}^-)\leq 
{\mathfrak A}\big(q_{\varepsilon,i}^+, q_{\varepsilon,i}^-\big),
\end{equation}
where 
${\mathfrak A}\big(q_{\varepsilon,i}^\pm\big)$ denotes the absolute turning angle of the polygonal arc $\Gamma$ at either of the vertices $q_{\varepsilon,i}^\pm\in\Gamma$, and, in terms of this, we have ${\mathfrak A}\big(q_{\varepsilon,i}^+,\,q_{\varepsilon,i}^-\big):=
{\mathfrak A}\big(q_{\varepsilon,i}^+\big)+{\mathfrak A}\big(q_{\varepsilon,i}^-\big)$ whenever the vertices $q_{\varepsilon,i}^+,\,q_{\varepsilon,i}^-\in\Gamma$ are distinct, whereas ${\mathfrak A}\big(q_{\varepsilon,i}^+, q_{\varepsilon,i}^-\big):={\mathfrak A}\big(q_{\varepsilon,i}^+\big)={\mathfrak A}\big(q_{\varepsilon,i}^-)$ in the case where the points $q_{\varepsilon,i}^\pm\in\Gamma$ coincide. Here, we also assume that ${\mathfrak A}(q_{\varepsilon,i}^\pm)=0$ if the point $q_{\varepsilon,i}^\pm\in\Gamma$ is actually not a vertex of $\Gamma$.
\vspace{.1in}

\noindent
By combining Eqs. (\ref{eqn 3.5.8})-(\ref{eqn 3.5.11}) for each $i=1, 2,\cdots, n(\varepsilon)$, and summing the terms, we see that 
\begin{equation}
\label{eqn 3.5.12}
K({\Gamma}_\varepsilon)\leq\big(1+A\,\varepsilon\big)\sum_{i=1}^{n(\varepsilon)} \Big(K(\dot{\gamma}_{\varepsilon,i})\,+K(\ddot{\gamma}_{\varepsilon,i})+{\mathfrak A}\big(q_{\varepsilon,i}^+, q_{\varepsilon,i}^-\big)\Big)\,+\,B\,\varepsilon\,\sum_{i=1}^{n(\varepsilon)} ||\gamma_{\varepsilon,i}||
\end{equation}
$$\leq \big(1+A\,\varepsilon\big)\,K(\Gamma)\,+\,B\,||\Gamma||\,\varepsilon,$$
where $A:={\rm max}\{C_4,C_6\}$ and $B:=\overline{C}_2$. Finally, we conclude from (\ref{eqn 3.5.12}) that the inequality (\ref{eqn 3.5.3}) holds uniformly for all $\varepsilon\in(0,\varepsilon_0)$ and all $P$-periodic (in $x$) polygonal arcs $\Gamma\in{\rm X}$, where the constants $A$ and $B$ are independent of $\varepsilon\in(0,\varepsilon_0)$ and the particular choice of the polygonal arc $\Gamma\in{\rm X}$. 
\vspace{.1in}

\noindent
It remains to extend the estimates (\ref{eqn 3.5.3}) and (\ref{eqn 3.5.4}) from the polygonal case discussed above
to the case of all $P$-periodic (in $x$), double-point-free smooth arcs $\Gamma\in\tilde{\bf X}_\varepsilon$, where $\Gamma_{\varepsilon}:=\Psi_{\varepsilon}(\Gamma)$ denotes either one of the arcs $\partial\,G_{\varepsilon,i}(\Gamma)$, $i=1,2$. To this end, given any smooth arc $\Gamma\in\tilde{{\rm X}}_\varepsilon$, we can choose a sequence $\big(\Gamma_n\big)_{n=1}^{\infty}$ of $P$-periodic polygonal arcs in $\Gamma_n\in{\rm X}$ such that 
\vspace{.1in}

\noindent
(${\bf i}$): for each $n\in N$, the vertices of the arc $\Gamma_n$ all lie in the arc $\Gamma$, and are positively-ordered relative to the natural ordering in both $\Gamma$ and $\Gamma_n$.
\vspace{.1in}

\noindent
In view of (${\bf i}$), the sequence  $\big(\Gamma_n\big)_{n=1}^\infty$ consists of double-point-free arcs such that 
\begin{equation}
\label{eqn 3.5.13}
||\Gamma_n||\leq ||\Gamma||\,\,\,{\rm and}\,\,\,K(\Gamma_n)\leq K(\Gamma).
\end{equation}
In view of (\ref{eqn 3.5.13}), it follows from (\ref{eqn 3.5.5}) and (\ref{eqn 3.5.12}) that for any $\varepsilon\in(0,\varepsilon_0)$ we have
\begin{equation}
\label{eqn 3.5.15}
\limsup_{n\rightarrow\infty}||\Psi_{\varepsilon}(\Gamma_n)||\leq||\Gamma||\,+\,C\,K(\Gamma)\varepsilon\,+\,\big[D\,||\Gamma||\,+\,E\,K(\Gamma)\big]\varepsilon^2,
\end{equation}
\begin{equation}
\label{eqn 3.5.16}
\limsup_{n\rightarrow\infty} K\big(\Psi_{\varepsilon}(\Gamma_n)\big)\leq\big(1+A\,{\varepsilon}\big)\,K(\Gamma)\,+\,B\,||\Gamma||\,\varepsilon.
\end{equation}
Now, we can also choose the sequence $\big(\Gamma_n\big)_{n=1}^\infty$, subject to (${\bf i}$), such that $\Gamma_n\rightarrow\Gamma\in\tilde{{\rm X}}_{\varepsilon}$ as $n\rightarrow\infty$. Therefore, there exist arc-sequences $\big(\Gamma_{n}^\pm\big)_{n=1}^\infty$ in $\tilde{\rm X}_{\varepsilon}$ such that (${\bf ii}$): $\Gamma_{n}^-<\Gamma,\,\Gamma_n<\Gamma_n^+$ for all $n\in N$, and such that (${\bf iii}$):  $\Gamma_n^+\downarrow\Gamma$ and $\Gamma_n^-\uparrow\Gamma$, both as $n\rightarrow\infty$, where $\uparrow$ and $\downarrow$ refer to monotone and uniform convergence. In view of Lem. \ref{lem 2.1.1}(d)(e) and Lem. \ref{lem 2.1.5}(a)(c), it follows from (${\bf ii}$) and (${\bf iii}$) that (${\bf iv}$): $\Psi_{\varepsilon}^-(\Gamma_n^-)\leq\Psi_{\varepsilon}^-(\Gamma_n)\leq\Psi_{\varepsilon}^+(\Gamma_n)\leq\Psi_{\varepsilon}^+(\Gamma_n^+)$ for all $n\in N$, and that (${\bf v}$): $\Psi_{\varepsilon}^-(\Gamma_n^-)\uparrow\Psi_{\varepsilon}^-(\Gamma)$ and $\Psi_{\varepsilon}^+(\Gamma_n^+)\downarrow\Psi_{\varepsilon}^+(\Gamma)$), both as $n\rightarrow\infty$, and it follows from (${\bf iv}$) and (${\bf v}$) that $\Psi_{\varepsilon}^\pm(\Gamma_n)\rightarrow\Psi_{\varepsilon}^\pm(\Gamma)$, both as $n\rightarrow\infty$, from which it in turn follows that (${\bf vi}$): $||\Psi_{\varepsilon}(\Gamma)||\leq\limsup_{n\rightarrow\infty}||\Psi_{\varepsilon}(\Gamma_n)||$ and (${\bf vii}$): $K\big(\Psi_{\varepsilon}(\Gamma)\big)\leq\limsup_{n\rightarrow\infty}K\big(\Psi_{\varepsilon}(\Gamma_n)\big)$, it follows from (\ref{eqn 3.5.15}) and (\ref{eqn 3.5.16}) that the assertions (\ref{eqn 3.5.3}) and (\ref{eqn 3.5.4}) hold in the cases where $\Psi_{\varepsilon}(\Gamma)=\Psi_{\varepsilon}^\pm(\Gamma)$, $\Gamma\in\tilde{{\rm X}}_{\varepsilon}$. 
\vspace{.1in}

\noindent
Finally, it remains to extend the estimates (\ref{eqn 3.5.3}) and (\ref{eqn 3.5.4}) to the remaining cases, in which, in terms of Def. \ref{def 2.1.4}, the notation $\Psi_{\varepsilon}(\Gamma)$ refers to either one of the mappings: $\hat{\Psi}_{\varepsilon}^-(\Gamma)=\partial\,H_{\varepsilon,1}(\Gamma):{\rm X}\rightarrow{\rm X}$ and $\hat{\Psi}_{\varepsilon}^+(\Gamma)=\partial\,H_{\varepsilon,2}(\Gamma):{\rm X}\rightarrow{\rm X}$. To accomplish this, we observe that $H_{\varepsilon,i}(\Gamma)=\bigcap_{\alpha\in(0,\varepsilon)}{\rm Cl}\big( G_{\alpha,i}(\Gamma)\big)$ for every $\varepsilon\in(0,\varepsilon_0)$, $i=1,2$, and $\Gamma\in\tilde{{\rm X}}_{\varepsilon}$. Therefore, for any $\varepsilon\in(0,\varepsilon_0)$ and $\Gamma\in\tilde{{\rm X}}_{\varepsilon}$, we have that  ${\Psi}_{\alpha}^\pm(\Gamma)\rightarrow\hat{\Psi}_{\varepsilon}^\pm(\Gamma)$ as $\alpha\uparrow\varepsilon$, from which it follows that $||\Psi_{\varepsilon}^\pm(\Gamma)||\leq\limsup||\Psi_{\alpha}^\pm(\Gamma)||$ and $K\big(\Psi_{\varepsilon}^\pm(\Gamma)\big)\leq\limsup K\big(\Psi_{\alpha}^\pm(\Gamma)\big)$, both as $\alpha\uparrow\varepsilon$.
In view of this, the asserted extensions of (\ref{eqn 3.5.3}) and (\ref{eqn 3.5.4}) to the cases $\Psi_{\varepsilon}(\Gamma)=\hat{\Psi}_{\varepsilon}^\pm(\Gamma)$, $\Gamma\in\tilde{{\rm X}}_{\varepsilon}$, follow from the previously proved cases. 
\vspace{.1in}

\noindent
{\bf Proof of Thm. \ref{thm 3.1.1}.} Concerning Part (a), we refer to Lems. \ref{lem 2.1.4} and \ref{lem 2.1.4*}. 
Concerning Part (b), we obtain the estimates (\ref{eqn 3.1.1}) and (\ref{eqn 3.1.2}) in more detail. By successive application of the estimates (\ref{eqn 3.5.1}), (\ref{eqn 3.5.2}), and (\ref{eqn 3.5.3}), we conclude that:
\begin{equation}
\label{eqn 3.6.1}
K(T_{\varepsilon,i}({\bf\Gamma}))=K(\Psi_{\varepsilon,i}(\Gamma_{\varepsilon,i}))\leq (1+A\varepsilon)K(\Gamma_{\varepsilon,i})+B\,||\Gamma_{\varepsilon,i}||\varepsilon
\end{equation}
$$\leq (1+A\,\varepsilon)\Big[K(\Gamma_i)+\big(K(\Gamma_{h,i})-K(\Gamma_i)\big)(\varepsilon/h)\Big]+B\Big[||\Gamma_i||+\big(||\Gamma_{h,i}||-||\Gamma_i||\big)(\varepsilon/h)\Big]\varepsilon$$
$$\leq K(\Gamma_i)+\Big[K(\Gamma_{h,i})+(Ah-1)K(\Gamma_i)+Bh||\Gamma_i||\Big](\varepsilon/h)$$
$$+\Big[A\,\big(K(\Gamma_{h,i})-K(\Gamma_i)\big)+B\,\big(||\Gamma_{h,i}||-||\Gamma_i||\big)\Big](\varepsilon^2/h),$$
for $i=1,2$ and $0<\varepsilon<\varepsilon_0(h):={\rm min}\{h,\varepsilon_0\}$, where we set $\Gamma_{\varepsilon,i}:=\Phi_{\varepsilon,i}({\bf\Gamma})$ and $\Gamma_{h,i}:=\Phi_{h,i}({\bf\Gamma})$. Similarly, one successively applies (\ref{eqn 3.5.1}), (\ref{eqn 3.5.2}), and (\ref{eqn 3.5.4}) to show that 
\begin{equation}
\label{eqn 3.6.2}
||T_{\varepsilon,i}({\bf\Gamma})||\leq ||\Gamma_i||+\Big[||\Gamma_{h,i}||-||\Gamma_i||+Ch\,K(\Gamma_i)\Big](\varepsilon/h)
\end{equation}
$$+\Big[C\,\big(K(\Gamma_{h,i})-K(\Gamma_i)\big)+Dh\,||\Gamma_i||+Eh\,K(\Gamma_i)\Big](\varepsilon^2/h)$$
$$+\Big[D\,\big(||\Gamma_{h,i}||-||\Gamma_i||\big)+E\,\big(K(\Gamma_{h,i})-K(\Gamma_i)\big)\Big](\varepsilon^3/h),$$
for $i=1,2$ and $0<\varepsilon<\varepsilon_0(h)$ and $\Gamma_{h,i}:=\Phi_{h,i}(\Gamma)$. The estimates (\ref{eqn 3.1.1}) and (\ref{eqn 3.1.2}) follow from this.
\vspace{.1in}

\noindent
Turning to the proof of to Thm. \ref{thm 3.1.1}(c), we set $F(\Gamma_i):=\lambda K(\Gamma_i)+\mu ||\Gamma_i||$, for $i=1,2$, and we combine the Eqs. (\ref{eqn 3.6.1}) and (\ref{eqn 3.6.2}) to see that
$$F(T_{\varepsilon,i}({\bf\Gamma}))\leq F(\Gamma_i)+\Big[F(\Gamma_{h,i})+(\lambda(Ah-1)+\mu Ch)K(\Gamma_i)+(\lambda Bh-\mu)||\Gamma_i||\Big](\varepsilon/h)$$
$$+\Big[(\lambda A+\mu C)+\mu E\varepsilon\Big]\Big[K(\Gamma_{h,i})-K(\Gamma_i)\Big](\varepsilon^2/h)+\Big[\lambda B+\mu D\varepsilon\Big]\Big[||\Gamma_{h,i}||-||\Gamma_i||\Big](\varepsilon^2/h)$$
\begin{equation}
\label{eqn 3.6.3}
+\mu \Big[D||\Gamma_i||+E\,K(\Gamma_i)\Big]\varepsilon^2,
\end{equation}
for $i=1,2$ and $0<\varepsilon<\varepsilon_0(h)$, which generalizes (\ref{eqn 3.1.3}).
\vspace{.1in}

\noindent
Finally, concerning Part (d), by choosing the values $\lambda,\mu>0$ such that $\lambda+\mu=1$ and $A+(\mu/\lambda)C=(\lambda/\mu)B$, we obtain the inequality
\begin{equation}
\label{eqn 3.6.4}
F(T_{\varepsilon,i}({\bf\Gamma}))\leq F(\Gamma_i)+\Big[F(\Gamma_{h,i})+(P_0 h-1)F(\Gamma_i)\Big](\varepsilon/h)
\end{equation}
$$+\Big[\lambda A+\mu(C+E)\Big]K(\Gamma_{h,i})(\varepsilon^2/h)+\Big[(\lambda B+\mu D\Big]||\Gamma_{h,i}||(\varepsilon^2/h)$$
$$+\mu\Big[D||\Gamma_i||+E\,K(\Gamma_i)\Big]\varepsilon^2,$$
valid for $i=1,2$ and $0<\varepsilon<\varepsilon_0(h)$, where $P_0:=((A+\sqrt{A^2+4BC})/2)$. It follows from (\ref{eqn  3.6.4}) that
\begin{equation}
\label{eqn 3.6.5}
F(T_{\varepsilon,i}({\bf\Gamma}))\leq  F(\Gamma_i)+\Big(F_{h,i}+(P_0 h-1)F(\Gamma_i)\Big)(\varepsilon/h)
\end{equation}
$$+\Big(M^*F_{h,i}+M h\,F(\Gamma_i)\Big)(\varepsilon^2/h)$$
for $i=1,2$ and $0<\varepsilon<\varepsilon_0(h)$, where $M^*:=[A+D+(\mu/\lambda)(C+E)+(\lambda/\mu)B]$, $M:=(D+(\mu/\lambda)E)$, and $F_{h,i}:=F(\Gamma_{h,i})=F(\Phi_{h,i}(\Gamma))$. This is the estimate (\ref{eqn 3.1.4}), given in greater detail.
\vspace{.1in}

\noindent
{\bf Proof of Thm. \ref{thm 3.1.2}.} 
By (\ref{eqn 3.6.5}), we have that
\begin{equation}
\label{eqn 3.6.6}
F(T_{\varepsilon,i}({\bf\Gamma}))\leq F(\Gamma_i)+\Big[(1+M^*\varepsilon)F(\Gamma_{h,i})+(Mh\varepsilon+P_0h-1)F(\Gamma_i)\Big](\varepsilon/h)
\end{equation}

\noindent+
for $i=1,2$, ${\bf\Gamma}\in {\bf X}$, and $\varepsilon\in(0,\varepsilon_0(h)]$. It follows from Eq. (\ref{eqn 3.6.6}) that Eq. (\ref{eqn 3.1.5}) holds for $i=1,2$, all $\varepsilon\in(0,\hat{\varepsilon}_0(h)]$, and all ${\bf \Gamma}\in{\bf X}$, where we define $\hat{\varepsilon}_0(h):={\rm min}\{\varepsilon_0(h),(P_0/M), (1/M^*)\}$. Thus Part (a) holds.
To prove Part (b), we let ${\bf\Gamma}_\varepsilon\in{\bf Y}$ denote a fixed point of ${\bf T}_\varepsilon$ which is obtained as the limit of a monotone sequence of operator iterates $({\bf\Gamma}_{\varepsilon,n})_{n=1}^\infty$, as in Thm. \ref{thm 2.1.2}. For fixed $0<h\leq 1/2$ and $\varepsilon\in(0,\hat{\varepsilon}_0(h)]$, let $\Gamma_{\varepsilon,h,n,i}:=\Phi_{h,i}({\bf\Gamma}_{\varepsilon,n})$ and $\Gamma_{\varepsilon,h,i}:=\Phi_{h,i}({\bf\Gamma}_{\varepsilon})$, 
$i=1,2$. It is easily seen that $F(\Gamma_{\varepsilon,h,n,i})$ is finite for each $n\in N$ and $i=1,2$, and that $F(\Gamma_{\varepsilon,h,n,i})\rightarrow F(\Gamma_{\varepsilon,h,i})<\infty$ as $n\rightarrow\infty$, $i=1,2$. 
Therefore, there exists a bound $M=M(\varepsilon,h)$ such that $F(\Gamma_{\varepsilon,h,n,i})\leq M${}
 uniformly for all $n\in N$ and $i=1,2$. In view of this, it follows from (\ref{eqn 3.1.5}) that for $i=1,2$ and for all $n\in N$, we have
\begin{equation}
\label{eqn 3.6.8}
F(\Gamma_{\varepsilon,n+1,i})\leq F(\Gamma_{\varepsilon,n,i})+(2M\varepsilon/h),\end{equation}
since $2hP_0\leq 1$, and 
\begin{equation}
\label{eqn 3.6.9}
F(\Gamma_{\varepsilon,n+1,i})\leq F(\Gamma_{\varepsilon,n,i})-(\alpha_n\varepsilon/h),
\end{equation}
where $\alpha_n\leq(1-2hP_0)F(\Gamma_{\varepsilon,n,i})-2M$. It follows from Eqs. (\ref{eqn 3.6.8}) and (\ref{eqn 3.6.9}) by induction that
\begin{equation}
\label{eqn 3.6.10}
F(\Gamma_{\varepsilon,n,i})\leq (2M/(1-2hP_0))+(2M\varepsilon/h)
\end{equation}
for $i=1,2$ and for all sufficiently large $n\in N$, where $M=M(\varepsilon,h)$. The assertion follows from Eq. (\ref{eqn 3.6.10}) in the limit as $n\rightarrow\infty$. Finally, the proof of Part (c), based on Parts (a) and (b) and Lem. \ref{lem 2.2.5}, is included in the assertion of Part (c).
\vspace{.1in}

\noindent
{\bf Proof of Cor. \ref{cor  3.1.1a}.}
 We have $A:={\rm max}\{C_4,C_6\}$, $B:=\overline{C}_2\leq \sqrt{2}C_2$, $C:=(1/\underline{a})$, $D:=C_3$, and $E:=(C_1^2/2)$. In view of this, the assertion follows from inequalities in Sections \ref{subsection 4.2} and \ref{subsection 4.3} by means of the inequalities: $A_1^2\leq 2\delta A_2$, $\varepsilon_0\leq (1/2)$, and $\varepsilon_0\leq(\underline{A}^2/2A_1)$. In fact we have that:
$$4C_1=C_6:=(8A_1/\underline{A}^2)\leq (8\sqrt{2\delta A_2}/\underline{A}^2),$$
$$C_2:=[(2\underline{A}A_2+4A_1^2)/\underline{A}^3]\leq [(2\underline{A}+8\delta)A_2)\big/\underline{A}^3],$$
$$C_3:=(A_1/\underline{A}^3)+(4A_1^2\varepsilon_0/\underline{A}^5)\leq(A_1/\underline{A}^3)+(8\delta A_2/\underline{A}^5)\varepsilon_0$$
$$C_5:=3(A_1/\underline{A}^3)+3\big[\big((\underline{A}+\underline{A}^2)A_1^2+
\underline{A}^3A_2\big)\big/\underline{A}^5\big]\varepsilon_0+(6A_1^2\big/\underline{A}^6)\varepsilon_0^2$$
$$C_4:=(1+C_5\varepsilon_0)(2A_1/\underline{A}^2)^2\varepsilon_0+C_5.$$

\subsection{Arc-length estimates for thin solutions}
\label{subsection 5.4}
\begin{remark}
\label{Z.2}
The present section is motivated by the search for upper bounds for the arc-length per $P$-period of solutions of Prob. \ref{prob 4I.+} at the pair ${\bm\lambda}\in\Re_+^2$, which are uniform in the most general possible sense as $\lambda_1,\lambda_2\rightarrow\infty$. Namely, this is a natural requirement for satisfying the estimate (\ref{eqn 2.7.1-1}). The first observation is that the fixed-point estimates in Thm. \ref{thm 3.1.2} are not helpful, because they do not hold uniformly relative to arbitrarily thin strips $\Omega({\bm\Gamma})$, as occur in the case where $\lambda_1,\lambda_2\rightarrow\infty$. This suggests an alternate arc-length study mostly restricted to the narrow-stream-limit, as introduced in Section 1.1. 
\end{remark}
\begin{lemma}
\label{lem Z.3}
{(Uniform self-separation property of arcs satisfying (\ref{eqn Z.1}))}
In the plane $\Re^2$, let be given a $P$-periodic region $G$ having exterior tangent balls of uniform radius at all boundary points, and a $P$-periodic (in $x$) function $a(p):\Re^2\rightarrow\Re_+$ in $\boldsymbol{\cal{A}}$ such that $|{\nabla}\,{\rm ln}\big(a(p)\big)|\leq B$ and $\Delta\,{\rm ln}\big(a(p)\big)$ $\geq H>0$, both uniformly in $G$. For any $\delta\geq 0$, let $\cal{M}(\delta)$ denote the family of all $P$-periodic (in $x$) directed $C^2$-arcs $\Gamma\in{\rm X}(G)$, each of which satisfies the condition:
\begin{equation}
\label{eqn Z.6}
\big|K(p)-(\partial\big/\partial\boldsymbol{\nu})\,{\rm ln}\big(a(p)\big)\big|\leq\delta
\end{equation}
for each $p\in\Gamma$ (see (\ref{eqn Z.2a})), where $K(p)$ (resp. $\boldsymbol{\nu}$) denotes the counter-clockwise-oriented curvature of (resp. the left-hand normal to) the arc $\Gamma$ at $p\in\Gamma$.
(Also, each $\Gamma$ does not cross itself, and therefore divides its complement $\Re^2\setminus\Gamma$ into two domains $D^\pm$ (as in Def. \ref{def 2.1.1}).) Then there exist uniform positive constants $\delta_0,\alpha_0>0$ such that 
$|p-q|\geq\alpha_0$, uniformly for all $\Gamma\in{\cal M}(\delta_0)$ and for all points $p,q\in\Gamma$ such that $||\gamma(p,q)||\geq\big(\pi\big/(B+\delta_0)\big)$, where $\gamma(p,q)$ denotes the shortest sub-arc of $\Gamma$ joining the points $p,q\in\Gamma$.
\end{lemma}
\vspace{.1in}

\noindent
{\bf Proof.} For any small value $\delta\geq 0$, let $\cal{M}(\delta)$ denote a family of all doubly-infinite, $P$-periodic (in $x$) directed smooth arcs $\Gamma$ satisfying (\ref{eqn Z.6}). 
By (\ref{eqn Z.6}), the numbers $K_0=K_0(\delta):=B+\delta$ and $R_0=R_0(\delta):=\big(1/(B+\delta)\big)$ (where $B:=\sup\big\{|{\nabla}\,{\rm ln}\big(a(p)\big)|:p\in G\big\}$) serve respectively as the uniform upper bound for the curvatures of the arcs $\Gamma\in {\cal M}(\delta)$, and uniform lower bounds for the corresponding radii of curvature of the same arcs. For any fixed $\Gamma\in{\cal M}(\delta)$, let $S$ denote the set of all ordered pairs $(p,q)\in\Gamma\times\Gamma$ such that $||\gamma||\geq \pi R_0$, where $\gamma$ denotes the shortest sub-arc of $\Gamma$ joining $p$ to $q$. Then $T:=\{(p,q)\in S:||\gamma||=\pi R_0\}\subset S$, where for any pair $(p,q)\in T$, we have that $|p-q|>2R_0$ unless $\gamma$ is exactly a half-circle of radius $R_0$. 
\vspace{.1in}

\noindent
We use $\alpha=\alpha(\Gamma)$ to denote the absolute minimum value of the (bounded, continuous) distance function $|p-q|:S\rightarrow\Re$. In the following study, we restrict our attention to the case where $|p-q|<2R_0$ for some pair $(p,q)\in S\setminus T$. Then it is clear from the above comments that the distance function $|p-q|:S\rightarrow\Re$ achieves an absolute minimum value $\alpha=\alpha(\Gamma)\in[0,2R_0)$ at a pair $(p,q)\in S$ such that $||\gamma||>\pi R_0$. In view of this, the shortest straight line-segment $L$ joining $p\in\Gamma$ to $q\in\Gamma$ is perpendicular to $\Gamma$ at both endpoints $p,q$.
\vspace{.1in}

\noindent

Given $\Gamma\in{\cal M}(\delta)$ and the "closest points" $p,q\in\Gamma$, we use $\omega$ to denote a simply-connected domain whose boundary is given by $\partial\omega=\gamma\cup L$, where $\gamma$ (resp. $L$) denotes the shortest sub-arc of $\Gamma$ (resp. the straight line-segment) joining $p$ and $q$. Since $L$ is perpendicular to $\Gamma$ at $p$ and $q$, and since $\Gamma$ does not cross itself, we have that $\int_{\gamma} K(p)ds=\pm\pi$. Also $$\int_{\partial\omega}\psi_{\boldsymbol{\nu}}(p)\,ds=\int\int_{\omega}\Delta\psi(p)\,dA\geq H\,||\omega||$$ by the divergence theorem, where $\psi(p):={\rm ln}\big(a(p)\big)$ and $\boldsymbol{\nu}$ is the exterior normal to $\omega$, and also $\int_{L}\big|\psi_{\boldsymbol{\nu}}\big|\,ds\leq B\,||L||$ if $L\subset G$. Finally, it follows directly	from (\ref{eqn Z.6}) and the above inequalities that
\begin{equation}
\label{eqn Z.7}
B\,||L||\geq \pi+H\,||\omega||-\delta\,||\gamma||,
\end{equation}  
provided that $L\subset G$, as is always true if $||L||\leq 2\,R_0$, and if the region $G$ has an exterior tangent ball of radius $(3R_0/2)$ at every point $p\in\partial G$.
Assuming that $({\bf i})$:  $\alpha=\alpha(\Gamma)<2R_0$, we define the double-point free strip-like domain $\omega_{(\alpha/2)}$, consisting of all points $p\in\omega$ such that ${\rm dist}(p,\Gamma)<(\alpha/2)$. One sees that $\omega_{(\alpha/2)}\subset\omega$, and therefore that 
\begin{equation}
\label{eqn Z.8}
||\omega||\geq||\omega_{(\alpha/2)}||=(\alpha/2)\int_{\gamma}\big(1+(\alpha/4)K(t)\big)\,dt\geq
(\alpha/2)||\gamma||-\pi(\alpha^2/8),
\end{equation} 
involving a counter-clockwise integral of the counter-clockwise oriented curvature function $K(t)$. Thus, assuming ${\bf(i)}$, it follows by (\ref{eqn Z.7}) and (\ref{eqn Z.8}) that 
\begin{equation}
\label{eqn Z.9}
B\alpha\geq\pi\big(1-(H\alpha^2/8)\big)+\big((H\alpha/2)-\delta\big)||\gamma||,
\end{equation}
and it follows from ${\bf(ii)}$:  $2\delta\leq H\alpha$ and (\ref{eqn Z.9}) that 
${\bf(iii)}$:  $B\alpha\geq\pi\big(1-(H\alpha^2/8)\big)$, without regard to the length of $\gamma$. Here, the inequality ${\bf(iii)}$ is equivalent to 
\begin{equation}
\label{eqn Z.10}
\alpha=\alpha(\Gamma)\geq E_0=E_0(B,H):=\frac{2\sqrt{2}\,\pi}{\sqrt{2}\,B+\sqrt{2\,B^2+\pi^2\,H}}
\end{equation}
In the case in where $E_0(B,H)\leq 2R_0$, it follows that $E_0(B,H)$ is a lower bound for $\alpha(\Gamma)$. If $E_0(B,H)>2R_0$, on the other hand, then the inequalities (\ref{eqn Z.10}) and ${\bf(i)}$ are mutually contradictory, implying that the assumption ${\bf(i)}$ is impossible and therefore that $\alpha(\Gamma)\geq 2R_0$. 
In other words, given a small, positive value $\delta_0>0$, we have
\begin{equation}
\label{eqn Z.11}
\alpha(\Gamma)\geq \alpha_0\,\,{\rm for}\,\,{\rm any}\,\,0\leq\delta\leq \delta_0\,\,{\rm and}\,\,\Gamma\in{\cal M}(\delta)\,\,{\rm such}\,\,{\rm that}\,\,H\alpha(\Gamma)\geq 2\delta,
\end{equation}
where we set $\alpha_0:=\min\big\{2R_0(\delta_0), E_0(B,H)\big\}$. It remains to prove the following: Given $\alpha_0$ defined above, there exists a value $\delta_1\in(0,\delta_0]$ so small that
\begin{equation}
\label{eqn Z.11a}
\alpha(\Gamma)\geq(\alpha_0\big/2)>0\,\,{\rm for}\,\,{\rm all}\,\,\Gamma\in{\cal M}(\delta_1).\,\,
\end{equation}
Toward the proof of (\ref{eqn Z.11a}) from (\ref{eqn Z.11}), we define $\phi(R_1,\delta):=\inf\big\{\alpha(\Gamma):\Gamma\in{\cal M}(\delta)$, $\pi R_0\leq||\gamma||\leq R_1\big\}$ for any $\delta\in[0,\delta_0]$ and $R_1\in[\pi R_0,\infty)$. Then ${\bf(iv)}$ $\phi(R_1,0)\geq\alpha_0$ for the fixed value $\alpha_0>0$ defined following (\ref{eqn Z.11}), and for any value $R_1\in(\pi R_0,\infty)$, as is seen by setting $\delta=0$ in (\ref{eqn Z.11}).
We claim that, given any fixed value $\alpha_1\in(0,\alpha_0)$, we have that ${\bf(v)}$:  $\phi(R_1,\delta)\geq \alpha_1$ for any $R_1\in(\pi R_0,\infty)$ and for any value $\delta\in(0,\delta_1]$, where the value $\delta_1=\delta_1(\alpha_1)\in(0,\delta_0]$ is sufficiently small, depending on the value $\alpha_1$. In fact if ${\bf(v)}$ is false, then there exist constants $\alpha_{0,1}\in(0,\alpha_0)$ and $R_1\in(\pi R_0,\infty)$, and a positive null-sequence $\big(\delta_n\big)_{n=1}^\infty$ of values in $(0,\delta_0]$ such that the corresponding sequence $\big(\Gamma_n\big)_{n=1}^\infty$ of $P$-periodic $C^{2}$-arcs such that $\Gamma_n\in{\cal M}(\delta_n)$ and the absolute curvature of $\Gamma_n$ is therefore uniformly bounded by $K_0(\delta_n):=B+\delta_n$, such that the related sequences $\big(\gamma_n\big)$, $\big(\omega_n\big)$, $\big(p_n\big)$, $\big(q_n\big)$, $\big(L_n\big)$ (in which $\gamma_n\subset\Gamma_n$ and the terms $\gamma_n, \omega_n, p_n, q_n, L_n$ are defined by analogy to $\gamma$, $\omega$, $p$, $q$, and $L$) such that $\pi R_0\leq||\gamma_n||\leq R_1$ a

nd $\alpha(\Gamma_n)=|p_n-q_n|\leq\alpha_{0,1}$, both for all $n\in N$. In view of the uniformly-boundedness of the absolute   curvatures of the arcs $\Gamma_n$, $n\in N$, it follows by the theorem or Ascoli-Arzela that there exists a subsequence (which we still index by $n\in N$) such that $\Gamma_n\rightarrow\Gamma\in{\cal M}(0)$, $||\gamma_n||\rightarrow||\gamma||\leq R_1$, and $\alpha(\Gamma_n)\rightarrow\alpha(\Gamma)\leq\alpha_{0,1}$, which contradicts ${\bf(iv)}$ $\phi(R_1,0)\geq\alpha_0$, completing the proof of ${\bf(v)}$. In view of ${\bf(v)}$, we have that ${\bf(vi)}$: $H\alpha(\Gamma)-2\delta\geq H\alpha_1-2\delta>0$ for all arcs $\Gamma\in{\cal M}(\delta_1)$, where we define $\alpha_1:=(\alpha_0/2)$ and $\delta_1:=\min\big\{\delta_1(\alpha_1), \big(H\alpha_1\big/2\big)\big\}$. In view of ${\bf(vi)}$, the assertion (\ref{eqn Z.11a}) follows from (\ref{eqn Z.11}).
\vspace{.1in}

\begin{lemma} 
{(Sequences of solution-arcs not having a uniform bound for their arc-lengths per $P$-period (in $x$))}
\label{lem Z.4}
In $\Re^2$, let be given an infinite sequence $\big(\Gamma_n\big)_{n=1}^\infty$ of $P$-periodic (in $x$) directed arcs in ${\rm X}(G)$ having uniformly bounded absolute curvature per $P$-period, and each having an arc-length parametrization $p_n(t):\Re\rightarrow {\rm Cl}(G)$ such that the mapping $p'_n(t):\Re\rightarrow\partial B_1(0)$ is $L_n$-periodic (where $L_n:=||\Gamma_n||=$ the length of one $P$-period (in $x$) of $\Gamma_n$). 
Assume each arc $\Gamma_n$ does not cross itself, and has absolute curvature $|p''_n(t)|$ never exceeding a uniform bound $K_0$. Then:
\vspace{.1in}

\noindent
(a) If there exists a closed, bounded ball $Q$ which contains one $P$-period (in $x$) of each arc $\Gamma_n$, and if, nevertheless, we have that $L_n\rightarrow\infty$ as $n\rightarrow\infty$, then
there exists a subsequence $\big(\Gamma_{n(i)}\big)_{i=1}^\infty$ and an arc $\Gamma\subset Q$ of bounded absolute curvature and infinite length such that $\Gamma$ has the arc-length parametrization $p(t):\Re\rightarrow Q$, and such that $p_{n(i)}(t)\rightarrow p(t)$ as $i\rightarrow\infty$, uniformly relative to any compact subset of $\Re$. Also, $\Gamma$ does not cross itself and does not have absolute curvature $|p''(t)|$ ever exceeding $K_0$. 
\vspace{.1in}

\noindent
(b) Let $\Gamma\subset Q$ denote the doubly-infinite directed arc of Part (a). Then either $\Gamma$ has infinite length or else it is (or contains) a closed (i.e. periodic) arc $\gamma$, because the function $p(t)$ either is periodic in $\Re$ or else is periodic for sufficiently positive or sufficiently negative $t\in\Re$. In the case where $\Gamma$ has infinite length as $t\rightarrow\infty$ (resp. $t\rightarrow -\infty$), the set of accumulation points $\gamma^+$ (resp. $\gamma^-$) of $\Gamma$ as $t\rightarrow\infty$ (resp. $t\rightarrow-\infty$) is a closed (i.e. periodic) $C^{2,1}$-arc having all the properties previously stated regarding $\Gamma$ in Part (a).
\vspace{.1in}

\noindent
(c) If $\Delta\,{\rm ln}\big(a(p)\big)\geq H\geq 0$ in $G$, then the infinite arc $\Gamma\subset Q$ with the properties stated in Parts (a) and (b) does not exist. Therefore, the sequence $\big(\Gamma_n\big)_{n=1}^\infty$ with properties stated in Part (a) does not exist.
\end{lemma}
\vspace{.1in}

\noindent
{\bf Proof of Part (a).} Let $p_n(t):\Re\rightarrow Q$, $n\in N$, (such that $p_n(t+L_n)=p_n(t)+(P,0)$ for all $t\in\Re$) denote the corresponding arc-length parametrizations of the arcs $\Gamma_n$, $n\in N$, such that for each $n\in N$, the related $L_n$-periodic functions $p'_n(t), p''_n(t))$ are such that $|p'_n(t)|=1$ and $|p''_n(t)|\leq K_0$, both for all $t\in\Re$. By applying the theorem of Ascoli-Arzela to the sequence of uniformly Lipschitz-continuous forward-tangent-vector functions $T_n(t)=p'_n(t):\Re\rightarrow \partial B_1(0)$, we pass to a subsequence $\big(T_{n(i)}\big)_{i=1}^{\infty}$ such that $T_{n(i)}(t)\rightarrow T(t)$ uniformly on any compact subset of $\Re$ as $i\rightarrow\infty$, where $T(t):\Re\rightarrow\partial B_1(0)$ denotes a $C^{1,1}$-function such that $|T(t)|=1$ and $|T'(t)|\leq K_0$, both for all $t\in\Re$. It also follows for the same subsequence that $p_n(t)\rightarrow p(t)$ on any compact subset of $\Re$ as $n\rightarrow\infty$, where the function $p(t):\Re\rightarrow Q$ is chosen such that $p(0)=\lim p_n(0)$ as $n\rightarrow\infty$ and $p'(t)=T(t)$ at all $t\in\Re$. Clearly $p(t):\Re\rightarrow Q$ is the arc-length parametrization of an arc $\Gamma\subset Q$ of infinite length and uniformly bounded curvature. Also, the arc $\Gamma$ does not cross itself, since the arcs $\Gamma_n$ don't.
\vspace{.1in}

\noindent
{\bf Proof of Part (b).} (See the proof of Lem. \ref{lem 4B.2}(a) for some details.) Let $p(t):\Re\rightarrow Q$ denote an arc-length parametrization of a directed arc $\Gamma\subset Q$ of infinite length, and of absolute curvature uniformly bounded by $K_0$. Then $|p(t)-p(\tau)|\geq(2/\pi)|t-\tau|$ for all $t,\tau\in\Re$ such that $|t-\tau|\leq\pi R_0$, where $R_0=(1/K_0)$.
Let $p_n:=p(t_n)\in Q$ for all $n\in N$, where the sequence $\big(t_n\big)_{n=1}^\infty$ is chosen such that $t_n<t_{n+1}<t_n+(\pi R_0/2)$ for all $n\in N$, so that $|p_{n+1}-p_n|\geq R_0$ for all $n$. Using the compactness of $Q$, we pass to a subsequence $\big(t_{n(i)}\big)_{i=1}^\infty$ such that $q_i:=p(t_{n(i)})\rightarrow p_0\in Q$ as $i\rightarrow\infty$ for some accumulation point $p_0$ of $\Gamma$ in $Q$. Then obviously $|q_{i+1}-q_{i}|\rightarrow 0$ as $i\rightarrow\infty$. Therefore ${\rm dist}(q_i,\Gamma_i)=|q_i-r_i|\rightarrow 0$ as $i\rightarrow\infty$, where $\Gamma_i:=\Gamma\setminus\hat{\Gamma}_i$, in which $\hat{\Gamma}_i$ denotes the arc-segment of $\Gamma$ associated with all $t\in\Re$ such that $|t-t_{n(i)}|\leq(\pi R_0/2)$, and where $r_i$ denotes a point in $\Gamma_i$ closest to $q_i$. Clearly, if $i$ is sufficiently large, then $r_i$ is not an endpoint of $\hat{\Gamma}_i$, and therefore the straight-line-segment $L_i$ joining $q_i$ to $r_i\in\Gamma_i$ is perpendicular to $\Gamma$ at the point $r_i$. In view of the bounded curvature of $\Gamma$, it follows that up to a small error for sufficiently large $i$, the directed arc-segment $\gamma_i$ of $\Gamma$ joining $q_i$ to $r_i$ (which does not cross itself) must have a turning angle of (b.1) $\,\,\pm \pi$ (half-turn either way) or (b.2) $\,\,\pm 2\pi$ (full turn either way). In either case, it follows that $||\omega_i||\geq \big(||\gamma_i||^2\big/4\pi\big)\geq(\pi R_0^2\big/4)$ up to a small error for sufficiently large $i$, where $\omega_i$ denotes the bounded connected region whose boundary consists of $\gamma_i$ and $L_i$ and $||\omega_i||$ is its Euclidean area.
The arc-segments of $\Gamma_i$ cannot cross the line-segment $L_i$, except possibly at its endpoints, $q_i, r_i\in\Gamma$, since otherwise $r_i$ would not minimize the distance from the point $q_i$ to the arc $\Gamma_i$. Since $\Gamma$ cannot cross itself, it follows that $\Gamma$ cannot cross the arc $\gamma_i$, and thus cannot pass from the domain $\omega_i$ to the interior of the complement of $\omega_i$, or visa versa, except possibly by passing through an end-point of $L_i$.
\vspace{.1in} 

\noindent
In case (b.1), the directed arc $\Gamma$ enters $\partial\omega_i$ at $q_i$, follows $\partial\omega_i$ to $r_i$ along $\gamma_i$, and then permanently exits ${\rm Cl}(\omega_i)$ at $r_i$. Therefore, given large $i\in N$, there exists $j\in N$ large enough such that $i<j$ and $q_j:=p(t_{n(j)})>q_i$ (in terms of the natural ordering in $\Gamma$). Therefore, there is a domain $\omega_j$ associated with $q_j$ which does not intersect $\omega_i$. Should case (b.1) occur infinitely often, it would follow that there exist infinitely many pairwise disjoint regions $\omega_i$ in $Q$, each having area exceeding $(\pi R_0^2\big/4)$. This contradiction of the compactness of $Q$ shows that case (b.2) must occur infinitely often. In the case (b.2), we observe that the directed arc $\Gamma$ enters ${\rm Cl}(\omega_i)$ at $q_i$, then follows $\partial \omega_i$ on $\gamma_i$ until it enters $\omega_i$ at the point at $r_i$, after which it must remain inside $\omega_i$ permanently as $i\rightarrow\infty$. It is important to distinguish two sub-cases of case (b.2), namely "left full turn" and "right full turn" (up to a small error for large $n$ in each case). In fact it is easily shown for sufficiently large $i$ that if a "left full turn" configuration at $i$ precedes a "right full turn" configuration (or visa versa) at $k>i$, then there exists a type (b.1) configuration at some $j$ with $i<j<k$. Therefore, the sequence cannot include infinitely-many  type (b.2) configurations of each sub-type, and there must exist a number $i_0\in N$ so large that all the type (b.2) configurations for $i\geq i_0$ are of the same sub-type. Thus, the infinite directed arc $\Gamma$ eventually becomes a counter-clockwise or clockwise spiral, and it has a closed directed arc $\gamma$ (which does not cross itself and satisfies the same absolute curvature bound as $\Gamma$) as its set of accumulation points as $t\rightarrow\infty$. Clearly $\Gamma$ does not cross $\gamma$, and is therefore located weakly inside or weakly outside $\gamma$; the convergence to $\gamma$ being eventually monotone in either case. 
\vspace{.1in}

\noindent
{\bf Proof of Part (c).} 
Under the assumptions of Part (a), there exists a doubly infinite spiral limit arc $\Gamma\subset Q$ whose sets $\gamma^\pm$ of accumulation points (as $t\rightarrow\pm\infty$) are simple closed curves each solving the limiting equation:
\begin{equation}
\label{eqn Z.2.1}
K(p)=\big(\partial\big/\partial\boldsymbol{\nu}\big)\,{\rm ln}\big(a(p)\big),
\end{equation}
where $\boldsymbol{\nu}$ (resp. $K(p)$) denotes the left normal (resp. left oriented curvature) at $p\in\gamma^\pm$. By integrating (\ref{eqn Z.2.1}) on $\gamma^\pm$, we conclude that
\begin{equation}
\label{eqn Z.2.2}
-2\pi=\int_{\gamma^\pm}K(p)ds=\int_{\gamma^\pm} \big(\partial\big/\partial\boldsymbol{\nu}\big)\,{\rm ln}\big(a(p)\big)ds
\end{equation}
$$=\int\int_{\omega^\pm}\Delta\,{\rm ln}\big(a(p)\big)dA\geq H\,||\omega^\pm||\geq 0,$$
\noindent
where $\gamma^\pm=\partial\omega^\pm$. This contradiction proves the assertion.

\begin{theorem} 
{(Uniform upper bound for the arc-length per $P$-period of periodic solution arcs)}
\label{thm Z.2}
Assume in the context of Prob. \ref{prob 4I.+} that $G$ is a $P$-periodic (in $x$) strip-like domain having interior and exterior tangent balls at all boundary points, and that the given $C^{2,1}$-function $a(p):\Re^2\rightarrow\Re_+$ is strictly logarithmically sub-harmonic throughout $G$, in fact we assume for constants $0<H\leq H_1$ that $H\leq\Delta\,{\rm ln}\big(a(p)\big)\leq H_1$ throughout $G$. For each $n\in N$, let ${\bf\Gamma}_n=(\Gamma_{n,1},\Gamma_{n,2})\in{\bf X}(G)$ denote a $P$-periodic (in $x$) classical solution of Prob. \ref{prob 4I.+} at ${\bm\lambda}_n=(\lambda_{n,1},\lambda_{n,2})\in\Re_+^2$, and let $\Gamma_n:=\{p\in\Omega_n:U_n(p)=1/2\}$ denote the corresponding $P$-periodic (in $x$) "center arc" of the solution (where $U_n(p):=U({\bm\Gamma}_n;p)$ in the closure of $\Omega_n:=\Omega({\bm\Gamma}_n)$). Here, we choose the vector sequence $\big({\bm\lambda}_n\big)_{n=1}^\infty$ in $\Re_+^2$ such that $\lambda_{n,i}\rightarrow\infty$ and $\lambda_{n,i}\mu_n\rightarrow 0$ both as $n\rightarrow\infty$ for $i=1,2$, where we set $\mu_n:={\rm ln}(\lambda_{n,2}\big/\lambda_{n,1})$ (compare to 
Thm. \ref{thm 4I.2}). Then there exists a uniform upper bound $M$ for the arc-lengths of the restrictions to a single $P$-period (in $x$) of the "center arcs" $\Gamma_n$ of the solutions ${\bf\Gamma}_n$, $n\in N$. The proof is based on the  curvature-estimate:
\begin{equation}
\label{eqn Z.2a}
\big|K_n(p)-\big(\partial\,{\rm ln}\big(a(p)\big)\big/\partial{\bm\nu}_n\big)\big|\leq \overline{\lambda}_n\,\,\overline{A}\,|\mu_n|+\big(2\overline{A}H_1\,{\rm exp}(\overline{\mu}))\big/\underline{A}^2\underline{\lambda}_n\big)
\end{equation}
valid at all points $p\in{\rm Cl}(\Omega_n)$, where $K_n(p)$ (resp. ${\bm \nu}_n$) denotes the left curvature of (resp. the left normal to) the level curve of $U_n$ through the point $p$, and where $\overline{\lambda}_n:={\rm max}\{\lambda_{n,1},\lambda_{n,2}\}$
and $\underline{\lambda}_n:={\rm min}\{\lambda_{n,1},\lambda_{n,2}\}$
(see Lem. \ref{lem Z.1}, Eq. (\ref{eqn Z.1})).
\end{theorem}

\begin{lemma}
\label{lem Z.5}
{(Toward the proof of Thm. \ref{thm Z.2}: Uniform upper bound for the horizontal spans of the single $P$-periods (in $x$) of all $P$-periodic solution arcs.)}         
Under the assumptions of Thm. \ref{thm Z.2}, there exists a uniform bound $M$ such that $x_n(t)-x_n(\tau)\leq M$ uniformly for all $n\in N$ and $t,\tau\in\Re$ with $t\leq\tau$, where $p_n(t)=(x_n(t),y_n(t)):\Re\rightarrow G$ denotes the arc-length parametrization of $\Gamma_n$. It follows that $|x_{n,1}-x_{n,2}|\leq (M+P)$ uniformly for the center curves of all solutions $\Gamma_n$, and for all points $p_{n,1}=(x_{n,1},y_{n,1})$ and $p_{n,2}=(x_{n,2},y_{n,2})$, both contained in the same $P$-period of $\Gamma_n$
\end{lemma}
\noindent
{\bf Proof.} Assume in Prob. \ref{prob 4I.+} that $G\subset\{|y|\leq b\}$ for a constant $b<\infty$. For each $n\in N$, let $p_n(t)=\big(x_n(t),y_n(t)\big):\Re\rightarrow G$ denote the arc-length parametrization of the $P$-periodic (in $x$) "center arc" $\Gamma_n$ of the solution ${\bf\Gamma}_n\in{\bf X}(G)$ of Prob. \ref{prob 4I.+} at ${\bm\lambda}_n$, and let $L_n$ denote the arc-length associated with one $P$-period of $\Gamma_n$. Suppose that no uniform upper bound $M:=\sup\big\{M_n:n\in N\big\}<\infty$ exists as described in the assertion, where we define $M_n:=\max\{x_n(t)-x_n(\tau):t\geq\tau\,\,{\rm in}\,\,\Re\}$ for each $n\in N$.
Then by passing to a subsequence (still indexed by $n\in N$), we can assume that $M_n\rightarrow\infty$ as $n\rightarrow\infty$. Obviously $L_n\geq M_n$ for any $n$. For any $n\in N$ and any integer $k$, the mapping $p_{n,k}(t):=\big(x_n(t)+kP, y_n(t)\big)=p_n(t+k\,L_n):\Re\rightarrow G$ is the arc-length parametrization of the horizontal $kP$-translation of $\Gamma_n$, which coincides with $\Gamma_n$. In view of this, we can translate the $x$ and $t$ variables to arrange things so that $x_n(\pm L_n/2)=\pm (M_n/2)$ for each $n\in N$. Thus, each arc $\Gamma_n$ is partitioned into the three arcs: $\hat{\Gamma}_n:=\{p_n(t):|t|\leq (L_n/2)\}$ and $\Gamma_n^\pm:=\{p_n(t):\pm t>(L_n/2)\}$, such that, in terms of the natural ordering in $\Gamma_n$, $\hat{\Gamma}_n$ passes from $\{x=+(M_n/2)\}\cap G$ to $\{x=-(M_n/2)\}\cap G$ without intersecting $\{|x|>(M_n/2)\}$, while $\Gamma_n^-$ passes from $\{x=-\infty, |y|\leq b\}$ to $\{x=(M_n/2)\}\cap G$ without intersecting $\{x>(M_n/2)\}$ and $\Gamma_n^+$ passes from $\{x=-(M_n/2)\}\cap G$ to $\{x=\infty, |y|\leq b\}$ without intersecting $\{x<-(M_n/2)\}$. We can assume, after applying the theorem of Ascoli-Arzela to pass to a subsequence, that the arcs $\hat{\Gamma}_n$ and $\Gamma_n^\pm$ converge respectively to arcs $\hat{\Gamma}$ (passing from $\{x=\infty, |y|\leq b\}$ to $\{x=-\infty, |y|\leq b\}$) and $\Gamma^\pm$ (both passing from  $\{x=-\infty, |y|\leq b\}$ to $\{x=\infty, |y|\leq b\}$). Moreover, the three arcs $\hat{\Gamma}$ and $\Gamma^\pm$ do not cross each other, because the original curves $\Gamma_n$, $n\in N$ do not cross themselves. In other words, we have
\begin{equation}
\label{eqn Z.3}
\Gamma^-\leq\hat{\Gamma}\leq\Gamma^+.
\end{equation}
It also follows from the uniform self-separation property of the arcs $\Gamma_n$, $n\in N$ (see Lem. \ref{lem Z.3}) that ${\rm dist}\big(\hat{\Gamma}_n,\Gamma^\pm_n\big)\geq(\alpha_0\big/2)$ for all sufficiently large $n\in N$, from which it follows in the limit that ${\rm dist}\big(\hat{\Gamma},\Gamma^\pm\big)\geq(\alpha_0\big/2)$, and therefore that the inequalities in (\ref{eqn Z.3}) are both strict. Also, it follows from (\ref{eqn Z.2a}) in the limit as $n\rightarrow \infty$ that the arcs $\Gamma^\pm$ satisfy the condition (\ref{eqn Z.2.1}) (namely $K(p)=(\partial/\partial{\bm\nu}){\rm ln}\big(a(p)\big)$, where $K(p)$ (resp. ${\bm\nu}$)
denotes the left-oriented curvature of (resp. left normal to)  $\Gamma^\pm$ at any point $p\in\Gamma^\pm$. By substituting (\ref{eqn Z.2.1}) into the divergence theorem relative to one $P$-period $\omega$ of the strip-like region bounded by the arcs $\Gamma^\pm$, we conclude that
$$\int\int_{\omega}\Delta\,{\rm ln}\big(a(p)\big)\,dA=0.$$
\noindent
But this is impossible, in view of our assumption that $\Delta\,{\rm ln}\big(a(p)\big)>0$ in $G$. This contradiction proves the assertion.
\vspace{.1in}

\noindent
{\bf Proof of Thm. \ref{thm Z.2}.} In the context of Thm. \ref{thm Z.2}, let be given a sequence $\big(\Gamma_n\big)_{n=1}^\infty$ of center-arcs $\Gamma_n\in{\rm X}(G)\cap C^{2}$ satisfying (\ref{eqn Z.2a}). 
By Lem. \ref{lem Z.5}(b), the horizontal span of one $P$-period of the $P$-periodic arc $\Gamma_n$ remains  uniformly bounded as $n\rightarrow\infty$. Therefore, there exists a compact ball $Q$ containing one $P$-period (in $x$) of each of the arcs $\Gamma_n$, $n\in N$. In view of this, the sequence $\big(L_n\big)_{n=1}^\infty$ remains uniformly bounded as $n\rightarrow\infty$ (where $L_n$ denotes the arc-length of one $P$-period of the center-arc $\Gamma_n$), since the unboundedness of the same sequence would contradict Lem. \ref{lem Z.4}.

\end{document}